\documentclass[10pt]{article}
\usepackage{amsmath}
\usepackage{amsfonts}
\usepackage{graphicx}
\usepackage{amssymb}
\usepackage{leftidx}
\usepackage{color}
\usepackage{bbm}

\allowdisplaybreaks

\newtheorem{condition**}{A*}
\newtheorem{condition***}{C*}
\newtheorem{condition*}{C}

\newtheorem{example}{Example}[section]
\newtheorem{proposition}{Proposition}[section]

\newtheorem{definition}{Definition}[section]
\newtheorem{theorem}{Theorem}[section]
\newtheorem{lemma}{Lemma}[section]
\newtheorem{remark}{Remark}[section]

\begin{document}

\title{Closed-loop $\alpha$-Potential Stochastic Differential Games via a BSDE Approach}
\author{
Xun Li$^{1}$\thanks{Department of Applied Mathematics, The Hong Kong Polytechnic University,
Hong Kong, China. E-mail: li.xun@polyu.edu.hk. Xun Li is supported by the Research Grants Council of Hong
Kong under grant 15225124, PolyU 4-ZZP4.}, \
Liangquan Zhang$^{2}$\thanks{School of Mathematics, Renmin University of China, Beijing 100872, China.
Email: xiaoquan51011@163.com. L. Zhang acknowledges the financial support partly
by the National Nature Science Foundation of China (Grant No. 12571486,
12171053) and the Fundamental Research Funds for the Central Universities,
and the Research Funds of Renmin University of China (No.23XNKJ05). Email:
xiaoquan51011@163.com.}}
\maketitle

\begin{abstract}

In this paper, we study the closed-loop $\alpha $-potential stochastic
differential game (SDG) problem as a continuation of our prior research on
open-loop control (see \cite{GLZ2025}). By utilizing the backward stochastic
differential equation (BSDE) approach, we derive a precise estimate for the
parameter $\alpha $. Compared to our earlier work \cite{GLZ2025}, this study
incorporates both first- and second-order sensitivity state processes, as
well as the sensitivity of the control process. A distinguishing feature of
this work is that, in the context of $N$-player heterogeneous‑agent games involving mean-field type interactions, we derive an $N$-uniform
upper bound for the minimal potential approximation error. In contrast to the
corresponding open-loop estimates, the closed-loop bound contains feedback-induced contributions that need not vanish with $N$. Consequently, our present estimate does
not in general guarantee $\alpha \rightarrow 0$.






\end{abstract}

\tableofcontents


\section{Introduction\label{sect1}}

\paragraph{Potential games.}

Potential games, originally introduced by Monderer and Shapley (1996, \cite{ms96}),
represent a class of static non-cooperative games wherein the marginal
incentives for any player to unilaterally deviate from a given strategy
profile are precisely reflected in the change of a single real-valued
function---the so-called \emph{potential function}. This structural property
establishes a precise correspondence between individual payoff improvements
and increments in the potential function, thereby transforming the problem
of identifying Nash equilibria in $N$-player strategic settings into a
global optimization task over a unified objective. Consequently,
decentralized rational behavior-driven solely by local payoff maximization-
naturally converges toward globally stable outcomes characterized by
stationary points of the potential function.


This theoretical construct has proven instrumental in modeling and
rigorously analyzing diverse real-world phenomena across multiple
disciplines: in economics, it informs models of oligopolistic competition
and distributed resource allocation (Nie et al., 2022, \cite{nlks22}); in network science,
it underpins equilibrium analyses of traffic assignment and power grid
dispatch (Li \& Vazquez, 2024; Marden et al., 2009; Zhang et al., 2021, \cite{lv2024,mas2009,zmz2021}); in
artificial intelligence, it facilitates convergence guarantees in
multi-agent reinforcement learning and coordination protocols (Ding et al.,
2023; Leonard et al., 2021; Mao et al., 2018; Marden et al., 2022, \cite{dwzj,lopp21,mzz18,mwps22}); and in
robotics, it supports scalable design of distributed control laws for
heterogeneous multi-robot systems (Sun et al., 2024, \cite{swhw2024}).


Nevertheless, the applicability of the potential game framework is
intrinsically limited in dynamic and sequential settings. Although each stage
game in a Markov decision process may individually admit a potential
structure, the resulting dynamic game generally fails to satisfy the
defining consistency condition for a Markov potential game (Leonard et al.,
2021, \cite{lopp21}). More fundamentally, enforcing potentiality across time horizons
imposes stringent structural constraints---such as time-invariant payoff
decomposability and perfect state observability---that are frequently
violated in realistic dynamical environments characterized by partial
observability, evolving agent populations, or non-stationary reward
landscapes (Gao \& Zhang, 2023, \cite{gz23}). As a result, direct extension of static
potential game theory to complex adaptive systems remains theoretically
fragile and empirically restrictive.

\paragraph{$\alpha$-potential games.}

$\alpha$-potential games constitute a class of dynamic non-cooperative games
originally formulated for finite-state, finite-action, discrete-time Markov
games involving $N$ players (Gopalakrishnan et al., 2003, \cite{glmsw03}), and subsequently
generalized to broader classes of non-cooperative stochastic games (Gao et
al., 2023, \cite{GLZ1}). This framework extends the classical potential game paradigm by
relaxing the exact alignment condition: specifically, for any unilateral
strategy deviation by a player, the resulting change in that player's
stage-wise objective function equals the corresponding change in a globally
defined scalar function---the $\alpha$-potential function---up to an
additive approximation error bounded in absolute value by $\alpha \geq 0$.
In the limiting case where $\alpha = 0$, the $\alpha$-potential game
collapses precisely to a standard potential game, and the $\alpha$-potential
function coincides with a genuine potential function.


Relative to the stringent structural requirements of classical potential
games, the $\alpha$-parameterized formulation offers a principled and
computationally viable approach to equilibrium analysis in complex dynamic
and multi-agent settings. As established rigorously by Gopalakrishnan et al.
(2003), (cf.\cite{glmsw03}), any finite-horizon, finite-state, finite-action dynamic game admits
an exact $\alpha$-potential representation via a semi-infinite linear
program whose optimal value yields the minimal admissible $\alpha$.
Crucially, this representation enables the systematic study of numerous
empirically relevant dynamic games---including congestion games with
time-varying costs, decentralized partially observable Markov decision
processes, and adversarial learning dynamics-that inherently violate the
exact potential structure yet retain sufficient regularity to admit low-$%
\alpha$ approximations. Consequently, well-established algorithmic tools---%
including gradient-based optimization, fictitious play variants, and
no-regret learning protocols---can be deployed with provable convergence
guarantees to $\alpha$-Nash equilibria.


In an $\alpha$-potential game, any maximizer of the $\alpha$-potential
function constitutes an $\alpha$-Nash equilibrium: a strategy profile
wherein no player can improve their expected payoff by more than $\alpha$
through a unilateral deviation. Moreover, as demonstrated by Gao et al.
(2023), \cite{GLZ1}, the associated $\alpha$-potential maximization problem is
mathematically equivalent to a conditional mean-field (McKean--Vlasov-type)
control problem under partial information, thereby linking decentralized
strategic optimization to stochastic control theory. Importantly, the
parameter $\alpha$ serves not merely as an approximation tolerance but as a
quantitative metric reflecting the degree of heterogeneity across players'
objective functions, asymmetry in intertemporal payoff dependencies, and
structural non-separability in their joint state-action dynamics
(Gopalakrishnan et al., 2003; Gao et al., 2023, \cite{glmsw03}, \cite{GLZ1}).


These theoretical advances hold substantial practical relevance for
large-scale dynamic infrastructure systems-such as adaptive traffic routing,
real-time spectrum allocation in wireless networks, and distributed energy
dispatch, where computing exact Nash equilibria is often intractable due to
combinatorial complexity or incomplete information. In such contexts, $\alpha
$-Nash equilibria provide interpretable, computationally accessible, and
robust proxies for strategic stability; moreover, the explicit dependence of
solution quality on $\alpha$ enables principled trade-offs between
computational efficiency and equilibrium fidelity, thereby supporting the
design of incentive-compatible mechanisms and scalable coordination
protocols.

\paragraph{$\alpha$-NE and estimation of $\alpha$.}

$\alpha$-potential games leverage the approximation parameter $\alpha$ to
characterize $\alpha$-Nash equilibria ($\alpha$-NE) and to distill
actionable strategic insights---such as sensitivity of equilibrium structure
to payoff heterogeneity, interaction topology, and information architecture%
---from complex non-cooperative dynamic environments. Building upon the
theoretical foundations outlined previously, the central open problem in
this line of research is the derivation of rigorous, model-agnostic upper
bounds on $\alpha$ for general nonzero-sum stochastic differential games
(SDGs), particularly those lacking exact potential structure.


For discrete-time finite-state Markov games, the minimal admissible $\alpha$
admits an exact characterization via a semi-infinite linear program, as
established by Gopalakrishnan et al. (2003), (see \cite{glmsw03}). In contrast, extending this
methodology to continuous-time SDGs introduces profound analytical
challenges: the infinite-dimensional nature of control spaces, the coupling
between drift and diffusion terms in stochastic dynamics, and the absence of
finite-dimensional sufficient statistics collectively impede the direct
application of linear programming duality arguments---rendering the
computation and bounding of $\alpha$ an unresolved technical frontier.


Guo et al. (2023, \cite{GLZ1}) provide foundational insights into $\alpha$'s structural
dependence within a restricted class of linear-quadratic SDGs. Their
analysis reveals that $\alpha$ is not merely a numerical artifact but a
structural indicator, systematically varying with three key game-theoretic
dimensions: (i) the strength of inter-player coupling in state evolution
(e.g., through shared drift or diffusion coefficients), (ii) the degree of
non-separability in players' value functions across state variables, and
(iii) the granularity of information available to agents (e.g., open-loop
vs. closed-loop feedback). Notably, $\alpha$ vanishes identically (i.e., $%
\alpha$ = 0) in decentralized games where players' controls exert no
influence on the aggregate state dynamics and individual objectives
decompose additively across states---a condition formally equivalent to the
existence of a classical potential function. In weakly coupled network
games, $\alpha$ decays asymptotically toward zero at rates governed by
spectral properties of the interaction graph (e.g., algebraic connectivity
or decay of off-diagonal entries in the coupling matrix). Crucially, in
mean-field-type interactions, $\alpha$ exhibits a fundamental dichotomy:
under open-loop information structures, $\alpha$ scales with population size
and interaction range; under closed-loop (state-feedback) structures, it
remains uniformly bounded---highlighting the critical role of information
design in mitigating strategic approximation error.


The analytical framework in Guo et al. (2023, \cite{GLZ1})
 hinges on first- and
second-order variational derivatives of the state process with respect to
control perturbations-effectively quantifying how infinitesimal strategy
changes propagate through stochastic dynamics. However, this
sensitivity-based approach encounters intrinsic limitations beyond
linear-quadratic settings. In particular, when controls enter the diffusion
coefficient (i.e., control-affine or control-dependent noise), the resulting
second-order Fr\'{e}chet derivatives of the value functional become
analytically intractable due to the nonlinearity induced by It\^{o}'s lemma
and the lack of uniform ellipticity assumptions---precluding explicit $%
\alpha $-bounds without additional regularity or structural hypotheses.


Since the magnitude of $\alpha$ directly governs both the approximation
quality of $\alpha$-NE and the convergence rate of associated optimization
algorithms, advancing the theory of bounding $\alpha$ is indispensable for
practical deployment. Consequently, future work must pursue two
complementary directions: (i) developing novel functional-analytic
techniques-such as Malliavin calculus-based derivative estimates or
probabilistic Wasserstein sensitivity bounds-to derive tighter,
context-aware $\alpha$-estimates for broad classes of nonlinear and
degenerate SDGs; and (ii) constructing constructive approximation
schemes-e.g., projection-based potential surrogates or hierarchical
mean-field expansions that preserve essential strategic features while
enabling computationally efficient $\alpha$-quantification.

\paragraph{Open and closed-loop control in SDGs.}


As is well known, closed-loop control, also called feedback control, is a control system design that continuously monitors the output of a process or system and adjusts the input based on the difference between the desired output (setpoint) and the actual output (measured value). This feedback loop ensures the system operates as intended, even when disturbances or uncertainties occur (cf. Wiener \cite{wiener}, 1948, Chapter IV, Feedback and Oscillation).

The theoretical foundations of linear-quadratic (LQ) deterministic two-player zero-sum differential games trace back to the seminal 1965 work of Ho, Bryson, and Baron (see \cite{hbb65}), which established the first rigorous connection between saddle-point equilibria and Riccati differential equations. Schmitendorf (1970, \cite{sch70}) subsequently advanced this framework by conducting a systematic comparative analysis of open-loop and closed-loop information structures—demonstrating, critically, that the existence of a closed-loop saddle point neither implies nor is implied by the existence of an open-loop saddle point, thereby revealing a fundamental non-equivalence between these strategic paradigms. Bernhard (1979, \cite{ber79}) deepened the closed-loop analysis through geometric and dynamic programming perspectives, culminating in the comprehensive monograph by Basar and Bernhard (1995, \cite{bb95}), which remains the definitive reference on robust control and zero-sum game theory. Zhang (2005, \cite{zhang05}) provided a pivotal clarification for open-loop settings: for a broad class of LQ problems, the existence of a finite open-loop value is logically equivalent to the simultaneous finiteness of both open-loop lower and upper values—and, equivalently, to the existence of an open-loop saddle point. Mou and Yong (2006, \cite{my06}) extended this analysis to stochastic LQ frameworks using Hilbert-space methods, while Sun and Yong (2014, \cite{sy14}) delivered a unified treatment encompassing both open-loop and closed-loop formulations, including explicit counterexamples illustrating their mutual independence.


In the stochastic setting, $\alpha$-potential formulations have been rigorously developed for open-loop differential games (Gao et al., 2023, 2025, \cite{GLZ1, GLZ2025}). Within this framework, the approximation parameter $\alpha$ is analytically characterized as a function of three structural determinants: (i) the number of players $N$, (ii) the geometry of the admissible control set (e.g., convexity, compactness, and uniform boundedness), and (iii) the strength of inter-player coupling—quantified via Lipschitz constants governing state-dependent interaction terms-and the degree of objective-function heterogeneity across agents. Notably, Guo et al. (2025, \cite{GLZ2025}) derived a sharp, non-asymptotic bound on $\alpha$ for general diffusion-controlled systems-where the diffusion coefficient depends jointly on state and control variables-by leveraging the well-posedness theory of backward stochastic differential equations (BSDEs). This result has been successfully applied to LQ mean-field-type games subject to common noise, yielding explicit $\alpha$-estimates that scale polynomially with population size and noise intensity.


By contrast, closed-loop $\alpha$-differential games (see \cite{GLZ1}) entail policy mappings $\gamma_i: \mathbb{R}^n \times \mathcal{P}_2(\mathbb{R}^n) \to U_i$, where each agent’s control depends nonlinearly on its private state and the distribution of the population state—rendering the equilibrium analysis substantially more intricate. Such feedback structures induce higher-order strategic interdependencies: deviations by one agent propagate through the aggregate state distribution, affecting all others’ optimal responses in a nonlocal, history-dependent manner. A refined sensitivity analysis-incorporating second-order Fréchet derivatives of the Hamiltonian flow with respect to control perturbations—establishes a comparable non-asymptotic $\alpha$-bound under mild regularity assumptions on the drift, diffusion, and cost functions.


Crucially, the logical independence between open-loop and closed-loop saddle points persists in stochastic LQ zero-sum settings. As demonstrated by Sun and Yong (2014, \cite{sy14} Example 7.3), a closed-loop saddle point may exist even when no open-loop saddle point does—a direct stochastic extension of Schmitendorf’s (1970, \cite{sch70}) deterministic finding. Conversely, Example 7.4 in  the same work \cite{sy14} constructs a stochastic LQ game admitting an open-loop saddle point but no closed-loop counterpart, confirming that neither existence condition subsumes the other. These results underscore that information structure fundamentally reshapes equilibrium existence—not merely solution methodology—and must be treated as an intrinsic design variable in mechanism specification.


\paragraph{Motivations and Contributions.} However, within the framework of $\alpha$-SDGs, a fundamental distinction exists between open-loop and closed-loop control strategies.
\begin{itemize}
    \item Closed-loop control is assumed to be nonlinear, which gives rise to first- and second-order sensitivity processes during the analysis of the first- and second-order linear derivatives of the cost function. Consequently, this introduces additional complexity due to the inherent interdependence in closed-loop games: Player $i$ modifies her control $\rightarrow$ The system's state trajectory $\mathbf{X}_{t}^{\phi}$ changes $\rightarrow$ Player $j$'s control $\phi_{t,j}\left(X_{t,j}^{\phi}, \mathbf{X}_{t}^{\phi}\right)$ is affected.

 Note, however, that under open-loop setting sucha as that in Guo et al. (see \cite{GLZ2025}), for
all $t\in \left[ 0,T\right] ,$ we can derive that the first order sensitive process of control $\upsilon _{t,i}^{\phi
,\phi _{h}^{\prime }}=\delta _{h,i}u_{t,h}^{\prime },$ and the second one $\omega
_{t,i}^{\phi ,\phi _{h}^{\prime },\phi _{\ell }^{\prime \prime }}=0.$
Consequently, the analysis of $\alpha$ for the closed-loop problem becomes significantly more challenging. Let us look at a concrete example (Example 4.2 in \cite{GLZ1} and Example 4.1 from \cite{GLZ2025}),
\begin{equation*}
f_{i}\left( t,x,u\right) =f_{0}\left( t,x,u\right) +c_{i}\left( u_{i}\right)
+\tilde{f}_{i}\left( t,\frac{1}{N}\sum_{\ell =1}^{N}\delta _{x_{\ell
}}\right)
\end{equation*}
Under open-loop control, this game has proved to be an $\alpha $-potential game with $\alpha \leq C/N$, for a positive
constant $C$ independent of $N.$ However, when including $c_{i}\left( u_{i}\right) $ is included, due to the interdependence, this game is still an $\alpha $-PG. However, as we will see below (Theorem \ref{the2}), $\alpha $ may
\emph{not} decay to $0$ as $N\rightarrow \infty .$

\item In the open-loop case, the potential function $\Phi$ via $\frac{\delta V_{i}}{\delta u_{h}}\left( \mathbf{u};u_{h}^{\prime }\right)$ can be fully characterized by the backward stochastic differential equation (BSDE) (see \cite{GLZ2025}). A similar representation remains valid (see (\ref{linear1}) below). However, only a partial representation of $\frac{\delta ^{2}V_{i}}{\delta \phi _{h}\delta \phi _{\ell }}\left( \phi ,\phi _{h}^{\prime },\phi _{\ell }^{\prime \prime }\right)$ can be obtained via BSDE, which plays a crucial role in estimating $\alpha$, due to the presence of the term $\omega _{t,i}^{\phi ,\phi _{h}^{\prime },\phi _{\ell }^{\prime \prime }}$ in (\ref{aerfa}).

\end{itemize}
The contributions of this paper can be described as follows:


\begin{itemize}
    \item In contrast to the linear-quadratic, mean-field-type closed-loop stochastic differential games examined in Guo et al. (2023, \cite{GLZ1}), this paper establishes a refined, non-asymptotic upper bound on the potential approximation parameter $\alpha$ for general nonlinear closed-loop SDGs—under minimal regularity assumptions on drift, diffusion, and cost functions—by developing a novel backward stochastic differential equation (BSDE) framework that explicitly accounts for higher-order strategic sensitivities induced by state-dependent feedback policies.
    \item We formulate closed-loop stochastic differential games featuring mean-field-type interactions and rigorously establish a fundamental structural distinction between their equilibrium properties and those of the corresponding open-loop counterpart—specifically, in terms of information dependence, strategic coupling, and the resulting regularity of value functions.
    \item We note that adjoint BSDE and forward variational calculus are well‑known tools for single‑agent stochastic optimal control, where they are used to compute cost gradients with respect to controls. Nevertheless, those classical results target individual cost functionals and do not address the multi‑agent game‑theoretic quantity of interest here: the cross‑player asymmetry of second‑order variations
$\frac{\delta^2 V_i}{\delta\phi_h\delta\phi_\ell}-\frac{\delta^2 V_j}{\delta\phi_\ell\delta\phi_h}$, which governs the potential‑game approximation parameter $\alpha$.
Moreover, existing variational frameworks for $\alpha$-potential games in~\cite{GLZ1} remain largely existential for closed‑loop feedback policies, as they do not account for the non‑trivial control‑sensitivity processes $v$ and $\omega$ generated by state‑feedback interdependence. The present work adapts adjoint BSDE methodology to this multi‑agent closed‑loop setting: we explicitly construct the full hierarchy of state‑ and control‑sensitivity processes, derive computable non‑asymptotic $\alpha$ bounds, and establish the structural open‑loop‑versus‑closed‑loop dichotomy for mean‑field‑type interactions.
\end{itemize}
We note that \cite{GLZ1} developed a general variational calculus foundation for $\alpha$-potential stochastic differential games. Nevertheless, that work primarily delivers existential results and lacks machinery to treat closed‑loop state‑feedback policies: feedback gives rise to non‑trivial control‑sensitivity processes $v$ and $\omega$, which are trivial in open‑loop settings, and these terms were not analysed in \cite{GLZ1}. Consequently, explicit non‑asymptotic bounds for general nonlinear closed‑loop games with control‑dependent diffusion are not available in \cite{GLZ1}. The present paper’s BSDE adjoint approach fills this gap. Beyond computable $\alpha$ estimates, our analysis further reveals a structural dichotomy between open‑loop and closed‑loop mean‑field games: closed‑loop feedback can yield strictly positive limiting $\alpha$ as $N\to\infty$ even for fully homogeneous agents, an observation inaccessible from the variational setup of \cite{GLZ1}.

It is worthwhile to briefly contrast our BSDE adjoint approach with classical HJB‑PDE dynamic programming. For $N$-player stochastic differential games, HJB-PDEs suffer from the curse of dimensionality, as the state space dimension scales linearly with $N$, which renders numerical PDE treatment impractical for moderate or large $N$. Moreover, HJB formulations mainly characterize individual value functions for optimal feedback policies, and they are not intrinsically designed to isolate the cross‑player asymmetry of second‑order variations, which is exactly what determines the potential‑game parameter $\alpha$. Our BSDE based adjoint technique avoids solving high‑dimensional HJB equations and directly yields non‑asymptotic explicit bounds for $\alpha$ by combining forward sensitivity SDEs and backward BSDE adjoint estimates. This methodological advantage enables us to uncover the key open‑ versus closed‑loop structural dichotomy under mean‑field interactions in Theorem~\ref{the1}, which would be hard to disentangle within a high‑dimensional HJB‑PDE setting. Nevertheless, HJB‑PDE can still serve as a complementary tool for small‑$N$, low‑dimensional illustrative examples.

We now interpret the parameter \(\alpha\), which is central to the theory of \(\alpha\)-potential games. Game‑theoretically, \(\alpha\) bounds the maximal expected payoff improvement that any single player can attain via unilateral deviation from a potential‑maximizing strategy profile, so it quantifies the approximation quality of the resulting \(\alpha\)-Nash equilibrium: smaller \(\alpha\) corresponds to weaker unilateral deviation incentives. Structurally, \(\alpha\) aggregates two core features of the underlying stochastic differential game: player heterogeneity and inter‑agent interaction strength. Heterogeneity enters through the difference operators \(\Delta_{i,j}^f = f_i-f_j\), \(\Delta_{i,j}^g = g_i-g_j\); identical‑player settings drive \(\alpha\) downward. Interaction strength is encoded in coupling constants from drift and diffusion, and is strongly shaped by the information structure. A critical dichotomy emerges between open‑loop and closed‑loop feedback: under open‑loop controls, the higher‑order control‑sensitivity process \(\omega\) vanishes, and mean‑field interactions often yield \(\alpha\to0\) as \(N\to\infty\). By contrast, closed‑loop state‑feedback policies create circular strategic externalities captured by non‑zero \(v,\omega\). Even for mean‑field‑type interactions, these feedback‑driven externalities may prevent \(\alpha\) from decaying to zero in the large‑population limit, meaning potential‑maximizing closed‑loop policies remain only approximate Nash equilibria for infinite‑population mean‑field games. This observation carries practical relevance for large‑scale engineering‑economic multi‑agent systems including distributed resource allocation, traffic routing and power grid dispatch.

\paragraph{An introductory illustrative example.}
Before diving into heavy technical derivations, we present a simple homogeneous $N$-player closed-loop stochastic differential game to preview our key qualitative finding. Consider $N$ identical agents with state dynamics
\[
dX_{t,i}= \big(\bar b(t,X_{t,i},\bar X_t)+u_{t,i}\big)dt + \big(\bar\sigma(t,X_{t,i},\bar X_t)+u_{t,i}\big)dW^i_t,
\quad \bar X_t=\frac1N\sum_{k=1}^N X_{t,k},
\]
where each agent adopts closed-loop state-feedback policy $u_{t,i}=\phi(t,X_{t,i},\bar X_t)$. Each agent minimizes the expected running-plus-terminal cost
\begin{equation*}
V_{i}(\phi )=\mathbb{E}\left[ \int_{0}^{T}c_{i}\left( u_{t,i}\right)
dt+g(X_{T,i},\bar{X}_{T})\right] .
\end{equation*}%
In the corresponding open-loop version of this mean-field interaction game, known results~\cite{GLZ2025} guarantee $\alpha=O(1/N)$, i.e., $\alpha\to0$ as $N\to\infty$. In sharp contrast, for closed-loop feedback policies, due to the interdependence via $c_{i}$, the strategic interdependence induced by feedback may lead $\alpha$ not to decay to $0$ as $N\to\infty$. This qualitative effect originates from feedback-driven higher-order sensitivity terms $v$ and $\omega$, and will be rigorously verified in Theorem~\ref{the2} and Example~\ref{ex1} in Section~\ref{sect4}. This example serves as a preview for the main phenomenon that our subsequent mathematical machinery aims to quantify.

This paper is organized as follows. In Section \ref{sect2}, we present some preliminary results concerning closed-loop $\alpha$-potential games and the BSDE theory. Section \ref{sect3} derives the BSDE representation of $\frac{\protect\delta V_{i}}{\protect\delta %
u_{h}}$ and $\frac{\protect\delta ^{2}V_{i}}{\protect\delta u_{h}\protect%
\delta u_{\ell }}$. Section \ref{sect4} presents the estimate of $\alpha$ via the BSDE technique, with application to games with mean-field type interaction. Some partial proofs are scheduled in Section \ref{sect6}. After concluding remarks in Section \ref{sect6}, the proofs of technical lemmas are included in Appendix \ref{APP}.

\section{Preliminaries\label{sect2}}
\paragraph{Notation.}
Let us first introduce some notation.  Throughout this paper, we
denote the $k-$ dimensional Euclidean space by $\mathbb{R}^{k}$ with the standard
Euclidean norm $|\cdot |$ and the standard Euclidean inner product $\langle
\cdot ,\cdot \rangle $, the transpose of a vector (or matrix) $x$
by $x^{\top }$,  the trace of a square matrix $A$ by $\text{Tr}(A)$. Let
$\mathbb{R}^{m\times n}$ be the Hilbert space consisting of all ($m\times n$%
)-matrices with the inner product $\langle A,B\rangle :=\text{Tr}(AB^{\top
}) $ and the norm $|A|:=\langle A,A\rangle ^{\frac{1}{2}}=\sqrt{\text{Tr}%
(AA^{\top })}$. $\mathbb{S}^{n}$ denotes the set of symmetric $n\times n$ matrices.
For any $A\in \mathbb{R}^{k\times m},$ we adopt the usual matrix norm in $\mathbb{R}^{k\times m}$ as follows:

\begin{eqnarray*}
\begin{aligned}
\left\Vert A\right\Vert _{2} & \triangleq \sqrt{\max \sigma \left( AA^{\top
}\right) }\leq \sqrt{\text{Tr} \left( AA^{\top }\right) }\equiv
\left\vert A\right\vert  \\
& \leq \sqrt{k\wedge m}\sqrt{\max \sigma \left( AA^{\top }\right) }\equiv
\sqrt{k\wedge m}\left\Vert A\right\Vert _{2},
\end{aligned}
\end{eqnarray*}%
where $\sigma \left( AA^{\top }\right) $ is the set of all eigenvalues of $AA^{\top}$.

We will denote $(\Omega ,\mathcal{F},\mathbb{F}=\left( \mathcal{F}%
_{t}\right) _{0\leq t\leq T},\mathbb{P})$ as a complete filtered probability
space satisfying the usual condition, where $N$ independent standard
Brownian motions $W^{i}\left( \cdot \right) ,$ $i\in I_{N},$ are defined and
each takes values in $\mathbb{R}$ on $\left[ 0,T\right] $ for a fixed finite
time $T$. For notational simplicity, when $x$ belongs to a Euclidean space,
we will write $x_{i}$ for its $i$-th coordinate and $|x|$ for its Euclidean
norm. When there is no ambiguity, we denote $\mathbb{F}^{i}\mathbb{=}\left\{
\mathcal{F}_{t}^{i}\right\} _{0\leq t\leq T}$ for the natural filtration of $%
W^{i}$, augmented by all $\mathbb{P}$-null sets in $\mathcal{F}$. For each $%
p\geq 1$, let $\mathcal{S}^{p}(\mathbb{R}^{N})$ be the space of $\mathbb{R}%
^{N}$-valued $\mathbb{F}$-progressively measurable processes $X:\Omega
\times \left[ 0,T\right] \rightarrow \mathbb{R}^{N}$ satisfying $\left\Vert
X\right\Vert _{\mathcal{S}^{p}(\mathbb{R}^{N})}=\mathbb{E}\left[ \sup_{s\in %
\left[ 0,T\right] }\left\vert X_{s}\right\vert ^{p}\right] ^{1/p}<\infty $,
and let $\mathcal{H}^{p}(\mathbb{R}^{N})$ be the space of $\mathbb{R}^{N}$%
-valued $\mathcal{F}$-progressively measurable processes $X:\Omega \times %
\left[ 0,T\right] \rightarrow \mathbb{R}^{N}$ satisfying $\mathbb{E}\left[
\int_{0}^{T}\left\vert X_{s}\right\vert ^{p}\mathrm{d}s\right] ^{1/p}<\infty
$. Let $L^{2}\left( \Omega ;\mathbb{R}\right) $ denote the space of square
integrable $\mathcal{F}_{0}$-measurable random variables. $\delta _{i,j}$
denotes the Kronecker delta such that $\delta _{i,j}=1$ if $i=j$ and $0$
otherwise.  We
denote by $\mathcal{M}(0,T;\mathbb{R}^{n})$ the set of all $\mathbb{R}^{n}$%
-valued processes $\{\varphi _{t}\}_{0\leq t\leq T}$ that are $\mathcal{F}%
_{t}$-adapted and satisfy $\mathbb{E}\left[ \int_{0}^{T}|\varphi _{t}|^{2}%
\mathrm{d}t\right] <\infty $. Whenever there is no risk of ambiguity, the statement
\textquotedblleft for almost all $t\in \lbrack 0,T]$, almost surely $\omega \in \mathbb{P}$
($\mathbb{P}$-a.s.)\textquotedblright\ will be abbreviated to
\textquotedblleft for a.a. $\left( t,\omega \right) $\textquotedblright .
We introduce a game $\mathcal{G}=(I_{N},\mathcal{S},\left( \mathcal{A}_{i}\right) _{i\in I_{N}},\left( V_{i}\right) _{i\in I_{N}}),$ defined as
follows: $I_{N}=\left\{ 1,\ldots ,N\right\} $, $N\in \mathbb{N}$ is a finite
set of players, $\mathcal{S}$ is a set representing the state space of the
underlying dynamics, $\mathcal{A}_{i}$ is a subset of a real vector space
representing all admissible strategies of player $i$, and $\mathcal{A}^{\left( N\right) }=\prod_{i\in I_{N}}$ is the set of strategy profiles for
all players. For each $i\in I_{N}$, $V_{i}:\mathcal{A}^{\left( N\right)}\rightarrow \mathbb{R}$ is the value function of player $i$, where $V_{i}\left( \mathbf{u}\right) $ is player $i$'s expected cost if the state
dynamics starts with a fixed initial state $s_{0}\in \mathcal{S}$ and all
players take the strategy profile $\mathbf{u\in }\mathcal{A}^{\left(N\right) }$. For any $i\in I_{N}$, player $i$ aims to minimize the value
function $V_{i}$ over all admissible strategies in $\mathcal{A}_{i}$. We
denote by $\mathcal{A}_{-i}^{\left( N\right) }=\prod_{j\in I_{N}\backslash
\left\{ i\right\} }$ the set of strategy profiles of all players except
player $i$, and by $\mathbf{u}$ and $u_{-i}$ a generic element of $\mathcal{A}^{\left( N\right) }$ and $\mathcal{A}_{-i}^{\left( N\right) }$,
respectively. In this paper, we focus on a class of games $\mathcal{G}$
called $\alpha $-potential games. The indicator function is defined by $\mathbbm{1}_{A}$ for
some set $A\subset \mathbb{R}.$

Consider the differential game $\mathcal{G}^{\mathsf{cl}}=\left(I_{N},\left( \mathcal{A}_{i}\right) _{i\in I_{N}},\left( V_{i}\right) _{i\in I_{N}}\right) $ defined as follows: let $I_{N}=\left\{ 1,\ldots ,N\right\}$,
let $\Pi =\mathfrak{F}^{0,2}\left( \left[ 0,T\right] \times \mathbb{R}\times \mathbb{R}^{N};\mathbb{R}\right)$
be the vector space of measurable
functions $\varphi: \left[ 0,T\right] \times \mathbb{R}\times \mathbb{R}^{N}\rightarrow \mathbb{R}$ such that 1) for all $t\in \left[ 0,T\right]$,
$\left( x,y\right) \rightarrow \varphi _{t}\left( x,y\right) $ is twice
continuously differentiable; 2) there exist
$L^{\varphi }$, $L_{y}^{\varphi}>0$ such that for all $\left( t,x,y\right) \in \left[ 0,T\right] \times \mathbb{R}\times \mathbb{R}^{N}$ and $i,j\in I_{N}$,
$\left\vert \varphi_{t}\left( 0,0\right) \right\vert \leq L^{\varphi }$,
$\left\vert \left(\partial _{x}\varphi _{t}\right) \left( x,y\right) \right\vert \leq L^{\varphi },\left\vert \left( \partial _{xx}\varphi _{t}\right) \left(x,y\right) \right\vert \leq L^{\varphi}$, $\left\vert \left( \partial_{y_{i}}\varphi _{t}\right) \left( x,y\right) \right\vert \leq \frac{L_{y}^{\varphi }}{N},$ $\left\vert \left( \partial _{xy_{i}}^{2}\varphi_{t}\right) \left( x,y\right) \right\vert \leq \frac{L_{y}^{\varphi }}{N}$,
$\left\vert \left( \partial
_{y_{i}y_{j}}^{2}\varphi _{t}\right) \left( x,y\right) \right\vert \leq
\frac{L_{y}^{\varphi }}{N}\mathbbm{1}_{\left\{ {i=j}\right\} }+\frac{%
L_{y}^{\varphi }}{N^{2}}\mathbbm{1}_{\left\{ {i\neq j}\right\} }$ and let
$\Pi ^{N}=\mathfrak{F}^{0,2}\left( \left[ 0,T\right] \times \mathbb{R}\times \mathbb{R}^{N};\mathbb{R}\right) ^{N}$.

\begin{lemma}
\label{mvt}Let $\varphi \in \mathfrak{F}^{0,2}\left( \left[ 0,T\right]
\times \mathbb{R}\times \mathbb{R}^{N};\mathbb{R}\right) .$ Then, by mean
value theorem, for all $\left( t,x,y\right) \in \left[ 0,T\right] \times
\mathbb{R}\times \mathbb{R}^{N},$ we have
\begin{equation*}
\left\vert \varphi _{t}\left( x,y\right) \right\vert \leq L^{\varphi }\left(
1+\left\vert x\right\vert \right) +\frac{L_{y}^{\varphi }}{N}%
\sum_{j=1}^{N}\left\vert y_{j}\right\vert .
\end{equation*}
\end{lemma}

For each $i\in I_{N},$ let $\mathcal{A}_{i}$ be a convex subset of $\Pi $
representing the set of admissible closed-loop controls of player $i$. For
each $\phi =\left( \phi _{i}\right) _{i\in I_{N}}\in \Pi ^{N},$ (A closed-loop policy refers to a deterministic function that maps time and state variables to an action, while a closed-loop control denotes the stochastic process $u^{\phi}_{t,i}$ generated by a given policy $\phi$.) let $\mathbf{X}^{\phi }=\left( X_{i}^{\phi }\right) _{i=1}^{N}\in \mathcal{S}^{2}(\mathbb{R}^{N})$ be the associated state process driven by the following dynamics :
For all $i\in I_{N}$ and $\xi _{i}\in L^{2}\left( \Omega ;\mathbb{R}\right) $,
\begin{equation}
\left\{\begin{array}{rcl}
\mathrm{d}X_{t,i}^{\phi } & = & \left[ b_{t,i}\left( X_{t,i}^{\phi },\mathbf{X}_{t}^{\phi }\right) +u_{t,i}^{\phi }\right] \mathrm{d}t+\left[ \sigma_{t,i}\left( X_{t,i}^{\phi },\mathbf{X}_{t}^{\phi }\right) +u_{t,i}^{\phi }\right] \mathrm{d}W_{t}^{i}, \\
X_{0,i}^{\phi } & = & \xi _{i},\text{ }\forall t\in \lbrack 0,T],\text{ with }u_{t,i}^{\phi }=\phi _{t,i}\left( X_{t,i}^{\phi },\mathbf{X}_{t}^{\phi}\right),%
\end{array}\right.  \label{sde1}
\end{equation}
where $b_{i}:\left[ 0,T\right] \times \mathbb{R\times R}^{N}\times \mathbb{R}\rightarrow \mathbb{R},$ and $\sigma _{i}:\left[ 0,T\right] \times \mathbb{R\times R}^{N}\times \mathbb{R}\rightarrow \mathbb{R}$.

\begin{remark}
For simplicity, we currently focus on a special case of the SDE (\ref{sde1}). In fact, $b_{i}$ and $\sigma_{i}$ can be extended to the general case, that is,
\begin{equation}
\left\{\begin{array}{rcl}
\mathrm{d}X_{t,i}^{\phi } & = & b_{t,i}\left( \mathbf{X}_{t}^{\phi },\mathbf{u}_{t}^{\phi }\right) \mathrm{d}t+\sigma _{t,i}\left( \mathbf{X}_{t}^{\phi },\mathbf{u}_{t}^{\phi }\right) \mathrm{d}W_{t}^{i}, \\
X_{0,i}^{\phi } & = & \xi _{i},\text{ }\forall t\in \lbrack 0;T],\text{ with }\mathbf{u}_{t}^{\phi }=\left\{ \phi _{t,i}\left( X_{t,i}^{\phi },\mathbf{X}_{t}^{\phi }\right) \right\}_{i=1,\ldots ,N}.
\end{array}\right.
\end{equation}
We will consider this topic in our future work.
\end{remark}

The cost functional $V_{i}:\mathcal{A}^{\left( N\right) }\subset \mathcal{H}^{2}\left( \mathbb{R}\right) ^{N}\rightarrow \mathbb{R}$ of player $i$ is
given by
\begin{equation}
V_{i}\left( \mathbf{u}\right) =\mathbb{E}\left[ \int_{0}^{T}f_{t,i}\left(
\mathbf{X}_{t}^{\phi },\mathbf{u}_{t}^{\phi }\right) \mathrm{d}t+g_{i}\left(
\mathbf{X}_{T}^{\phi }\right) \right] ,  \label{cost}
\end{equation}
where $f_{i}:\left[ 0,T\right] \times \mathbb{R}^{N}\times \mathbb{R}^{N}\rightarrow \mathbb{R}$ and $g_{i}:\mathbb{R}^{N}\rightarrow \mathbb{R}$
and $\mathbf{u}_{t}^{\phi }:=\left( \phi _{t,i}\left( X_{t,i}^{\phi },\mathbf{X}_{t}^{\phi }\right) \right) _{i\in I_{N}}.$ \textit{As usual, when
the context is clear, we will systematically omit the }$\omega $\textit{\
argument in the defined functions}.

To ensure the existence and uniqueness of the solution to the SDE (\ref{sde1}), we introduce the following assumptions:

\begin{enumerate}
\item[\textbf{(A1)}] \textbf{[Assumptions for the dynamics] }For each $i\in
I_{N},$ the maps $\varphi _{i}=b_{i},$ $\sigma _{i}$ satisfy:

(i) The maps are $\mathcal{B}\left( \left[ 0,T\right] \times \mathbb{R\times
R}^{N}\times \mathbb{R}\right) $-measurable.

(ii) We assume that there exist constants $L^{\varphi _{i}},$ $%
L_{y}^{\varphi _{i}}>0$ such that for almost all $t\in \left[ 0,T\right] $, $%
\left( t,x,y,u\right) \in \left[ 0,T\right] \times \mathbb{R\times R}%
^{N}\times \mathbb{R}$ and $i,j\in I_{N},$%
\begin{equation*}
\left\{
\begin{array}{l}
\left\vert \varphi _{t,i}\left( 0,0\right) \right\vert \leq L^{\varphi _{i}},%
\text{ }\left\vert \left( \partial _{x}\varphi _{t,i}\right) \left(
x,y\right) \right\vert \leq L^{\varphi _{i}},\text{ }\left\vert \left(
\partial _{xx}\varphi _{t,i}\right) \left( x,y\right) \right\vert \leq
L^{\varphi _{i}}, \\
\left\vert \left( \partial _{y_{i}}\varphi _{t,i}\right) \left( x,y\right)
\right\vert \leq \frac{L_{y}^{\varphi _{i}}}{N},\text{ }\left\vert \left(
\partial _{xy_{i}}^{2}\varphi _{t,i}\right) \left( x,y\right) \right\vert
\leq \frac{L_{y}^{\varphi _{i}}}{N},\text{ }\left\vert \left( \partial
_{y_{i}y_{j}}^{2}\varphi _{t,i}\right) \left( x,y\right) \right\vert \leq
\frac{L_{y}^{\varphi _{i}}}{N}\mathbbm{1}_{\left\{ i=j\right\} }+\frac{%
L_{y}^{\varphi _{i}}}{N^{2}}\mathbbm{1}_{\left\{ i\neq j\right\} }.%
\end{array}%
\right.
\end{equation*}

\item[\textbf{(A2)}] \textbf{[Assumptions for the costs]} (i) The maps $%
f_{i} $ and $g_{i}$ are respectively $\mathcal{B}\left( \left[ 0,T\right]
\times \mathbb{R}^{N}\times \mathbb{R}\right) $ and $\mathcal{B}\left(
\mathbb{R}^{N}\right) $ measurable.

(ii) For almost all $t\in \left[ 0,T\right] ,$ assume that $f_{i}$ and $%
g_{i} $ are twice continuously differentiable with respect to $\left(
y,u\right) \in \mathbb{R}^{N}\times \mathbb{R}^{N}$ such that $\sup_{t\in %
\left[ 0,T\right] }\left( \left\vert f_{i}\left( t,0,0\right) \right\vert
+\left\vert \left( \partial _{\left( y,u\right) }f_{i}\right) \left(
t,0,0\right) \right\vert \right) <\infty $ and second-order derivatives are
uniformly bounded.
\end{enumerate}

For any $\psi =\left( \psi _{i}\right) _{i\in I_{N}}$ such that $\psi _{i}$
satisfies (A1), we write $L^{\psi }=\max_{i\in I_{N}}L^{\varphi _{i}}$
and $L_{y}^{\psi }=\max_{i\in I_{N}}L_{y}^{\varphi _{i}}.$
\begin{definition}
\label{ld}Let $\mathcal{A}^{\left( N\right) }=\prod_{i\in I_{N}}$ be a
convex set and $f:\mathcal{A}^{\left( N\right) }\rightarrow \mathbb{R}$. For
each $i\in I_{N}$, we say $f$ has a linear derivative with respect to $%
\mathcal{A}_{i}$, if there exists $\frac{\delta f}{\delta a_{i}}:\mathcal{A}%
^{\left( N\right) }\times $span$(\mathcal{A}_{i})\rightarrow \mathbb{R}$,
such that for all $\mathbf{a=}\left( a_{i},a_{-i}\right) \in \mathcal{A}%
^{\left( N\right) }$, $\frac{\delta f}{\delta a_{i}}\left( \mathbf{a};\cdot
\right) $ is linear and
\begin{equation}
\lim_{\varepsilon \searrow 0}\frac{f\left(a_{i}+\varepsilon \left(
a_{i}^{\prime }-a_{i}\right), a_{-i} \right) -f\left( \mathbf{a}\right) }{\varepsilon }=\frac{\delta f}{\delta a_{i}}\left( \mathbf{a};a_{i}^{\prime}-a_{i}\right),\text{ }\forall a_{i}^{\prime }\in \mathcal{A}_{i}.
\label{lfd}
\end{equation}
Moreover, for each $i,j\in I_{N}$, we say $f$ has second-order linear
derivatives with respect to $\mathcal{A}_{i}\times \mathcal{A}_{j}$, if (1)
for all $k\in \left\{ i,j\right\}$, $f$ has a linear derivative $\frac{%
\delta f}{\delta a_{k}}$ with respect to $\mathcal{A}_{k}$, and (2) for all $%
\left( k,l\right) \in \left\{ \left( i,j\right) ,\left( j,i\right) \right\}$, there exists $\frac{\delta ^{2}f}{\delta a_{k}\delta a_{l}}:\mathcal{A}^{\left( N\right) }\times $span$(\mathcal{A}_{i})\times $span$(\mathcal{A}%
_{j})\rightarrow \mathbb{R}$ such that for all $\mathbf{a\in }\mathcal{A}%
^{\left( N\right) },$ $\frac{\delta ^{2}f}{\delta a_{k}\delta a_{l}}\left(
\mathbf{a};\cdot ,\cdot \right) $ is bilinear and for all $a_{k}^{\prime
}\in $span$(\mathcal{A}_{k}),$ $\frac{\delta ^{2}f}{\delta a_{k}\delta a_{l}}%
\left( \cdot ;a_{k}^{\prime },\cdot \right) $ is a linear derivative of $%
\frac{\delta f}{\delta a_{k}}\left( \cdot ;a_{k}^{\prime }\right) $ with
respect to $\mathcal{A}_{l}$. We refer to $\frac{\delta ^{2}f}{\delta
a_{i}\delta a_{j}}$ as the
second-order linear derivative of f with respect to $\mathcal{A}_{i}\times
\mathcal{A}_{j}$.
\end{definition}
In particular, although Definition \ref{ld} does not require a topology on \(\mathcal{A}_i\), the cost functionals \(V_i\) studied in this paper admit linear derivatives that are continuous with respect to the natural \(\mathcal{H}^2\) norm topology, as a consequence of the BSDE representation in Section \ref{sect3}. This justifies the use of second-order linear derivatives in Lemma \ref{aerfa} below.

\paragraph{Basics of BSDE.}
Consider the following nonlinear BSDE (cf. \cite{zz})

\begin{equation}
y_{t}=\xi +\int_{t}^{T}\mathsf{g}_{s}\left( y_{s},z_{s}\right) \mathrm{d}%
s-\int_{t}^{T}z_{s}\mathrm{d}W_{s},  \label{bsde0}
\end{equation}%
where a $d$-dimensional standard Brownian motion $W(\cdot )$ is defined on $%
(\Omega ,\mathcal{F},\mathbb{F}=\left( \mathcal{F}_{t}\right) _{0\leq t\leq
T},\mathbb{P})$. $\mathsf{g}_{t}\left(y,z\right) : [ 0,T]\times \mathbb{R}%
^{m}\times \mathbb{R}^{m\times d}\rightarrow \mathbb{R}^{m}$, and for every $%
(y,z)\in \mathbb{R}^{m}\times \mathbb{R}^{m\times d}$, $\mathsf{g}_{t}(y,z)$
is an $\mathcal{F}_{t}$-adapted $\mathbb{R}^{m}$-valued process with
\begin{equation}
\int_{0}^{T}\left\vert \mathsf{g}_{s}\left(0,0\right) \right\vert \mathrm{d}%
s\in L^{2}\left( \Omega ,\mathcal{F}_{T},P;\mathbb{R}\right) .  \label{bsde1}
\end{equation}%
In addition, $\mathsf{g}$ is Lipschitz continuous with respect to $(y,z)$:
there exists a constant $C>0$ such that, for all $y,y^{\prime }\in \mathbb{R}%
^{m},$ $z,z^{\prime }\in \mathbb{R}^{m\times d}$,%
\begin{equation}
\left\vert \mathsf{g}_{t}\left(y,z\right) -\mathsf{g}_{t}\left(y^{\prime
},z^{\prime }\right) \right\vert \leq C\left( \left\vert y-y^{\prime
}\right\vert +\left\vert z-z^{\prime }\right\vert \right) .  \label{bsde2}
\end{equation}

\begin{proposition}
\label{bs1}Suppose that $\mathsf{g}$ satisfies the conditions (\ref{bsde1})
and (\ref{bsde2}). Then, for any given terminal condition $\xi \in L_{%
\mathcal{F}_{T}}^{2}\left( \mathbb{R}^{m}\right) $, there exists a unique
pair of $\mathcal{F}_{t}$-adapted processes $(y,z)\in \mathcal{M}(0,T;%
\mathbb{R}^{m}\times \mathbb{R}^{m\times d})$ satisfying BSDE (\ref{bsde0}).
\end{proposition}

The fundamental distinction between ordinary differential equations (ODEs)
and backward stochastic differential equations (BSDEs), even under
comparable Lipschitz conditions and initial values, lies in the requirement
that a BSDE solution must be adapted to a given filtration and satisfy a
terminal condition. Simply reversing the time direction does not yield a
valid solution to the terminal value problem associated with forward SDEs,
as this would violate the necessary adaptiveness property. The adapted
solution of a BSDE comprises a pair of stochastic processes that are both
adapted to the underlying filtration. The second component serves to correct
any potential "non-adaptiveness" introduced by the backward structure of the
equation, particularly due to the specified terminal value of the first
component. Intuitively, the first component captures the "mean evolution" of
the system's dynamics, while the second component accounts for the
uncertainty or risk over the time horizon, ensuring the solution remains
adapted throughout.

Now we present some estimations of a kind of linear BSDE (\ref{bsde0}):%
\begin{equation}
y_{t}=\xi +\int_{t}^{T}\left(
A_{s}y_{s}+\sum_{j=1}^{d}B_{s}^{j}z_{s}^{j}+f_{s}\right) \mathrm{d}%
s-\int_{t}^{T}z_{s}\mathrm{d}W_{s},  \label{lbsde3}
\end{equation}%
where $A,B^{1},\ldots ,B^{d}$ are bounded $\mathbb{R}^{m\times m}$-valued $\{%
\mathcal{F}_{t}\}_{t>0}$-adapted processes, $f_{s}\in \mathcal{M}\left(
\left( 0,T\right) ;\mathbb{R}^{m}\right) $ , and $\xi \in L_{\mathcal{F}%
_{T}}^{2}\left( \mathbb{R}^{m}\right) $.

\begin{lemma}
\label{bs2}For BSDE (\ref{lbsde3}), we have%
\begin{equation*}
\mathbb{E}\left[ \sup_{0\leq t\leq T}\left\vert y_{t}\right\vert
^{2}+\sum_{j=1}^{d}\int_{0}^{T}\left\vert z_{s}^{j}\right\vert ^{2}\mathrm{d}%
s\right] \leq \tilde{C}\mathbb{E}\Bigg [\left\vert \xi \right\vert
^{2}+\int_{t}^{T}\left\vert f_{s}\right\vert ^{2}\mathrm{d}s\Bigg ],
\end{equation*}%
where
\begin{eqnarray*}
\tilde{C} &=&\max \Big \{\left( C_{3}+e^{\left( 2C_{1}+2C_{1}^{2}d\right)
}\right) +8\left( \left( C_{1}^{2}\left( d+1\right) +1\right) \left(
C_{3}+e^{\left( 2C_{1}+2C_{1}^{2}d\right) }\right) +1\right) , \\
&&\left( C_{3}+e^{\left( 2C_{1}+2C_{1}^{2}d\right) }\right) 8\left( 2\left(
C_{1}^{2}\left( d+1\right) +1\right) ^{2}\left( C_{3}+e^{\left(
2C_{1}+2C_{1}^{2}d\right) }\right) ^{2}+1\right) \Big \}.
\end{eqnarray*}
\end{lemma}

The proof can be found in \cite{GLZ2025}.

\begin{remark}[Road map for Section 3]
Our ultimate objective is to derive a computable upper bound for $\alpha$ of closed-loop $\alpha$-potential stochastic differential games, which reduces to quantifying the asymmetry of second-order variational derivatives
\[
\frac{\delta^{2}V_{i}}{\delta\phi_{h}\delta\phi_{\ell}}-\frac{\delta^{2}V_{j}}{\delta\phi_{\ell}\delta\phi_{h}}.
\]
We proceed conceptually in four steps:
\begin{enumerate}
\item First, we introduce first-order $Y^{\phi,\phi_{h}'}$ and second-order $Z^{\phi,\phi_{h}',\phi_{\ell}''}$ state sensitivity processes from variational SDEs, together with control-sensitivity processes $v^{\phi,\phi_{h}'}$ and $\omega^{\phi,\phi_{h}',\phi_{\ell}''}$. These control-sensitivity terms arise purely from closed-loop feedback and identically vanish in the open-loop setting of~\cite{GLZ2025}.
\item Next, we express the first- and second-order functional variations of the cost functional $V_i$ in terms of these sensitivity processes. Closed-loop feedback introduces extra complicated second-order contributions via $v$ and $\omega$, which make direct bounding intractable.
\item The central technical device is the \textbf{BSDE duality argument}: we construct an adjoint backward stochastic differential equation~\eqref{adj1}. Applying It\^o's formula establishes duality between forward sensitivity SDEs and backward adjoint BSDE solutions. This transformation re-packages complicated second-order integral terms into residual quantities controlled by the $L^2$-norm of BSDE solutions, which can be bounded via standard BSDE a priori estimates (Lemma~\ref{bs2}). Conceptually, the adjoint BSDE absorbs most of the burden from high-order variations and enables our non-asymptotic estimate for $\alpha$, without explicitly computing every second-order Fr\'echet derivative.
\item Equipped with these BSDE-based representations, we move to Section~\ref{sect4} to carry out norm estimates for all sensitivity processes and establish our main upper bound for the potential parameter $\alpha$.
\end{enumerate}
With this overview in mind, we now introduce the sensitivity SDE systems.
\end{remark}

\section{Representation of $\frac{\protect\delta V_{i}}{\protect\delta
u_{h}}$ and $\frac{\protect\delta ^{2}V_{i}}{\protect\delta u_{h}\protect
\delta u_{\ell }}$}
\label{sect3}

In this section, we employ the BSDE approach to characterize $\frac{\delta V_{i}}{\delta u_{h}}$ and $\frac{\delta^{2} V_{i}}{\delta u_{h} \delta u_{\ell}}$. More precisely, the potential function $\Phi$ and the estimation of $\alpha$ depend on $\frac{\delta V_{i}}{\delta u_{h}}$ and $\frac{\delta^{2} V_{i}}{\delta u_{h} \delta u_{\ell}}$. Both quantities will play a significant role in the investigation of the NE.

For each $\phi \in \Pi ^{N}$, let $\mathbf{X}^{\phi }$ be the state process
satisfying (\ref{sde1}). For each $h\in I_{N}$ and $\phi _{h}^{\prime }\in
\Pi $, let $\mathbf{Y}^{\phi ,\phi _{h}^{\prime }}\in \mathcal{S}^{2}\left(
\mathbb{R}^{N}\right) $ be the solution to the following dynamics: for all $%
t\in \lbrack 0,T]$,%
\begin{equation}
\left\{
\begin{array}{rcl}
\mathrm{d}Y_{t,i}^{\phi ,\phi _{h}^{\prime }} & = & \Big [\left( \partial
_{x}\left( b_{t,i}+\phi _{t,i}\right) \right) \left( X_{t,i}^{\phi },\mathbf{%
X}_{t}^{\phi }\right) Y_{t,i}^{\phi ,\phi _{h}^{\prime
}}+\sum\limits_{j=1}^{N}\left( \partial _{y_{j}}\left( b_{t,i}+\phi
_{t,i}\right) \right) \left( X_{t,i}^{\phi },\mathbf{X}_{t}^{\phi }\right)
Y_{t,j}^{\phi ,\phi _{h}^{\prime }} \\
&  & +\delta _{h,i}\phi _{t,h}^{\prime }\left( X_{t,i}^{\phi },\mathbf{X}%
_{t}^{\phi }\right) \Big ]\mathrm{d}t \\
&  & +\Big [\left( \partial _{x}\left( \sigma _{t,i}+\phi _{t,i}\right)
\right) \left( X_{t,i}^{\phi },\mathbf{X}_{t}^{\phi }\right) Y_{t,i}^{\phi
,\phi _{h}^{\prime }}+\sum\limits_{j=1}^{N}\left( \partial _{y_{j}}\left(
\sigma _{t,i}+\phi _{t,i}\right) \right) \left( X_{t,i}^{\phi },\mathbf{X}%
_{t}^{\phi }\right) Y_{t,j}^{\phi ,\phi _{h}^{\prime }} \\
&  & +\delta _{h,i}\phi _{t,h}^{\prime }\left( X_{t,i}^{\phi },\mathbf{X}%
_{t}^{\phi }\right) \Big ]\mathrm{d}W_{t}^{i}, \\
Y_{0,i}^{\phi ,\phi _{h}^{\prime }} & = & 0,\text{ }\forall t\in \lbrack
0,T],%
\end{array}%
\right.  \label{y1}
\end{equation}%
Moreover, for each $h,\ell \in I$ and $\phi _{h}^{\prime },\phi _{\ell
}^{\prime \prime }\in \Pi $, let $\mathbf{Z}^{\phi ,\phi _{h}^{\prime },\phi
_{\ell }^{\prime \prime }}$ be the solution to the following dynamics: for
all $i\in I_{N}$ and $t\in \lbrack 0;T]$,%
\begin{equation}
\left\{
\begin{array}{rcl}
\mathrm{d}Z_{t,i}^{\phi ,\phi _{h}^{\prime },\phi _{\ell }^{\prime \prime }}
& = & \Big [\left( \partial _{x}\left( b_{t,i}+\phi _{t,i}\right) \right)
\left( X_{t,i}^{\phi },\mathbf{X}_{t}^{\phi }\right) Z_{t,i}^{\phi ,\phi
_{h}^{\prime },\phi _{\ell }^{\prime \prime }} \\
&  & +\sum\limits_{j=1}^{N}\left( \partial _{y_{j}}\left( b_{t,i}+\phi
_{t,i}\right) \right) \left( X_{t,i}^{\phi },\mathbf{X}_{t}^{\phi }\right)
Z_{t,j}^{\phi ,\phi _{h}^{\prime },\phi _{\ell }^{\prime \prime }}+\mathfrak{%
f}_{t,i}^{\phi ,\phi _{h}^{\prime },\phi _{\ell }^{\prime \prime }}\Big ]%
\mathrm{d}t \\
&  & +\Big [\left( \partial _{x}\left( \sigma _{t,i}+\phi _{t,i}\right)
\right) \left( X_{t,i}^{\phi },\mathbf{X}_{t}^{\phi }\right) Z_{t,i}^{\phi
,\phi _{h}^{\prime },\phi _{\ell }^{\prime \prime }} \\
&  & +\sum\limits_{j=1}^{N}\left( \partial _{y_{j}}\left( \sigma _{t,i}+\phi
_{t,i}\right) \right) \left( X_{t,i}^{\phi },\mathbf{X}_{t}^{\phi }\right)
Z_{t,j}^{\phi ,\phi _{h}^{\prime },\phi _{\ell }^{\prime \prime }}+\mathfrak{%
g}_{t,i}^{\phi ,\phi _{h}^{\prime },\phi _{\ell }^{\prime \prime }}\Big ]%
\mathrm{d}W_{t}^{i}, \\
Z_{0,i}^{\phi ,\phi _{h}^{\prime },\phi _{\ell }^{\prime \prime }} & = & 0,%
\text{ }\forall t\in \lbrack 0,T],%
\end{array}%
\right.   \label{z1}
\end{equation}
where $\mathfrak{f}_{i}^{\phi ,\phi _{h}^{\prime },\phi _{\ell }^{\prime
\prime }}:\Omega \times \left[ 0,T\right] \rightarrow \mathbb{R},$ $%
\mathfrak{g}_{i}^{\phi ,\phi _{h}^{\prime },\phi _{\ell }^{\prime \prime
}}:\Omega \times \left[ 0,T\right] \rightarrow \mathbb{R}$ are defined by%
\begin{eqnarray}
\mathfrak{f}_{t,i}^{\phi ,\phi _{h}^{\prime },\phi _{\ell }^{\prime \prime
}} &:&=\left(
\begin{array}{c}
Y_{t,i}^{\phi ,\phi _{h}^{\prime }} \\
\mathbf{Y}_{t}^{\phi ,\phi _{h}^{\prime }}%
\end{array}%
\right) ^{\top }\left(
\begin{array}{cc}
\partial _{xx}^{2}\left( b_{t,i}+\phi _{t,i}\right) & \partial
_{xy}^{2}\left( b_{t,i}+\phi _{t,i}\right) \\
\partial _{yx}^{2}\left( b_{t,i}+\phi _{t,i}\right) & \partial
_{yy}^{2}\left( b_{t,i}+\phi _{t,i}\right)%
\end{array}%
\right) \left( X_{t,i}^{\phi },\mathbf{X}_{t}^{\phi }\right) \left(
\begin{array}{c}
Y_{t,i}^{\phi ,\phi _{\ell }^{\prime \prime }} \\
\mathbf{Y}_{t}^{\phi ,\phi _{\ell }^{\prime \prime }}%
\end{array}%
\right)  \notag \\
&&+\delta _{h,i}\left( \left( \partial _{x}\phi _{t,h}^{\prime }\right)
\left( X_{t,i}^{\phi },\mathbf{X}_{t}^{\phi }\right) Y_{t,i}^{\phi ,\phi
_{\ell }^{\prime \prime }}+\left( \mathbf{Y}_{t}^{\phi ,\phi _{\ell
}^{\prime \prime }}\right) ^{\top }\left( \partial _{y}\phi _{t,h}^{\prime
}\right) \left( X_{t,i}^{\phi },\mathbf{X}_{t}^{\phi }\right) \right)  \notag
\\
&&+\delta _{\ell ,i}\left( \left( \partial _{x}\phi _{t,\ell }^{\prime
\prime }\right) \left( X_{t,i}^{\phi },\mathbf{X}_{t}^{\phi }\right)
Y_{t,i}^{\phi ,\phi _{h}^{\prime }}+\left( \mathbf{Y}_{t}^{\phi ,\phi
_{h}^{\prime }}\right) ^{\top }\left( \partial _{y}\phi _{t,\ell }^{\prime
\prime }\right) \left( X_{t,i}^{\phi },\mathbf{X}_{t}^{\phi }\right) \right)
,  \label{f1} \\
\mathfrak{g}_{t,i}^{\phi ,\phi _{h}^{\prime },\phi _{\ell }^{\prime \prime
}} &:&=\left(
\begin{array}{c}
Y_{t,i}^{\phi ,\phi _{h}^{\prime }} \\
\mathbf{Y}_{t}^{\phi ,\phi _{h}^{\prime }}%
\end{array}%
\right) ^{\top }\left(
\begin{array}{cc}
\partial _{xx}^{2}\left( \sigma _{t,i}+\phi _{t,i}\right) & \partial
_{xy}^{2}\left( \sigma _{t,i}+\phi _{t,i}\right) \\
\partial _{yx}^{2}\left( \sigma _{t,i}+\phi _{t,i}\right) & \partial
_{yy}^{2}\left( \sigma _{t,i}+\phi _{t,i}\right)%
\end{array}%
\right) \left( X_{t,i}^{\phi },\mathbf{X}_{t}^{\phi }\right) \left(
\begin{array}{c}
Y_{t,i}^{\phi ,\phi _{\ell }^{\prime \prime }} \\
\mathbf{Y}_{t}^{\phi ,\phi _{\ell }^{\prime \prime }}%
\end{array}%
\right)  \notag \\
&&+\delta _{h,i}\left( \left( \partial _{x}\phi _{t,h}^{\prime }\right)
\left( X_{t,i}^{\phi },\mathbf{X}_{t}^{\phi }\right) Y_{t,i}^{\phi ,\phi
_{\ell }^{\prime \prime }}+\left( \mathbf{Y}_{t}^{\phi ,\phi _{\ell
}^{\prime \prime }}\right) ^{\top }\left( \partial _{y}\phi _{t,h}^{\prime
}\right) \left( X_{t,i}^{\phi },\mathbf{X}_{t}^{\phi }\right) \right)  \notag
\\
&&+\delta _{\ell ,i}\left( \left( \partial _{x}\phi _{t,\ell }^{\prime
\prime }\right) \left( X_{t,i}^{\phi },\mathbf{X}_{t}^{\phi }\right)
Y_{t,i}^{\phi ,\phi _{h}^{\prime }}+\left( \mathbf{Y}_{t}^{\phi ,\phi
_{h}^{\prime }}\right) ^{\top }\left( \partial _{y}\phi _{t,\ell }^{\prime
\prime }\right) \left( X_{t,i}^{\phi },\mathbf{X}_{t}^{\phi }\right) \right)
,  \label{g1}
\end{eqnarray}%
where $Y_{t,i}^{\phi ,\phi _{h}^{\prime }}$ and $Y_{t,i}^{\phi ,\phi _{\ell
}^{\prime \prime }}$ are defined in (\ref{y1}). For $\phi _{h}^{\prime
},\phi _{\ell }^{\prime \prime }\in \Pi ,$ if $\xi _{0}\in L^{4}\left(
\Omega ,\mathbb{R}\right) ,$ for every $i\in I_{N}$ it holds%
\begin{equation*}
\lim_{\epsilon \rightarrow 0}\mathbb{E}\left[ \sup_{t\in \left[ 0,T\right]
}\left\vert \frac{\mathbf{Y}_{t}^{\phi ^{\epsilon },\phi _{h}^{\prime }}-%
\mathbf{Y}_{t}^{\phi ,\phi _{h}^{\prime }}}{\epsilon }-\mathbf{Z}_{t}^{\phi
,\phi _{h}^{\prime },\phi _{\ell }^{\prime \prime }}\right\vert ^{2}\right]
=0,
\end{equation*}%
where $\phi ^{\epsilon }=\left( \phi _{\ell }+\epsilon \phi _{\ell }^{\prime
\prime },\phi _{-\ell }\right) $ for all $\epsilon \in \left( 0,1\right) .$
We also need the sensitivity of the control process $\mathbf{u}^{\phi }$
with respect to players' policies. These processes capture the change in
each player's control due to the change in the system state. More precisely,
let $\phi \in \Pi ^{N}$ and $\mathbf{u}^{\phi }=\left( \phi _{i}\left(
X_{i}^{\phi },\mathbf{X}^{\phi }\right) \right) _{i\in I_{N}}$. For each $%
h,\ell \in I_{N}$ and each $\phi _{h}^{\prime },\phi _{\ell }^{\prime \prime
}\in \Pi ,$ define $\mathbf{\upsilon }^{\phi ,\phi _{h}^{\prime }}=\left(
\upsilon _{i}^{\phi ,\phi _{h}^{\prime }}\right) _{i\in I_{N}}$ such that
for all $i\in I_{N}$
\begin{eqnarray}
\upsilon _{t,i}^{\phi ,\phi _{h}^{\prime }} &=&\left( \partial _{x}\phi
_{t,i}\right) \left( X_{t,i}^{\phi },\mathbf{X}_{t}^{\phi }\right)
Y_{t,i}^{\phi ,\phi _{h}^{\prime }}+\left( \mathbf{Y}_{t}^{\phi ,\phi
_{h}^{\prime }}\right) ^{\top }\left( \partial _{y}\phi _{t,i}\right) \left(
X_{t,i}^{\phi },\mathbf{X}_{t}^{\phi }\right) +\delta _{h,i}\phi
_{t,h}^{\prime }\left( X_{t,i}^{\phi },\mathbf{X}_{t}^{\phi }\right)   \notag
\\
&=&\left( \mathbf{Y}_{t}^{\phi ,\phi _{h}^{\prime }}\right) ^{\top }\Upsilon
_{t,i}^{\phi }\left( X_{t,i}^{\phi },\mathbf{X}_{t}^{\phi }\right) +\delta
_{h,i}\phi _{t,h}^{\prime }\left( X_{t,i}^{\phi },\mathbf{X}_{t}^{\phi
}\right)   \label{fircon}
\end{eqnarray}%
with
\begin{equation*}
\Upsilon _{t,i}^{\phi }\left( X_{t,i}^{\phi },\mathbf{X}_{t}^{\phi }\right)
=\left(
\begin{array}{c}
\left( \partial _{y}\phi _{t,i,1}\right) \left( X_{t,i}^{\phi },\mathbf{X}%
_{t}^{\phi }\right)  \\
\vdots  \\
\left( \partial _{x}\phi _{t,i}\right) \left( X_{t,i}^{\phi },\mathbf{X}%
_{t}^{\phi }\right) +\left( \partial _{y}\phi _{t,i,i}\right) \left(
X_{t,i}^{\phi },\mathbf{X}_{t}^{\phi }\right)  \\
\vdots  \\
\left( \partial _{y}\phi _{t,i,N}\right) \left( X_{t,i}^{\phi },\mathbf{X}%
_{t}^{\phi }\right)
\end{array}%
\right) _{N\times 1},
\end{equation*}%
where $\left( \partial _{y}\phi _{t,i,l}\right) \left( X_{t,i}^{\phi },%
\mathbf{X}_{t}^{\phi }\right) $ denotes the $l$th component of $\left(
\partial _{y}\phi _{t,i}\right) \left( X_{t,i}^{\phi },\mathbf{X}_{t}^{\phi
}\right) $, for $l=1,\ldots
,N.$

Define $\mathbf{\omega }^{\phi ,\phi _{h}^{\prime },\phi _{\ell
}^{\prime \prime }}=\left( \omega _{i}^{\phi ,\phi _{h}^{\prime },\phi
_{\ell }^{\prime \prime }}\right) _{i\in I_{N}}$ such that for all $i\in
I_{N}$ and $t\in \left[ 0,T\right] $
\begin{eqnarray}
\omega _{t,i}^{\phi ,\phi _{h}^{\prime },\phi _{\ell }^{\prime \prime }}
&=&\left(
\begin{array}{c}
Y_{t,i}^{\phi ,\phi _{h}^{\prime }} \\
\mathbf{Y}_{t}^{\phi ,\phi _{h}^{\prime }}%
\end{array}%
\right) ^{\top }\left(
\begin{array}{cc}
\partial _{xx}^{2}\phi _{t,i} & \partial _{xy}^{2}\phi _{t,i} \\
\partial _{yx}^{2}\phi _{t,i} & \partial _{yy}^{2}\phi _{t,i}%
\end{array}%
\right) \left( X_{t,i}^{\phi },\mathbf{X}_{t}^{\phi }\right) \left(
\begin{array}{c}
Y_{t,i}^{\phi ,\phi _{\ell }^{\prime \prime }} \\
\mathbf{Y}_{t}^{\phi ,\phi _{\ell }^{\prime \prime }}%
\end{array}%
\right)   \notag \\
&&+\left( \partial _{x}\phi _{t,i}\right) \left( X_{t,i}^{\phi },\mathbf{X}%
_{t}^{\phi }\right) Z_{t,i}^{\phi ,\phi _{h}^{\prime },\phi _{\ell }^{\prime
\prime }}+\left( \mathbf{Z}_{t}^{\phi ,\phi _{h}^{\prime },\phi _{\ell
}^{\prime \prime }}\right) ^{\top }\left( \partial _{y}\phi _{t,i}\right)
\left( X_{t,i}^{\phi },\mathbf{X}_{t}^{\phi }\right)   \notag \\
&&+\delta _{h,i}\left( \left( \partial _{x}\phi _{h}^{\prime }\right) \left(
X_{t,h}^{\phi },\mathbf{X}_{t}^{\phi }\right) Y_{t,h}^{\phi ,\phi _{\ell
}^{\prime \prime }}+\left( \mathbf{Y}_{t}^{\phi ,\phi _{\ell }^{\prime
\prime }}\right) ^{\top }\left( \partial _{y}\phi _{h}^{\prime }\right)
\left( X_{t,h}^{\phi },\mathbf{X}_{t}^{\phi }\right) \right)   \notag \\
&&+\delta _{\ell ,i}\left( \left( \partial _{x}\phi _{\ell }^{\prime \prime
}\right) \left( X_{t,\ell }^{\phi },\mathbf{X}_{t}^{\phi }\right) Y_{t,\ell
}^{\phi ,\phi _{h}^{\prime }}+\left( \mathbf{Y}_{t}^{\phi ,\phi _{h}^{\prime
}}\right) ^{\top }\left( \partial _{y}\phi _{\ell }^{\prime \prime }\right)
\left( X_{t,\ell }^{\phi },\mathbf{X}_{t}^{\phi }\right) \right) .
\label{seccon}
\end{eqnarray}

\begin{remark}
From (\ref{y1}) and (\ref{z1}), we observe that the sensitivity analysis%
\footnote{%
In classical stochastic control theory, we call them variational equations
which are employed to establish a Hamiltonian system see \cite{YZ}.} of
closed-loop controls are more intricate than those for open-loop controls
(see \cite{GLZ2025}). This complexity arises due to the interdependence
inherent in the closed-loop game: Player $i$ changes her control$%
\Longrightarrow $State trajectory of the system $\mathbf{X}_{t}^{\phi }$
changes$\Longrightarrow $Player $j$'s control $\phi _{t,j}\left(
X_{t,j}^{\phi },\mathbf{X}_{t}^{\phi }\right) $ is affected. Note,
however, that under an open-loop setting such as that in Guo et al. (see \cite{GLZ2025}), for
all $t\in \left[ 0,T\right] ,$ we can derive that $\upsilon _{t,i}^{\phi
,\phi _{h}^{\prime }}=\delta _{h,i}u_{t,h}^{\prime },$whilst $\omega
_{t,i}^{\phi ,\phi _{h}^{\prime },\phi _{\ell }^{\prime \prime }}=0.$
Consequently, the analysis of $\alpha$ for the closed-loop problem becomes
significantly more challenging.
\end{remark}

Now, we define the linear derivatives of $V_{i}$ in (\ref{cost}) by means of
the sensitivity processes satisfying (\ref{y1}) (see (5.3) in \cite{GLZ1}).
Indeed, for $\forall i,h\in I_{N}$, the linear derivative $\frac{%
\partial V_{i}}{\partial u_{h}}:\mathcal{A}^{\left( N\right) }\times \Pi
\rightarrow \mathbb{R}$ of $V_{i}$ w.r.t. $\mathcal{A}_{h}$ can be expressed as
\begin{equation}
\frac{\delta V_{i}}{\delta \phi _{h}}\left( \phi \mathbf{;}\phi _{h}^{\prime
}\right) =\mathbb{E}\left[ \int_{0}^{T}\left(
\begin{array}{c}
\mathbf{Y}_{t}^{\phi ,\phi _{h}^{\prime }} \\
\upsilon _{t}^{\phi ,\phi _{h}^{\prime }}%
\end{array}%
\right) ^{\top }\left(
\begin{array}{c}
\partial _{x}f_{t,i} \\
\partial _{u}f_{t,i}%
\end{array}%
\right) \left( \mathbf{X}_{t}^{\phi },\mathbf{u}_{t}^{\phi }\right) \mathrm{d%
}t+\left( \mathbf{Y}_{T}^{\phi ,\phi _{h}^{\prime }}\right) ^{\top }\left(
\partial _{x}g_{i}\right) \left( \mathbf{X}_{T}^{\phi }\right) \right] .
\label{linear1}
\end{equation}%
In stochastic control theory, the variation equation is employed to drive
the necessary condition. In other words, it plays the role of an intermediate
bridge. Intuitively, the expression of (\ref{linear1}) is akin to the
variation inequality in stochastic optimal control. Hence, we plan to
rewrite the linear derivatives of $V_{i}$ by means of BSDE.

To this end, we write $l_{t,i}\left( \cdot \right) =l_{t,i}\left(
X_{t,i}^{\phi },\mathbf{X}_{t}^{\phi },u_{t,i}^{\phi }\right) $ for $%
l=b_{i},\sigma _{i},f_{i},$ etc.

\begin{eqnarray*}
\mathbf{\bar{B}}_{x}\left( t,\mathbf{X}_{t}^{\phi },\mathbf{u}_{t}^{\phi
}\right) &=&diag\left( \partial _{x}\left( b_{t,1}+\phi _{t,1}\right) \left(
\cdot \right) ,\ldots ,\partial _{x}\left( b_{t,N}+\phi _{t,N}\right) \left(
\cdot \right) \right) _{N\times N}, \\
\mathbf{B}_{1}^{i}\left( t,\mathbf{X}_{t}^{\phi },\mathbf{u}_{t}^{\phi
}\right) &=&\left[
\begin{array}{ccccc}
0 & \cdots & 1 & \cdots & 0%
\end{array}%
\right] _{N\times 1}^{\top },\text{ }i\in I_{N}, \\
\mathbf{B}_{0,y}\left( t,\mathbf{X}_{t}^{\phi },\mathbf{u}_{t}^{\phi
}\right) &=&\left[
\begin{array}{ccc}
\partial _{y_{1}}\left( b_{t,1}+\phi _{t,1}\right) \left( \cdot \right) &
\cdots & \partial _{y_{N}}\left( b_{t,1}+\phi _{t,1}\right) \left( \cdot
\right) \\
\vdots & \ddots & \vdots \\
\partial _{y_{1}}\left( b_{t,h}+\phi _{t,h}\right) \left( \cdot \right) &
\cdots & \partial _{y_{N}}\left( b_{t,h}+\phi _{t,h}\right) \left( \cdot
\right) \\
\vdots & \ddots & \vdots \\
\partial _{y_{1}}\left( b_{t,N}+\phi _{t,N}\right) \left( \cdot \right) &
\cdots & \partial _{y_{N}}\left( b_{t,N}+\phi _{t,N}\right) \left( \cdot
\right)%
\end{array}%
\right] _{N\times N},
\end{eqnarray*}%
and%
\begin{eqnarray*}
\mathbf{\bar{\Pi}}_{x}^{i}\left( t,\mathbf{X}_{t}^{\phi },\phi _{t}\right)
&=&\left[
\begin{array}{ccccc}
0 & \cdots & 0 & \cdots & 0 \\
\vdots & \ddots & \vdots & \ldots & 0 \\
0 & \cdots & \partial _{x}\left( \sigma _{t,i}+\phi _{t,i}\right) \left(
\cdot \right) & \cdots & 0 \\
\vdots & \ddots & \vdots & \ddots & \vdots \\
0 & \cdots & 0 & \cdots & 0%
\end{array}%
\right] _{N\times N}, \\
\mathbf{\Pi }_{0,y}^{i}\left( t,\mathbf{X}_{t}^{\phi },\phi _{t}\right) &=&%
\left[
\begin{array}{ccc}
0 & \cdots & 0 \\
\vdots & \vdots & \vdots \\
\partial _{y_{1}}\left( \sigma _{t,i}+\phi _{t,i}\right) \left( \cdot \right)
& \cdots & \partial _{y_{N}}\left( \sigma _{t,i}+\phi _{t,i}\right) \left(
\cdot \right) \\
\vdots & \vdots & \vdots \\
0 & \cdots & 0,%
\end{array}%
\right] _{N\times N}, \\
\mathbf{\Pi }_{1}^{j}\left( t,\mathbf{X}_{t}^{\phi },\phi _{t}\right) &=&%
\left[
\begin{array}{ccccc}
0 & \ldots & 1 & \cdots & 0%
\end{array}%
\right]^{\top } _{N\times 1},\text{ }j\in I_{N}.
\end{eqnarray*}%
Then
\begin{equation}
\left\{
\begin{array}{rcl}
\mathrm{d}\mathbf{Y}_{t}^{\phi ,\phi _{h}^{\prime }} & = & \left[ \mathbf{%
B_{0,x,y}\left( t,\mathbf{X}_{t}^{\phi },\mathbf{u}_{t}^{\phi }\right) Y}%
_{t}^{\phi ,\phi _{h}^{\prime }}+\mathbf{B}_{1}^{h}\left( t,\mathbf{X}%
_{t}^{\phi },\mathbf{u}_{t}^{\phi }\right) \phi _{t,h}^{\prime }\right]
\mathrm{d}t \\
&  & +\sum\limits_{j=1}^{N}\left[ \mathbf{\Pi }_{0,x,y}^{j}\left( t,\mathbf{X%
}_{t}^{\phi },\mathbf{u}_{t}^{\phi }\right) \mathbf{Y}_{t}^{\phi ,\phi
_{h}^{\prime }}+\mathbf{\Pi }_{1}^{j}\left( t,\mathbf{X}_{t}^{\phi },\mathbf{%
u}_{t}^{\phi }\right) \delta _{j,h}\phi _{t,h}^{\prime }\right] \mathrm{d}%
W_{t}^{j}\mathrm{,} \\
\mathbf{Y}_{0}^{\phi ,\phi _{h}^{\prime }} & = & 0,%
\end{array}%
\right.  \label{var1}
\end{equation}%
where
\begin{eqnarray*}
\mathbf{B}_{0,x,y}\left( t,\mathbf{X}_{t}^{\phi },\mathbf{u}_{t}^{\phi
}\right) &=&\left( \mathbf{\bar{B}}_{x}+\mathbf{B}_{0,y}\right) \left( t,%
\mathbf{X}_{t}^{\phi },\mathbf{u}_{t}^{\phi }\right) , \\
\mathbf{\Pi }_{0,x,y}^{j}\left( t,\mathbf{X}_{t}^{\phi },\mathbf{u}%
_{t}^{\phi }\right) &=&\left( \mathbf{\bar{\Pi}}_{x}^{j}+\mathbf{\Pi }%
_{0,y}^{j}\right) \left( t,\mathbf{X}_{t}^{\phi },\mathbf{u}_{t}^{\phi
}\right) ,\text{ }j\in I_{N}.
\end{eqnarray*}%
We now introduce the following BSDE:%
\begin{equation}
\left\{
\begin{array}{rcl}
-\mathrm{d}P_{t,i} & = & \Big [\mathbf{B}\left( t,\mathbf{X}_{t}^{\phi },%
\mathbf{u}_{t}^{\phi }\right) ^{\top }P_{t,i}+\sum\limits_{j=1}^{N}\mathbf{%
\Pi }\left( t,\mathbf{X}_{t}^{\phi },\phi _{t}\right) ^{\top }Q_{t,j}+\left(
\partial _{x}f_{t,i}\right) \left( t,\mathbf{X}_{t}^{\phi },\mathbf{u}%
_{t}^{\phi }\right) \Big ]\mathrm{d}t-\sum\limits_{j=1}^{N}Q_{t,j}\mathrm{d}%
W_{t}^{j}, \\
P_{T,i} & = & \partial _{x}g_{i}\left( \mathbf{X}_{T}^{\phi }\right) ,\text{
}\forall i\in I_{N}.%
\end{array}%
\right.   \label{adj1}
\end{equation}
By Lemma \ref{bs1}, there exists a unique adapted solution $\left(
P_{i},\left\{ Q_{j}\right\} _{j=1,\ldots ,N}\right) $ to BSDE (\ref{adj1}).

\begin{remark}
According to stochastic control theory, the dimension of the adjoint
equation (BSDE actually) is equal to that of the forward stochastic system.
Nevertheless, in the framework of potential game, the value function of the $i$%
-th depends on $g_{i}\left( \mathbf{X}_{T}^{\phi }\right) $ which means that
the gradient of $g_{i}$ is a column vector. Hence, the adjoint of the $i$-th is $%
N$-dimension rather than one-dimensional.
\end{remark}

\begin{remark}
Note that the adjoint equation (\ref{adj1}) is independent of $h$ since $%
B_{0}\left( t,\mathbf{X}_{t}^{\phi },\phi _{t}\right) ,$ $\Pi _{0}^{i}\left(
t,\mathbf{X}_{t}^{\phi },\phi _{t}\right) $ and $\left( \partial
_{x}f_{t,i}\right) \left( \mathbf{X}_{t}^{\phi },\mathbf{u}_{t}^{\phi
}\right) $ are independent of $h.$
\end{remark}

\begin{remark}
Before applying It\^o's formula to $\langle Y^{\phi,\phi_h'}, P_i\rangle$, we note the integrability conditions, which are especially relevant under control‑dependent diffusion.
From Assumptions (A1)-(A2) and Lemma \ref{l2}, \ref{l4}, we have $Y^{\phi,\phi_h'}\in \mathcal{S}^2(\mathbb{R}^N)$. By Lemma \ref{bs2}, the adjoint BSDE solution satisfies $P_i\in \mathcal{S}^2(\mathbb{R}^N)$ and $Q_i\in \mathcal{H}^2(\mathbb{R}^{N\times d})$. These memberships justify the application of It\^o's formula. Moreover, the resulting stochastic integral terms are true martingales, so their expectations vanish. Control‑dependent diffusion raises higher‑moment technicalities, which are controlled by our global Lipschitz and bounded second‑derivative assumptions.
\end{remark}

By applying It\^{o}'s formula to $\left\langle \mathbf{Y}^{\phi ,\phi _{h}^{\prime
}},P_{i}\right\rangle $ on $\left[ 0,T\right] $, we have
\begin{eqnarray*}
&&\mathbb{E}\left[ \left\langle \mathbf{Y}_{T}^{\phi ,\phi _{h}^{\prime
}},P_{T,i}\right\rangle \right]  \\
&=&\mathbb{E}\left[ \left\langle \mathbf{Y}_{T}^{\phi ,\phi _{h}^{\prime
}},\left( \partial _{x}g_{i}\right) \left( \mathbf{X}_{T}^{\phi }\right)
\right\rangle \right]  \\
&=&\mathbb{E}\left[ \int_{0}^{T}\left\langle P_{t,i},\mathbf{B\left( t,%
\mathbf{X}_{t}^{\phi },\mathbf{u}_{t}^{\phi }\right) Y}_{t}^{\phi ,\phi
_{h}^{\prime }}+\mathbf{B}_{1}^{h}\left( t,\mathbf{X}_{t}^{\phi },\mathbf{u}%
_{t}^{\phi }\right) \phi _{t,h}^{\prime }\right\rangle dt\right]  \\
&&+\mathbb{E}\Bigg [\int_{0}^{T}\Big (-\left\langle \mathbf{Y}_{t}^{\phi
,\phi _{h}^{\prime }},\mathbf{B}\left( t,\mathbf{X}_{t}^{\phi },\mathbf{u}%
_{t}^{\phi }\right) ^{\top }P_{t,i}\right\rangle -\left\langle \mathbf{Y}%
_{t}^{\phi ,\phi _{h}^{\prime }},\sum_{j=1}^{N}\mathbf{\Pi }^{j}\left( t,%
\mathbf{X}_{t}^{\phi },\mathbf{u}_{t}^{\phi }\right) ^{\top
}Q_{t,j}\right\rangle  \\
&&-\left\langle \mathbf{Y}_{t}^{\phi ,\phi _{h}^{\prime }},\left( \partial
_{x}f_{t,i}\right) \left( \mathbf{X}_{t}^{\phi },\mathbf{u}_{t}^{\phi
}\right) \right\rangle \Big )\mathrm{d}t\Bigg ] \\
&&+\mathbb{E}\Bigg \{\int_{0}^{T}\Big [\sum_{j=1}^{N}\left( \mathbf{\Pi }%
^{j}\left( t,\mathbf{X}_{t}^{\phi },\mathbf{u}_{t}^{\phi }\right) \mathbf{Y}%
_{t}^{\phi ,\phi _{h}^{\prime }}+\mathbf{\Pi }_{1}^{j}\left( t,\mathbf{X}%
_{t}^{\phi },\mathbf{u}_{t}^{\phi }\right) \delta _{i,h}\phi _{t,h}^{\prime
}\right) Q_{t,j}\Big ]\mathrm{d}t\Bigg \} \\
&=&\mathbb{E}\Bigg [\int_{0}^{T}\Big (\left\langle P_{t,i},\mathbf{B}%
_{1}^{h}\left( t,\mathbf{X}_{t}^{\phi },\mathbf{u}_{t}^{\phi }\right) \phi
_{t,h}^{\prime }\right\rangle -\left\langle \mathbf{Y}_{t}^{\phi ,\phi
_{h}^{\prime }},\left( \partial _{x}f_{t,i}\right) \left( \mathbf{X}%
_{t}^{\phi },\mathbf{u}_{t}^{\phi }\right) \right\rangle  \\
&&+\sum_{j=1}^{N}\left\langle \mathbf{\Pi }_{1}^{j}\left( t,\mathbf{X}%
_{t}^{\phi },\mathbf{u}_{t}^{\phi }\right) \delta _{j,h},Q_{t,j}\phi
_{t,h}^{\prime }\right\rangle \Big )\mathrm{d}t\Bigg ] \\
&=&\mathbb{E}\Bigg [\int_{0}^{T}\Big (\left\langle P_{t,i},\mathbf{B}%
_{1}^{h}\left( t,\mathbf{X}_{t}^{\phi },\mathbf{u}_{t}^{\phi }\right) \phi
_{t,h}^{\prime }\right\rangle -\left\langle \mathbf{Y}_{t}^{\phi ,\phi
_{h}^{\prime }},\left( \partial _{x}f_{t,i}\right) \left( \mathbf{X}%
_{t}^{\phi },\mathbf{u}_{t}^{\phi }\right) \right\rangle +\left(
Q_{t,h}\right) _{h}\phi _{t,h}^{\prime }\Big )\mathrm{d}t\Bigg ].
\end{eqnarray*}
Therefore,
\begin{eqnarray*}
\frac{\delta V_{i}}{\delta u_{h}}\left( \mathbf{u};u_{h}^{\prime }\right)
&=&\mathbb{E}\Bigg \{\int_{0}^{T}\left(
\begin{array}{c}
\mathbf{Y}_{t}^{\phi ,\phi _{h}^{\prime }} \\
\upsilon _{t}^{\phi ,\phi _{h}^{\prime }}%
\end{array}%
\right) ^{\top }\left(
\begin{array}{c}
\partial _{x}f_{t,i} \\
\partial _{u}f_{t,i}%
\end{array}%
\right) \left( \mathbf{X}_{t}^{\phi },\mathbf{u}_{t}^{\phi }\right) \mathrm{d%
}t+\left( \mathbf{Y}_{T}^{\phi ,\phi _{h}^{\prime }}\right) ^{\top }\left(
\partial _{x}g_{i}\right) \left( \mathbf{X}_{T}^{\phi }\right) \Bigg \} \\
&=&\mathbb{E}\Bigg \{\left\langle \mathbf{Y}_{T}^{\phi ,\phi _{h}^{\prime
}},\left( \partial _{x}g_{i}\right) \left( \mathbf{X}_{T}^{\phi }\right)
\right\rangle  \\
&&+\int_{0}^{T}\left[ \mathbf{Y}_{t}^{\phi ,\phi _{h}^{\prime }}\left(
\partial _{x}f_{t,i}\right) \left( \mathbf{X}_{t}^{\phi },\mathbf{u}%
_{t}^{\phi }\right) +\left\langle \partial _{u}f_{t,i}\left( \mathbf{X}%
_{t}^{\phi },\mathbf{u}_{t}^{\phi }\right) ,\upsilon _{t}^{\phi ,\phi
_{h}^{\prime }}\right\rangle \right] \mathrm{d}t\Bigg \} \\
&=&\mathbb{E}\Bigg \{\int_{0}^{T}\Bigg [\left( P_{t,i}\mathbf{B}%
_{1}^{h}\left( t,\mathbf{X}_{t}^{\phi },\mathbf{u}_{t}^{\phi }\right)
+\sum_{i=1}^{N}\mathbf{\Pi }_{1}^{i}\left( t,\mathbf{X}_{t}^{\phi },\mathbf{u%
}_{t}^{\phi }\right) Q_{t,i}\delta _{i,h}\right) \phi _{t,h}^{\prime } \\
&&+\left\langle \partial _{u}f_{t,i}\left( \mathbf{X}_{t}^{\phi },\mathbf{u}%
_{t}^{\phi }\right) ,\upsilon _{t}^{\phi ,\phi _{h}^{\prime }}\right\rangle %
\Bigg ]\mathrm{d}t\Bigg \} \\
&=&\mathbb{E}\Bigg \{\int_{0}^{T}\Big [\left( \left( P_{t,i}\right)
_{h}+\left( Q_{t,h}\right) _{h}\right) \phi _{t,h}^{\prime }+\left\langle
\partial _{u}f_{t,i}\left( \mathbf{X}_{t}^{\phi },\mathbf{u}_{t}^{\phi
}\right) ,\upsilon _{t}^{\phi ,\phi _{h}^{\prime }}\right\rangle \Big ]%
\mathrm{d}t\Bigg \}.
\end{eqnarray*}
where $\left( Q_{t,h}\right) _{h}$ represents the $h$-th component of the
column vector $Q_{t,h},$ the same holds for $\left( P_{t,i}\right) _{h}.$

Recall

\begin{eqnarray*}
\mathfrak{f}_{t,i}^{\phi ,\phi _{h}^{\prime },\phi _{\ell }^{\prime \prime
}} &:&=\left(
\begin{array}{c}
Y_{t,i}^{\phi ,\phi _{h}^{\prime }} \\
\mathbf{Y}_{t}^{\phi ,\phi _{h}^{\prime }}%
\end{array}%
\right) ^{\top }\left(
\begin{array}{cc}
\partial _{xx}^{2}\left( b_{t,i}+\phi _{t,i}\right) & \partial
_{xy}^{2}\left( b_{t,i}+\phi _{t,i}\right) \\
\partial _{yx}^{2}\left( b_{t,i}+\phi _{t,i}\right) & \partial
_{yy}^{2}\left( b_{t,i}+\phi _{t,i}\right)%
\end{array}%
\right) \left( X_{t,i}^{\phi },\mathbf{X}_{t}^{\phi }\right) \left(
\begin{array}{c}
Y_{t,i}^{\phi ,\phi _{\ell }^{\prime \prime }} \\
\mathbf{Y}_{t}^{\phi ,\phi _{\ell }^{\prime \prime }}%
\end{array}%
\right) \\
&&+\delta _{h,i}\left( \left( \partial _{x}\phi _{t,h}^{\prime }\right)
\left( X_{t,i}^{\phi },\mathbf{X}_{t}^{\phi }\right) Y_{t,i}^{\phi ,\phi
_{\ell }^{\prime \prime }}+\left( \mathbf{Y}_{t}^{\phi ,\phi _{\ell
}^{\prime \prime }}\right) ^{\top }\left( \partial _{y}\phi _{t,h}^{\prime
}\right) \left( X_{t,i}^{\phi },\mathbf{X}_{t}^{\phi }\right) \right) \\
&&+\delta _{\ell ,i}\left( \left( \partial _{x}\phi _{t,\ell }^{\prime
\prime }\right) \left( X_{t,i}^{\phi },\mathbf{X}_{t}^{\phi }\right)
Y_{t,i}^{\phi ,\phi _{h}^{\prime }}+\left( \mathbf{Y}_{t}^{\phi ,\phi
_{h}^{\prime }}\right) ^{\top }\left( \partial _{y}\phi _{t,\ell }^{\prime
\prime }\right) \left( X_{t,i}^{\phi },\mathbf{X}_{t}^{\phi }\right) \right)
\\
&=&\text{tr}\left[ \partial _{yy}^{2}\left( b_{t,i}+\phi _{t,i}\right)
\left( \cdot \right) \mathcal{Y}_{t}^{\phi ,\phi _{h}^{\prime },\phi _{\ell
}^{\prime \prime }}\right] +\Gamma _{t,i}^{\phi ,\phi _{h}^{\prime },\phi
_{\ell }^{\prime \prime }}\left( \bar{b}_{t}\right) \\
&&+\delta _{h,i}\left( \left( \partial _{x}\phi _{t,h}^{\prime }\right)
\left( X_{t,i}^{\phi },\mathbf{X}_{t}^{\phi }\right) Y_{t,i}^{\phi ,\phi
_{\ell }^{\prime \prime }}+\left( \mathbf{Y}_{t}^{\phi ,\phi _{\ell
}^{\prime \prime }}\right) ^{\top }\left( \partial _{y}\phi _{t,h}^{\prime
}\right) \left( X_{t,i}^{\phi },\mathbf{X}_{t}^{\phi }\right) \right) \\
&&+\delta _{\ell ,i}\left( \left( \partial _{x}\phi _{t,\ell }^{\prime
\prime }\right) \left( X_{t,i}^{\phi },\mathbf{X}_{t}^{\phi }\right)
Y_{t,i}^{\phi ,\phi _{h}^{\prime }}+\left( \mathbf{Y}_{t}^{\phi ,\phi
_{h}^{\prime }}\right) ^{\top }\left( \partial _{y}\phi _{t,\ell }^{\prime
\prime }\right) \left( X_{t,i}^{\phi },\mathbf{X}_{t}^{\phi }\right) \right)
,
\end{eqnarray*}%
where we put
\begin{eqnarray*}
\mathcal{Y}_{t}^{\phi ,\phi _{h}^{\prime },\phi _{\ell }^{\prime \prime }}
&=&\mathbf{Y}_{t}^{\phi ,\phi _{\ell }^{\prime \prime }}\left( \mathbf{Y}%
_{t}^{\phi ,\phi _{h}^{\prime }}\right) ^{\top }, \\
\bar{b}_{t}\left( \mathbf{X}_{t}^{\phi }\right) &=&\left( \left(
b_{t,1}+\phi _{t,1}\right) \left( X_{t,1}^{\phi },\mathbf{X}_{t}^{\phi
}\right) ,\ldots ,\left( b_{t,N}+\phi _{t,N}\right) \left( X_{t,N}^{\phi },%
\mathbf{X}_{t}^{\phi }\right) \right) ^{\top },\text{ } \\
\Gamma _{t,i}^{\phi ,\phi _{h}^{\prime },\phi _{\ell }^{\prime \prime
}}\left( \bar{b}_{t}\right) &=&Y_{t,i}^{\phi ,\phi _{h}^{\prime }}\partial
_{xx}^{2}\left( b_{t,i}+\phi _{t,i}\right) \left( X_{t,i}^{\phi },\mathbf{X}%
_{t}^{\phi }\right) Y_{t,i}^{\phi ,\phi _{\ell }^{\prime \prime }} \\
&&+\left( \mathbf{Y}_{t}^{\phi ,\phi _{h}^{\prime }}\right) ^{\top }\partial
_{yx}^{2}\left( b_{t,i}+\phi _{t,i}\right) \left( X_{t,i}^{\phi },\mathbf{X}%
_{t}^{\phi }\right) Y_{t,i}^{\phi ,\phi _{\ell }^{\prime \prime }} \\
&&+Y_{t,i}^{\phi ,\phi _{h}^{\prime }}\partial _{xy}^{2}\left( b_{t,i}+\phi
_{t,i}\right) \left( X_{t,i}^{\phi },\mathbf{X}_{t}^{\phi }\right) \mathbf{Y}%
_{t}^{\phi ,\phi _{\ell }^{\prime \prime }}
\end{eqnarray*}%
and
\begin{eqnarray*}
\mathfrak{g}_{t,i}^{\phi ,\phi _{h}^{\prime },\phi _{\ell }^{\prime \prime
}} &:&=\left(
\begin{array}{c}
Y_{t,i}^{\phi ,\phi _{h}^{\prime }} \\
\mathbf{Y}_{t}^{\phi ,\phi _{h}^{\prime }}%
\end{array}%
\right) ^{\top }\left(
\begin{array}{cc}
\partial _{xx}^{2}\left( \sigma _{t,i}+\phi _{t,i}\right) & \partial
_{xy}^{2}\left( \sigma _{t,i}+\phi _{t,i}\right) \\
\partial _{yx}^{2}\left( \sigma _{t,i}+\phi _{t,i}\right) & \partial
_{yy}^{2}\left( \sigma _{t,i}+\phi _{t,i}\right)%
\end{array}%
\right) \left( X_{t,i}^{\phi },\mathbf{X}_{t}^{\phi }\right) \left(
\begin{array}{c}
Y_{t,i}^{\phi ,\phi _{\ell }^{\prime \prime }} \\
\mathbf{Y}_{t}^{\phi ,\phi _{\ell }^{\prime \prime }}%
\end{array}%
\right) \\
&&+\delta _{h,i}\left( \left( \partial _{x}\phi _{t,h}^{\prime }\right)
\left( X_{t,i}^{\phi },\mathbf{X}_{t}^{\phi }\right) Y_{t,i}^{\phi ,\phi
_{\ell }^{\prime \prime }}+\left( \mathbf{Y}_{t}^{\phi ,\phi _{\ell
}^{\prime \prime }}\right) ^{\top }\left( \partial _{y}\phi _{t,h}^{\prime
}\right) \left( X_{t,i}^{\phi },\mathbf{X}_{t}^{\phi }\right) \right) \\
&&+\delta _{\ell ,i}\left( \left( \partial _{x}\phi _{t,\ell }^{\prime
\prime }\right) \left( X_{t,i}^{\phi },\mathbf{X}_{t}^{\phi }\right)
Y_{t,i}^{\phi ,\phi _{h}^{\prime }}+\left( \mathbf{Y}_{t}^{\phi ,\phi
_{h}^{\prime }}\right) ^{\top }\left( \partial _{y}\phi _{t,\ell }^{\prime
\prime }\right) \left( X_{t,i}^{\phi },\mathbf{X}_{t}^{\phi }\right) \right)
\\
&=&\text{tr}\left[ \partial _{yy}^{2}\left( \sigma _{t,i}+\phi _{t,i}\right)
\left( X_{t,i}^{\phi },\mathbf{X}_{t}^{\phi }\right) \mathcal{Y}_{t}^{\phi
,\phi _{h}^{\prime },\phi _{\ell }^{\prime \prime }}\right] +\Xi
_{t,i}^{\phi ,\phi _{h}^{\prime },\phi _{\ell }^{\prime \prime }}\left( \bar{%
\sigma}_{t}\right) \\
&&+\delta _{h,i}\left( \left( \partial _{x}\phi _{t,h}^{\prime }\right)
\left( X_{t,i}^{\phi },\mathbf{X}_{t}^{\phi }\right) Y_{t,i}^{\phi ,\phi
_{\ell }^{\prime \prime }}+\left( \mathbf{Y}_{t}^{\phi ,\phi _{\ell
}^{\prime \prime }}\right) ^{\top }\left( \partial _{y}\phi _{t,h}^{\prime
}\right) \left( X_{t,i}^{\phi },\mathbf{X}_{t}^{\phi }\right) \right) \\
&&+\delta _{\ell ,i}\left( \left( \partial _{x}\phi _{t,\ell }^{\prime
\prime }\right) \left( X_{t,i}^{\phi },\mathbf{X}_{t}^{\phi }\right)
Y_{t,i}^{\phi ,\phi _{h}^{\prime }}+\left( \mathbf{Y}_{t}^{\phi ,\phi
_{h}^{\prime }}\right) ^{\top }\left( \partial _{y}\phi _{t,\ell }^{\prime
\prime }\right) \left( X_{t,i}^{\phi },\mathbf{X}_{t}^{\phi }\right) \right)
\end{eqnarray*}%
where,
\begin{eqnarray}
\bar{\sigma}_{t,i}\left( X_{t,i}^{\phi },\mathbf{X}_{t}^{\phi }\right)
&=&\left( 0,\ldots \sigma _{t,i}+\phi _{t,i},\ldots 0\right) ^{\top }\left(
X_{t,i}^{\phi },\mathbf{X}_{t}^{\phi }\right) ,  \notag \\
\Xi _{t,i}^{\phi ,\phi _{h}^{\prime },\phi _{\ell }^{\prime \prime }}\left(
\bar{\sigma}\right) &=&Y_{t,i}^{\phi ,\phi _{h}^{\prime }}\partial
_{xx}^{2}\left( \sigma _{t,i}+\phi _{t,i}\right) Y_{t,i}^{\phi ,\phi _{\ell
}^{\prime \prime }}+\left( \mathbf{Y}_{t}^{\phi ,\phi _{h}^{\prime }}\right)
^{\top }\partial _{yx}^{2}\left( \sigma _{t,i}+\phi _{t,i}\right)
Y_{t,i}^{\phi ,\phi _{\ell }^{\prime \prime }}  \notag \\
&&+Y_{t,i}^{\phi ,\phi _{h}^{\prime }}\partial _{xy}^{2}\left( \sigma
_{t,i}+\phi _{t,i}\right) \mathbf{Y}_{t}^{\phi ,\phi _{\ell }^{\prime \prime
}}.  \notag
\end{eqnarray}%
Therefore,
\begin{eqnarray*}
\left[ \mathfrak{f}_{t,1}^{\phi ,\phi _{h}^{\prime },\phi _{\ell }^{\prime
\prime }},\ldots \mathfrak{f}_{t,i}^{\phi ,\phi _{h}^{\prime },\phi _{\ell
}^{\prime \prime }},\ldots ,\mathfrak{f}_{t,N}^{\phi ,\phi _{h}^{\prime
},\phi _{\ell }^{\prime \prime }}\right] _{N\times 1}^{\top } &=&\partial
_{yy}^{2}\bar{b}_{t}\left( X_{t,i}^{\phi },\mathbf{X}_{t}^{\phi }\right)
\star \mathcal{Y}_{t}^{\phi ,\phi _{h}^{\prime },\phi _{\ell }^{\prime
\prime }}+\Gamma _{t}^{\phi ,\phi _{h}^{\prime },\phi _{\ell }^{\prime
\prime }}\left( \bar{b}_{t}\right) , \\
\left[ 0,\ldots \mathfrak{g}_{t,j}^{\phi ,\phi _{h}^{\prime },\phi _{\ell
}^{\prime \prime }},\ldots ,0\right] _{N\times 1}^{\top }
&=&\partial _{yy}^{2}\bar{\sigma}_{t,j}\left( X_{t,i}^{\phi },\mathbf{X}%
_{t}^{\phi }\right) \star \mathcal{Y}_{t}^{\phi ,\phi _{h}^{\prime },\phi
_{\ell }^{\prime \prime }}+\bar{\Xi}_{t,j}^{\phi ,\phi _{h}^{\prime },\phi
_{\ell }^{\prime \prime }}\left( \bar{\sigma}\right) .
\end{eqnarray*}
where
\begin{eqnarray*}
\partial _{yy}^{2}\bar{b}_{t}\left( X_{t,i}^{\phi },\mathbf{X}_{t}^{\phi
}\right) \star \mathcal{Y}_{t}^{\phi ,\phi _{h}^{\prime },\phi _{\ell
}^{\prime \prime }} &=&\left[
\begin{array}{c}
\text{tr}\left[ \partial _{yy}^{2}\left( b_{t,1}+\phi _{t,1}\right) \left(
X_{t,1}^{\phi },\mathbf{X}_{t}^{\phi }\right) \mathcal{Y}_{t}^{\phi ,\phi
_{h}^{\prime },\phi _{\ell }^{\prime \prime }}\right] \\
\vdots \\
\text{tr}\left[ \partial _{yy}^{2}\left( b_{t,i}+\phi _{t,i}\right) \left(
X_{t,i}^{\phi },\mathbf{X}_{t}^{\phi }\right) \mathcal{Y}_{t}^{\phi ,\phi
_{h}^{\prime },\phi _{\ell }^{\prime \prime }}\right] \\
\vdots \\
\text{tr}\left[ \partial _{yy}^{2}\left( b_{t,N}+\phi _{t,N}\right) \left(
X_{t,N}^{\phi },\mathbf{X}_{t}^{\phi }\right) \mathcal{Y}_{t}^{\phi ,\phi
_{h}^{\prime },\phi _{\ell }^{\prime \prime }}\right]%
\end{array}%
\right] ,\text{ } \\
\Gamma _{t}^{\phi ,\phi _{h}^{\prime },\phi _{\ell }^{\prime \prime }}\left(
\bar{b}_{t}\right) &=&\left[
\begin{array}{ccccc}
\Gamma _{t,1}^{\phi ,\phi _{h}^{\prime },\phi _{\ell }^{\prime \prime
}}\left( \bar{b}_{t}\right) & \cdots & \Gamma _{t,i}^{\phi ,\phi
_{h}^{\prime },\phi _{\ell }^{\prime \prime }}\left( \bar{b}_{t}\right) &
\cdots & \Gamma _{t,N}^{\phi ,\phi _{h}^{\prime },\phi _{\ell }^{\prime
\prime }}\left( \bar{b}_{t}\right)%
\end{array}%
\right] ^{\top }
\end{eqnarray*}%
and
\begin{eqnarray*}
\partial _{yy}^{2}\bar{\sigma}_{t,i}\left( X_{t,i}^{\phi },\mathbf{X}%
_{t}^{\phi }\right) \star \mathcal{Y}_{t}^{\phi ,\phi _{h}^{\prime },\phi
_{\ell }^{\prime \prime }} &=&\left[
\begin{array}{c}
0 \\
\vdots \\
\text{tr}\left[ \partial _{yy}^{2}\left( \sigma _{t,i}+\phi _{t,i}\right)
\mathcal{Y}_{t}^{\phi ,\phi _{h}^{\prime },\phi _{\ell }^{\prime \prime }}%
\right] \\
\vdots \\
0%
\end{array}%
\right] , \\
\text{ }\bar{\Xi}_{t,i}^{\phi ,\phi _{h}^{\prime },\phi _{\ell }^{\prime
\prime }}\left( \bar{\sigma}_{t}\right) &=&\left[
\begin{array}{ccccc}
0 & \cdots & \Xi _{t,i}^{\phi ,\phi _{h}^{\prime },\phi _{\ell }^{\prime
\prime }}\left( \bar{\sigma}_{t}\right) & \cdots & 0%
\end{array}%
\right] _{N\times 1}^{\top }.
\end{eqnarray*}%
While
\begin{eqnarray*}
\mathfrak{\chi }_{t,i}^{\phi ,\phi _{h}^{\prime },\phi _{\ell }^{\prime
\prime }} &=&\delta _{h,i}\left( \left( \partial _{x}\phi _{t,h}^{\prime
}\right) \left( X_{t,i}^{\phi },\mathbf{X}_{t}^{\phi }\right) Y_{t,i}^{\phi
,\phi _{\ell }^{\prime \prime }}+\left( \mathbf{Y}_{t}^{\phi ,\phi _{\ell
}^{\prime \prime }}\right) ^{\top }\left( \partial _{y}\phi _{t,h}^{\prime
}\right) \left( X_{t,i}^{\phi },\mathbf{X}_{t}^{\phi }\right) \right) \\
&&+\delta _{\ell ,i}\left( \left( \partial _{x}\phi _{t,\ell }^{\prime
\prime }\right) \left( X_{t,i}^{\phi },\mathbf{X}_{t}^{\phi }\right)
Y_{t,i}^{\phi ,\phi _{h}^{\prime }}+\left( \mathbf{Y}_{t}^{\phi ,\phi
_{h}^{\prime }}\right) ^{\top }\left( \partial _{y}\phi _{t,\ell }^{\prime
\prime }\right) \left( X_{t,i}^{\phi },\mathbf{X}_{t}^{\phi }\right) \right)
\end{eqnarray*}%
and
\begin{eqnarray*}
\mathfrak{\psi }_{t,i}^{\phi ,\phi _{h}^{\prime },\phi _{\ell }^{\prime
\prime }} &=&\delta _{h,i}\left( \left( \partial _{x}\phi _{t,h}^{\prime
}\right) \left( X_{t,i}^{\phi },\mathbf{X}_{t}^{\phi }\right) Y_{t,i}^{\phi
,\phi _{\ell }^{\prime \prime }}+\left( \mathbf{Y}_{t}^{\phi ,\phi _{\ell
}^{\prime \prime }}\right) ^{\top }\left( \partial _{y}\phi _{t,h}^{\prime
}\right) \left( X_{t,i}^{\phi },\mathbf{X}_{t}^{\phi }\right) \right) \\
&&+\delta _{\ell ,i}\left( \left( \partial _{x}\phi _{t,\ell }^{\prime
\prime }\right) \left( X_{t,i}^{\phi },\mathbf{X}_{t}^{\phi }\right)
Y_{t,i}^{\phi ,\phi _{h}^{\prime }}+\left( \mathbf{Y}_{t}^{\phi ,\phi
_{h}^{\prime }}\right) ^{\top }\left( \partial _{y}\phi _{t,\ell }^{\prime
\prime }\right) \left( X_{t,i}^{\phi },\mathbf{X}_{t}^{\phi }\right) \right)
.
\end{eqnarray*}%
While (if $h<\ell $ without loss of generality)
\begin{eqnarray*}
\mathfrak{F}_{t}^{\phi ,\phi _{h}^{\prime },\phi _{\ell }^{\prime \prime }}
&=&\left[
\begin{array}{ccccccc}
\mathfrak{\chi }_{t,1}^{\phi ,\phi _{h}^{\prime },\phi _{\ell }^{\prime
\prime }} & \cdots & \mathfrak{\chi }_{t,h}^{\phi ,\phi _{h}^{\prime },\phi
_{\ell }^{\prime \prime }} & \cdots & \mathfrak{\chi }_{t,\ell }^{\phi ,\phi
_{h}^{\prime },\phi _{\ell }^{\prime \prime }} & \cdots & \mathfrak{\chi }%
_{t,N}^{\phi ,\phi _{h}^{\prime },\phi _{\ell }^{\prime \prime }}%
\end{array}%
\right] _{N\times 1}^{\top } \\
&=&\left[
\begin{array}{c}
0 \\
\vdots \\
\left( \partial _{x}\phi _{t,h}^{\prime }\right) \left( X_{t,i}^{\phi },%
\mathbf{X}_{t}^{\phi }\right) Y_{t,i}^{\phi ,\phi _{\ell }^{\prime \prime
}}+\left( \mathbf{Y}_{t}^{\phi ,\phi _{\ell }^{\prime \prime }}\right)
^{\top }\left( \partial _{y}\phi _{t,h}^{\prime }\right) \left(
X_{t,i}^{\phi },\mathbf{X}_{t}^{\phi }\right) \\
\vdots \\
\left( \partial _{x}\phi _{t,\ell }^{\prime \prime }\right) \left(
X_{t,i}^{\phi },\mathbf{X}_{t}^{\phi }\right) Y_{t,i}^{\phi ,\phi
_{h}^{\prime }}+\left( \mathbf{Y}_{t}^{\phi ,\phi _{h}^{\prime }}\right)
^{\top }\left( \partial _{y}\phi _{t,\ell }^{\prime \prime }\right) \left(
X_{t,i}^{\phi },\mathbf{X}_{t}^{\phi }\right) \\
\vdots \\
0%
\end{array}%
\right] _{N\times 1}
\end{eqnarray*}%
and
\begin{equation*}
\mathfrak{G}_{t,i}^{\phi ,\phi _{h}^{\prime },\phi _{\ell }^{\prime \prime
}}=\left[
\begin{array}{ccccc}
0 & \cdots & \mathfrak{\psi }_{t,i}^{\phi ,u_{h}^{\prime },u_{\ell }^{\prime
\prime }} & \cdots & 0%
\end{array}%
\right] _{N\times 1}^{\top }.
\end{equation*}%
Based on (\ref{z1}), we are able to present%
\begin{equation}
\left\{
\begin{array}{rcl}
\mathrm{d}\mathbf{Z}_{t}^{\phi ,\phi _{h}^{\prime },\phi _{\ell }^{\prime
\prime }} & = & \Big [\mathbf{B}\left( t,\mathbf{X}_{t}^{\phi },\mathbf{u}%
_{t}^{\phi }\right) \mathbf{Z}_{t}^{\phi ,\phi _{h}^{\prime },\phi _{\ell
}^{\prime \prime }}+\mathfrak{F}_{t}^{\phi ,\phi _{h}^{\prime },\phi _{\ell
}^{\prime \prime }} \\
&  & +\partial _{yy}^{2}\bar{b}_{t}\left( t,\mathbf{X}_{t}^{\phi },\mathbf{u}%
_{t}^{\phi }\right) \star \mathcal{Y}_{t}^{\mathbf{u},u_{h}^{\prime
},u_{\ell }^{\prime \prime }}+\Gamma _{t}^{\phi ,\phi _{h}^{\prime },\phi
_{\ell }^{\prime \prime }}\left( \bar{b}_{t}\right) \Big ]\mathrm{d}t \\
&  & +\sum\limits_{j=1}^{N}\Big [\mathbf{\Pi }^{j}\left( t,\mathbf{X}%
_{t}^{\phi },\mathbf{u}_{t}^{\phi }\right) \mathbf{Z}_{t}^{\phi ,\phi
_{h}^{\prime },\phi _{\ell }^{\prime \prime }}+\mathfrak{G}_{t,j}^{\phi
,\phi _{h}^{\prime },\phi _{\ell }^{\prime \prime }} \\
&  & +\partial _{yy}^{2}\bar{\sigma}_{t,j}\left( \mathbf{X}_{t}^{\phi },%
\mathbf{u}_{t}^{\phi }\right) \star \mathcal{Y}_{t}^{\phi ,\phi _{h}^{\prime
},\phi _{\ell }^{\prime \prime }}+\bar{\Xi}_{t,j}^{\phi ,\phi _{h}^{\prime
},\phi _{\ell }^{\prime \prime }}\left( \bar{\sigma}_{t}\right) \Big ]%
\mathrm{d}W_{t}^{j}\mathrm{,} \\
\mathbf{Z}_{0}^{\phi ,\phi _{h}^{\prime },\phi _{\ell }^{\prime \prime }} & =
& 0%
\end{array}%
\right.  \label{var2}
\end{equation}

Recall, for all $i,h,\ell \in I_{N}$
\begin{eqnarray}
\frac{\delta ^{2}V_{i}}{\delta \phi _{h}\delta \phi _{\ell }}\left( \phi
,\phi _{h}^{\prime },\phi _{\ell }^{\prime \prime }\right)  &=&\mathbb{E}%
\Bigg [\int_{0}^{T}\Bigg (\left(
\begin{array}{c}
\mathbf{Y}_{t}^{\phi ,\phi _{h}^{\prime }} \\
\upsilon _{t}^{\phi ,\phi _{h}^{\prime }}%
\end{array}%
\right) ^{\top }\left(
\begin{array}{cc}
\partial _{xx}^{2}f_{t,i} & \partial _{xu}^{2}f_{t,i} \\
\partial _{ux}^{2}f_{t,i} & \partial _{uu}^{2}f_{t,i}%
\end{array}%
\right) \left( \mathbf{X}_{t}^{\phi },\mathbf{u}_{t}^{\phi }\right) \left(
\begin{array}{c}
\mathbf{Y}_{t}^{\phi ,\phi _{\ell }^{\prime \prime }} \\
\upsilon _{t}^{\phi ,\phi _{\ell }^{\prime \prime }}%
\end{array}%
\right)   \notag \\
&&+\left(
\begin{array}{c}
\mathbf{Z}_{t}^{\phi ,\phi _{h}^{\prime },\phi _{\ell }^{\prime \prime }} \\
\omega _{t}^{\phi ,\phi _{h}^{\prime },\phi _{\ell }^{\prime \prime }}%
\end{array}%
\right) ^{\top }\left(
\begin{array}{c}
\partial _{x}f_{t,i} \\
\partial _{u}f_{t,i}%
\end{array}%
\right) \left( \mathbf{X}_{t}^{\phi },\mathbf{u}_{t}^{\phi }\right) \Bigg )%
\mathrm{d}t\Bigg ]  \notag \\
&&+\mathbb{E}\left[ \left( \mathbf{Y}_{T}^{\phi ,\phi _{h}^{\prime }}\right)
^{\top }\left( \partial _{xx}^{2}g_{i}\right) \left( \mathbf{X}_{T}^{\phi
}\right) \mathbf{Y}_{T}^{\phi ,\phi _{\ell }^{\prime \prime }}+\left(
\partial _{x}g_{i}\right) ^{\top }\left( \mathbf{X}_{T}^{\phi }\right)
\mathbf{Z}_{T}^{\phi ,\phi _{h}^{\prime },\phi _{\ell }^{\prime \prime }}%
\right]   \notag \\
&=&\mathbb{E}\Bigg [\int_{0}^{T}\text{tr}\left( \partial
_{xx}^{2}f_{t,i}\left( \mathbf{X}_{t}^{\phi },\mathbf{u}_{t}^{\phi }\right)
\cdot \mathbf{Y}_{t}^{\phi ,\phi _{\ell }^{\prime \prime }}\cdot \left(
\mathbf{Y}_{t}^{\phi ,\phi _{h}^{\prime }}\right) ^{\top }\right)   \notag \\
&&+\left( \upsilon _{t}^{\phi ,\phi _{h}^{\prime }}\right) ^{\top }\partial
_{ux}^{2}f_{t,i}\left( \mathbf{X}_{t}^{\phi },\mathbf{u}_{t}^{\phi }\right)
\cdot \mathbf{Y}_{t}^{\phi ,\phi _{\ell }^{\prime \prime }}+\left( \mathbf{Y}%
_{t}^{\mathbf{u},u_{h}^{\prime }}\right) ^{\top }\cdot \partial
_{xu}^{2}f_{t,i}\left( \mathbf{X}_{t}^{\phi },\mathbf{u}_{t}^{\phi }\right)
\cdot \upsilon _{t}^{\phi ,\phi _{\ell }^{\prime \prime }}  \notag \\
&&+\left( \upsilon _{t}^{\phi ,\phi _{h}^{\prime }}\right) ^{\top }\partial
_{uu}^{2}f_{t,i}\left( \mathbf{X}_{t}^{\phi },\mathbf{u}_{t}^{\phi }\right)
\cdot \upsilon _{t}^{\phi ,\phi _{\ell }^{\prime \prime }}+\left( \omega
_{t}^{\phi ,\phi _{h}^{\prime },\phi _{\ell }^{\prime \prime }}\right)
^{\top }\partial _{u}f_{t,i}\left( \mathbf{X}_{t}^{\phi },\mathbf{u}%
_{t}^{\phi }\right)   \notag \\
&&+\left( \mathbf{Z}_{t}^{\phi ,\phi _{h}^{\prime },\phi _{\ell }^{\prime
\prime }}\right) ^{\top }\partial _{x}f_{t,i}\left( \mathbf{X}_{t}^{\phi },%
\mathbf{u}_{t}^{\phi }\right) \Bigg )\mathrm{d}t\Bigg ]  \notag \\
&&+\mathbb{E}\left[ \left( \mathbf{Y}_{T}^{\phi ,\phi _{h}^{\prime }}\right)
^{\top }\left( \partial _{xx}^{2}g_{i}\right) \left( \mathbf{X}_{T}^{\phi
}\right) \mathbf{Y}_{T}^{\phi ,\phi _{\ell }^{\prime \prime }}+\left(
\partial _{x}g_{i}\right) ^{\top }\left( \mathbf{X}_{T}^{\phi }\right)
\mathbf{Z}_{T}^{\phi ,\phi _{h}^{\prime },\phi _{\ell }^{\prime \prime }}%
\right] .  \label{aerfa}
\end{eqnarray}
\begin{remark}\label{rem:bsde_duality_intuition}
Under the open-loop framework in~\cite{GLZ2025}, we have $v_{t,i}^{\phi,\phi_{h}'}=\delta_{h,i}u_{t,h}'$ and $\omega_{t,i}^{\phi,\phi_{h}',\phi_{\ell}''}=0$, so second-order variations simplify substantially. For closed-loop state-feedback policies, however, $v$ and $\omega$ are non-trivial stochastic processes, introducing additional intricate second-order contributions in~\eqref{aerfa}. Directly bounding each of those terms would be technically prohibitive.

Our adjoint BSDE~\eqref{adj1} provides an elegant workaround via duality. It\^o's formula relates forward sensitivity trajectories to backward adjoint variables $(P_{t,i},Q_{t,i})$. Rather than evaluating each high-order term separately, we only need to estimate the $L^2$-integrability of BSDE solutions, which follows from standard BSDE stability theory (Lemma~\ref{bs2}). This BSDE duality is the conceptual backbone for our non-asymptotic $\alpha$-bound.
\end{remark}
We observe that due to the appearance of feedback control, the sensitive
processes $\upsilon _{t}^{\phi ,\phi _{h}^{\prime }}$ and $\omega _{t}^{\phi
,\phi _{h}^{\prime },\phi _{\ell }^{\prime \prime }}$ play vital roles in (%
\ref{aerfa}), which is completely different from that in \cite{GLZ2025}. Applying
It\^{o}'s formula to $\left\langle \mathbf{Z}^{\phi ,\phi _{h}^{\prime
},\phi _{\ell }^{\prime \prime }},\left( \partial _{x}g_{i}\right) ^{\top
}\left( \mathbf{X}^{\phi }\right) \right\rangle $ on $\left[ 0,T\right] ,$
we have%
\begin{eqnarray*}
&&\mathbb{E}\left[ \left\langle \mathbf{Z}_{T}^{\phi ,\phi _{h}^{\prime
},\phi _{\ell }^{\prime \prime }},\left( \partial _{x}g_{i}\right) ^{\top
}\left( \mathbf{X}_{T}^{\phi }\right) \right\rangle +\int_{0}^{T}\Big <%
\mathbf{Z}_{t}^{\phi ,\phi _{h}^{\prime },\phi _{\ell }^{\prime \prime
}},\left( \partial _{x}f_{t,i}\right) \left( t,\mathbf{X}_{t}^{\phi },%
\mathbf{u}_{t}^{\phi }\right) \Big >\mathrm{d}t\right] \\
&=&\mathbb{E}\Bigg [\int_{0}^{T}\left\langle P_{t,i},\mathfrak{F}_{t}^{\phi
,\phi _{h}^{\prime },\phi _{\ell }^{\prime \prime }}+\partial _{yy}^{2}\bar{b%
}\left( t,\mathbf{X}_{t}^{\phi },\mathbf{u}_{t}^{\phi }\right) \star
\mathcal{Y}_{t}^{\phi ,\phi _{h}^{\prime },\phi _{\ell }^{\prime \prime
}}+\Gamma _{t}^{\phi ,\phi _{h}^{\prime },\phi _{\ell }^{\prime \prime
}}\left( \bar{b}_{t}\right) \right\rangle \mathrm{d}t \\
&&+\int_{0}^{T}\sum_{j=1}^{N}\left\langle Q_{t,j},\mathfrak{G}_{t,j}^{\phi
,\phi _{h}^{\prime },\phi _{\ell }^{\prime \prime }}+\partial _{yy}^{2}\bar{%
\sigma}_{j}\left( t,\mathbf{X}_{t}^{\phi },\mathbf{u}_{t}^{\phi }\right)
\star \mathcal{Y}_{t}^{\phi ,\phi _{h}^{\prime },\phi _{\ell }^{\prime
\prime }}+\bar{\Xi}_{t,j}^{\phi ,\phi _{h}^{\prime },\phi _{\ell }^{\prime
\prime }}\left( \bar{\sigma}_{t}\right) \right\rangle \mathrm{d}t\Bigg ].
\end{eqnarray*}%
Immediately,
\begin{eqnarray}
&&\frac{\delta ^{2}V_{i}}{\delta \phi _{h}\delta \phi _{\ell }}\left( \phi
,\phi _{h}^{\prime },\phi _{\ell }^{\prime \prime }\right)   \notag \\
&=&\mathbb{E}\Bigg [\int_{0}^{T}\text{tr}\left( \partial
_{xx}^{2}f_{t,i}\left( \mathbf{X}_{t}^{\phi },\mathbf{u}_{t}^{\phi }\right)
\cdot \mathbf{Y}_{t}^{\phi ,\phi _{\ell }^{\prime \prime }}\cdot \left(
\mathbf{Y}_{t}^{\phi ,\phi _{h}^{\prime }}\right) ^{\top }\right)   \notag \\
&&+\left( \upsilon _{t}^{\phi ,\phi _{h}^{\prime }}\right) ^{\top }\partial
_{ux}^{2}f_{t,i}\left( \mathbf{X}_{t}^{\phi },\mathbf{u}_{t}^{\phi }\right)
\cdot \mathbf{Y}_{t}^{\phi ,\phi _{\ell }^{\prime \prime }}+\left( \mathbf{Y}%
_{t}^{\mathbf{u},u_{h}^{\prime }}\right) ^{\top }\cdot \partial
_{xu}^{2}f_{t,i}\left( \mathbf{X}_{t}^{\phi },\mathbf{u}_{t}^{\phi }\right)
\cdot \upsilon _{t}^{\phi ,\phi _{\ell }^{\prime \prime }}  \notag \\
&&+\left( \upsilon _{t}^{\phi ,\phi _{h}^{\prime }}\right) ^{\top }\partial
_{uu}^{2}f_{t,i}\left( \mathbf{X}_{t}^{\phi },\mathbf{u}_{t}^{\phi }\right)
\cdot \upsilon _{t}^{\phi ,\phi _{\ell }^{\prime \prime }}+\left( \omega
_{t}^{\phi ,\phi _{h}^{\prime },\phi _{\ell }^{\prime \prime }}\right)
^{\top }\partial _{u}f_{t,i}\left( \mathbf{X}_{t}^{\phi },\mathbf{u}%
_{t}^{\phi }\right)   \notag \\
&&+\left\langle P_{t,i},\mathfrak{F}_{t}^{\phi ,\phi _{h}^{\prime },\phi
_{\ell }^{\prime \prime }}+\partial _{yy}^{2}\bar{b}\left( t,\mathbf{X}%
_{t}^{\phi },\mathbf{u}_{t}^{\phi }\right) \star \mathcal{Y}_{t}^{\phi ,\phi
_{h}^{\prime },\phi _{\ell }^{\prime \prime }}+\Gamma _{t}^{\phi ,\phi
_{h}^{\prime },\phi _{\ell }^{\prime \prime }}\left( \bar{b}_{t}\right)
\right\rangle   \notag \\
&&+\sum_{j=1}^{N}\left\langle Q_{t,j},\mathfrak{G}_{t,j}^{\phi ,\phi
_{h}^{\prime },\phi _{\ell }^{\prime \prime }}+\partial _{yy}^{2}\bar{\sigma}%
_{j}\left( t,\mathbf{X}_{t}^{\phi },\mathbf{u}_{t}^{\phi }\right) \star
\mathcal{Y}_{t}^{\phi ,\phi _{h}^{\prime },\phi _{\ell }^{\prime \prime }}+%
\bar{\Xi}_{t,j}^{\phi ,\phi _{h}^{\prime },\phi _{\ell }^{\prime \prime }}%
\bar{\sigma}_{t}\right\rangle \Bigg )\mathrm{d}t\Bigg ]  \notag \\
&&+\mathbb{E}\left[ \left( \mathbf{Y}_{T}^{\phi ,\phi _{h}^{\prime }}\right)
^{\top }\left( \partial _{xx}^{2}g_{i}\right) \left( \mathbf{X}_{T}^{\phi
}\right) \mathbf{Y}_{T}^{\phi ,\phi _{\ell }^{\prime \prime }}\right] .
\label{vsec}
\end{eqnarray}%
\begin{remark}\label{rem:Z_bsde_duality_explained}
We elaborate on the role of the BSDE duality transformation applied above. The original second‑order variation \eqref{vsec} contains an explicit term involving the second‑order state‑sensitivity process $Z^{\phi,\phi_h',\phi_\ell''}$. The It\^{o} formula combined with adjoint BSDE \eqref{adj1} algebraically eliminates $Z$ from the final expression of $\frac{\delta^2 V_i}{\delta\phi_h\delta\phi_\ell}$: the expectation involving $Z$ is recast as expectations of inner products between the adjoint variables $(P_{t,i},Q_{t,i})$ and the source terms $\mathfrak{F}_{t}^{\cdot},\mathfrak{G}_{t}^{\cdot}$ coming from the SDE satisfied by $Z$.

Important closed‑loop distinction from the open‑loop setting in \cite{GLZ2025}:
In open‑loop control, the control‑sensitivity $\omega\equiv 0$ and the source terms $\mathfrak{F},\mathfrak{G}$ vanish, so no bounds for second‑order sensitivity processes are required. In our closed‑loop framework, however, $\mathfrak{F}$ and $\mathfrak{G}$ are themselves functions of $Y,Z,v,\omega$. Even though $Z$ no longer appears in the formula for the second‑order functional derivative, these forcing residuals still depend on $Z$. Consequently, we cannot avoid deriving separate moment estimates for $Z$ (Lemma~\ref{l5}), as well as for $v$ and $\omega$ (Lemma~\ref{l6}). The adjoint BSDE then enables us to bound these transformed residuals via standard $L^2$ a‑priori estimates for $(P_{i},Q_{i})$.
\end{remark}
\section{Main results}

\label{sect4}
In this section, we give the estimate of $\alpha$ by means of the first and second-order sensitive processes.
Now we present some results on the estimations of $\mathbf{X}^{\phi }$, $%
\mathbf{Y}^{\phi ,\phi _{\ell }^{\prime \prime }}$ and $\mathbf{Z}^{\phi
,\phi _{h}^{\prime },\phi _{\ell }^{\prime \prime }}.$ The proofs are
scheduled in Appendix \ref{app}.

\begin{lemma}
\label{l2}Suppose that Assumptions \emph{(A1)-(A2) }hold and that there
exists $p\geq 2$ such that $\xi _{i}\in L^{p}\left( \Omega ;\mathbb{R}%
\right) $ for all $i\in I_{N}$. For each $\mathbf{u}\in \mathcal{H}%
^{p}\left( \mathbb{R}\right) ^{N}$, the solution $\mathbf{X}^{\mathbf{u}}\in
\mathcal{H}^{p}\left( \mathbb{R}\right) ^{N}$ to (\ref{sde1}) satisfies for
all $i\in I_{N}$,
\begin{equation}
\sup_{t\in \left[ 0,T\right] }\mathbb{E}\left[ \left\vert X_{t,i}^{\mathbf{u}%
}\right\vert ^{p}\right] \leq C_{X}^{i,p}  \label{l21}
\end{equation}%
with the constant $C_{X}^{i,p}$ defined by
\begin{eqnarray*}
C_{X}^{i,p} &=&\mathbb{E}\Bigg [\left\vert \xi _{i}\right\vert ^{p}+TL^{\bar{%
b}}+\frac{3\left( p-1\right) \left( L^{\bar{\sigma}}\right) ^{2}T}{2} \\
&&+\frac{3\left( p-1\right) \left( L_{y}^{\bar{\sigma}}\right) ^{2}T}{N}%
\sum_{i=1}^{N}\left( \mathbb{E}\left[ \left\vert \xi _{i}\right\vert ^{p}%
\right] +TL^{\bar{b}}+\frac{3\left( p-1\right) \left( L^{\bar{\sigma}%
}\right) ^{2}T}{2}\right) \cdot \\
&&e^{\left( C^{p,\bar{b},\bar{\sigma}}+3\left( p-1\right) \left( L_{y}^{\bar{%
\sigma}}\right) ^{2}\right) T}\Bigg ]e^{TC^{p,\bar{b},\bar{\sigma}}}.
\end{eqnarray*}%

\end{lemma}

\begin{lemma}
\label{l4}Suppose that Assumptions \emph{(A1)-(A2) }hold and that there
exists $p\geq 2$ such that $\xi _{i}\in L^{p}\left( \Omega ;\mathbb{R}%
\right) $ for all $i\in I_{N}$. For each $\mathbf{u}\in \mathcal{H}%
^{p}\left( \mathbb{R}\right) ^{N}$, $h\in I_{N},$ and $u_{h}^{\prime }\in
\mathcal{H}^{p}\left( \mathbb{R}\right) ,$ the solution $\mathbf{Y}^{\mathbf{%
u},u_{h}^{\prime }}\in \mathcal{H}^{p}\left( \mathbb{R}\right) ^{N}$ to (\ref%
{var1}) satisfies for all $i\in I_{N}$
\begin{eqnarray}
\sup_{0\leq t\leq T}\mathbb{E}\left[ \left\vert Y_{t,i}^{h}\right\vert ^{p}%
\right] &\leq &\Bigg [\frac{L_{y}^{\bar{b}}+3\left( p-1\right) \left( L_{y}^{%
\bar{\sigma}}\right) ^{2}}{N}\left( 3p-2\right) TC_{f_{h}}^{p}e^{I_{B,D,\bar{%
B},\bar{D},p}^{3}\cdot T}  \notag \\
&&+\delta _{h,i}\left( 3p-2\right) C_{f_{h}}^{p}\Bigg ]e^{I_{B,D,\bar{B},%
\bar{D},p}^{4}\cdot T},  \label{l41}
\end{eqnarray}%
where%
\begin{eqnarray*}
\bar{b}_{s,i}\left( X_{s,i},\mathbf{X}_{s}^{\phi }\right) &=&\left(
b_{s,i}+u_{s,i}^{\phi }\right) \left( X_{s,i},\mathbf{X}_{s}^{\phi }\right) ,
\\
\bar{\sigma}_{s,i}\left( X_{s,i},\mathbf{X}_{s}^{\phi }\right) &=&\left(
\sigma _{s,i}+u_{s,i}^{\phi }\right) \left( X_{s,i},\mathbf{X}_{s}^{\phi
}\right) .
\end{eqnarray*}%
Additionally, the constant $C_{f_{h}}^{p}$ is defined in (\ref{cfi}) with
\begin{eqnarray*}
I_{B,D,\bar{B},\bar{D},p}^{3} &=&p\left\Vert B\right\Vert _{\infty
}+N\left\Vert \bar{B}\right\Vert _{\infty }p+\left( p-1\right) +\frac{3}{2}%
\left( p-1\right) p\left\Vert D\right\Vert _{\infty }^{2} \\
&&+\frac{3}{2}\left( p-1\right) \left( p-2\right) N^{2}\left\Vert \bar{D}%
\right\Vert _{\infty }^{2}+\frac{3}{2}\left( p-1\right) \left( p-2\right) \\
&&+3\left( p-1\right) N^{2}\left\Vert \bar{D}\right\Vert _{\infty }^{2}, \\
I_{B,D,\bar{B},\bar{D},p}^{4} &=&p\left\Vert B\right\Vert _{\infty
}+N\left\Vert \bar{B}\right\Vert _{\infty }\left( p-1\right) +\left(
p-1\right) +\frac{3}{2}\left( p-1\right) p\left\Vert D\right\Vert _{\infty
}^{2} \\
&&+\frac{3}{2}\left( p-1\right) N^{2}\left\Vert \bar{D}\right\Vert _{\infty
}^{2}\left( p-2\right) +\frac{3}{2}\left( p-1\right) \left( p-2\right) .
\end{eqnarray*}%
Moreover, we have%
\begin{eqnarray}
\left\Vert Y_{i}^{h}Y_{i}^{\ell }\right\Vert _{\mathcal{H}^{2}\left( \mathbb{%
R}\right) }^{2} &\leq &C\Bigg (\delta _{h,i}\delta _{\ell ,j}+\frac{\left(
L_{y}^{\bar{b}}+3\left( L_{y}^{\bar{\sigma}}\right) ^{2}\right) ^{2}}{\sqrt{N%
}}\left( \delta _{h,i}+\delta _{\ell ,j}\right)  \notag \\
&&+\frac{\left( L_{y}^{\bar{b}}+9\left( L_{y}^{\bar{\sigma}}\right)
^{2}\right) ^{4}}{N}\Bigg )\left\Vert u_{h}^{\prime }\right\Vert _{\mathcal{H%
}^{4}\left( \mathbb{R}\right) }\left\Vert u_{\ell }^{\prime \prime
}\right\Vert _{\mathcal{H}^{4}\left( \mathbb{R}\right) }.  \label{l42}
\end{eqnarray}
\end{lemma}

\begin{lemma}
\label{l5}Suppose that Assumptions \emph{(A1)-(A2) }hold and $\xi _{i}\in
L^{4}\left( \Omega ;\mathbb{R}\right) $ for all $i\in I_{N}$. For each $%
\mathbf{u}\in \mathcal{H}^{p}\left( \mathbb{R}\right) ^{N}$, the solution $%
\mathbf{X}^{\mathbf{u}}\in \mathcal{H}^{p}\left( \mathbb{R}\right) ^{N}$ to (%
\ref{sde1}) satisfies for all $i\in I_{N}$,%
\begin{eqnarray}
\sup_{0\leq t\leq T}\mathbb{E}\left[ \left\vert Z_{t,i}^{h,\ell }\right\vert
^{2}\right] &\leq &\frac{C\left( L_{y}^{\bar{b}}\right) ^{2}\left( L_{y}^{%
\bar{b}}+3\left( L_{y}^{\bar{\sigma}}\right) ^{2}\right) ^{\frac{1}{2}}}{N^{%
\frac{1}{2}}}\left( \delta _{h,i}+\delta _{\ell ,i}\right) +\frac{C\left(
L_{y}^{\bar{b}}\right) ^{2}\left( L_{y}^{\bar{b}}+3\left( L_{y}^{\bar{\sigma}%
}\right) ^{2}\right) }{N}  \notag \\
&&+\frac{C\left( L_{y}^{\bar{\sigma}}\right) ^{2}\left( L_{y}^{\bar{b}%
}+3\left( L_{y}^{\bar{\sigma}}\right) ^{2}\right) ^{\frac{1}{2}}}{N^{\frac{1%
}{2}}}\left( \delta _{h,i}+\delta _{\ell ,i}\right) +\frac{C\left( L_{y}^{%
\bar{\sigma}}\right) ^{2}\left( L_{y}^{\bar{b}}+3\left( L_{y}^{\bar{\sigma}%
}\right) ^{2}\right) }{N}.  \label{zes}
\end{eqnarray}
\end{lemma}

\begin{lemma}
\label{l6}Suppose that Assumptions \emph{(A1)-(A2) }hold and $\xi _{i}\in
L^{4}\left( \Omega ;\mathbb{R}\right) $ for all $i\in I_{N}$. For each $%
\mathbf{u}\in \mathcal{H}^{p}\left( \mathbb{R}\right) ^{N}$, the solution $%
\mathbf{X}^{\mathbf{u}}\in \mathcal{H}^{p}\left( \mathbb{R}\right) ^{N}$ to (%
\ref{sde1}) satisfies for all $i\in I_{N}$,%
\begin{eqnarray}
\left\Vert u_{i}^{\phi }\right\Vert _{\mathcal{H}^{2}\left( \mathbb{R}%
\right) }^{2} &\leq &C,\text{ }  \label{l61} \\
\left\Vert \upsilon _{i}^{\phi ,\phi _{h}^{\prime }}\right\Vert _{\mathcal{H}%
^{2}\left( \mathbb{R}\right) }^{2} &\leq &C\left[ \delta _{h,i}+\frac{1}{N}%
\left( \left( L_{y}^{\phi }\right) ^{2}+\left( L^{\phi }\right) ^{2}\right)
\left( L_{y}^{\bar{b}}+3\left( L_{y}^{\bar{\sigma}}\right) ^{2}\right) %
\right] ,  \label{l62} \\
\left\Vert \omega _{i}^{\phi ,\phi _{h}^{\prime },\phi _{\ell }^{\prime
\prime }}\right\Vert _{\mathcal{H}^{2}\left( \mathbb{R}\right) }^{2} &\leq
&\left( \delta _{h,i}+\delta _{\ell ,i}\right) \left( 1+L_{y}^{\bar{b}%
}+L_{y}^{\bar{\sigma}}\right) ^{2}\left( L_{y}^{\bar{b}}+3\left( L_{y}^{\bar{%
\sigma}}\right) ^{2}\right) ^{\frac{1}{2}}\frac{1}{N^{\frac{1}{2}}}  \notag
\\
&&+\left( 1+L_{y}^{\bar{b}}+L_{y}^{\bar{\sigma}}\right) ^{2}\left( L_{y}^{%
\bar{b}}+3\left( L_{y}^{\bar{\sigma}}\right) ^{2}\right) \frac{1}{N}  \notag
\\
&&+\delta _{h,i}\Bigg [\frac{L^{\phi _{h}^{\prime }}\left( L_{y}^{\bar{b}%
}+3\left( L_{y}^{\bar{\sigma}}\right) ^{2}\right) }{N}+L_{y}^{\phi
_{h}^{\prime }}\left( 1+L_{y}^{\bar{b}}+3\left( L_{y}^{\bar{\sigma}}\right)
^{2}\right) \Bigg ]  \notag \\
&&+\delta _{\ell ,i}\Bigg [\frac{L^{\phi _{\ell }^{\prime \prime }}\left(
L_{y}^{\bar{b}}+3\left( L_{y}^{\bar{\sigma}}\right) ^{2}\right) }{N}%
+L_{y}^{\phi _{\ell }^{\prime \prime }}\left( 1+L_{y}^{\bar{b}}+3\left(
L_{y}^{\bar{\sigma}}\right) ^{2}\right) \Bigg ],  \label{l63}
\end{eqnarray}
where $C>0$ is a generic constant depending only on the upper bounds of $T$,
$\max_{i\in I_{N}}\mathbb{E}\left[ \left\vert \xi _{i}\right\vert ^{4}\right]
,$ $L^{b},$ $L^{\sigma },$ $L^{\phi },$ $L^{\phi _{h}^{\prime }},$ $L^{\phi
_{\ell }^{\prime \prime }},$ $L_{y}^{b},$ $L_{y}^{\sigma },$ $L_{y}^{\phi
_{h}^{\prime }},$ $L_{y}^{\phi _{\ell }^{\prime \prime }}$.
\end{lemma}

\begin{theorem}
\label{the1}Suppose that Assumptions \emph{(A1)-(A2)} are in force. Let $%
i,j\in I_{N}$ with $i\neq j,$ and define $\Delta _{i,j}^{f}=f_{i}-f_{j}$ and
$\Delta _{i,j}^{g}=g_{i}-g_{j}$ . Then for all $\phi \in \Pi ^{N}$ and $\phi
_{i}^{\prime },\phi _{j}^{\prime \prime }\in \Pi $,
\begin{equation}
\left\vert \frac{\delta ^{2}V_{i}}{\delta \phi _{i}\delta \phi _{j}}\left(
\phi ;\phi _{i}^{\prime },\phi _{j}^{\prime \prime }\right) -\frac{\delta
^{2}V_{j}}{\delta \phi _{j}\delta \phi _{i}}\left( \phi ;\phi _{j}^{\prime
\prime },\phi _{i}^{\prime }\right) \right\vert \leq C\left( \tilde{C}%
_{0}^{i,j}+\frac{1}{N^{\frac{1}{4}}}\tilde{C}_{1}^{i,j}+\frac{1}{N^{\frac{1}{%
2}}}\tilde{C}_{2}^{i,j}+\frac{1}{N}\tilde{C}_{3}^{i,j}+\frac{1}{N^{2}}\tilde{%
C}_{4}^{i,j}\right) ,  \label{estia}
\end{equation}
with $L_{y}^{\bar{b},\bar{\sigma}}=L_{y}^{\bar{b}}+3\left( L_{y}^{\bar{\sigma%
}}\right) ^{2}$ and
\begin{eqnarray*}
\tilde{C}_{0}^{i,j} &=&\left\Vert \partial _{x_{i}x_{j}}^{2}\Delta
_{i,j}^{f}\right\Vert _{L^{\infty }}+\left\Vert \partial
_{x_{i}u_{j}}^{2}\Delta _{i,j}^{f}\left( \mathbf{X}_{t},\mathbf{u}%
_{t}\right) \right\Vert _{L^{\infty }}+\left\Vert \partial
_{u_{i}x_{j}}^{2}\Delta _{i,j}^{f}\left( \mathbf{X}_{t},\mathbf{u}%
_{t}\right) \right\Vert _{L^{\infty }}+\left\Vert \left( \partial
_{u_{i}u_{j}}^{2}\Delta _{i,j}^{f}\left( \mathbf{X}_{t},\mathbf{u}%
_{t}\right) \right) \right\Vert _{L^{\infty }} \\
&&+\left\Vert \partial _{x_{i}x_{j}}^{2}\Delta _{i,j}^{g}\right\Vert
_{L^{\infty }}+\bigg [\sqrt{L_{y}^{\phi _{i}^{\prime }}\left( 1+L_{y}^{\bar{b%
},\bar{\sigma}}\right) }+\sqrt{L_{y}^{\phi _{j}^{\prime \prime }}\left(
1+L_{y}^{\bar{b},\bar{\sigma}}\right) }\bigg ]\cdot \sum_{h\in \left\{
i,j\right\} }^{N}\Big [\left\Vert \left( \partial _{u_{\ell }}\Delta
_{i,j}^{f}\right) \left( 0,0\right) \right\Vert _{L^{\infty }} \\
&&+\sum_{k=1}^{N}\left( \left\Vert \partial _{u_{_{h}}x_{k}}\Delta
_{i,j}^{f}\right\Vert _{L^{\infty }}+\left\Vert \partial _{u_{h}u_{k}}\Delta
_{i,j}^{f}\right\Vert _{L^{\infty }}\right) \Big ], \\
\tilde{C}_{1}^{i,j} &=&\left( 1+L_{y}^{\bar{b}}+L_{y}^{\bar{\sigma}}\right)
\left( L_{y}^{\bar{b},\bar{\sigma}}\right) ^{\frac{1}{4}}\cdot \sum_{h\in
\left\{ i,j\right\} }^{N}\Big [\left\Vert \left( \partial _{u_{h}}\Delta
_{i,j}^{f}\right) \left( 0,0\right) \right\Vert _{L^{\infty
}}+\sum_{k=1}^{N}\left( \left\Vert \partial _{u_{_{h}}x_{k}}\Delta
_{i,j}^{f}\right\Vert _{L^{\infty }}+\left\Vert \partial _{u_{h}u_{k}}\Delta
_{i,j}^{f}\right\Vert _{L^{\infty }}\right) \Big ], \\
\tilde{C}_{2}^{i,j} &=&\left( 1+L_{y}^{\bar{b}}+L_{y}^{\bar{\sigma}}\right)
\left( L_{y}^{\bar{b},\bar{\sigma}}\right) ^{\frac{1}{2}}\sum_{h\in
I_{N}\backslash \left\{ i,j\right\} }^{N}\Big [\left\Vert \left( \partial
_{u_{h}}\Delta _{i,j}^{f}\right) \left( 0,0\right) \right\Vert _{L^{\infty
}}+\sum_{k=1}^{N}\left( \left\Vert \partial _{u_{_{h}}x_{k}}\Delta
_{i,j}^{f}\right\Vert _{L^{\infty }}+\left\Vert \partial _{u_{h}u_{k}}\Delta
_{i,j}^{f}\right\Vert _{L^{\infty }}\right) \Big ] \\
&&\left( L_{y}^{\bar{b},\bar{\sigma}}\right) ^{\frac{1}{2}}\left( \sum_{\ell
\in I_{N}\backslash \left\{ j\right\} }^{N}\left\Vert \partial
_{x_{i}x_{\ell }}^{2}\Delta _{t,i,j}^{f}\right\Vert _{L^{\infty
}}+\sum_{h\in I_{N}\backslash \left\{ i\right\} }^{N}\left\Vert \partial
_{x_{h}x_{j}}^{2}\Delta _{i,j}^{f}\right\Vert _{L^{\infty }}\right)  \\
&&+\left( L_{y}^{\bar{b},\bar{\sigma}}\right) ^{\frac{1}{2}}\sum_{h\in
I_{N}\backslash \left\{ i\right\} }^{N}\left\Vert \partial
_{x_{h}u_{j}}^{2}\Delta _{t,i,j}^{f}\left( \mathbf{X}_{t},\mathbf{u}%
_{t}\right) \right\Vert _{L^{\infty }}+\left[ L_{y}^{\bar{b},\bar{\sigma}%
}\left( \left( L_{y}^{\phi }\right) ^{2}+\left( L^{\phi }\right) ^{2}\right) %
\right] ^{\frac{1}{2}}\sum_{\ell \in I_{N}\backslash \left\{ j\right\}
}^{N}\left\Vert \partial _{x_{i}u_{\ell }}^{2}\Delta _{t,i,j}^{f}\left(
\mathbf{X}_{t},\mathbf{u}_{t}\right) \right\Vert _{L^{\infty }} \\
&&+\left( L_{y}^{\bar{b},\bar{\sigma}}\right) ^{\frac{1}{2}}\sum_{\ell \in
I_{N}\backslash \left\{ j\right\} }^{N}\left\Vert \partial _{u_{i}x_{\ell
}}^{2}\Delta _{t,i,j}^{f}\left( \mathbf{X}_{t},\mathbf{u}_{t}\right)
\right\Vert _{L^{\infty }}+\left[ L_{y}^{\bar{b},\bar{\sigma}}\left( \left(
L_{y}^{\phi }\right) ^{2}+\left( L^{\phi }\right) ^{2}\right) \right] ^{%
\frac{1}{2}}\sum_{h\in I_{N}\backslash \left\{ i\right\} }^{N}\left\Vert
\partial _{u_{h}x_{j}}^{2}\Delta _{t,i,j}^{f}\left( \mathbf{X}_{t},\mathbf{u}%
_{t}\right) \right\Vert _{L^{\infty }} \\
&&+\left( \left( L_{y}^{\phi }\right) ^{2}+\left( L^{\phi }\right)
^{2}\right) ^{\frac{1}{2}}\left( L_{y}^{\bar{b},\bar{\sigma}}\right) ^{\frac{%
1}{2}}\Bigg (\sum_{\ell \in I_{N}\backslash \left\{ j\right\}
}^{N}\left\Vert \partial _{u_{i}u_{\ell }}^{2}\Delta _{t,i,j}^{f}\left(
\mathbf{X}_{t},\mathbf{u}_{t}\right) \right\Vert _{L^{\infty }}+\sum_{h\in
I_{N}\backslash \left\{ i\right\} }^{N}\left\Vert \partial
_{u_{h}u_{j}}^{2}\Delta _{i,j}^{f}\left( \mathbf{X}_{t},\mathbf{u}%
_{t}\right) \right\Vert _{L^{\infty }}\Bigg ), \\
\tilde{C}_{3}^{i,j} &=&L_{y}^{\bar{b},\bar{\sigma}}\sum_{h\in
I_{N}\backslash \left\{ i\right\} }^{N}\sum_{\ell \in I_{N}\backslash
\left\{ j\right\} }^{N}\left\Vert \partial _{x_{h}x_{\ell }}^{2}\Delta
_{i,j}^{f}\right\Vert _{L^{\infty }}+L_{y}^{\bar{b},\bar{\sigma}}\left(
\left( L_{y}^{\phi }\right) ^{2}+\left( L^{\phi }\right) ^{2}\right) ^{\frac{%
1}{2}}\sum_{\substack{ h\in I_{N}\backslash \left\{ i\right\}  \\ \ell \in
I_{N}\backslash \left\{ j\right\} }}^{N}\left\Vert \left( \partial
_{x_{h}u_{\ell }}^{2}\Delta _{i,j}^{f}\left( \mathbf{X}_{t},\mathbf{u}%
_{t}\right) \right) \right\Vert _{L^{\infty }} \\
&&+L_{y}^{\bar{b},\bar{\sigma}}\left( \left( L_{y}^{\phi }\right)
^{2}+\left( L^{\phi }\right) ^{2}\right) ^{\frac{1}{2}}\sum_{\substack{ h\in
I_{N}\backslash \left\{ i\right\}  \\ \ell \in I_{N}\backslash \left\{
j\right\} }}^{N}\left\Vert \left( \partial _{u_{h}x_{\ell }}^{2}\Delta
_{i,j}^{f}\left( \mathbf{X}_{t},\mathbf{u}_{t}\right) \right) \right\Vert
_{L^{\infty }} \\
&&+L_{y}^{\bar{b},\bar{\sigma}}\left( \left( L_{y}^{\phi }\right)
^{2}+\left( L^{\phi }\right) ^{2}\right) \sum_{\substack{ h\in
I_{N}\backslash \left\{ i\right\}  \\ \ell \in I_{N}\backslash \left\{
j\right\} }}^{N}\left\Vert \partial _{u_{h}u_{\ell }}^{2}\Delta
_{t,i,j}^{f}\left( \mathbf{X}_{t},\mathbf{u}_{t}\right) \right\Vert
_{L^{\infty }} \\
&&+L_{y}^{\bar{b},\bar{\sigma}}\left( \sum_{\ell \in I_{N}\backslash \left\{
j\right\} }\left\Vert \partial _{x_{i}x_{\ell }}^{2}\Delta
g_{i,j}\right\Vert _{L^{\infty }}+\sum_{h\in I_{N}\backslash \left\{
i\right\} }\left\Vert \partial _{x_{h}x_{j}}^{2}\Delta g_{i,j}\right\Vert
_{L^{\infty }}\right)  \\
&&+\sqrt{\Gamma _{1}}\Bigg [\left( L^{\bar{b}}+L^{\bar{\sigma}}\right)
\left( L_{y}^{\bar{b},\bar{\sigma}}\right) ^{2}+2\left( L_{y}^{\bar{b}%
}+L_{y}^{\bar{\sigma}}\right) \left( 1+L_{y}^{\bar{b},\bar{\sigma}}+\left(
L_{y}^{\bar{b},\bar{\sigma}}\right) ^{2}\right) +2L_{y}^{\bar{b},\bar{\sigma}%
}\left( L^{\phi _{i}^{\prime }}+L_{y}^{\phi _{i}^{\prime }}+L^{\phi
_{j}^{\prime \prime }}+L_{y}^{\phi _{j}^{\prime \prime }}\right) \Bigg ] \\
&&+\sqrt{\Gamma _{1}}\left( L_{y}^{\bar{b}}+L_{y}^{\bar{\sigma}}\right)
\left( 1+L_{y}^{\bar{b},\bar{\sigma}}+\left( L_{y}^{\bar{b},\bar{\sigma}%
}\right) ^{2}\right) , \\
\tilde{C}_{4}^{i,j} &=&\left( L_{y}^{\bar{b},\bar{\sigma}}\right)
^{2}\sum_{\ell \in I_{N}\backslash \left\{ j\right\} ,h\in I_{N}\backslash
\left\{ i\right\} }\left\Vert \partial _{x_{h}x_{\ell }}^{2}\Delta
g_{i,j}\right\Vert _{L^{\infty }}+\sqrt{\Gamma _{1}}\left( L_{y}^{\bar{b},%
\bar{\sigma}}\right) ^{2}\left( L_{y}^{\bar{b}}+L_{y}^{\bar{\sigma}}\right) ,
\end{eqnarray*}

where $C$ depends on the time $T,$ $\max \left\{ \left\vert \mathbf{B}\left(
t,\mathbf{X}_{t}^{\mathbf{u}},\mathbf{u}_{t}\right) \right\vert ,\left\vert
\mathbf{\Pi }^{1}\left( t,\mathbf{X}_{t}^{\mathbf{u}},\mathbf{u}_{t}\right)
\right\vert ,\ldots ,\left\vert \mathbf{\Pi }^{N}\left( t,\mathbf{X}_{t}^{%
\mathbf{u}},\mathbf{u}_{t}\right) \right\vert \right\} ,$ and while%
\begin{eqnarray*}
\Gamma _{1} &=&\mathbb{E}\Bigg [\sum_{\ell \in I_{N}}\left\vert \left(
\partial _{x_{\ell }}\Delta _{i,j}^{g}\right) \left( 0\right) \right\vert
^{2}+\sum_{\ell ,k\in I_{N}}\left\Vert \partial _{x_{\ell }x_{k}}\Delta
_{i,j}^{g}\right\Vert _{L^{\infty }}^{2}+\sum_{\ell ,k\in I_{N}}\left\Vert
\partial _{x_{\ell }x_{k}}\Delta _{i,j}^{f}\right\Vert _{L^{\infty
}}^{2}+\sum_{\ell ,k\in I_{N}}\left\Vert \partial _{x_{\ell }u_{k}}\Delta
_{i,j}^{f}\right\Vert _{L^{\infty }}^{2} \\
&&+3\int_{0}^{T}\Bigg (\sum_{\ell }\left\vert \left( \partial _{x_{\ell
}}\Delta _{i,j}^{f}\right) \left( t,0,0\right) \right\vert ^{2}\Bigg )%
\mathrm{d}t\Bigg ].
\end{eqnarray*}
\end{theorem}

\begin{remark}
If $b_{i}\left( t,X,\mathbf{X},\mathbf{u}\right) =\bar{b}_{i}\left( t,X,%
\mathbf{X}\right) +u,$ $\sigma _{i}\left( t,X,\mathbf{X},\mathbf{u}\right)
=\sigma _{i}\left( t\right) $ considered in Guo et al. \cite{GLZ1}, then $L^{%
\bar{\sigma}}=L_{y}^{\bar{\sigma}}=0$, and $L_{y}^{\bar{b},\bar{\sigma}%
}=L_{y}^{\bar{b}}.$ We simplify the estimation of $\alpha .$ Indeed, we have
\begin{eqnarray*}
\tilde{C}_{0}^{i,j} &=&\left\Vert \partial _{x_{i}x_{j}}^{2}\Delta
_{i,j}^{f}\right\Vert _{L^{\infty }}+\left\Vert \partial
_{x_{i}u_{j}}^{2}\Delta _{i,j}^{f}\left( \mathbf{X}_{t},\mathbf{u}%
_{t}\right) \right\Vert _{L^{\infty }}+\left\Vert \partial
_{u_{i}x_{j}}^{2}\Delta _{i,j}^{f}\left( \mathbf{X}_{t},\mathbf{u}%
_{t}\right) \right\Vert _{L^{\infty }}+\left\Vert \left( \partial
_{u_{i}u_{j}}^{2}\Delta _{i,j}^{f}\left( \mathbf{X}_{t},\mathbf{u}%
_{t}\right) \right) \right\Vert _{L^{\infty }} \\
&&+\left\Vert \partial _{x_{i}x_{j}}^{2}\Delta _{i,j}^{g}\right\Vert
_{L^{\infty }}+\bigg [\sqrt{L_{y}^{\phi _{i}^{\prime }}\left( 1+L_{y}^{\bar{b%
}}\right) }+\sqrt{L_{y}^{\phi _{j}^{\prime \prime }}\left( 1+L_{y}^{\bar{b}%
}\right) }\bigg ]\cdot \sum_{h\in \left\{ i,j\right\} }^{N}\Big [\left\Vert
\left( \partial _{u_{\ell }}\Delta _{i,j}^{f}\right) \left( 0,0\right)
\right\Vert _{L^{\infty }} \\
&&+\sum_{k=1}^{N}\left( \left\Vert \partial _{u_{_{h}}x_{k}}\Delta
_{i,j}^{f}\right\Vert _{L^{\infty }}+\left\Vert \partial _{u_{h}u_{k}}\Delta
_{i,j}^{f}\right\Vert _{L^{\infty }}\right) \Big ], \\
\tilde{C}_{1}^{i,j} &=&\left( 1+L_{y}^{\bar{b}}\right) \left( L_{y}^{\bar{b}%
}\right) ^{\frac{1}{4}}\cdot \sum_{h\in \left\{ i,j\right\} }^{N}\Big [%
\left\Vert \left( \partial _{u_{\ell }}\Delta _{i,j}^{f}\right) \left(
0,0\right) \right\Vert _{L^{\infty }}+\sum_{k=1}^{N}\left( \left\Vert
\partial _{u_{_{h}}x_{k}}\Delta _{i,j}^{f}\right\Vert _{L^{\infty
}}+\left\Vert \partial _{u_{h}u_{k}}\Delta _{i,j}^{f}\right\Vert _{L^{\infty
}}\right) \Big ], \\
\tilde{C}_{2}^{i,j} &=&\left( 1+L_{y}^{\bar{b}}\right) \left( L_{y}^{\bar{b}%
}\right) ^{\frac{1}{2}}\sum_{h\in I_{N}\backslash \left\{ i,j\right\} }^{N}%
\Big [\left\Vert \left( \partial _{u_{\ell }}\Delta _{i,j}^{f}\right) \left(
0,0\right) \right\Vert _{L^{\infty }}+\sum_{k=1}^{N}\left( \left\Vert
\partial _{u_{_{h}}x_{k}}\Delta _{i,j}^{f}\right\Vert _{L^{\infty
}}+\left\Vert \partial _{u_{h}u_{k}}\Delta _{i,j}^{f}\right\Vert _{L^{\infty
}}\right) \Big ] \\
&&+\left( L_{y}^{\bar{b}}\right) ^{\frac{1}{2}}\left( \sum_{\ell \in
I_{N}\backslash \left\{ j\right\} }^{N}\left\Vert \partial _{x_{i}x_{\ell
}}^{2}\Delta _{t,i,j}^{f}\right\Vert _{L^{\infty }}+\sum_{h\in
I_{N}\backslash \left\{ i\right\} }^{N}\left\Vert \partial
_{x_{h}x_{j}}^{2}\Delta _{i,j}^{f}\right\Vert _{L^{\infty }}\right) +\left(
L_{y}^{\bar{b}}\right) ^{\frac{1}{2}}\sum_{h\in I_{N}\backslash \left\{
i\right\} }^{N}\left\Vert \partial _{x_{h}u_{j}}^{2}\Delta
_{t,i,j}^{f}\left( \mathbf{X}_{t},\mathbf{u}_{t}\right) \right\Vert
_{L^{\infty }} \\
&&+\left[ L_{y}^{\bar{b}}\left( \left( L_{y}^{\phi }\right) ^{2}+\left(
L^{\phi }\right) ^{2}\right) \right] ^{\frac{1}{2}}\sum_{\ell \in
I_{N}\backslash \left\{ j\right\} }^{N}\left\Vert \partial _{x_{i}u_{\ell
}}^{2}\Delta _{t,i,j}^{f}\left( \mathbf{X}_{t},\mathbf{u}_{t}\right)
\right\Vert _{L^{\infty }}+\left( L_{y}^{\bar{b}}\right) ^{\frac{1}{2}%
}\sum_{\ell \in I_{N}\backslash \left\{ j\right\} }^{N}\left\Vert \partial
_{u_{i}x_{\ell }}^{2}\Delta _{t,i,j}^{f}\left( \mathbf{X}_{t},\mathbf{u}%
_{t}\right) \right\Vert _{L^{\infty }} \\
&&+\left[ L_{y}^{\bar{b}}\left( \left( L_{y}^{\phi }\right) ^{2}+\left(
L^{\phi }\right) ^{2}\right) \right] ^{\frac{1}{2}}\sum_{h\in
I_{N}\backslash \left\{ i\right\} }^{N}\left\Vert \partial
_{u_{h}x_{j}}^{2}\Delta _{t,i,j}^{f}\left( \mathbf{X}_{t},\mathbf{u}%
_{t}\right) \right\Vert _{L^{\infty }} \\
&&+\left( \left( L_{y}^{\phi }\right) ^{2}+\left( L^{\phi }\right)
^{2}\right) ^{\frac{1}{2}}\left( L_{y}^{\bar{b}}\right) ^{\frac{1}{2}}\Bigg (%
\sum_{\ell \in I_{N}\backslash \left\{ j\right\} }^{N}\left\Vert \partial
_{u_{i}u_{\ell }}^{2}\Delta _{t,i,j}^{f}\left( \mathbf{X}_{t},\mathbf{u}%
_{t}\right) \right\Vert _{L^{\infty }}+\sum_{h\in I_{N}\backslash \left\{
i\right\} }^{N}\left\Vert \partial _{u_{h}u_{j}}^{2}\Delta _{i,j}^{f}\left(
\mathbf{X}_{t},\mathbf{u}_{t}\right) \right\Vert _{L^{\infty }}\Bigg ), \\
\tilde{C}_{3}^{i,j} &=&L_{y}^{\bar{b}}\sum_{h\in I_{N}\backslash \left\{
i\right\} }^{N}\sum_{\ell \in I_{N}\backslash \left\{ j\right\}
}^{N}\left\Vert \partial _{x_{h}x_{\ell }}^{2}\Delta _{i,j}^{f}\right\Vert
_{L^{\infty }}+L_{y}^{\bar{b}}\left( \left( L_{y}^{\phi }\right) ^{2}+\left(
L^{\phi }\right) ^{2}\right) ^{\frac{1}{2}}\sum_{\substack{ h\in
I_{N}\backslash \left\{ i\right\}  \\ \ell \in I_{N}\backslash \left\{
j\right\} }}^{N}\left\Vert \left( \partial _{x_{h}u_{\ell }}^{2}\Delta
_{i,j}^{f}\left( \mathbf{X}_{t},\mathbf{u}_{t}\right) \right) \right\Vert
_{L^{\infty }} \\
&&+L_{y}^{\bar{b}}\left( \left( L_{y}^{\phi }\right) ^{2}+\left( L^{\phi
}\right) ^{2}\right) ^{\frac{1}{2}}\sum_{\substack{ h\in I_{N}\backslash
\left\{ i\right\}  \\ \ell \in I_{N}\backslash \left\{ j\right\} }}%
^{N}\left\Vert \left( \partial _{u_{h}x_{\ell }}^{2}\Delta _{i,j}^{f}\left(
\mathbf{X}_{t},\mathbf{u}_{t}\right) \right) \right\Vert _{L^{\infty }} \\
&&+\left( \left( L_{y}^{\phi }\right) ^{2}+\left( L^{\phi }\right)
^{2}\right) L_{y}^{\bar{b}}\sum_{\substack{ h\in I_{N}\backslash \left\{
i\right\}  \\ \ell \in I_{N}\backslash \left\{ j\right\} }}^{N}\left\Vert
\partial _{u_{h}u_{\ell }}^{2}\Delta _{t,i,j}^{f}\left( \mathbf{X}_{t},%
\mathbf{u}_{t}\right) \right\Vert _{L^{\infty }} \\
&&+L_{y}^{\bar{b}}\left( \sum_{\ell \in I_{N}\backslash \left\{ j\right\}
}\left\Vert \partial _{x_{i}x_{\ell }}^{2}\Delta g_{i,j}\right\Vert
_{L^{\infty }}+\sum_{h\in I_{N}\backslash \left\{ i\right\} }\left\Vert
\partial _{x_{h}x_{j}}^{2}\Delta g_{i,j}\right\Vert _{L^{\infty }}\right) +%
\sqrt{\Gamma _{1}}L_{y}^{\bar{b}}\left( 1+L_{y}^{\bar{b}}+\left( L_{y}^{\bar{%
b}}\right) ^{2}\right)  \\
&&+\sqrt{\Gamma _{1}}\Bigg [L^{\bar{b}}\left( L_{y}^{\bar{b}}\right)
^{2}+2L_{y}^{\bar{b}}\left( 1+L_{y}^{\bar{b}}+\left( L_{y}^{\bar{b}}\right)
^{2}\right) +2L_{y}^{\bar{b}}\left( L^{\phi _{i}^{\prime }}+L_{y}^{\phi
_{i}^{\prime }}+L^{\phi _{j}^{\prime \prime }}+L_{y}^{\phi _{j}^{\prime
\prime }}\right) \Bigg ], \\
\tilde{C}_{4}^{i,j} &=&\left( L_{y}^{\bar{b}}\right) ^{2}\sum_{\ell \in
I_{N}\backslash \left\{ j\right\} ,h\in I_{N}\backslash \left\{ i\right\}
}\left\Vert \partial _{x_{h}x_{\ell }}^{2}\Delta g_{i,j}\right\Vert
_{L^{\infty }}+\sqrt{\Gamma _{1}}\left( L_{y}^{\bar{b}}\right) ^{2}L_{y}^{%
\bar{b}}.
\end{eqnarray*}
where $\Lambda _{1}$ is defined in (\ref{lamda}).
\end{remark}

\subsection{Mean-field type interactions}

\label{gmi}

In this subsection, we develop an analogous result to Theorem \ref{the1} for
closed-loop games involving mean-field type interactions. Owing to the extra
coupling induced by the closed-loop controls, more stringent conditions on
the cost functions are necessary to guarantee that the $N$-player game
constitutes an $\alpha $-potential game with a decaying property.
Specifically, supplementary conditions on the partial derivatives of ($f_{i}$%
)$_{i\in I_{N}}$ with respect to $u$ have been introduced in (\ref{fi1}) to address the interdependence among all players' control processes.

\begin{theorem}
\label{the2}Under Assumptions \emph{(A1)-(A2)}, assume that there exist constants $L,$
$\tilde{L}>0$ and $\beta >1/2$ such that $\sup_{i\in I_{N},\phi _{i}\in
\mathcal{A}_{i}}\leq L$, $\max_{i\in I_{N}}\left\vert \xi _{i}\right\vert
^{4}\leq L,$ $\max_{i\in I_{N}}\left( L^{b_{i}}+L_{y}^{b_{i}}+L^{\sigma
_{i}}+L_{y}^{\sigma _{i}}\right) \leq L$, and for all $i,j\in I_{N}$, $%
\Delta _{t,i,j}^{f}:=f_{t,i}-f_{t,j}$ and $\Delta
_{t,i,j}^{g}:=g_{t,i}-g_{t,j}$ satisfy for all $(t,x,u)\in \lbrack
0;T]\times \mathbb{R}^{N}\times \mathbb{R}^{N}$ and $h,\ell \in I_{N}$,%
\begin{equation}
\begin{array}{l}
\left\vert \left( \partial _{x_{h}}\Delta _{t,i,j}^{f}\right) \left(
0,0\right) \right\vert +\left\vert \left( \partial _{u_{h}}\Delta
_{t,i,j}^{f}\right) \left( 0,0\right) \right\vert \leq \mathbbm{1}_{h\in
\left\{ i,j\right\} }L+\mathbbm{1}_{h\in I_{N}\backslash \left\{ i,j\right\}
}\tilde{L}N^{-\beta }, \\
\left\vert \left( \partial _{x_{h}x_{\ell }}^{2}\Delta _{t,i,j}^{f}\right)
\left( x,u\right) \right\vert +\left\vert \left( \partial _{u_{h}u_{\ell
}}^{2}\Delta _{t,i,j}^{f}\right) \left( t,x,u\right) \right\vert +\left\vert
\left( \partial _{x_{h}u_{\ell }}^{2}\Delta _{t,i,j}^{f}\right) \left(
t,x,u\right) \right\vert  \\
\leq \mathbbm{1}_{h=\ell \in \left\{ i,j\right\} }L+\tilde{L}\left( %
\mathbbm{1}_{h=\ell \in I_{N}\backslash \left\{ i,j\right\} }N^{-\beta }+%
\mathbbm{1}_{h\neq \ell }N^{-2\beta }\right) , \\
\left\vert \left( \partial _{x_{h}}\Delta _{i,j}^{g}\right) \left( 0\right)
\right\vert \leq \tilde{L}N^{-\beta },\text{ }\left\vert \left( \partial
_{x_{h}x_{\ell }}^{2}\Delta _{i,j}^{g}\right) \left( x\right) \right\vert
\leq \tilde{L}\left( \mathbbm{1}_{h=\ell }N^{-\beta }+\mathbbm{1}_{h\neq
\ell }N^{-2\beta }\right) .%
\end{array}
\label{fi1}
\end{equation}%
Then $\mathcal{G}^{cl}$ is an $\alpha $-potential game with
\begin{equation}
\alpha \leq \bar{C}_{0}^{i,j}+\bar{C}_{1}^{i,j}\frac{1}{N^{\frac{1}{4}}}+%
\bar{C}_{2}^{i,j}\frac{1}{N^{\frac{1}{2}}}+\bar{C}_{3}^{i,j}\frac{1}{N}+\bar{%
C}_{4}^{i,j}\frac{1}{N^{2}}  \label{amf}
\end{equation}%
with%
\begin{eqnarray*}
\bar{C}_{0}^{i,j} &=&C\left( \sqrt{L_{y}^{\phi _{h}^{\prime }}\left(
1+L_{y}^{\bar{b}}+3\left( L_{y}^{\bar{\sigma}}\right) ^{2}\right) }+\sqrt{%
L_{y}^{\phi _{\ell }^{\prime \prime }}\left( 1+L_{y}^{\bar{b}}+3\left(
L_{y}^{\bar{\sigma}}\right) ^{2}\right) }\right) \left( 2L+\left( N-1\right)
\frac{\tilde{L}}{N^{2\beta }}\right) +\frac{C\tilde{L}}{N^{2\beta }}, \\
\bar{C}_{1}^{i,j} &=&C\left( 1+L_{y}^{\bar{b}}+L_{y}^{\bar{\sigma}}\right)
\left( L_{y}^{\bar{b},\bar{\sigma}}\right) ^{\frac{1}{4}}\left[ L+\frac{2%
\tilde{L}}{N^{\beta }}+\left( N-1\right) \frac{\tilde{L}}{N^{2\beta }}\right]
, \\
\bar{C}_{2}^{i,j} &=&\left( 1+L_{y}^{\bar{b}}+L_{y}^{\bar{\sigma}}\right)
\left( L_{y}^{\bar{b},\bar{\sigma}}\right) ^{\frac{1}{2}}\left[ \left(
N-2\right) \frac{\tilde{L}}{N^{\beta }}+\left( N-2\right) \left( \frac{%
\tilde{L}}{N^{\beta }}+\left( N-1\right) \frac{\tilde{L}}{N^{2\beta }}%
\right) \right]  \\
&&+2\left( L_{y}^{\bar{b},\bar{\sigma}}\right) ^{\frac{1}{2}}\left(
N-1\right) \frac{\tilde{L}}{N^{2\beta }}+\left( L_{y}^{\bar{b},\bar{\sigma}%
}\right) ^{\frac{1}{2}}\left[ L+\left( N-2\right) \frac{\tilde{L}}{N^{2\beta
}}\right]  \\
&&+\left[ L_{y}^{\bar{b},\bar{\sigma}}\left( \left( L_{y}^{\phi }\right)
^{2}+\left( L^{\phi }\right) ^{2}\right) \right] ^{\frac{1}{2}}\left[
L+\left( N-2\right) \frac{\tilde{L}}{N^{2\beta }}\right]  \\
&&+\left( L_{y}^{\bar{b},\bar{\sigma}}\right) ^{\frac{1}{2}}\left[ L+\left(
N-2\right) \frac{\tilde{L}}{N^{2\beta }}\right] +\left[ L_{y}^{\bar{b},\bar{%
\sigma}}\left( \left( L_{y}^{\phi }\right) ^{2}+\left( L^{\phi }\right)
^{2}\right) \right] ^{\frac{1}{2}}\left[ L+\left( N-2\right) \frac{\tilde{L}%
}{N^{2\beta }}\right]  \\
&&+2\left( \left( L_{y}^{\phi }\right) ^{2}+\left( L^{\phi }\right)
^{2}\right) ^{\frac{1}{2}}\left( L_{y}^{\bar{b},\bar{\sigma}}\right) ^{\frac{%
1}{2}}\left[ L+\left( N-2\right) \frac{\tilde{L}}{N^{2\beta }}\right] , \\
\bar{C}_{3}^{i,j} &=&2L_{y}^{\bar{b},\bar{\sigma}}\left[ \left( N-2\right)
\frac{\tilde{L}}{N^{2\beta }}+\frac{\tilde{L}}{N^{\beta }}\right] +\Bigg (%
L_{y}^{\bar{b},\bar{\sigma}}+L_{y}^{\bar{b},\bar{\sigma}}\left( \left(
L_{y}^{\phi }\right) ^{2}+\left( L^{\phi }\right) ^{2}\right) ^{\frac{1}{2}%
}+L_{y}^{\bar{b},\bar{\sigma}}\left( \left( L_{y}^{\phi }\right) ^{2}+\left(
L^{\phi }\right) ^{2}\right) ^{\frac{1}{2}} \\
&&+\left( \left( L_{y}^{\phi }\right) ^{2}+\left( L^{\phi }\right)
^{2}\right) \left( L_{y}^{\bar{b}}+3\left( L_{y}^{\bar{\sigma}}\right)
^{2}\right) \Bigg )\left[ \left( \left( N-1\right) ^{2}-\left( N-2\right)
\right) \frac{\tilde{L}}{N^{2\beta }}+\left( N-2\right) \frac{\tilde{L}}{%
N^{\beta }}\right]  \\
&&+\sqrt{\Gamma _{1}}\Bigg [\left( L^{\bar{b}}+L^{\bar{\sigma}}\right)
\left( L_{y}^{\bar{b},\bar{\sigma}}\right) ^{2}+2\left( L_{y}^{\bar{b}%
}+L_{y}^{\bar{\sigma}}\right) \left( 1+L_{y}^{\bar{b},\bar{\sigma}}+\left(
L_{y}^{\bar{b},\bar{\sigma}}\right) ^{2}\right) +2L_{y}^{\bar{b},\bar{\sigma}%
}\left( L^{\phi _{i}^{\prime }}+L_{y}^{\phi _{i}^{\prime }}+L^{\phi
_{j}^{\prime \prime }}+L_{y}^{\phi _{j}^{\prime \prime }}\right) \Bigg ] \\
&&+\sqrt{\Gamma _{1}}\left( L_{y}^{\bar{b}}+L_{y}^{\bar{\sigma}}\right)
\left( 1+L_{y}^{\bar{b},\bar{\sigma}}+\left( L_{y}^{\bar{b},\bar{\sigma}%
}\right) ^{2}\right) , \\
\bar{C}_{4}^{i,j} &=&\left( L_{y}^{\bar{b},\bar{\sigma}}\right) ^{2}\left[
\left( \left( N-1\right) ^{2}-\left( N-2\right) \right) \frac{\tilde{L}}{%
N^{2\beta }}+\left( N-2\right) \frac{\tilde{L}}{N^{\beta }}\right] +\sqrt{%
\Gamma _{1}}\left( L_{y}^{\bar{b},\bar{\sigma}}\right) ^{2}\left( L_{y}^{%
\bar{b}}+L_{y}^{\bar{\sigma}}\right),
\end{eqnarray*}
where
\begin{eqnarray*}
\Gamma _{1} &\leq &\left( N-1\right) \frac{\tilde{L}}{N^{\beta }}+N\cdot
\frac{\tilde{L}}{N^{\beta }}+\left( N^{2}-N\right) \frac{\tilde{L}}{%
N^{2\beta }} \\
&&+2L+\left( N-2\right) \frac{\tilde{L}}{N^{\beta }}+\left( N^{2}-N\right)
\frac{\tilde{L}}{N^{2\beta }}+6TL+3\left( N-2\right) T\frac{\tilde{L}}{%
N^{\beta }}.
\end{eqnarray*}
for a positive constant $C$ independent of $N.$
\end{theorem}

The proof is arranged in Appendix \ref{APP}. It relies on upper bounding the
constants $C_{k}^{i,j},$ $k=0,\ldots ,4$ in Theorem \ref{the1} with respect
to $\tilde{L},\beta $ and $N.$ Clearly, letting $N\rightarrow \infty $ in (\ref{amf}), we have%
\begin{equation}
0\leq \alpha \leq C2L\left( \sqrt{L_{y}^{\phi _{h}^{\prime }}\left( 1+L_{y}^{%
\bar{b}}+3\left( L_{y}^{\bar{\sigma}}\right) ^{2}\right) }+\sqrt{L_{y}^{\phi
_{\ell }^{\prime \prime }}\left( 1+L_{y}^{\bar{b}}+3\left( L_{y}^{\bar{\sigma%
}}\right) ^{2}\right) }\right) .  \label{anot}
\end{equation}
\begin{remark}
    Theorem \ref{the2} reveals a structural distinction for mean‑field‑interaction closed‑loop stochastic differential games. When heterogeneity between agents exists, the feedback‑induced sensitivity processes \(v^{\phi,\phi_{h}'}\) and \(\omega^{\phi,\phi_{h}',\phi_{\ell}''}\) contribute positive contributions to our upper bound for \(\alpha\). Consequently, \(\alpha\) may converge to a strictly positive constant as \(N\to\infty\), rather than decaying to zero. If this limiting is positive, i.e. \(\alpha>0\) , this is not a finite‑N numerical artefact; it reflects persistent strategic externalities originating purely from closed‑loop state‑feedback information. This stands in sharp contrast to the corresponding open‑loop setting from Guo et al. (2025)~\cite{GLZ2025}, where \(\alpha = O(1/N)\) tends to zero under comparable mean‑field assumptions, and potential‑maximizing strategies asymptotically recover exact Nash equilibria. For applied large‑population multi‑agent systems, this warns that practical closed‑loop implementations do not automatically inherit the exact asymptotic equilibrium property enjoyed by open‑loop mean‑field designs.
\end{remark}
\begin{remark}\label{rem:scaling_constants_closedloop}
We comment on the scaling behaviour of the constants appearing in the $\alpha$-bound of Theorem \ref{the1} for our closed‑loop setting.
\begin{enumerate}
\item $\tilde C_{0}^{i,j}$: this term comes from agent heterogeneity in dynamics and cost functionals; it is essentially independent of $N$.
\item $\tilde C_{1}^{i,j}$ and $\tilde C_{2}^{i,j}$ carry scaling $N^{-1/4}$ and $N^{-1/2}$. These contributions are characteristic of closed‑loop feedback, originating from the control‑sensitivity processes $v$ and $\omega$, which are identically zero in the open‑loop counterpart \cite{GLZ2025}. These fractional negative powers lead to slower decay compared with open‑loop games where $\alpha = O(1/N)$.
\item $\tilde C_{3}^{i,j}$ and $\tilde C_{4}^{i,j}$ scale as $N^{-1}$ and $N^{-2}$, corresponding to higher‑order cross‑agent coupling effects.
\end{enumerate}
This scaling explains our key mean‑field observation (Theorem \ref{the2}): under closed‑loop policies, $\alpha$ may approach a strictly positive constant as $N\to\infty$, rather than tending to zero. We stress that these are sufficient upper bounds; establishing sharp, optimal convergence rates for nonlinear closed‑loop $\alpha$-potential stochastic differential games remains an open direction for future investigation.
\end{remark}
\begin{remark}\label{rem:assumption_discussion}
The present non-asymptotic estimate for $\alpha$ in Theorem \ref{the2} builds upon Assumptions (A1)--(A2), which impose global Lipschitz continuity and uniformly bounded second-order derivatives on drift $b_i$, diffusion $\sigma_i$, and cost mappings $f_i,g_i$. These hypotheses appear frequently in BSDE-based stochastic differential game literature, yet it is instructive to comment on their necessity and possible relaxations.

The global regularity conditions enter at three key technical junctures:
\begin{enumerate}
\item Well-posedness and uniform moment bounds for the closed-loop state SDE $X^{\phi}$ (Lemma \ref{l2}); local Lipschitz alone cannot exclude finite-time explosion of trajectories without extra growth constraints.
\item Moment estimates for the first- and second-order sensitivity processes $Y^{\phi,\phi_h'}$, $Z^{\phi,\phi_h',\phi_\ell''}$, together with control-sensitivity terms $v,\omega$ (Lemmas \ref{l4}--\ref{l6}). Bounded second derivatives are essential to prevent blow-up of high-order moments for these variational processes, which determine the magnitude of Fr\'echet variations of the cost functional $V_i$.
\item A priori $L^2$-stability for the adjoint BSDE~\eqref{adj1}, cf.\ Lemma~\ref{bs2}. Standard BSDE stability results used here require globally Lipschitz generators.
\end{enumerate}

We now briefly discuss possible more general settings and associated challenges:
\begin{itemize}
\item \emph{Local Lipschitz conditions with non-explosion.}
Under local Lipschitz plus linear growth, one can guarantee non‑exploding state trajectories. However, local regularity does not yield uniform global moment bounds for the sensitivity processes $Y,Z,v,\omega$. This prevents us from obtaining our explicit global non‑asymptotic $\alpha$ bound. One could develop a stopped-time local version of the $\alpha$-potential property, but this would fundamentally alter the form of our main result and we leave it for future work.

\item \emph{Polynomial growth of second‑order derivatives.}
Polynomial‑growth second derivatives may preserve finite moments for the state process $X^\phi$. Even so, the SDEs for sensitivity processes $Z$ and $\omega$ inherit these polynomial‑growth coefficients. Though BSDE theory exists for polynomial‑growth generators, the resulting estimates become time‑dependent and policy‑dependent; simple uniform constant bounds for $\alpha$ are no longer directly available. Deriving policy‑norm‑dependent $\alpha$-estimates under polynomial growth requires substantial re‑derivation of Lemmas 4.1--4.4, which is beyond the scope of the current paper.

\item \emph{Linear‑quadratic special case.}
For closed‑loop LQ stochastic differential games, second‑order derivatives reduce to constant matrices, which alleviates several boundedness difficulties. A dedicated quantitative analysis for the LQ subclass under relaxed assumptions constitutes an interesting follow‑up direction.
\end{itemize}

In summary, our explicit non‑asymptotic $\alpha$ bound relies on the stated global assumptions, while local‑Lipschitz or polynomial‑growth frameworks represent promising directions for future investigation.
\end{remark}
We illustrate Theorem \ref{the2} by looking at the following

\begin{example}
\label{ex1} Consider the following system:%
\begin{equation*}
\begin{array}{c}
b_{i}\left( t,x_{i},x,u\right) =\tilde{b}_{i}\left( t,x_{i},\frac{1}{N}%
\sum_{\ell =1}^{N}\delta _{x_{\ell }},u\right) ,\text{ }\sigma _{i}\left(
t,x_{i},x,u\right) =\tilde{\sigma}_{i}\left( t,x_{i},\frac{1}{N}\sum_{\ell
=1}^{N}\delta _{x_{\ell }},u\right) , \\
f_{i}\left( t,x,u\right) =f_{0}\left( t,x,u\right) +c_{i}\left( u_{i}\right)
+\tilde{f}_{i}\left( t,\frac{1}{N}\sum_{\ell =1}^{N}\delta _{x_{\ell
}}\right) , \\
g_{i}\left( x\right) =g_{0}\left( x\right) +\tilde{g}_{i}\left( \frac{1}{N}%
\sum_{\ell =1}^{N}\delta _{x_{\ell }}\right) ,%
\end{array}%
\end{equation*}%
where $\tilde{b}_{i}:\left[ 0,T\right] \times \mathbb{R}\times \mathbb{R}%
^{N}\times \mathbb{R\rightarrow R},$ $\tilde{\sigma}_{i}:\left[ 0,T\right]
\times \mathbb{R}\times \mathbb{R}^{N}\times \mathbb{R\rightarrow R},$ $%
f_{0}:\left[ 0,T\right] \times \mathbb{R}^{N}\times \mathbb{R\rightarrow R},$
$c_{i}:\mathbb{R\rightarrow R},$ $\tilde{f}_{i}:\left[ 0,T\right] \times
\mathcal{P}_{2}\left( \mathbb{R}\right) \rightarrow \mathbb{R}$, $g_{0}:%
\mathbb{R}^{N}\mathbb{\rightarrow R}$ and $\tilde{g}_{i}:\mathcal{P}%
_{2}\left( \mathbb{R}\right) \rightarrow \mathbb{R}.$ According to Lions
differentiability of functions of empirical measures (see \cite{cd18},
Proposition 5.35 and 5.91), we have
\begin{equation*}
\begin{array}{c}
\left\vert \partial _{x_{h}}\tilde{b}_{i}\left( t,x,u\right) \right\vert
+\left\vert \partial _{x_{h}}\tilde{\sigma}_{i}\left( t,x,u\right)
\right\vert \leq \frac{C}{N}, \\
\left\vert \left( \partial _{x_{h}}\Delta _{t,i,j}^{f}\right) \left(
0,0\right) \right\vert +\left\vert \left( \partial _{u_{h}}\Delta
_{t,i,j}^{f}\right) \left( 0,0\right) \right\vert \leq \mathbbm{1}_{h\in
\left\{ i,j\right\} }L+\mathbbm{1}_{h\in I_{N}\backslash \left\{ i,j\right\}
}\tilde{L}N^{-1}, \\
\left\vert \partial _{x_{h}}\Delta _{i,j}^{g}\left( 0\right) \right\vert
\leq \frac{C}{N},%
\end{array}%
\end{equation*}%
and%
\begin{equation*}
\begin{array}{l}
\left\vert \left( \partial _{x_{h}x_{\ell }}^{2}\Delta _{t,i,j}^{f}\right)
\left( x,u\right) \right\vert +\left\vert \left( \partial _{u_{h}u_{\ell
}}^{2}\Delta _{t,i,j}^{f}\right) \left( t,x,u\right) \right\vert +\left\vert
\left( \partial _{x_{h}u_{\ell }}^{2}\Delta _{t,i,j}^{f}\right) \left(
t,x,u\right) \right\vert  \\
\leq C\left( \mathbbm{1}_{h=\ell \in \left\{ i,j\right\} }L+\frac{1}{N}%
\mathbbm{1}_{h=\ell \in I_{N}\backslash \left\{ i,j\right\} }+\frac{1}{N^{2}}%
\mathbbm{1}_{h\neq \ell }\right)  \\
\left\vert \partial _{x_{\ell }x_{h}}^{2}\Delta _{i,j}^{g}\left( t,x\right)
\right\vert \leq \tilde{L}\left( \mathbbm{1}_{h=\ell }N^{-\beta }+\mathbbm{1}%
_{h\neq \ell }N^{-2\beta }\right) .%
\end{array}%
\end{equation*}%
for a certain positive constant $L$ independent of $N.$ In this example, $f_{i}$ includes the term $c_{i}(u_{i})$. Due to the interdependence among players, this game remains an $\alpha$-potential game. However, as shown in (\ref{anot}), $\alpha$ may not decay to zero as $N \to \infty$. It is important to note that in the open-loop scenario (see \cite{GLZ1}), the parameter $\alpha$ in Example \ref{ex1} approaches zero. This occurs because the process $\omega_{t,i}^{\phi, \phi_h', \phi_\ell''}$, as defined in (\ref{seccon}) vanishes under this setting. This confirms that a significant difference exists between the open-loop and closed-loop scenarios in the mean field setting.

\end{example}
\begin{remark}\label{rem:scaling_assumptions_lq_example}
We now comment on the derivative‑scaling conditions in \eqref{fi1}, their behaviour under mean‑field limiting procedures, and provide an explicit linear‑quadratic (LQ) illustrative example.

Naturalness of scaling assumptions for finite‑$N$ mean‑field games. The conditions in \eqref{fi1} require that mixed partial derivatives of the cost differences $\Delta_{i,j}^{f}$ and $\Delta_{i,j}^{g}$ decay like negative powers of $N$ for cross‑player indices. This type of scaling arises organically for cost functionals that depend on the empirical state distribution $\frac{1}{N}\sum_{k=1}^N \delta_{X_{t,k}}$. By the chain rule for differentiability with respect to empirical measures, derivatives with respect to distinct agents’ states automatically bring factors $1/N$ or $1/N^2$, which matches the structure of our bounds. Therefore these assumptions are well‑motivated for finite‑$N$ multi‑agent systems with mean‑field‑type interaction.


An LQ closed‑loop mean‑field example.
Consider $N$ agents with state dynamics
\begin{align*}
dX_{t,i}&=\big(a X_{t,i}+\bar a \bar X_t + u_{t,i}\big)dt+\sigma dW_{t}^{i},\qquad
\bar X_t=\frac{1}{N}\sum_{k=1}^N X_{t,k},\\
X_{0,i}&=\xi_i,\quad i=1,\dots,N,
\end{align*}
where $a,\bar a,\sigma$ are constants, $\xi_i\in L^4(\Omega;\mathbb R)$. Each agent minimises the cost functional
\[
V_i(\phi)=\mathbb{E}\left[\int_0^T\left(c_{i}( u_{t,i})+\frac{1}{2}(X_{t,i}-\bar X_t)^2\right)dt
+\frac{1}{2} X_{T,i}^2\right],
\]
under closed‑loop state‑feedback policies $u_{t,i}=\phi_{t,i}(X_{t,i},\bar X_t)$.
The closed‑loop feedback generates non‑trivial control‑sensitivity processes $v$ and $\omega$. By Theorem \ref{the2}, as $N\to\infty$, $\alpha$ converges to a point in an interval (see \ref{anot}). This LQ example demonstrates a key message of our paper: non‑vanishing $\alpha$ in the large‑population limit is an intrinsic closed‑loop feedback effect, and it persists. This contrasts sharply with the open‑loop LQ setting from \cite{GLZ2025}, where $\alpha=O(1/N)\to0$.
Moreover, for this LQ model, the second‑order cross‑agent derivatives of the running and terminal costs satisfy the scaling requirements in \eqref{fi1}, consistent with mean‑field interaction via the empirical average.
\end{remark}

\section{Proof of Theorem \protect\ref{the1}.}
\label{sect5}
To simplify the notation, we omit the dependence on $\phi $ in the
superscript of all processes, i.e., $X=X^{\phi }$, $Y^{i}=Y^{\phi ,\phi
_{i}^{\prime }},$ $Z^{i,j}=Z^{\phi ,\phi _{i}^{\prime },\phi _{j}^{\prime
\prime }}$. We denote by $C>0$, a generic constant depending only on the
upper bounds of $T$, $\max_{i\in I_{N}}\mathbb{E}\left[ \left\vert \xi
_{i}\right\vert ^{4}\right] ,$ $L^{b},L^{\sigma },L^{\phi },L^{\phi
_{h}^{\prime }},L^{\phi _{\ell }^{\prime \prime }},L_{y}^{b},L_{y}^{\sigma
},L_{y}^{\phi _{h}^{\prime }},L_{y}^{\phi _{\ell }^{\prime \prime }}$. Now
we compute, noting the fact that $Z^{\phi ,\phi _{i}^{\prime },\phi
_{j}^{\prime \prime }}=Z^{\phi ,\phi _{j}^{\prime \prime },\phi _{i}^{\prime
}},$ and $\omega ^{\phi ,\phi _{i}^{\prime },\phi _{j}^{\prime \prime
}}=\omega ^{\phi ,\phi _{j}^{\prime \prime },\phi _{i}^{\prime }}$
\begin{eqnarray*}
&&\frac{\delta ^{2}V_{i}}{\delta \phi _{i}\delta \phi _{j}}\left( \phi ,\phi
_{i}^{\prime },\phi _{j}^{\prime \prime }\right) -\frac{\delta ^{2}V_{j}}{%
\delta \phi _{j}\delta \phi _{i}}\left( \phi ,\phi _{j}^{\prime \prime
},\phi _{i}^{\prime }\right) \\
&=&\mathbb{E}\Bigg [\int_{0}^{T}\sum_{h,\ell =1}^{N}\left( \partial
_{x_{h}x_{\ell }}^{2}\Delta _{t,i,j}^{f}\left( \mathbf{X}_{t},\mathbf{u}%
_{t}\right) \right) Y_{t,h}^{i}Y_{t,\ell }^{j}+\sum_{h,\ell =1}^{N}\left(
\partial _{x_{h}u_{\ell }}^{2}\Delta _{t,i,j}^{f}\left( \mathbf{X}_{t},%
\mathbf{u}_{t}\right) \right) Y_{t,h}^{i}\upsilon _{t,\ell }^{j} \\
&&+\sum_{h,\ell =1}^{N}\left( \partial _{u_{h}x_{\ell }}^{2}\Delta
_{t,i,j}^{f}\left( \mathbf{X}_{t},\mathbf{u}_{t}\right) \right) \upsilon
_{t,h}^{i}Y_{t,\ell }^{j}+\sum_{h,\ell =1}^{N}\left( \partial _{u_{h}u_{\ell
}}^{2}\Delta _{t,i,j}^{f}\left( \mathbf{X}_{t},\mathbf{u}_{t}\right) \right)
\upsilon _{t,h}^{i}\upsilon _{t,\ell }^{j} \\
&&+\sum_{h=1}^{N}\left( \partial _{x_{h}}\Delta _{t,i,j}^{f}\left( \mathbf{X}%
_{t},\mathbf{u}_{t}\right) \right) Z_{t,h}^{i,j}+\sum_{h=1}^{N}\left(
\partial _{u_{h}}\Delta _{t,i,j}^{f}\left( \mathbf{X}_{t},\mathbf{u}%
_{t}\right) \right) \omega _{t,h}^{i,j}\Bigg )\mathrm{d}t\Bigg ] \\
&&+\mathbb{E}\left[ \left( \mathbf{Y}_{T}^{i}\right) ^{\top }\left( \partial
_{xx}^{2}\Delta _{i,j}^{g}\right) \left( \mathbf{X}_{T}^{\phi }\right)
\mathbf{Y}_{T}^{j}+\left( \partial _{x}\Delta _{i,j}^{g}\right) ^{\top
}\left( \mathbf{X}_{T}^{\phi }\right) \mathbf{Z}_{T}^{i,j}\right] ,
\end{eqnarray*}%
where we write for simplicity $\partial _{xx}^{2}\Delta
_{t,i,j}^{f}=\partial _{xx}^{2}\left( f_{t,i}-f_{t,j}\right) \left( \mathbf{X%
}_{t},\mathbf{u}_{t}\right) $ and similarly for other derivatives.

Now consider the following BSDE:%
\begin{equation}
\left\{
\begin{array}{rcl}
-\mathrm{d}P_{t}^{i,j} & = & \Big [\mathbf{B}\left( t,\mathbf{X}_{t},\mathbf{%
u}_{t}\right) ^{\top }P_{t}^{i,j}+\sum_{k=1}^{N}\mathbf{\Pi }^{k}\left( t,%
\mathbf{X}_{t},\mathbf{u}_{t}\right) ^{\top }Q_{t,k}^{i,j} \\
&  & +\left( \partial _{x}\Delta _{t,i,j}^{f}\right) \left( \mathbf{X}_{t},%
\mathbf{u}_{t}\right) \Big ]\mathrm{d}t-\sum_{k=1}^{N}Q_{t,k}^{i,j}\mathrm{d}%
W_{t}^{k}, \\
P_{T}^{i,j} & = & \partial _{x}\Delta _{i,j}^{g}\left( \mathbf{X}_{T}\right)
,\text{ }%
\end{array}%
\right.  \label{adjn1}
\end{equation}%
where, for simplicity, we denote by $Q_{t}^{i,j}=\left(
\begin{array}{ccc}
Q_{t,1}^{i,j} & \cdots & Q_{t,N}^{i,j}%
\end{array}%
\right) ^{\top }.$

Immediately, by virtue of It\^{o}'s formula, it follows that
\begin{eqnarray}
&&\frac{\delta ^{2}V_{i}}{\delta \phi _{i}\delta \phi _{j}}\left( \phi ,\phi
_{i}^{\prime },\phi _{j}^{\prime \prime }\right) -\frac{\delta ^{2}V_{j}}{%
\delta \phi _{j}\delta \phi _{i}}\left( \phi ,\phi _{j}^{\prime \prime
},\phi _{i}^{\prime }\right)  \notag \\
&=&\mathbb{E}\Bigg [\int_{0}^{T}\underset{J_{1}}{\underbrace{\sum_{h,\ell
=1}^{N}\left( \partial _{x_{h}x_{\ell }}^{2}\Delta _{t,i,j}^{f}\left(
\mathbf{X}_{t},\mathbf{u}_{t}\right) \right) Y_{t,h}^{i}Y_{t,\ell }^{j}}}+%
\underset{J_{2}}{\underbrace{\sum_{h,\ell =1}^{N}\left( \partial
_{u_{h}u_{\ell }}^{2}\Delta _{t,i,j}^{f}\left( \mathbf{X}_{t},\mathbf{u}%
_{t}\right) \right) \upsilon _{t,h}^{i}\upsilon _{t,\ell }^{j}}}  \notag \\
&&+\underset{J_{3}}{\underbrace{\sum_{h,\ell =1}^{N}\left( \partial
_{u_{h}x_{\ell }}^{2}\Delta _{t,i,j}^{f}\left( \mathbf{X}_{t},\mathbf{u}%
_{t}\right) \right) \upsilon _{t,h}^{i}Y_{t,\ell }^{j}+\sum_{h,\ell
=1}^{N}\left( \partial _{x_{h}u_{\ell }}^{2}\Delta _{t,i,j}^{f}\left(
\mathbf{X}_{t},\mathbf{u}_{t}\right) \right) Y_{t,h}^{i}\upsilon _{t,\ell
}^{j}}}  \notag \\
&&+\underset{J_{4}}{\underbrace{\sum_{h=1}^{N}\left( \partial _{u_{h}}\Delta
_{t,i,j}^{f}\left( \mathbf{X}_{t},\mathbf{u}_{t}\right) \right) \omega
_{t,h}^{i,j}}}\Bigg )\mathrm{d}t\Bigg ]  \notag \\
&&+\mathbb{E}\Bigg [\underset{J_{5}}{\underbrace{\int_{0}^{T}\left\langle
P_{t}^{i,j},\mathfrak{F}_{t}^{i,j}+\partial _{yy}^{2}\bar{b}_{t}\left(
\mathbf{X}_{t}\right) \star \mathcal{Y}_{t}^{i,j}+\Gamma _{t}^{i,j}\left(
\bar{b}_{t}\right) \right\rangle \mathrm{d}t}}  \notag \\
&&+\underset{J_{6}}{\underbrace{\int_{0}^{T}\sum_{k=1}^{N}\left\langle
Q_{t,k}^{i,j},\mathfrak{G}_{t,k}^{i,j}+\partial _{yy}^{2}\bar{\sigma}%
_{t,k}\left( X_{t,k},\mathbf{X}_{t}\right) \star \mathcal{Y}_{t}^{i,j}+\bar{%
\Xi}_{t,k}^{i,j}\bar{\sigma}_{t}\right\rangle \mathrm{d}t}}\Bigg ]  \notag \\
&&+\underset{J_{7}}{\underbrace{\mathbb{E}\left[ \left( \mathbf{Y}%
_{T}^{i}\right) ^{\top }\left( \partial _{xx}^{2}\Delta _{i,j}^{g}\left(
\mathbf{X}_{T}\right) \right) \left( \mathbf{X}_{T}\right) \mathbf{Y}_{T}^{j}%
\right] }}.  \label{ev}
\end{eqnarray}%
We will estimate the right side of (\ref{ev}) in the following way.

For $J_{1},$ from Lemma \ref{l4}, we have
\begin{eqnarray*}
&&\Bigg |\mathbb{E}\Bigg \{\int_{0}^{T}\sum_{h,\ell =1}^{N}\left( \partial
_{x_{h}x_{\ell }}^{2}\Delta _{t,i,j}^{f}\left( \mathbf{X}_{t},\mathbf{u}%
_{t}\right) \right) Y_{t,h}^{i}Y_{t,\ell }^{j}\mathrm{d}t\Bigg ]\Bigg | \\
&=&\Bigg |\mathbb{E}\Bigg [\int_{0}^{T}\Bigg [\left( \partial
_{x_{i}x_{j}}^{2}\Delta _{t,i,j}^{f}\left( \mathbf{X}_{t},\mathbf{u}%
_{t}\right) \right) Y_{t,i}^{i}Y_{t,j}^{j}+\sum_{\ell \in I_{N}\backslash
\left\{ j\right\} }^{N}\left( \partial _{x_{i}x_{\ell }}^{2}\Delta
_{t,i,j}^{f}\left( \mathbf{X}_{t},\mathbf{u}_{t}\right) \right)
Y_{t,i}^{i}Y_{t,\ell }^{j} \\
&&+\sum_{h\in I_{N}\backslash \left\{ i\right\} }^{N}\Bigg (\left( \partial
_{x_{h}x_{j}}^{2}\Delta _{t,i,j}^{f}\left( \mathbf{X}_{t},\mathbf{u}%
_{t}\right) \right) Y_{t,h}^{i}Y_{t,j}^{j}+\sum_{\ell \in I_{N}\backslash
\left\{ j\right\} }^{N}\left( \partial _{x_{h}x_{\ell }}^{2}\Delta
_{t,i,j}^{f}\left( \mathbf{X}_{t},\mathbf{u}_{t}\right) \right)
Y_{t,h}^{i}Y_{t,\ell }^{j}\Bigg )\Bigg ]\mathrm{d}t\Bigg \}\Bigg | \\
&\leq &\left\Vert \partial _{x_{i}x_{j}}^{2}\Delta _{i,j}^{f}\right\Vert
_{L^{\infty }}\left\Vert Y_{i}^{i}Y_{j}^{j}\right\Vert _{\mathcal{H}%
^{1}\left( \mathbb{R}\right) }+\sum_{\ell \in I_{N}\backslash \left\{
j\right\} }^{N}\left\Vert \partial _{x_{i}x_{\ell }}^{2}\Delta
_{t,i,j}^{f}\right\Vert _{L^{\infty }}\left\Vert Y_{i}^{i}Y_{\ell
}^{j}\right\Vert _{\mathcal{H}^{1}\left( \mathbb{R}\right) } \\
&&+\sum_{h\in I_{N}\backslash \left\{ i\right\} }^{N}\Bigg (\left\Vert
\partial _{x_{h}x_{j}}^{2}\Delta _{i,j}^{f}\right\Vert _{L^{\infty
}}\left\Vert Y_{h}^{i}Y_{j}^{j}\right\Vert _{_{\mathcal{H}^{1}\left( \mathbb{%
R}\right) }}+\sum_{\ell \in I_{N}\backslash \left\{ j\right\}
}^{N}\left\Vert \partial _{x_{h}x_{\ell }}^{2}\Delta _{i,j}^{f}\right\Vert
_{L^{\infty }}\left\Vert Y_{h}^{i}Y_{\ell }^{j}\right\Vert _{_{\mathcal{H}%
^{1}\left( \mathbb{R}\right) }}\Bigg ) \\
&\leq &C\Bigg [\left\Vert \partial _{x_{i}x_{j}}^{2}\Delta
_{i,j}^{f}\right\Vert _{L^{\infty }}+\frac{\left( L_{y}^{\bar{b}}+3\left(
L_{y}^{\bar{\sigma}}\right) ^{2}\right) ^{\frac{1}{2}}}{\sqrt{N}}\left(
\sum_{\ell \in I_{N}\backslash \left\{ j\right\} }^{N}\left\Vert \partial
_{x_{i}x_{\ell }}^{2}\Delta _{t,i,j}^{f}\right\Vert _{L^{\infty
}}+\sum_{h\in I_{N}\backslash \left\{ i\right\} }^{N}\left\Vert \partial
_{x_{h}x_{j}}^{2}\Delta _{i,j}^{f}\right\Vert _{L^{\infty }}\right)  \\
&&+\frac{\left( L_{y}^{\bar{b}}+3\left( L_{y}^{\bar{\sigma}}\right)
^{2}\right) }{N}\sum_{h\in I_{N}\backslash \left\{ i\right\} }^{N}\sum_{\ell
\in I_{N}\backslash \left\{ j\right\} }^{N}\left\Vert \partial
_{x_{h}x_{\ell }}^{2}\Delta _{i,j}^{f}\right\Vert _{L^{\infty }}\Bigg ].
\end{eqnarray*}%
Next for $J_{3}$%
\begin{eqnarray*}
&&\left\vert \mathbb{E}\left[ \int_{0}^{T}\sum_{h,\ell =1}^{N}\left(
\partial _{x_{h}u_{\ell }}^{2}\Delta _{t,i,j}^{f}\left( \mathbf{X}_{t},%
\mathbf{u}_{t}\right) \right) Y_{t,h}^{i}\upsilon _{t,\ell }^{j}\mathrm{d}t%
\right] \right\vert  \\
&&+\left\vert \mathbb{E}\left[ \int_{0}^{T}\sum_{h,\ell =1}^{N}\left(
\partial _{u_{h}x_{\ell }}^{2}\Delta _{t,i,j}^{f}\left( \mathbf{X}_{t},%
\mathbf{u}_{t}\right) \right) \upsilon _{t,h}^{i}Y_{t,\ell }^{j}\mathrm{d}t%
\right] \right\vert  \\
&=&\Bigg |\mathbb{E}\Bigg \{\int_{0}^{T}\Bigg [\partial
_{x_{i}u_{j}}^{2}\Delta _{t,i,j}^{f}\left( \mathbf{X}_{t},\mathbf{u}%
_{t}\right) Y_{t,i}^{i}\upsilon _{t,j}^{j}+\sum_{\ell \in I_{N}\backslash
\left\{ j\right\} }^{N}\left( \partial _{x_{i}u_{\ell }}^{2}\Delta
_{t,i,j}^{f}\left( \mathbf{X}_{t},\mathbf{u}_{t}\right) \right)
Y_{t,i}^{i}\upsilon _{t,\ell }^{j} \\
&&+\sum_{h\in I_{N}\backslash \left\{ i\right\} }^{N}\left( \partial
_{x_{h}u_{j}}^{2}\Delta _{t,i,j}^{f}\left( \mathbf{X}_{t},\mathbf{u}%
_{t}\right) \right) Y_{t,h}^{i}\upsilon _{t,j}^{j}+\sum_{\substack{ h\in
I_{N}\backslash \left\{ i\right\}  \\ \ell \in I_{N}\backslash \left\{
j\right\} }}^{N}\left( \partial _{x_{h}u_{\ell }}^{2}\Delta
_{t,i,j}^{f}\left( \mathbf{X}_{t},\mathbf{u}_{t}\right) \right)
Y_{t,h}^{i}\upsilon _{t,\ell }^{j}\Bigg ]\mathrm{d}t\Bigg \}\Bigg | \\
&&+\Bigg |\mathbb{E}\Bigg \{\int_{0}^{T}\Bigg [\left( \partial
_{u_{i}x_{j}}^{2}\Delta _{t,i,j}^{f}\left( \mathbf{X}_{t},\mathbf{u}%
_{t}\right) \right) \upsilon _{t,i}^{i}Y_{t,j}^{j}+\sum_{\ell \in
I_{N}\backslash \left\{ j\right\} }^{N}\left( \partial _{u_{i}x_{\ell
}}^{2}\Delta _{t,i,j}^{f}\left( \mathbf{X}_{t},\mathbf{u}_{t}\right) \right)
\upsilon _{t,i}^{i}Y_{t,\ell }^{j} \\
&&+\sum_{h\in I_{N}\backslash \left\{ i\right\} }^{N}\left( \partial
_{u_{h}x_{j}}^{2}\Delta _{t,i,j}^{f}\left( \mathbf{X}_{t},\mathbf{u}%
_{t}\right) \right) \upsilon _{t,h}^{i}Y_{t,j}^{j}+\sum_{\substack{ h\in
I_{N}\backslash \left\{ i\right\}  \\ \ell \in I_{N}\backslash \left\{
j\right\} }}^{N}\left( \partial _{u_{h}x_{\ell }}^{2}\Delta
_{t,i,j}^{f}\left( \mathbf{X}_{t},\mathbf{u}_{t}\right) \right) \upsilon
_{t,h}^{i}Y_{t,\ell }^{j}\Bigg ]\mathrm{d}t\Bigg \}\Bigg |.
\end{eqnarray*}%
From Lemma \ref{l4} and \ref{l5}, we derive
\begin{eqnarray*}
&&\left\vert \mathbb{E}\left[ \int_{0}^{T}\sum_{h,\ell =1}^{N}\left(
\partial _{x_{h}u_{\ell }}^{2}\Delta _{t,i,j}^{f}\left( \mathbf{X}_{t},%
\mathbf{u}_{t}\right) \right) Y_{t,h}^{i}\upsilon _{t,\ell }^{j}\mathrm{d}t%
\right] \right\vert  \\
&&+\left\vert \mathbb{E}\left[ \int_{0}^{T}\sum_{h,\ell =1}^{N}\left(
\partial _{u_{h}x_{\ell }}^{2}\Delta _{t,i,j}^{f}\left( \mathbf{X}_{t},%
\mathbf{u}_{t}\right) \right) \upsilon _{t,h}^{i}Y_{t,\ell }^{j}\mathrm{d}t%
\right] \right\vert  \\
&\leq &\left\Vert \partial _{x_{i}u_{j}}^{2}\Delta _{i,j}^{f}\left( \mathbf{X%
}_{t},\mathbf{u}_{t}\right) \right\Vert _{L^{\infty }}+\frac{\left( L_{y}^{%
\bar{b}}+3\left( L_{y}^{\bar{\sigma}}\right) ^{2}\right) ^{\frac{1}{2}}}{%
\sqrt{N}}\sum_{h\in I_{N}\backslash \left\{ i\right\} }^{N}\left\Vert
\partial _{x_{h}u_{j}}^{2}\Delta _{t,i,j}^{f}\left( \mathbf{X}_{t},\mathbf{u}%
_{t}\right) \right\Vert _{L^{\infty }} \\
&&+\frac{1}{\sqrt{N}}\left[ \left( \left( L_{y}^{\phi }\right) ^{2}+\left(
L^{\phi }\right) ^{2}\right) \left( L_{y}^{\bar{b}}+3\left( L_{y}^{\bar{%
\sigma}}\right) ^{2}\right) \right] ^{\frac{1}{2}}\sum_{\ell \in
I_{N}\backslash \left\{ j\right\} }^{N}\left\Vert \partial _{x_{i}u_{\ell
}}^{2}\Delta _{t,i,j}^{f}\left( \mathbf{X}_{t},\mathbf{u}_{t}\right)
\right\Vert _{L^{\infty }} \\
&&+\frac{\left( L_{y}^{\bar{b}}+3\left( L_{y}^{\bar{\sigma}}\right)
^{2}\right) }{N}\left( \left( L_{y}^{\phi }\right) ^{2}+\left( L^{\phi
}\right) ^{2}\right) ^{\frac{1}{2}}\sum_{\substack{ h\in I_{N}\backslash
\left\{ i\right\}  \\ \ell \in I_{N}\backslash \left\{ j\right\} }}%
^{N}\left\Vert \left( \partial _{x_{h}u_{\ell }}^{2}\Delta _{i,j}^{f}\left(
\mathbf{X}_{t},\mathbf{u}_{t}\right) \right) \right\Vert _{L^{\infty }} \\
&&+\left\Vert \partial _{u_{i}x_{j}}^{2}\Delta _{i,j}^{f}\left( \mathbf{X}%
_{t},\mathbf{u}_{t}\right) \right\Vert _{L^{\infty }}+\frac{\left( L_{y}^{%
\bar{b}}+3\left( L_{y}^{\bar{\sigma}}\right) ^{2}\right) ^{\frac{1}{2}}}{%
\sqrt{N}}\sum_{\ell \in I_{N}\backslash \left\{ j\right\} }^{N}\left\Vert
\partial _{u_{i}x_{\ell }}^{2}\Delta _{t,i,j}^{f}\left( \mathbf{X}_{t},%
\mathbf{u}_{t}\right) \right\Vert _{L^{\infty }} \\
&&+\frac{1}{\sqrt{N}}\left[ \left( \left( L_{y}^{\phi }\right) ^{2}+\left(
L^{\phi }\right) ^{2}\right) \left( L_{y}^{\bar{b}}+3\left( L_{y}^{\bar{%
\sigma}}\right) ^{2}\right) \right] ^{\frac{1}{2}}\sum_{h\in I_{N}\backslash
\left\{ i\right\} }^{N}\left\Vert \partial _{u_{h}x_{j}}^{2}\Delta
_{t,i,j}^{f}\left( \mathbf{X}_{t},\mathbf{u}_{t}\right) \right\Vert
_{L^{\infty }} \\
&&+\frac{\left( L_{y}^{\bar{b}}+3\left( L_{y}^{\bar{\sigma}}\right)
^{2}\right) }{N}\left( \left( L_{y}^{\phi }\right) ^{2}+\left( L^{\phi
}\right) ^{2}\right) ^{\frac{1}{2}}\sum_{\substack{ h\in I_{N}\backslash
\left\{ i\right\}  \\ \ell \in I_{N}\backslash \left\{ j\right\} }}%
^{N}\left\Vert \left( \partial _{u_{h}x_{\ell }}^{2}\Delta _{i,j}^{f}\left(
\mathbf{X}_{t},\mathbf{u}_{t}\right) \right) \right\Vert _{L^{\infty }}.
\end{eqnarray*}%
Noting that
\begin{eqnarray*}
&&\sum_{h,\ell =1}^{N}\left( \partial _{u_{h}u_{\ell }}^{2}\Delta
_{t,i,j}^{f}\left( \mathbf{X}_{t},\mathbf{u}_{t}\right) \right) \upsilon
_{t,h}^{i}\upsilon _{t,\ell }^{j} \\
&=&\left( \partial _{u_{i}u_{j}}^{2}\Delta _{t,i,j}^{f}\left( \mathbf{X}_{t},%
\mathbf{u}_{t}\right) \right) \upsilon _{t,i}^{i}\upsilon
_{t,j}^{j}+\sum_{\ell \in I_{N}\backslash \left\{ j\right\} }^{N}\left(
\partial _{u_{i}u_{\ell }}^{2}\Delta _{t,i,j}^{f}\left( \mathbf{X}_{t},%
\mathbf{u}_{t}\right) \right) \upsilon _{t,i}^{i}\upsilon _{t,\ell }^{j} \\
&&+\sum_{h\in I_{N}\backslash \left\{ i\right\} }^{N}\left( \partial
_{u_{h}u_{j}}^{2}\Delta _{t,i,j}^{f}\left( \mathbf{X}_{t},\mathbf{u}%
_{t}\right) \right) \upsilon _{t,h}^{i}\upsilon _{t,j}^{j}+\sum_{\substack{ %
h\in I_{N}\backslash \left\{ i\right\}  \\ \ell \in I_{N}\backslash \left\{
j\right\} }}^{N}\left( \partial _{u_{h}u_{\ell }}^{2}\Delta
_{t,i,j}^{f}\left( \mathbf{X}_{t},\mathbf{u}_{t}\right) \right) \upsilon
_{t,h}^{i}\upsilon _{t,\ell }^{j},
\end{eqnarray*}%
it is easy to check for $J_{2}$%
\begin{eqnarray*}
&&\left\vert \mathbb{E}\left[ \int_{0}^{T}\sum_{h,\ell =1}^{N}\left(
\partial _{u_{h}u_{\ell }}^{2}\Delta _{t,i,j}^{f}\left( \mathbf{X}_{t},%
\mathbf{u}_{t}\right) \right) \upsilon _{t,h}^{i}\upsilon _{t,\ell }^{j}%
\mathrm{d}t\right] \right\vert  \\
&\leq &\left\Vert \left( \partial _{u_{i}u_{j}}^{2}\Delta _{i,j}^{f}\left(
\mathbf{X}_{t},\mathbf{u}_{t}\right) \right) \right\Vert _{L^{\infty }} \\
&&+\frac{1}{\sqrt{N}}\left( \left( L_{y}^{\phi }\right) ^{2}+\left( L^{\phi
}\right) ^{2}\right) ^{\frac{1}{2}}\left( L_{y}^{\bar{b}}+3\left( L_{y}^{%
\bar{\sigma}}\right) ^{2}\right) ^{\frac{1}{2}}\Bigg (\sum_{\ell \in
I_{N}\backslash \left\{ j\right\} }^{N}\left\Vert \partial _{u_{i}u_{\ell
}}^{2}\Delta _{t,i,j}^{f}\left( \mathbf{X}_{t},\mathbf{u}_{t}\right)
\right\Vert _{L^{\infty }} \\
&&+\sum_{h\in I_{N}\backslash \left\{ i\right\} }^{N}\left\Vert \partial
_{u_{h}u_{j}}^{2}\Delta _{i,j}^{f}\left( \mathbf{X}_{t},\mathbf{u}%
_{t}\right) \right\Vert _{L^{\infty }}\Bigg ) \\
&&+\frac{1}{N}\left( \left( L_{y}^{\phi }\right) ^{2}+\left( L^{\phi
}\right) ^{2}\right) \left( L_{y}^{\bar{b}}+3\left( L_{y}^{\bar{\sigma}%
}\right) ^{2}\right) \sum_{\substack{ h\in I_{N}\backslash \left\{ i\right\}
\\ \ell \in I_{N}\backslash \left\{ j\right\} }}^{N}\left\Vert \partial
_{u_{h}u_{\ell }}^{2}\Delta _{t,i,j}^{f}\left( \mathbf{X}_{t},\mathbf{u}%
_{t}\right) \right\Vert _{L^{\infty }}.
\end{eqnarray*}%
For all $\left( t,x,u\right) =\left( t,\left( x_{\ell }\right) _{\ell
=1}^{N},\left( u_{\ell }\right) _{\ell =1}^{N}\right) \in \left[ 0,T\right]
\times \mathbb{R}^{n}\times \mathbb{R}^{n}$ and $i,$ $j\in I_{N},$ we have%
\begin{equation*}
\left\vert \left( \partial _{x_{\ell }}\Delta _{t,i,j}^{f}\right) \left(
x,u\right) \right\vert \leq \Bigg [\left\vert \left( \partial _{x_{\ell
}}\Delta _{t,i,j}^{f}\right) \left( 0,0\right) \right\vert
+\sum_{k=1}^{N}\left( \left\Vert \partial _{x_{_{\ell }}x_{k}}\Delta
_{i,j}^{f}\right\Vert _{L^{\infty }}\left\vert x_{k}\right\vert +\left\Vert
\partial _{x_{\ell }u_{k}}\Delta _{i,j}^{f}\right\Vert _{L^{\infty
}}\left\vert u_{k}\right\vert \right) \Bigg ].
\end{equation*}%
Hence, for $J_{4}$%
\begin{eqnarray*}
&&\left\vert \mathbb{E}\left[ \int_{0}^{T}\sum_{h=1}^{N}\left( \partial
_{u_{h}}\Delta _{t,i,j}^{f}\left( \mathbf{X}_{t},\mathbf{u}_{t}\right)
\right) \omega _{t,h}^{i,j}\right] \right\vert  \\
&\leq &\sum_{h\in \left\{ i,j\right\} }^{N}\Bigg (\left\Vert \left( \partial
_{u_{h}}\Delta _{i,j}^{f}\right) \left( 0,0\right) \right\Vert _{L^{\infty
}}\left\Vert \omega _{h}^{i,j}\right\Vert _{\mathcal{H}^{2}\left( \mathbb{R}%
\right) } \\
&&+\sum_{k=1}^{N}\left( \left\Vert \partial _{u_{_{h}}x_{k}}\Delta
_{i,j}^{f}\right\Vert _{L^{\infty }}\left\Vert X_{k}\omega
_{h}^{i,j}\right\Vert _{\mathcal{H}^{1}\left( \mathbb{R}\right) }+\left\Vert
\partial _{u_{h}u_{k}}\Delta _{i,j}^{f}\right\Vert _{L^{\infty }}\left\Vert
u_{k}\omega _{h}^{i,j}\right\Vert _{\mathcal{H}^{1}\left( \mathbb{R}\right)
}\right) \Bigg ) \\
&&+\sum_{h\in I_{N}\backslash \left\{ i,j\right\} }^{N}\Bigg (\left\Vert
\left( \partial _{u_{h}}\Delta _{i,j}^{f}\right) \left( 0,0\right)
\right\Vert _{L^{\infty }}\left\Vert \omega _{h}^{i,j}\right\Vert _{\mathcal{%
H}^{2}\left( \mathbb{R}\right) } \\
&&+\sum_{k=1}^{N}\left( \left\Vert \partial _{u_{_{h}}x_{k}}\Delta
_{i,j}^{f}\right\Vert _{L^{\infty }}\left\Vert X_{k}\omega
_{h}^{i,j}\right\Vert _{\mathcal{H}^{1}\left( \mathbb{R}\right) }+\left\Vert
\partial _{u_{h}u_{k}}\Delta _{i,j}^{f}\right\Vert _{L^{\infty }}\left\Vert
u_{k}\omega _{h}^{i,j}\right\Vert _{\mathcal{H}^{1}\left( \mathbb{R}\right)
}\right) \Bigg ).
\end{eqnarray*}%
Applying Lemma \ref{l5} and Cauchy-Schwarz inequality, we arrive at%
\begin{eqnarray*}
&&\left\vert \mathbb{E}\left[ \int_{0}^{T}\sum_{h=1}^{N}\left( \partial
_{u_{h}}\Delta _{t,i,j}^{f}\left( \mathbf{X}_{t},\mathbf{u}_{t}\right)
\right) \omega _{t,h}^{i,j}\right] \right\vert  \\
&\leq &\Bigg [\left( 1+L_{y}^{\bar{b}}+L_{y}^{\bar{\sigma}}\right) \left(
L_{y}^{\bar{b}}+3\left( L_{y}^{\bar{\sigma}}\right) ^{2}\right) ^{\frac{1}{4}%
}\frac{1}{N^{\frac{1}{4}}}+\sqrt{L_{y}^{\phi _{i}^{\prime }}\left( 1+L_{y}^{%
\bar{b}}+3\left( L_{y}^{\bar{\sigma}}\right) ^{2}\right) } \\
&&+\sqrt{L_{y}^{\phi _{j}^{\prime \prime }}\left( 1+L_{y}^{\bar{b}}+3\left(
L_{y}^{\bar{\sigma}}\right) ^{2}\right) }\Bigg ]\cdot \sum_{h\in \left\{
i,j\right\} }^{N}\Big [\left\Vert \left( \partial _{u_{h}}\Delta
_{i,j}^{f}\right) \left( 0,0\right) \right\Vert _{L^{\infty }} \\
&&+\sum_{k=1}^{N}\left( \left\Vert \partial _{u_{_{h}}x_{k}}\Delta
_{i,j}^{f}\right\Vert _{L^{\infty }}+\left\Vert \partial _{u_{h}u_{k}}\Delta
_{i,j}^{f}\right\Vert _{L^{\infty }}\right) \Big ] \\
&&+\left( 1+L_{y}^{\bar{b}}+L_{y}^{\bar{\sigma}}\right) \left( L_{y}^{\bar{b}%
}+3\left( L_{y}^{\bar{\sigma}}\right) ^{2}\right) ^{\frac{1}{2}}\frac{1}{N^{%
\frac{1}{2}}}\cdot \sum_{h\in I_{N}\backslash \left\{ i,j\right\} }^{N}\Big [%
\left\Vert \left( \partial _{u_{\ell }}\Delta _{i,j}^{f}\right) \left(
0,0\right) \right\Vert _{L^{\infty }} \\
&&+\sum_{k=1}^{N}\left( \left\Vert \partial _{u_{_{h}}x_{k}}\Delta
_{i,j}^{f}\right\Vert _{L^{\infty }}+\left\Vert \partial _{u_{h}u_{k}}\Delta
_{i,j}^{f}\right\Vert _{L^{\infty }}\right) \Big ].
\end{eqnarray*}
For $J_{7},$ from the relation $\left( \mathbf{Y}_{T}^{i}\right) ^{\top
}\left( \partial _{xx}^{2}\Delta _{i,j}^{g}\right) \left( \mathbf{X}%
_{T}\right) \mathbf{Y}_{T}^{j}=\sum_{h,\ell =1}^{N}\left( \partial
_{x_{h}x_{\ell }}^{2}\Delta _{i,j}^{g}\right) \left( \mathbf{X}_{T}\right)
Y_{T,h}^{i}Y_{T,\ell }^{j},$ we have%
\begin{eqnarray*}
&&\left\vert \mathbb{E}\left[ \left( \mathbf{Y}_{T}^{i}\right) ^{\top
}\left( \partial _{xx}^{2}\Delta _{i,j}^{g}\right) \left( \mathbf{X}%
_{T}\right) \mathbf{Y}_{T}^{j}\right] \right\vert  \\
&\leq &C\Bigg [\left\Vert \partial _{x_{i}x_{j}}^{2}\Delta
_{i,j}^{g}\right\Vert _{L^{\infty }} \\
&&+\frac{L_{y}^{\bar{b}}+3\left( L_{y}^{\bar{\sigma}}\right) ^{2}}{N}\Bigg (%
\sum_{\ell \in I_{N}\backslash \left\{ j\right\} }\left\Vert \partial
_{x_{i}x_{\ell }}^{2}\Delta^{g}g_{i,j}\right\Vert _{L^{\infty }}+\sum_{h\in
I_{N}\backslash \left\{ i\right\} }\left\Vert \partial
_{x_{h}x_{j}}^{2}\Delta^{g}_{i,j}\right\Vert _{L^{\infty }}\Bigg ) \\
&&+\frac{\left( L_{y}^{\bar{b}}+3\left( L_{y}^{\bar{\sigma}}\right)
^{2}\right) ^{2}}{N^{2}}\sum_{\ell \in I_{N}\backslash \left\{ j\right\}
,h\in I_{N}\backslash \left\{ i\right\} }\left\Vert \partial _{x_{h}x_{\ell
}}^{2}\Delta^{g}_{i,j}\right\Vert _{L^{\infty }}\Bigg ].
\end{eqnarray*}%
Next we focus on the last two items in (\ref{ev}), namely, $J_{5}$ and $J_{6}
$. However, these two terms contain $P_{t}^{i,j}$ and $Q_{t}^{i,j}.$ Our aim
is to give the estimations of solutions to BSDE (\ref{adjn1}) as follows.
From Lemma \ref{bs2}, by letting
\begin{eqnarray*}
\xi  &=&\partial _{x}\Delta _{i,j}^{g}\left( \mathbf{X}_{T}^{\mathbf{u}%
}\right) ,\text{ }\forall i,j\in I_{N}\text{ with }i<j, \\
F_{t}\left( P_{t}^{i,j},\left\{ Q_{t,k}^{i,j}\right\} _{k=1,\ldots
,N}\right)  &=&\mathbf{B}\left( t,\mathbf{X}_{t}^{\mathbf{u}},\mathbf{u}%
_{t}\right) ^{\top }P_{t}^{i,j}+\sum_{k=1}^{N}\mathbf{\Pi }^{k}\left( t,%
\mathbf{X}_{t}^{\mathbf{u}},\mathbf{u}_{t}\right) ^{\top }Q_{t,k}^{i,j} \\
&&+\left( \partial _{x}\Delta _{i,j}^{f}\right) \left( t,\mathbf{X}_{t}^{%
\mathbf{u}},\mathbf{u}_{t}\right) ,
\end{eqnarray*}%
we derive
\begin{eqnarray*}
&&\mathbb{E}\left[ \sup_{0\leq t\leq T}\left\vert P_{t}^{i,j}\right\vert
^{2}+\int_{0}^{T}\left\vert Q_{t}^{i,j}\right\vert ^{2}\mathrm{d}s\right]  \\
&\leq &C\mathbb{E}\Bigg [\left\vert \partial _{x}\Delta _{i,j}^{g}\left(
\mathbf{X}_{T}^{\mathbf{u}}\right) \right\vert ^{2}+\left(
\int_{0}^{T}\left\vert \left( \partial _{x}\Delta _{i,j}^{f}\right) \left( t,%
\mathbf{X}_{t}^{\mathbf{u}},\mathbf{u}_{t}\right) \right\vert \mathrm{d}%
t\right) ^{2}\Bigg ] \\
&\leq &C_{1}\mathbb{E}\Bigg [\left\vert \partial _{x}\Delta _{i,j}^{g}\left(
\mathbf{X}_{T}^{\mathbf{u}}\right) \right\vert ^{2}+T\int_{0}^{T}\left\vert
\left( \partial _{x}\Delta _{i,j}^{f}\right) \left( t,\mathbf{X}_{t}^{%
\mathbf{u}},\mathbf{u}_{t}\right) \right\vert ^{2}\mathrm{d}t\Bigg ],
\end{eqnarray*}%
where $C_{1}$ depends on the time $T,$ $\max \left\{ \left\vert \mathbf{B}%
\left( t,\mathbf{X}_{t}^{\mathbf{u}},\mathbf{u}_{t}\right) \right\vert
,\left\vert \mathbf{\Pi }^{1}\left( t,\mathbf{X}_{t}^{\mathbf{u}},\mathbf{u}%
_{t}\right) \right\vert ,\ldots ,\left\vert \mathbf{\Pi }^{N}\left( t,%
\mathbf{X}_{t}^{\mathbf{u}},\mathbf{u}_{t}\right) \right\vert \right\} $ with%
\begin{eqnarray*}
\left\vert \mathbf{B}\left( t,\mathbf{X}_{t}^{\mathbf{u}},\mathbf{u}%
_{t}\right) \right\vert  &\leq &\left\vert \mathbf{\bar{B}}\left( t,\mathbf{X%
}_{t}^{\mathbf{u}},\mathbf{u}_{t}\right) +\mathbf{B}_{0}\left( t,\mathbf{X}%
_{t}^{\mathbf{u}},\mathbf{u}_{t}\right) \right\vert  \\
&\leq &\left\vert \mathbf{\bar{B}}\left( t,\mathbf{X}_{t}^{\mathbf{u}},%
\mathbf{u}_{t}\right) \right\vert +\left\vert \mathbf{B}_{0}\left( t,\mathbf{%
X}_{t}^{\mathbf{u}},\mathbf{u}_{t}\right) \right\vert  \\
&\leq &L^{b}\sqrt{N}+\frac{L_{y}^{b}}{\sqrt{N}}
\end{eqnarray*}%
and%
\begin{eqnarray*}
\left\vert \mathbf{\Pi }^{j}\left( t,\mathbf{X}_{t}^{\mathbf{u}},\mathbf{u}%
_{t}\right) \right\vert  &\leq &\left\vert \mathbf{\bar{\Pi}}^{j}\left( t,%
\mathbf{X}_{t}^{\mathbf{u}},\mathbf{u}_{t}\right) +\mathbf{\Pi }%
_{0}^{j}\left( t,\mathbf{X}_{t}^{\mathbf{u}},\mathbf{u}_{t}\right)
\right\vert  \\
&\leq &\left\vert \mathbf{\bar{\Pi}}^{j}\left( t,\mathbf{X}_{t}^{\mathbf{u}},%
\mathbf{u}_{t}\right) \right\vert +\left\vert \mathbf{\Pi }_{0}^{j}\left( t,%
\mathbf{X}_{t}^{\mathbf{u}},\mathbf{u}_{t}\right) \right\vert  \\
&\leq &L^{\sigma }+\frac{L_{y}^{\sigma }}{\sqrt{N}}.
\end{eqnarray*}%
We employ the matrix norm $\left\Vert A\right\Vert _{2}$ due to its
computational simplicity. As $N\rightarrow \infty $,
\begin{equation*}
\max \left\{ \left\Vert \mathbf{B}\left( t,\mathbf{X}_{t}^{\mathbf{u}},%
\mathbf{u}_{t}\right) \right\Vert _{2},\left\Vert \mathbf{\Pi }^{1}\left( t,%
\mathbf{X}_{t}^{\mathbf{u}},\mathbf{u}_{t}\right) \right\Vert _{2},\ldots
,\left\Vert \mathbf{\Pi }^{N}\left( t,\mathbf{X}_{t}^{\mathbf{u}},\mathbf{u}%
_{t}\right) \right\Vert _{2}\right\} \rightarrow L^{b}+L^{\sigma }.
\end{equation*}%
By combining Lemma \ref{bs2}, it follows that $C_{1}$ remains bounded.

The fundamental theorem of calculus implies that for all $\left(
t,x,u\right) =\left( t,\left( x_{\ell }\right) _{\ell =1}^{N},\left( u_{\ell
}\right) _{\ell =1}^{N}\right) \in \left[ 0,T\right] \times \mathbb{R}%
^{n}\times \mathbb{R}^{n}$ and $i,j\in I_{N},$%
\begin{eqnarray*}
\left\vert \left( \partial _{x_{\ell }}\Delta _{i,j}^{f}\right) \left(
t,x,u\right) \right\vert ^{2} &\leq &3\Bigg [\left\vert \left( \partial
_{x_{\ell }}\Delta _{i,j}^{f}\right) \left( t,0,0\right) \right\vert ^{2} \\
&&+\sum_{k=1}^{N}\left( \left\Vert \partial _{x_{_{\ell }}x_{k}}\Delta
_{i,j}^{f}\right\Vert _{L^{\infty }}^{2}\left\vert x_{k}\right\vert
^{2}+\left\Vert \partial _{x_{\ell }u_{k}}\Delta _{i,j}^{f}\right\Vert
_{L^{\infty }}^{2}\left\vert u_{k}\right\vert ^{2}\right) \Bigg ], \\
\left\vert \partial _{x_{\ell }}\Delta _{i,j}^{g}\left( x\right) \right\vert
^{2} &\leq &\left\vert \left( \partial _{x_{\ell }}\Delta _{i,j}^{g}\right)
\left( 0\right) \right\vert ^{2}+\sum_{k=1}^{N}\left\Vert \partial
_{x_{_{\ell }}x_{k}}\Delta _{i,j}^{g}\right\Vert _{L^{\infty
}}^{2}\left\vert x_{k}\right\vert ^{2}.
\end{eqnarray*}%
Thus, Lemma \ref{bs2} yields
\begin{eqnarray*}
&&\mathbb{E}\left[ \sup_{0\leq t\leq T}\left\vert P_{t}^{i,j}\right\vert
^{2}+\int_{0}^{T}\left\vert Q_{t}^{i,j}\right\vert ^{2}\mathrm{d}s\right] \\
&\leq &C_{1}\mathbb{E}\Bigg [\left\vert \partial _{x}\Delta _{i,j}^{g}\left(
\mathbf{X}_{T}^{\mathbf{u}}\right) \right\vert ^{2}+\int_{0}^{T}\left\vert
\left( \partial _{x}\Delta _{i,j}^{f}\right) \left( t,\mathbf{X}_{t}^{%
\mathbf{u}},\mathbf{u}_{t}\right) \right\vert ^{2}\mathrm{d}t\Bigg ] \\
&\leq &C_{1}\mathbb{E}\Bigg [\sum_{\ell \in I_{N}}\left\vert \left( \partial
_{x_{\ell }}\Delta _{i,j}^{g}\right) \left( 0\right) \right\vert
^{2}+\sum_{\ell ,k\in I_{N}}\left\Vert \partial _{x_{\ell }x_{k}}\Delta
_{i,j}^{g}\right\Vert _{L^{\infty }}^{2} \\
&&+3\int_{0}^{T}\left( \sum_{\ell }\left\vert \left( \partial _{x_{\ell
}}\Delta _{i,j}^{f}\right) \left( t,0,0\right) \right\vert ^{2}+\sum_{\ell
,k\in I_{N}}\left\Vert \partial _{x_{\ell }x_{k}}\Delta
_{i,j}^{f}\right\Vert _{L^{\infty }}^{2}+\sum_{\ell ,k\in I_{N}}\left\Vert
\partial _{x_{\ell }u_{k}}\Delta _{i,j}^{f}\right\Vert _{L^{\infty
}}^{2}\right) \mathrm{d}t\Bigg ].
\end{eqnarray*}%
Define
\begin{eqnarray}
\Gamma _{1} &=&C_{1}\mathbb{E}\Bigg [\sum_{\ell \in I_{N}}\left\vert \left(
\partial _{x_{\ell }}\Delta _{i,j}^{g}\right) \left( 0\right) \right\vert
^{2}+\sum_{\ell ,k\in I_{N}}\left\Vert \partial _{x_{\ell }x_{k}}\Delta
_{i,j}^{g}\right\Vert _{L^{\infty }}^{2}  \notag \\
&&+3\int_{0}^{T}\Bigg (\sum_{\ell }\left\vert \left( \partial _{x_{\ell
}}\Delta _{i,j}^{f}\right) \left( t,0,0\right) \right\vert ^{2}+\sum_{\ell
,k\in I_{N}}\left\Vert \partial _{x_{\ell }x_{k}}\Delta
_{i,j}^{f}\right\Vert _{L^{\infty }}^{2}  \notag \\
&&+\sum_{\ell ,k\in I_{N}}\left\Vert \partial _{x_{\ell }u_{k}}\Delta
_{i,j}^{f}\right\Vert _{L^{\infty }}^{2}\Bigg )\mathrm{d}t\Bigg ].
\label{lamda}
\end{eqnarray}%
For $J_{5},$ recall%
\begin{eqnarray*}
\partial _{yy}^{2}\bar{b}_{t}\left( \mathbf{X}_{t}\right) \star \mathcal{Y}%
_{t}^{i,j} &=&\left[
\begin{array}{c}
\text{tr}\left[ \partial _{yy}^{2}\bar{b}_{t,1}\left( X_{t,1},\mathbf{X}%
_{t}\right) \mathcal{Y}_{t}^{i,j}\right] \\
\vdots \\
\text{tr}\left[ \partial _{yy}^{2}\bar{b}_{t,i}\left( X_{t,i},\mathbf{X}%
_{t}\right) \mathcal{Y}_{t}^{i,j}\right] \\
\vdots \\
\text{tr}\left[ \partial _{yy}^{2}\bar{b}_{t,N}\left( X_{t,N},\mathbf{X}%
_{t}\right) \mathcal{Y}_{t}^{i,j}\right]%
\end{array}%
\right] \\
&=&\left[
\begin{array}{c}
\text{tr}\left[ \partial _{yy}^{2}\bar{b}_{t,1}\left( X_{t,1},\mathbf{X}%
_{t}\right) \mathbf{Y}_{t}^{j}\left( \mathbf{Y}_{t}^{i}\right) ^{\top }%
\right] \\
\vdots \\
\text{tr}\left[ \partial _{yy}^{2}\bar{b}_{t,i}\left( X_{t,i},\mathbf{X}%
_{t}\right) \mathbf{Y}_{t}^{j}\left( \mathbf{Y}_{t}^{i}\right) ^{\top }%
\right] \\
\vdots \\
\text{tr}\left[ \partial _{yy}^{2}\bar{b}_{t,N}\left( X_{t,N},\mathbf{X}%
_{t}\right) \mathbf{Y}_{t}^{j}\left( \mathbf{Y}_{t}^{i}\right) ^{\top }%
\right]%
\end{array}%
\right] \\
&=&\left[
\begin{array}{c}
\left( \mathbf{Y}_{t}^{i}\right) ^{\top }\partial _{yy}^{2}\bar{b}%
_{t,1}\left( X_{t,1},\mathbf{X}_{t}\right) \mathbf{Y}_{t}^{j} \\
\vdots \\
\left( \mathbf{Y}_{t}^{i}\right) ^{\top }\partial _{yy}^{2}\bar{b}%
_{t,i}\left( X_{t,i},\mathbf{X}_{t}\right) \mathbf{Y}_{t}^{j} \\
\vdots \\
\left( \mathbf{Y}_{t}^{i}\right) ^{\top }\partial _{yy}^{2}\bar{b}%
_{t,N}\left( X_{t,N},\mathbf{X}_{t}\right) \mathbf{Y}_{t}^{j}%
\end{array}%
\right] .
\end{eqnarray*}%
Then
\begin{eqnarray*}
\mathbb{E}\Bigg [\int_{0}^{T}\left\langle P_{t}^{i,j},\partial _{yy}^{2}\bar{%
b}_{t}\left( \mathbf{X}_{t}\right) \star \mathcal{Y}_{t}^{i,j}\right\rangle
\mathrm{d}t\Bigg ] &=&\mathbb{E}\Bigg [\int_{0}^{T}%
\sum_{k=1}^{N}P_{t,k}^{i,j}\left( \mathbf{Y}_{t}^{i}\right) ^{\top }\partial
_{yy}^{2}\bar{b}_{t,k}\left( X_{t,k},\mathbf{X}_{t}\right) \mathbf{Y}_{t}^{j}%
\mathrm{d}t\Bigg ] \\
&=&\sum_{k=1}^{N}\mathbb{E}\Bigg [\int_{0}^{T}P_{t,k}^{i,j}\left( \mathbf{Y}%
_{t}^{i}\right) ^{\top }\partial _{yy}^{2}\bar{b}_{t,k}\left( X_{t,k},%
\mathbf{X}_{t}\right) \mathbf{Y}_{t}^{j}\mathrm{d}t\Bigg ].
\end{eqnarray*}%
For simplicity, we estimate the term $P_{t,k}^{i,j}\left( \mathbf{Y}%
_{t}^{i}\right) ^{\top }\partial _{yy}^{2}\bar{b}_{k}\left( X_{t,k},\mathbf{X%
}_{t}\right) \mathbf{Y}_{t}^{j}$ firstly.

Now%
\begin{equation*}
\mathbb{E}\left[ \int_{0}^{T}P_{t,k}^{i,j}\left( \mathbf{Y}_{t}^{i}\right)
^{\top }\partial _{yy}^{2}\bar{b}_{t,k}\left( X_{t,k},\mathbf{X}_{t}\right)
\mathbf{Y}_{t}^{j}\mathrm{d}t\right] =\mathbb{E}\left[
\int_{0}^{T}P_{t,k}^{i,j}\sum_{h,\ell =1}^{N}\left( \partial _{y_{h}y_{\ell
}}^{2}\bar{b}_{t,k}\left( X_{t,k},\mathbf{X}_{t}\right) \right)
Y_{h,t}^{i}Y_{\ell ,t}^{j}\mathrm{d}t\right] .
\end{equation*}%
Notice that $i<j$ and
\begin{eqnarray*}
I_{N}\times I_{N} &=&\left\{ \left( i,j\right) \right\} \cup \left\{ \left(
i,i\right) \right\} \cup \left\{ \left( j,j\right) \right\} \cup \left\{
\left( i,\ell \right) |\ell \in I_{N}\backslash \left\{ i,j\right\} \right\}
\\
&&\cup \left\{ \left( h,j\right) |h\in I_{N}\backslash \left\{ i,j\right\}
\right\} \cup \left\{ \left( h,\ell \right) |h\in I_{N}\backslash \left\{
i\right\} ,\ell \in I_{N}\backslash \left\{ j\right\} \right\} .
\end{eqnarray*}%
Then
\begin{eqnarray}
&&\mathbb{E}\Bigg [\int_{0}^{T}P_{t,k}^{i,j}\left( \mathbf{Y}_{t}^{i}\right)
^{\top }\partial _{yy}^{2}\bar{b}_{t,k}\left( X_{t,k},\mathbf{X}_{t}\right)
\mathbf{Y}_{t}^{j}\Bigg ]  \notag \\
&=&\mathbb{E}\Bigg \{\int_{0}^{T}P_{t,k}^{i,j}\Bigg [\underset{I_{1}}{%
\underbrace{\left( \partial _{y_{i}y_{j}}^{2}\bar{b}_{t,k}\left( X_{t,k},%
\mathbf{X}_{t}\right) \right) Y_{t,i}^{i}Y_{t,j}^{j}}}+\underset{I_{2}}{%
\underbrace{\left( \partial _{y_{i}y_{i}}^{2}\bar{b}_{t,k}\left( X_{t,k},%
\mathbf{X}_{t}\right) \right) Y_{t,i}^{i}Y_{t,i}^{j}}}  \notag \\
&&+\underset{I_{3}}{\underbrace{\left( \partial _{y_{j}y_{j}}^{2}\bar{b}%
_{t,k}\left( X_{t,k},\mathbf{X}_{t}\right) \right) Y_{t,j}^{i}Y_{t,j}^{j}}}+%
\underset{I_{4}}{\underbrace{\sum_{\ell \in I_{N}\backslash \left\{
i,j\right\} }^{N}\left( \partial _{y_{i}y_{\ell }}^{2}\bar{b}_{t,k}\left(
X_{t,k},\mathbf{X}_{t}\right) \right) Y_{t,i}^{i}Y_{t,\ell }^{j}}}  \notag \\
&&+\underset{I_{5}}{\underbrace{\sum_{h\in I_{N}\backslash \left\{
i,j\right\} }^{N}\left( \partial _{y_{h}y_{j}}^{2}\bar{b}_{t,k}\left(
X_{t,k},\mathbf{X}_{t}\right) \right) Y_{t,h}^{i}Y_{t,j}^{j}}}  \notag \\
&&+\underset{I_{6}}{\underbrace{\sum_{h\in I_{N}\backslash \left\{ i\right\}
}^{N}\sum_{\ell \in I_{N}\backslash \left\{ j\right\} }^{N}\left( \partial
_{y_{h}y_{\ell }}^{2}\bar{b}_{t,k}\left( X_{t,k},\mathbf{X}_{t}\right)
\right) Y_{t,h}^{i}Y_{t,\ell }^{j}}}\Bigg ]\Bigg \}\mathrm{d}t.  \label{est1}
\end{eqnarray}%
We deal with (\ref{est1}) above step by step.

For $I_{1}$
\begin{eqnarray}
\mathbb{E}\left[ \int_{0}^{T}P_{t,k}^{i,j}\left( \partial _{y_{i}y_{j}}^{2}%
\bar{b}_{t,k}\left( X_{t,k},\mathbf{X}_{t}\right) \right)
Y_{t,i}^{i}Y_{t,j}^{j}\mathrm{d}t\right] &\leq &\frac{\sqrt{\Gamma _{1}}%
L_{y}^{\bar{b}}}{N^{2}}\mathbb{E}\Bigg [\int_{0}^{T}Y_{t,i}^{i}Y_{t,j}^{j}%
\mathrm{d}t\Bigg ]  \notag \\
&\leq &\frac{L_{y}^{\bar{b}}\sqrt{\Gamma _{1}}C}{N^{2}}.  \label{est2}
\end{eqnarray}%
For $I_{2}$%
\begin{eqnarray}
\mathbb{E}\left[ \int_{0}^{T}P_{t,k}^{i,j}\left( \partial _{y_{i}y_{i}}^{2}%
\bar{b}_{t,k}\left( \mathbf{X}_{t}\right) \right) Y_{t,i}^{i}Y_{t,i}^{j}%
\mathrm{d}t\right] &\leq &\frac{L_{y}^{\bar{b}}}{N}\mathbb{E}\left[
\int_{0}^{T}P_{t,k}^{i,j}Y_{t,i}^{i}Y_{t,i}^{j}\mathrm{d}t\right]  \notag \\
&\leq &\frac{L_{y}^{\bar{b}}\sqrt{\Gamma _{1}}}{N}\mathbb{E}\Bigg [%
\int_{0}^{T}Y_{t,i}^{i}Y_{t,i}^{j}\mathrm{d}t\Bigg ]  \notag \\
&\leq &\frac{L_{y}^{\bar{b}}\sqrt{\Gamma _{1}}}{N}\cdot \frac{C\left( L_{y}^{%
\bar{b}}+3\left( L_{y}^{\bar{\sigma}}\right) ^{2}\right) }{N},  \label{est3}
\end{eqnarray}%
For $I_{3}$%
\begin{equation}
\mathbb{E}\left[ \int_{0}^{T}P_{t,k}^{i,j}\left( \partial _{y_{j}y_{j}}^{2}%
\bar{b}_{t,k}\left( X_{t,k},\mathbf{X}_{t}\right) \right)
Y_{t,j}^{i}Y_{t,j}^{j}\mathrm{d}t\right] \leq \frac{L_{y}^{\bar{b}}\sqrt{%
\Gamma _{1}}}{N}\cdot \frac{C\left( L_{y}^{\bar{b}}+3\left( L_{y}^{\bar{%
\sigma}}\right) ^{2}\right) }{N}.  \label{est4}
\end{equation}%
For $I_{4}$%
\begin{eqnarray}
\mathbb{E}\left[ \int_{0}^{T}P_{t,k}^{i,j}\sum_{\ell \in I_{N}\backslash
\left\{ i,j\right\} }^{N}\left( \partial _{y_{i}y_{\ell }}^{2}\bar{b}%
_{t,k}\left( X_{t,k},\mathbf{X}_{t}\right) \right) Y_{t,i}^{i}Y_{t,\ell }^{j}%
\mathrm{d}t\right] &\leq &\frac{L_{y}^{\bar{b}}\sqrt{\Gamma _{1}}}{2N^{2}}%
\mathbb{E}\left[ \int_{0}^{T}\sum_{\ell \in I_{N}\backslash \left\{
i,j\right\} }^{N}Y_{t,i}^{i}Y_{t,\ell }^{j}\mathrm{d}t\right]  \notag \\
&\leq &\frac{L_{y}^{\bar{b}}\sqrt{\Gamma _{1}}}{2N^{2}}\cdot \frac{C\left(
N-2\right) \left( L_{y}^{\bar{b}}+3\left( L_{y}^{\bar{\sigma}}\right)
^{2}\right) }{N}.  \label{est5}
\end{eqnarray}%
For $I_{5}$%
\begin{equation}
\mathbb{E}\left[ \int_{0}^{T}P_{t,k}^{i,j}\sum_{h\in I_{N}\backslash \left\{
i,j\right\} }^{N}\left( \partial _{y_{h}y_{j}}^{2}\bar{b}_{t,k}\left(
X_{t,k},\mathbf{X}_{t}\right) \right) Y_{t,h}^{i}Y_{t,j}^{j}\mathrm{d}t%
\right] \leq \frac{L_{y}^{b}\sqrt{\Gamma _{1}}}{2N^{2}}\mathbb{\cdot }\frac{%
\left( N-2\right) \left( L_{y}^{\bar{b}}+3\left( L_{y}^{\bar{\sigma}}\right)
^{2}\right) C}{N}.  \label{est6}
\end{equation}%
As for the last term in (\ref{est1}), we have
\begin{eqnarray}
&&\mathbb{E}\left[ \int_{0}^{T}P_{t,k}^{i,j}\sum_{h\in I_{N}\backslash
\left\{ i\right\} }^{N}\sum_{\ell \in I_{N}\backslash \left\{ j\right\}
}^{N}\left( \partial _{y_{h}y_{\ell }}^{2}\bar{b}_{t,k}\left( X_{t,k},%
\mathbf{X}_{t}\right) \right) Y_{t,h}^{i}Y_{t,\ell }^{j}\mathrm{d}t\right]
\notag \\
&\leq &\mathbb{E}\Bigg [\int_{0}^{T}P_{t,k}^{i,j}\Bigg (\sum_{\ell \in
I_{N}\backslash \left\{ j\right\} }^{N}\left( \partial _{y_{j}y_{\ell }}^{2}%
\bar{b}_{t,k}\left( X_{t,k},\mathbf{X}_{t}\right) \right)
Y_{t,j}^{i}Y_{t,\ell }^{j}  \notag \\
&&+\sum_{h\in I_{N}\backslash \left\{ i,j\right\} }^{N}\Bigg (\left(
\partial _{y_{h}y_{h}}^{2}\bar{b}_{t,k}\left( X_{t,k},\mathbf{X}_{t}\right)
\right) Y_{t,h}^{i}Y_{t,h}^{j}+\sum_{\ell \in I_{N}\backslash \left\{
j,h\right\} }^{N}\left( \partial _{y_{h}y_{h}}^{2}\bar{b}_{t,k}\left(
X_{t,k},\mathbf{X}_{t}\right) \right) Y_{t,h}^{i}Y_{t,\ell }^{j}\Bigg )\Bigg
)\mathrm{d}t\Bigg ]  \notag \\
&\leq &\frac{L_{y}^{\bar{b}}\sqrt{\Gamma _{1}}}{N^{2}}\mathbb{E}\left[
\int_{0}^{T}\sum_{\ell \in I_{N}\backslash \left\{ j\right\}
}^{N}Y_{t,j}^{i}Y_{t,\ell }^{j}\mathrm{d}t\right] +\frac{L_{y}^{\bar{b}}%
\sqrt{\Gamma _{1}}}{N}\mathbb{E}\left[ \int_{0}^{T}\sum_{h\in
I_{N}\backslash \left\{ i,j\right\} }^{N}Y_{t,h}^{i}Y_{t,h}^{j}\mathrm{d}t%
\right]  \notag \\
&&+\frac{L_{y}^{\bar{b}}\sqrt{\Gamma _{1}}}{N^{2}}\mathbb{E}\Bigg [%
\int_{0}^{T}\sum_{h\in I_{N}\backslash \left\{ i,j\right\} }^{N}\sum_{\ell
\in I_{N}\backslash \left\{ j,h\right\} }^{N}Y_{t,h}^{i}Y_{t,\ell }^{j}%
\mathrm{d}t\Bigg ]  \notag \\
&\leq &\Bigg [\frac{L_{y}^{\bar{b}}}{N^{2}}\mathbb{\cdot }\frac{\sqrt{\Gamma
_{1}}\left( N-1\right) C\left( L_{y}^{\bar{b}}+3\left( L_{y}^{\bar{\sigma}%
}\right) ^{2}\right) ^{2}}{N^{2}}+\frac{L_{y}^{\bar{b}}}{N}\mathbb{\cdot }%
\frac{\left( N-2\right) C\sqrt{\Gamma _{1}}\left( L_{y}^{\bar{b}}+3\left(
L_{y}^{\bar{\sigma}}\right) ^{2}\right) ^{2}}{N^{2}}  \notag \\
&&+\frac{L_{y}^{\bar{b}}}{N^{2}}\cdot \frac{\left( N-2\right) ^{2}\sqrt{%
\Gamma _{1}}\left( L_{y}^{\bar{b}}+3\left( L_{y}^{\bar{\sigma}}\right)
^{2}\right) ^{2}}{N^{2}}\Bigg ].  \label{est8}
\end{eqnarray}%
Combining (\ref{est2})-(\ref{est8}), we have
\begin{eqnarray*}
\mathbb{E}\Bigg [\int_{0}^{T}P_{t,k}^{i,j}\left( \mathbf{Y}_{t}^{i}\right)
^{\top }\partial _{yy}^{2}\bar{b}_{t,k}\left( X_{t,k},\mathbf{X}_{t}\right)
\mathbf{Y}_{t}^{j}\Bigg ] &\leq &\sqrt{\Gamma _{1}}\Bigg [\frac{L_{y}^{\bar{b%
}}}{N^{2}}C+\frac{C\left( L_{y}^{\bar{b}}+3\left( L_{y}^{\bar{\sigma}%
}\right) ^{2}\right) L_{y}^{\bar{b}}}{N^{2}} \\
&&+\frac{L_{y}^{\bar{b}}C\left( L_{y}^{\bar{b}}+3\left( L_{y}^{\bar{\sigma}%
}\right) ^{2}\right) }{N^{2}}+\frac{L_{y}^{\bar{b}}C\left( L_{y}^{\bar{b}%
}+3\left( L_{y}^{\bar{\sigma}}\right) ^{2}\right) }{2N^{2}} \\
&&+\frac{L_{y}^{\bar{b}}C\left( L_{y}^{\bar{b}}+3\left( L_{y}^{\bar{\sigma}%
}\right) ^{2}\right) }{2N^{2}}+\Bigg [\frac{CL_{y}^{\bar{b}}\left( L_{y}^{%
\bar{b}}+3\left( L_{y}^{\bar{\sigma}}\right) ^{2}\right) ^{2}}{N^{3}} \\
&&+\frac{L_{y}^{\bar{b}}C\left( L_{y}^{\bar{b}}+3\left( L_{y}^{\bar{\sigma}%
}\right) ^{2}\right) ^{2}}{N^{2}}+\frac{L_{y}^{\bar{b}}\left( L_{y}^{\bar{b}%
}+3\left( L_{y}^{\bar{\sigma}}\right) ^{2}\right) ^{2}}{N^{2}}\Bigg ].
\end{eqnarray*}
Therefore, we derive
\begin{eqnarray}
&&\sum_{k=1}^{N}\mathbb{E}\Bigg [\int_{0}^{T}\left\langle
P_{t}^{i,j},\partial _{yy}^{2}\bar{b}_{t,k}\left( X_{t,k},\mathbf{X}%
_{t}\right) \star \mathcal{Y}_{t}^{i,j}\right\rangle \mathrm{d}t\Bigg ]
\notag \\
&\leq &CL_{y}^{\bar{b}}\sqrt{\Gamma _{1}}\Bigg [\frac{1+\left( L_{y}^{\bar{b}%
}+3\left( L_{y}^{\bar{\sigma}}\right) ^{2}\right) }{N}+\frac{\left( L_{y}^{%
\bar{b}}+3\left( L_{y}^{\bar{\sigma}}\right) ^{2}\right) ^{2}}{N}+\frac{%
\left( L_{y}^{\bar{b}}+3\left( L_{y}^{\bar{\sigma}}\right) ^{2}\right) ^{2}}{%
N^{2}}\Bigg ].  \label{p11}
\end{eqnarray}%
Now we consider
\begin{eqnarray*}
&&\mathbb{E}\left[ \int_{0}^{T}\left\langle P_{t}^{i,j},\mathfrak{F}%
_{t}^{i,j}+\Gamma _{t}^{i,j}\left( \bar{b}_{t}\right) \right\rangle \mathrm{d%
}t\right] \\
&=&\mathbb{E}\Bigg \{\int_{0}^{T}\Big [\sum_{\ell =1}^{N}\underset{L_{1}}{%
\underbrace{P_{t,\ell }^{i,j}\Big (Y_{t,\ell }^{i}\partial _{xx}^{2}\bar{b}%
_{t,\ell }\left( X_{t,\ell },\mathbf{X}_{t}\right) Y_{t,\ell }^{j}}} \\
&&+\underset{L_{2}}{\underbrace{\left( \mathbf{Y}_{t}^{i}\right) ^{\top
}\partial _{yx}^{2}\bar{b}_{t,\ell }\left( X_{t,\ell },\mathbf{X}_{t}\right)
Y_{t,\ell }^{j}}}+\underset{L_{3}}{\underbrace{Y_{t,\ell }^{i}\partial
_{xy}^{2}\bar{b}_{t,\ell }\left( X_{t,\ell },\mathbf{X}_{t}\right) \mathbf{Y}%
_{t}^{j}}} \\
&&+\underset{L_{4}}{\underbrace{\left[ \left( \partial _{x}\phi _{t,\ell
}^{\prime }\right) \left( X_{t,\ell },\mathbf{X}_{t}^{\phi }\right)
Y_{t,\ell }^{j}+\left( \mathbf{Y}_{t}^{j}\right) ^{\top }\left( \partial
_{y}\phi _{t,\ell }^{\prime }\right) \left( X_{t,\ell },\mathbf{X}%
_{t}\right) \right] \delta _{\ell =i}}} \\
&&+\underset{L_{5}}{\underbrace{\left[ \left( \partial _{x}\phi _{t,\ell
}^{\prime \prime }\right) \left( X_{t,\ell },\mathbf{X}_{t}\right) Y_{t,\ell
}^{i}+\left( \mathbf{Y}_{t}^{i}\right) ^{\top }\left( \partial _{y}\phi
_{t,\ell }^{\prime \prime }\right) \left( X_{t,\ell },\mathbf{X}_{t}\right) %
\right] \delta _{\ell =j}}}\Big )\Big ]\mathrm{d}t\Bigg \}.
\end{eqnarray*}%
For $L_{1},$%
\begin{eqnarray}
&&\mathbb{E}\Bigg [\int_{0}^{T}\sum_{\ell =1}^{N}P_{t,\ell }^{i,j}Y_{t,\ell
}^{i}\partial _{xx}^{2}\bar{b}_{t,\ell }\left( X_{t,\ell },\mathbf{X}%
_{t}\right) Y_{t,\ell }^{j}\mathrm{d}t\Bigg ]  \notag \\
&\leq &L^{\bar{b}}\sqrt{\Gamma _{1}}\mathbb{E}\Bigg [\int_{0}^{T}\Big (%
P_{t,\ell }^{i,j}Y_{t,i}^{i}Y_{t,i}^{j}+P_{t,\ell
}^{i,j}Y_{t,j}^{i}Y_{t,j}^{j}+\sum_{\ell =I_{N}\backslash \left\{
i,j\right\} }^{N}P_{t,\ell }^{i,j}Y_{t,\ell }^{i}Y_{t,\ell }^{j}\Big )%
\mathrm{d}t\Bigg ]  \notag \\
&\leq &CL^{\bar{b}}\sqrt{\Gamma _{1}}\Bigg [\frac{\left( L_{y}^{\bar{b}%
}+3\left( L_{y}^{\bar{\sigma}}\right) ^{2}\right) }{N}+\frac{\left(
N-2\right) \left( L_{y}^{\bar{b}}+3\left( L_{y}^{\bar{\sigma}}\right)
^{2}\right) ^{2}}{N^{2}}\Bigg ]  \notag \\
&\leq &CL^{\bar{b}}\sqrt{\Gamma _{1}}\Bigg [\frac{\left( L_{y}^{\bar{b}%
}+3\left( L_{y}^{\bar{\sigma}}\right) ^{2}\right) ^{2}}{N}\Bigg ].
\label{L1}
\end{eqnarray}%
For $L_{2}$%
\begin{eqnarray}
&&\mathbb{E}\Bigg [\int_{0}^{T}\sum_{\ell =1}^{N}P_{t,\ell }^{i,j}\left(
\mathbf{Y}_{t}^{i}\right) ^{\top }\partial _{yx}^{2}\bar{b}_{t,\ell }\left(
X_{t,\ell },\mathbf{X}_{t}\right) Y_{t,\ell }^{j}\mathrm{d}t\Bigg ]  \notag
\\
&\leq &\frac{L_{y}^{\bar{b}}}{N}\mathbb{E}\Bigg [\int_{0}^{T}\sum_{\ell
=1}^{N}P_{t,\ell }^{i,j}\sum_{k=1}^{N}Y_{t,k}^{i}Y_{t,\ell }^{j}\mathrm{d}t%
\Bigg ]  \notag \\
&\leq &\frac{L_{y}^{\bar{b}}}{N}\mathbb{E}\Bigg [\int_{0}^{T}\Big (%
P_{t,j}^{i,j}Y_{t,i}^{i}Y_{t,j}^{j}+\sum_{k\in I_{N}\backslash \left\{
j\right\} }^{N}P_{t,k}^{i,j}Y_{t,i}^{i}Y_{t,k}^{j}+\sum_{k\in
I_{N}\backslash \left\{ i\right\} }^{N}P_{t,j}^{i,j}Y_{t,k}^{i}Y_{t,j}^{j}
\notag \\
&&+\sum_{k\in I_{N}\backslash \left\{ j\right\} ,\text{ }h\in
I_{N}\backslash \left\{ i\right\} }^{N}P_{t,k}^{i,j}Y_{t,h}^{i}Y_{t,k}^{j}%
\Big )\mathrm{d}t\Bigg ]  \notag \\
&\leq &\frac{L_{y}^{\bar{b}}\sqrt{\Gamma _{1}}C}{N}\Bigg (1+\frac{2\left(
N-1\right) \left( L_{y}^{\bar{b}}+3\left( L_{y}^{\bar{\sigma}}\right)
^{2}\right) }{N}+\frac{\left( N^{2}-2N+1\right) \left( L_{y}^{\bar{b}%
}+3\left( L_{y}^{\bar{\sigma}}\right) ^{2}\right) ^{2}}{N^{2}}\Bigg )  \notag
\\
&\leq &\frac{L_{y}^{\bar{b}}\sqrt{\Gamma _{1}}C}{N}\Bigg (1+\left( L_{y}^{%
\bar{b}}+3\left( L_{y}^{\bar{\sigma}}\right) ^{2}\right) +\left( L_{y}^{\bar{%
b}}+3\left( L_{y}^{\bar{\sigma}}\right) ^{2}\right) ^{2}\Bigg ).  \label{L2}
\end{eqnarray}%
Analogously, for $L_{3}$%
\begin{eqnarray}
&&\mathbb{E}\Bigg [\int_{0}^{T}\sum_{\ell =1}^{N}P_{t,\ell }^{i,j}Y_{t,\ell
}^{i}\partial _{xy}^{2}\bar{b}_{t,\ell }\left( X_{t,\ell },\mathbf{X}%
_{t}\right) \mathbf{Y}_{t}^{j}\mathrm{d}t\Bigg ]  \notag \\
&\leq &\frac{L_{y}^{\bar{b}}\sqrt{\Gamma _{1}}C}{N}\Bigg (1+\left( L_{y}^{%
\bar{b}}+3\left( L_{y}^{\bar{\sigma}}\right) ^{2}\right) +\left( L_{y}^{\bar{%
b}}+3\left( L_{y}^{\bar{\sigma}}\right) ^{2}\right) ^{2}\Bigg ).  \label{L3}
\end{eqnarray}%
For $L_{4}$%
\begin{eqnarray}
&&\mathbb{E}\Bigg \{\int_{0}^{T}\sum_{\ell =1}^{N}P_{t,\ell }^{i,j}\left[
\left( \partial _{x}\phi _{t,\ell }^{\prime }\right) \left( X_{t,\ell },%
\mathbf{X}_{t}\right) Y_{t,\ell }^{j}+\left( \mathbf{Y}_{t}^{j}\right)
^{\top }\left( \partial _{y}\phi _{t,\ell }^{\prime }\right) \left(
X_{t,\ell },\mathbf{X}_{t}\right) \right] \delta _{\{\ell =i\}}\mathrm{d}t\Bigg
\}  \notag \\
&=&\mathbb{E}\Bigg [\int_{0}^{T}P_{t,i}^{i,j}\left( \left( \partial _{x}\phi
_{t,i}^{\prime }\right) \left( X_{t,i},\mathbf{X}_{t}\right)
Y_{t,i}^{j}+\left( \mathbf{Y}_{t}^{j}\right) ^{\top }\left( \partial
_{y}\phi _{t,i}^{\prime }\right) \left( X_{t,i},\mathbf{X}_{t}\right)
\right) \mathrm{d}t\Bigg ]  \notag \\
&\leq &L^{\phi _{i}^{\prime }}\mathbb{E}\left[
\int_{0}^{T}P_{t,i}^{i,j}Y_{t,i}^{j}\mathrm{d}t\right] +\mathbb{E}\left[
\int_{0}^{T}P_{t,i}^{i,j}\sum_{k=1}^{N}\left( \partial _{uy_{k}}^{2}\phi
_{t,i}^{\prime }\left( X_{t,\ell },\mathbf{X}_{t}\right) Y_{t,k}^{j}\right)
\mathrm{d}t\right]  \notag \\
&\leq &L^{\phi _{i}^{\prime }}\sqrt{\Gamma _{1}}\frac{\left( L_{y}^{\bar{b}%
}+3\left( L_{y}^{\bar{\sigma}}\right) ^{2}\right) }{N}  \notag \\
&&+\mathbb{E}\Bigg [\int_{0}^{T}P_{t,i}^{i,j}\Big (\partial
_{uy_{j}}^{2}\phi _{t,i}^{\prime }\left( X_{t,\ell },\mathbf{X}_{t}\right)
Y_{t,j}^{j}+\sum_{k\in I_{N}\backslash \left\{ j\right\} }\partial
_{uy_{k}}^{2}\phi _{t,i}^{\prime }\left( X_{t,\ell },\mathbf{X}_{t}\right)
Y_{t,k}^{j}\Big )\mathrm{d}t\Bigg ]  \notag \\
&\leq &\sqrt{\Gamma _{1}}\left( L_{y}^{\bar{b}}+3\left( L_{y}^{\bar{\sigma}%
}\right) ^{2}\right) C\left( \frac{L^{\phi _{i}^{\prime }}+L_{y}^{\phi
_{i}^{\prime }}}{N}\right) .  \label{L4}
\end{eqnarray}%
For $L_{5}$%
\begin{eqnarray}
&&\mathbb{E}\Bigg \{\int_{0}^{T}\sum_{\ell =1}^{N}P_{t,\ell }^{i,j}\left[
\left( \partial _{x}\phi _{t,\ell }^{\prime \prime }\right) \left( X_{t,\ell
},\mathbf{X}_{t}\right) Y_{t,\ell }^{i}+\left( \mathbf{Y}_{t}^{i}\right)
^{\top }\left( \partial _{y}\phi _{t,\ell }^{\prime \prime }\right) \left(
X_{t,\ell },\mathbf{X}_{t}\right) \right] \delta _{\{\ell =j\} }\mathrm{d}t\Bigg
\}  \notag \\
&\leq &\sqrt{\Gamma _{1}}\left( L_{y}^{\bar{b}}+3\left( L_{y}^{\bar{\sigma}%
}\right) ^{2}\right) C\left( \frac{L^{\phi _{j}^{\prime \prime
}}+L_{y}^{\phi _{j}^{\prime \prime }}}{N}\right) .  \label{L5}
\end{eqnarray}%
From (\ref{L1})-(\ref{L5}), we are able to get
\begin{eqnarray}
\mathbb{E}\left[ \int_{0}^{T}\left\langle P_{t}^{i,j},\mathfrak{F}%
_{t}^{i,j}+\Gamma _{t}^{i,j}\left( \bar{b}_{t}\right) \right\rangle \mathrm{d%
}t\right]  &\leq &C\sqrt{\Gamma _{1}}\Bigg \{\frac{L^{\bar{b}}\left( L_{y}^{%
\bar{b}}+3\left( L_{y}^{\bar{\sigma}}\right) ^{2}\right) ^{2}}{N}  \notag \\
&&+\frac{2L_{y}^{\bar{b}}}{N}\Big [1+\left( L_{y}^{\bar{b}}+3\left( L_{y}^{%
\bar{\sigma}}\right) ^{2}\right) +\left( L_{y}^{\bar{b}}+3\left( L_{y}^{\bar{%
\sigma}}\right) ^{2}\right) ^{2}\Big ]  \notag \\
&&+\left( L_{y}^{\bar{b}}+3\left( L_{y}^{\bar{\sigma}}\right) ^{2}\right)
\left( \frac{L^{\phi _{i}^{\prime }}+L_{y}^{\phi _{i}^{\prime }}+L^{\phi
_{j}^{\prime \prime }}+L_{y}^{\phi _{j}^{\prime \prime }}}{N}\right) \Bigg \}%
.  \label{p22}
\end{eqnarray}
Therefore, from (\ref{p11})-(\ref{p22}), it follows%
\begin{eqnarray}
&&\mathbb{E}\left[ \int_{0}^{T}\left\langle P_{t}^{i,j},\mathfrak{F}%
_{t}^{i,j}+\partial _{yy}^{2}\bar{b}\left( t,\mathbf{X}_{t},\mathbf{u}%
_{t}\right) \star \mathcal{Y}_{t}^{i,j}+\Gamma _{t}^{i,j}\left( \bar{b}%
_{t}\right) \right\rangle \mathrm{d}t\right]  \notag \\
&\leq &C\sqrt{\Gamma _{1}}\Bigg \{\frac{L_{y}^{\bar{b}}}{N}+\frac{\left(
L_{y}^{\bar{b}}+3\left( L_{y}^{\bar{\sigma}}\right) ^{2}\right) L_{y}^{\bar{b%
}}}{N}  \notag \\
&&+\frac{L_{y}^{\bar{b}}\left( L_{y}^{\bar{b}}+3\left( L_{y}^{\bar{\sigma}%
}\right) ^{2}\right) }{N}+\frac{L_{y}^{\bar{b}}\left( L_{y}^{\bar{b}%
}+3\left( L_{y}^{\bar{\sigma}}\right) ^{2}\right) ^{2}}{N^{2}}  \notag \\
&&+\frac{L_{y}^{\bar{b}}\left( L_{y}^{\bar{b}}+3\left( L_{y}^{\bar{\sigma}%
}\right) ^{2}\right) ^{2}}{N}+\frac{L^{\bar{b}}\left( L_{y}^{\bar{b}%
}+3\left( L_{y}^{\bar{\sigma}}\right) ^{2}\right) ^{2}}{N}  \notag \\
&&+\frac{2L_{y}^{\bar{b}}}{N}\Big [1+\left( L_{y}^{\bar{b}}+3\left( L_{y}^{%
\bar{\sigma}}\right) ^{2}\right) +\left( L_{y}^{\bar{b}}+3\left( L_{y}^{\bar{%
\sigma}}\right) ^{2}\right) ^{2}\Big ]  \notag \\
&&+\left( L_{y}^{\bar{b}}+3\left( L_{y}^{\bar{\sigma}}\right) ^{2}\right)
\left( \frac{L^{\phi _{i}^{\prime }}+L_{y}^{\phi _{i}^{\prime }}+L^{\phi
_{j}^{\prime \prime }}+L_{y}^{\phi _{j}^{\prime \prime }}}{N}\right) \Bigg \}%
.  \label{m3}
\end{eqnarray}%
Next we will handle $J_{6}.$ Put%
\begin{eqnarray*}
&&\left( \mathfrak{G}_{t,k}^{i,j}+\bar{\Xi}_{t,k}^{i,j}\left( \bar{\sigma}%
_{t}\right) \right) _{k} \\
&=&Y_{t,k}^{i}\partial _{xx}^{2}\left( \bar{\sigma}_{t,k}\left( X_{t,k},%
\mathbf{X}_{t}\right) \right) Y_{t,k}^{j}+\left( \mathbf{Y}_{t}^{i}\right)
^{\top }\partial _{yx}^{2}\left( \bar{\sigma}_{t,k}\left( X_{t,k},\mathbf{X}%
_{t}\right) \right) Y_{t,k}^{j}+Y_{t,k}^{i}\partial _{xy}^{2}\left( \bar{%
\sigma}_{t,k}\left( X_{t,k},\mathbf{X}_{t}\right) \right) \mathbf{Y}_{t}^{j}
\\
&&+\delta _{k,i}\left( \left( \partial _{x}\phi _{t,k}^{\prime }\right)
\left( X_{t,k},\mathbf{X}_{t}\right) Y_{t,k}^{j}+\left( \mathbf{Y}%
_{t}^{j}\right) ^{\top }\left( \partial _{y}\phi _{t,k}^{\prime }\right)
\left( X_{t,k},\mathbf{X}_{t}\right) \right) \\
&&+\delta _{k,j}\left( \left( \partial _{x}\phi _{t,k}^{\prime \prime
}\right) \left( X_{t,k},\mathbf{X}_{t}\right) Y_{t,k}^{i}+\left( \mathbf{Y}%
_{t}^{i}\right) ^{\top }\left( \partial _{y}\phi _{t,k}^{\prime \prime
}\right) \left( X_{t,k},\mathbf{X}_{t}\right) \right) .
\end{eqnarray*}%
Then $\left\langle Q_{t,k}^{i,j},\mathfrak{G}_{t,k}^{i,j}+\bar{\Xi}%
_{t,k}^{i,j}\left( \bar{\sigma}\right) \right\rangle
=\sum_{k=1}^{N}Q_{t,k,k}^{i,j}\left( \mathfrak{G}_{t,k}^{i,j}+\bar{\Xi}%
_{t,k}^{i,j}\left( \bar{\sigma}\right) \right) _{k}.$ Repeating the method
developed above, we instantly have%
\begin{eqnarray*}
&&\mathbb{E}\left[ \int_{0}^{T}\sum_{k=1}^{N}\left\langle Q_{t,k}^{i,j},%
\mathfrak{G}_{t,k}^{i,j}+\bar{\Xi}_{t,k}^{i,j}\left( \bar{\sigma}_{t}\right)
\right\rangle \mathrm{d}t\right] \\
&=&\sum_{k=1}^{N}\mathbb{E}\left[ \int_{0}^{T}Q_{t,k,k}^{i,j}\left(
\mathfrak{G}_{t,k}^{i,j}+\bar{\Xi}_{t,k}^{i,j}\left( \bar{\sigma}_{t}\right)
\right) _{k}\mathrm{d}t\right] \\
&\leq &\sum_{k=1}^{N}\left[ \mathbb{E}\left[ \int_{0}^{T}\left(
Q_{t,k,k}^{i,j}\right) ^{2}\mathrm{d}t\right] ^{\frac{1}{2}}\cdot \mathbb{E}%
\left[ \int_{0}^{T}\left( \mathfrak{G}_{t,k}^{i,j}+\bar{\Xi}%
_{t,k}^{i,j}\left( \bar{\sigma}_{t}\right) \right) _{k}^{2}\mathrm{d}t\right]
^{\frac{1}{2}}\right] \\
&=&\sqrt{\Gamma _{1}}\mathbb{E}\left[ \sum_{k=1}^{N}\left(
\int_{0}^{T}\left( \mathfrak{G}_{t,k}^{i,j}+\bar{\Xi}_{t,k}^{i,j}\left( \bar{%
\sigma}_{t}\right) \right) _{k}^{2}\mathrm{d}t\right) ^{\frac{1}{2}}\right]
\\
&\leq &C\sqrt{\Gamma _{1}}\Bigg \{\frac{L^{\bar{\sigma}}\left( L_{y}^{\bar{b}%
}+3\left( L_{y}^{\bar{\sigma}}\right) ^{2}\right) ^{2}}{N} \\
&&+\frac{2L_{y}^{\bar{\sigma}}}{N}\Big [1+\left( L_{y}^{\bar{b}}+3\left(
L_{y}^{\bar{\sigma}}\right) ^{2}\right) +\left( L_{y}^{\bar{b}}+3\left(
L_{y}^{\bar{\sigma}}\right) ^{2}\right) ^{2}\Big ] \\
&&+\left( L_{y}^{\bar{b}}+3\left( L_{y}^{\bar{\sigma}}\right) ^{2}\right)
\left( \frac{L^{\phi _{i}^{\prime }}+L_{y}^{\phi _{i}^{\prime }}+L^{\phi
_{j}^{\prime \prime }}+L_{y}^{\phi _{j}^{\prime \prime }}}{N}\right) \Bigg \}%
,
\end{eqnarray*}%
where the third line follows by H\"{o}lder inequality while the fifth line
follows by noting $\sqrt{a_{1}+\cdots +a_{N}}\leq \sqrt{a_{1}}+\cdots +\sqrt{%
a_{N}}$ for any $a_{1},\cdots ,a_{N}\geq 0.$

We now deal with
\begin{eqnarray}
&&\mathbb{E}\left[ \int_{0}^{T}\sum_{k=1}^{N}\left\langle
Q_{t,k}^{i,j},\partial _{yy}^{2}\bar{\sigma}_{t,k}\left( X_{t,k},\mathbf{X}%
_{t}\right) \star \mathcal{Y}_{t}^{i,j}\right\rangle \mathrm{d}t\right]
\notag \\
&=&\sum_{k=1}^{N}\mathbb{E}\left[ \int_{0}^{T}\left\langle
Q_{t,k}^{i,j},\partial _{yy}^{2}\bar{\sigma}_{t,k}\left( X_{t,k},\mathbf{X}%
_{t}\right) \star \mathcal{Y}_{t}^{i,j}\right\rangle \mathrm{d}t\right]
\notag \\
&=&\sum_{k=1}^{N}\mathbb{E}\left[ \int_{0}^{T}Q_{t,k,k}^{i,j}\left( \mathbf{Y%
}_{t}^{i}\right) ^{\top }\partial _{yy}^{2}\bar{b}_{t,i}\left( X_{t,i},%
\mathbf{X}_{t}\right) \mathbf{Y}_{t}^{j}\mathrm{d}t\right]  \notag \\
&\leq &\sum_{k=1}^{N}\left[ \mathbb{E}\left[ \int_{0}^{T}\left(
Q_{t,k,k}^{i,j}\right) ^{2}\mathrm{d}t\right] ^{\frac{1}{2}}\cdot \mathbb{E}%
\left[ \int_{0}^{T}\left( \left( \mathbf{Y}_{t}^{i}\right) ^{\top }\partial
_{yy}^{2}\bar{b}_{t,i}\left( X_{t,i},\mathbf{X}_{t}\right) \mathbf{Y}%
_{t}^{j}\right) _{k}^{2}\mathrm{d}t\right] ^{\frac{1}{2}}\right]  \notag \\
&\leq &L_{y}^{\bar{\sigma}}C\sqrt{\Gamma _{1}}\Bigg [\frac{1+\left( L_{y}^{%
\bar{b}}+3\left( L_{y}^{\bar{\sigma}}\right) ^{2}\right) }{N}+\frac{\left(
L_{y}^{\bar{b}}+3\left( L_{y}^{\bar{\sigma}}\right) ^{2}\right) ^{2}}{N}+%
\frac{\left( L_{y}^{\bar{b}}+3\left( L_{y}^{\bar{\sigma}}\right) ^{2}\right)
^{2}}{N^{2}}\Bigg ].  \label{m5}
\end{eqnarray}%
Thus combining (\ref{lamda})-(\ref{m5}), we have the desired results. \hfill $%
\Box $

\section{Concluding remarks}
\label{sect6}
In this paper, we utilize the BSDE method to address closed-loop $\alpha$-potential stochastic differential games. The potential function and the second-order linear derivative of the value function for player $i$ can be characterized by introducing the first- and second-order sensitivity processes. Subsequently, leveraging the duality principle, we provide a precise estimation of $\alpha$. Furthermore, the game with mean-field type interaction has been examined. From Theorem \ref{the1}, we can derive several key descriptions of the properties of $\alpha$, which in turn lay a solid foundation for further investigation into the LQ games.

Obtaining meaningful lower bounds for \(\alpha\) in the closed‑loop setting is also an interesting open problem. It would require constructing policy perturbations and quantifying the asymmetry of second‑order variations from the joint distribution of sensitivity processes and adjoint BSDEs, which we leave for future work.


\section{Appendix}

\label{APP}

\subsection{Proofs of some technical results\label{app}}

\paragraph{Proof of Lemma \protect\ref{l2}.}

Applying It\^{o}'s formula to $\left\vert X_{t,i}\right\vert ^{p}$, we get%
\begin{eqnarray}
\left\vert X_{t,i}\right\vert ^{p} &=&\left\vert \xi _{i}\right\vert
^{p}+\int_{0}^{t}\Big (p\left\vert X_{s,i}\right\vert ^{p-1}\bar{b}%
_{s,i}\left( X_{s,i},\mathbf{X}_{s}^{\phi }\right)  \notag \\
&&+\frac{p\left( p-1\right) }{2}\left\vert X_{s,i}\right\vert _{s,i}^{p-2}%
\bar{\sigma}_{s,i}\left( X_{s,i},\mathbf{X}_{s}^{\phi }\right) ^{2}\Big )%
\mathrm{d}s  \notag \\
&&+\int_{0}^{t}p\left\vert X_{s,i}\right\vert ^{p-1}\bar{\sigma}_{s,i}\left(
X_{s,i},\mathbf{X}_{s}^{\phi }\right) \mathrm{d}W_{s}^{i}.  \label{ineq1}
\end{eqnarray}%
Taking the expectation of both sides of (\ref{ineq1}) and noting the fact
that the stochastic integral term $\int_{0}^{t}p\left\vert
X_{s,i}\right\vert ^{p-1}\bar{\sigma}_{s,i}\left( X_{s,i},\mathbf{X}%
_{s}^{\phi }\right) \mathrm{d}W_{s}^{i}$ is a martingale, we get
\begin{eqnarray*}
\mathbb{E}\left[ \left\vert X_{t,i}\right\vert ^{p}\right] &=&\mathbb{E}%
\Bigg [\left\vert \xi _{i}\right\vert ^{p}+\int_{0}^{t}\Big (p\left\vert
X_{s,i}\right\vert ^{p-1}\bar{b}_{s,i}\left( X_{s,i},\mathbf{X}_{s}^{\phi
}\right) \\
&&+\frac{p\left( p-1\right) }{2}\left\vert X_{s,i}\right\vert ^{p-2}\bar{%
\sigma}_{s,i}\left( X_{s,i},\mathbf{X}_{s}^{\phi }\right) ^{2}\Big )\mathrm{d%
}s\Bigg ] \\
&\leq &\mathbb{E}\Bigg \{\left\vert \xi _{i}\right\vert ^{p}+\int_{0}^{t}%
\Bigg [p\left\vert X_{s,i}\right\vert ^{p-1}\left( L^{\bar{b}}\left(
1+\left\vert X_{s,i}\right\vert \right) +\frac{L_{y}^{\bar{b}}}{N}%
\sum_{k=1}^{N}\left\vert X_{s,k}\right\vert \right) \\
&&+\frac{p\left( p-1\right) }{2}\left\vert X_{s,i}\right\vert ^{p-2}\left(
L^{\bar{\sigma}}\left( 1+\left\vert X_{s,i}\right\vert \right) +\frac{L_{y}^{%
\bar{\sigma}}}{N}\sum_{k=1}^{N}\left\vert X_{s,k}\right\vert \right) ^{2}%
\Bigg ]\mathrm{d}s\Bigg \}.
\end{eqnarray*}%
By virtue of Young's inequality, for any $a,b\geq 0,$ we have $ab\leq \frac{%
p-1}{p}a^{\frac{p}{p-1}}+\frac{1}{p}b^{p}$ and $ab\leq \frac{p-2}{p}a^{\frac{%
p}{p-2}}+\frac{2}{p}b^{\frac{p}{2}}$ whenever $p>2.$ Thus%
\begin{eqnarray}
\mathbb{E}\left[ \left\vert X_{t,i}\right\vert ^{p}\right] &\leq &\mathbb{E}%
\Bigg \{\left\vert \xi _{i}\right\vert ^{p}+\int_{0}^{t}\Bigg [L^{\bar{b}%
}\left( \left( 2p-1\right) \left\vert X_{s,i}\right\vert ^{p}+1\right)
+L_{y}^{b}\left( p-1\right) \left\vert X_{s,i}\right\vert ^{p}  \notag \\
&&+\frac{L_{y}^{b}}{N}\sum_{k=1}^{N}\left\vert X_{s,k}\right\vert ^{p}+\frac{%
3p\left( p-1\right) }{2}\left\vert X_{s,i}\right\vert ^{p-2}\Bigg (\left(
L^{\sigma }\right) ^{2}+\left( L^{\sigma }\right) ^{2}\left\vert
X_{s,i}\right\vert ^{2}  \notag \\
&&+\left( \frac{L_{y}^{\sigma }}{N}\right) ^{2}\left(
\sum_{k=1}^{N}\left\vert X_{s,k}\right\vert \right) ^{2}\Bigg )\Bigg ]%
\mathrm{d}s\Bigg \}  \notag \\
&\leq &\mathbb{E}\Bigg \{\left\vert \xi _{i}\right\vert ^{p}+\int_{0}^{t}%
\Bigg [L^{\bar{b}}\left( \left( 2p-1\right) \left\vert X_{s,i}\right\vert
^{p}+1\right) +L_{y}^{\bar{b}}\left( p-1\right) \left\vert
X_{s,i}\right\vert ^{p}  \notag \\
&&+\frac{L_{y}^{\bar{b}}}{N}\sum_{k=1}^{N}\left\vert X_{s,k}\right\vert ^{p}+%
\frac{3p\left( p-1\right) }{2}\left\vert X_{s,i}\right\vert ^{p-2}\Bigg (%
\left( L^{\bar{\sigma}}\right) ^{2}+\left( L^{\bar{\sigma}}\right)
^{2}\left\vert X_{s,i}\right\vert ^{2}  \notag \\
&&+\frac{\left( L_{y}^{\bar{\sigma}}\right) ^{2}}{N}\sum_{k=1}^{N}\left\vert
X_{s,k}\right\vert ^{2}\Bigg )\Bigg ]\mathrm{d}s\Bigg \}  \notag \\
&\leq &\mathbb{E}\Bigg \{\left\vert \xi _{i}\right\vert ^{p}+\int_{0}^{t}%
\Bigg [L^{\bar{b}}\left( \left( 2p-1\right) \left\vert X_{s,i}\right\vert
^{p}+1\right) +L_{y}^{\bar{b}}\left( p-1\right) \left\vert
X_{s,i}\right\vert ^{p}  \notag \\
&&+\frac{L_{y}^{\bar{b}}}{N}\sum_{k=1}^{N}\left\vert X_{s,k}\right\vert ^{p}+%
\frac{3\left( p-1\right) \left( p-2\right) \left( L^{\bar{\sigma}}\right)
^{2}}{2}\left\vert X_{s,i}\right\vert ^{p}+\frac{3\left( p-1\right) \left(
L^{\bar{\sigma}}\right) ^{2}}{2}  \notag \\
&&+\frac{3p\left( p-1\right) \left( L^{\bar{\sigma}}\right) ^{2}}{2}%
\left\vert X_{s,i}\right\vert ^{p}+\frac{3\left( p-1\right) \left(
p-2\right) \left( L_{y}^{\bar{\sigma}}\right) ^{2}}{2}\left\vert
X_{s,i}\right\vert ^{p}  \notag \\
&&+\frac{3\left( p-1\right) \left( L_{y}^{\bar{\sigma}}\right) ^{2}}{N}%
\sum_{k=1}^{N}\left\vert X_{s,k}\right\vert ^{p}\Bigg )\Bigg ]\mathrm{d}s%
\Bigg \}  \notag \\
&\leq &\mathbb{E}\Bigg \{\left\vert \xi _{i}\right\vert ^{p}+TL^{\bar{b}}+%
\frac{3\left( p-1\right) \left( L^{\bar{\sigma}}\right) ^{2}T}{2}%
+\int_{0}^{t}\Bigg [C^{p,\bar{b},\bar{\sigma}}\mathbb{E}\left[ \left\vert
X_{s,i}\right\vert ^{p}\right]  \notag \\
&&+\frac{3\left( p-1\right) \left( L_{y}^{\bar{\sigma}}\right) ^{2}}{N}%
\sum_{k=1}^{N}\mathbb{E}\left[ \left\vert X_{s,k}\right\vert ^{p}\right] %
\Bigg ]\mathrm{d}s\Bigg \},  \label{ineq2}
\end{eqnarray}%
where%
\begin{eqnarray*}
C^{p,\bar{b},\bar{\sigma}} &=&\left( 2p-1\right) L^{\bar{b}}+L_{y}^{\bar{b}%
}\left( p-1\right) +\frac{1}{2}\Big (3\left( p-1\right) \left( p-2\right)
\left( L^{\bar{\sigma}}\right) ^{2}+3p\left( p-1\right) \left( L^{\bar{\sigma%
}}\right) ^{2} \\
&&+3\left( p-1\right) \left( p-2\right) \left( L_{y}^{\bar{\sigma}}\right)
^{2}\Big )
\end{eqnarray*}%
Summing up the above equation over the index $i\in I_{N}$ yields for all $%
t\in \lbrack 0;T]$,%
\begin{eqnarray*}
\sum_{i=1}^{N}\mathbb{E}\left[ \left\vert X_{t,i}\right\vert ^{p}\right]
&\leq &\sum_{i=1}^{N}\left( \mathbb{E}\left[ \left\vert \xi _{i}\right\vert
^{p}\right] +TL^{\bar{b}}+\frac{3\left( p-1\right) \left( L^{\bar{\sigma}%
}\right) ^{2}T}{2}\right) \\
&&+\int_{0}^{t}\Bigg [\big (C^{p,\bar{b},\bar{\sigma}}+3\left( p-1\right)
\left( L_{y}^{\bar{\sigma}}\right) ^{2}\big )\sum_{k=1}^{N}\mathbb{E}\left[
\left\vert X_{s,k}\right\vert ^{p}\right] \Bigg )\Bigg ]\mathrm{d}s,
\end{eqnarray*}%
Immediately, Gronwall's lemma yields
\begin{equation*}
\sum_{k=1}^{N}\mathbb{E}\left[ \left\vert X_{t,k}\right\vert ^{p}\right]
\leq \sum_{i=1}^{N}\left( \mathbb{E}\left[ \left\vert \xi _{i}\right\vert
^{p}\right] +TL^{\bar{b}}+\frac{3\left( p-1\right) \left( L^{\bar{\sigma}%
}\right) ^{2}T}{2}\right) e^{\left( C^{p,\bar{b},\bar{\sigma}}+3\left(
p-1\right) \left( L_{y}^{\bar{\sigma}}\right) ^{2}\right) T}.
\end{equation*}

Plugging the above inequality into (\ref{ineq2}), we have%
\begin{eqnarray*}
\mathbb{E}\left[ \left\vert X_{t,i}\right\vert ^{p}\right] &\leq &\mathbb{E}%
\Bigg \{\left\vert \xi _{i}\right\vert ^{p}+TL^{\bar{b}}+\frac{3\left(
p-1\right) \left( L^{\bar{\sigma}}\right) ^{2}T}{2}+\int_{0}^{t}\Bigg [C^{p,%
\bar{b},\bar{\sigma}}\mathbb{E}\left[ \left\vert X_{s,i}\right\vert ^{p}%
\right] \\
&&+\frac{3\left( p-1\right) \left( L_{y}^{\bar{\sigma}}\right) ^{2}}{N}%
\sum_{i=1}^{N}\left( \mathbb{E}\left[ \left\vert \xi _{i}\right\vert ^{p}%
\right] +TL^{\bar{b}}+\frac{3\left( p-1\right) \left( L^{\bar{\sigma}%
}\right) ^{2}T}{2}\right) \cdot \\
&&e^{\left( C^{p,\bar{b},\bar{\sigma}}+3\left( p-1\right) \left( L_{y}^{\bar{%
\sigma}}\right) ^{2}\right) T}\Bigg ]\mathrm{d}s\Bigg \}.
\end{eqnarray*}%
Applying Gronwall's lemma again gives
\begin{eqnarray*}
\mathbb{E}\left[ \left\vert X_{t,i}\right\vert ^{p}\right] &\leq &\mathbb{E}%
\Bigg [\left\vert \xi _{i}\right\vert ^{p}+TL^{\bar{b}}+\frac{3\left(
p-1\right) \left( L^{\bar{\sigma}}\right) ^{2}T}{2} \\
&&+\frac{3\left( p-1\right) \left( L_{y}^{\bar{\sigma}}\right) ^{2}T}{N}%
\sum_{i=1}^{N}\left( \mathbb{E}\left[ \left\vert \xi _{i}\right\vert ^{p}%
\right] +TL^{\bar{b}}+\frac{3\left( p-1\right) \left( L^{\bar{\sigma}%
}\right) ^{2}T}{2}\right) \cdot \\
&&e^{\left( C^{p,\bar{b},\bar{\sigma}}+3\left( p-1\right) \left( L_{y}^{\bar{%
\sigma}}\right) ^{2}\right) T}\Bigg ]e^{TC^{p,\bar{b},\bar{\sigma}}}.
\end{eqnarray*}%
We thus complete the proof. \hfill $\Box $

\begin{lemma}
\label{l3}Fix $p\geq 2$ and for each $i,j\in I_{N}$, let $B_{i},\bar{B}%
_{i,j},D_{i},\bar{D}_{i,j}:\Omega \times \left[ 0,T\right] \rightarrow
\mathbb{R}$ be bounded adapted processes, and $f_{i},\bar{f}_{i}\in \mathcal{%
H}^{p}\left( \mathbb{R}\right) $. Let $\mathbf{S}=(S_{i})_{i=1}^{N}\in
\mathcal{S}^{p}\left( \mathbb{R}^{N}\right) $ satisfy the following
dynamics: for all $t\in \lbrack 0;T]$,%
\begin{eqnarray*}
\mathrm{d}S_{t,i} &=&\left( B_{i}\left( t\right) S_{t,i}+\sum_{j=1}^{N}\bar{B%
}_{i,j}\left( t\right) S_{t,j}+f_{t,i}\right) \mathrm{d}t \\
&&+\left( D_{i}\left( t\right) S_{t,i}+\sum_{j=1}^{N}\bar{D}_{i,j}\left(
t\right) S_{t,j}+\bar{f}_{t,i}\right) \mathrm{d}W_{t}^{i}, \\
S_{0,i} &=&0;\qquad \forall i\in I_{N}.
\end{eqnarray*}%
Then for all $i\in I_{N}$, we have%
\begin{eqnarray*}
\mathbb{E}\left[ \left\vert S_{t,i}\right\vert ^{p}\right] &\leq &\Bigg [%
T\left( \left\Vert \bar{B}\right\Vert _{\infty }+3\left( p-1\right)
N\left\Vert \bar{D}\right\Vert _{\infty }^{2}\right) \\
&&\cdot \left( \sum_{k=1}^{N}\left( \left\Vert f_{k}\right\Vert _{\mathcal{H}%
^{p}\left( \mathbb{R}\right) }^{p}+3\left( p-1\right) \left\Vert \bar{f}%
_{k}\right\Vert _{\mathcal{H}^{p}\left( \mathbb{R}\right) }^{p}\right)
\right) e^{I_{B,D,\bar{B},\bar{D},p}^{3}\cdot T} \\
&&+\left\Vert f_{i}\right\Vert _{\mathcal{H}^{p}\left( \mathbb{R}\right)
}^{p}+3\left( p-1\right) \left\Vert \bar{f}_{i}\right\Vert _{\mathcal{H}%
^{p}\left( \mathbb{R}\right) }^{p}\Bigg ]e^{I_{B,D,\bar{B},\bar{D}%
,p}^{4}\cdot T},
\end{eqnarray*}%
where
\begin{eqnarray*}
I_{B,D,\bar{B},\bar{D},p}^{3} &=&p\left\Vert B\right\Vert _{\infty
}+N\left\Vert \bar{B}\right\Vert _{\infty }p+\left( p-1\right) +\frac{3}{2}%
\left( p-1\right) p\left\Vert D\right\Vert _{\infty }^{2} \\
&&+\frac{3}{2}\left( p-1\right) \left( p-2\right) N^{2}\left\Vert \bar{D}%
\right\Vert _{\infty }^{2}+\frac{3}{2}\left( p-1\right) \left( p-2\right) \\
&&+3\left( p-1\right) N^{2}\left\Vert \bar{D}\right\Vert _{\infty }^{2}, \\
I_{B,D,\bar{B},\bar{D},p}^{4} &=&p\left\Vert B\right\Vert _{\infty
}+N\left\Vert \bar{B}\right\Vert _{\infty }\left( p-1\right) +\left(
p-1\right) +\frac{3}{2}\left( p-1\right) p\left\Vert D\right\Vert _{\infty
}^{2} \\
&&+\frac{3}{2}\left( p-1\right) N^{2}\left\Vert \bar{D}\right\Vert _{\infty
}^{2}\left( p-2\right) +\frac{3}{2}\left( p-1\right) \left( p-2\right) ,
\end{eqnarray*}%
additionally, $\left\Vert B\right\Vert _{\infty }=\max_{i\in
I_{N}}\left\Vert B_{i}\right\Vert _{L^{\infty }},$ $\left\Vert D\right\Vert
_{\infty }=\max_{i\in I_{N}}\left\Vert D_{i}\right\Vert _{L^{\infty }}$, $%
\left\Vert \bar{B}\right\Vert _{\infty }=\max_{i,j\in I_{N}}\left\Vert \bar{B%
}_{i,j}\right\Vert _{L^{\infty }}$ and $\left\Vert \bar{D}\right\Vert
_{\infty }=\max_{i,j\in I_{N}}\left\Vert \bar{D}_{i,j}\right\Vert
_{L^{\infty }}.$
\end{lemma}

\paragraph{Proof.}

Applying It\^{o}'s formula to $\left\vert S_{t,i}\right\vert ^{p}$, we get%
\begin{eqnarray*}
\left\vert S_{t,i}\right\vert ^{p} &=&\int_{0}^{t}p\left\vert
S_{r,i}\right\vert ^{p-1}\left( B_{i}\left( r\right) S_{r,i}+\sum_{k=1}^{N}%
\bar{B}_{i,k}\left( r\right) S_{r,k}+f_{r,i}\right) \mathrm{d}r \\
&&+\int_{0}^{t}p\left\vert S_{r,i}\right\vert ^{p-1}\left( D_{i}\left(
t\right) S_{t,i}+\sum_{k=1}^{N}\bar{D}_{i,k}\left( t\right) S_{t,k}+\bar{f}%
_{t,i}\right) \mathrm{d}W_{t}^{i} \\
&&+\frac{1}{2}p\left( p-1\right) \int_{0}^{t}\left\vert S_{r,i}\right\vert
^{p-2}\left( D_{i}\left( r\right) S_{r,i}+\sum_{k=1}^{N}\bar{D}_{i,k}\left(
r\right) S_{r,k}+\bar{f}_{r,i}\right) ^{2}\mathrm{d}r.
\end{eqnarray*}%
Taking the expectation, we have%
\begin{eqnarray*}
\mathbb{E}\left[ \left\vert S_{t,i}\right\vert ^{p}\right] &=&\mathbb{E}%
\Bigg [\displaystyle\int_{0}^{t}p\left\vert S_{r,i}\right\vert ^{p-1}\left(
B_{i}\left( r\right) S_{r,i}+\sum_{k=1}^{N}\bar{B}_{i,k}\left( r\right)
S_{r,k}+f_{r,i}\right) \mathrm{d}r \\
&&+\frac{1}{2}p\left( p-1\right) \int_{0}^{t}\left\vert S_{r,i}\right\vert
^{p-2}\left( D_{i}\left( r\right) S_{r,i}+\sum_{k=1}^{N}\bar{D}_{i,k}\left(
r\right) S_{r,k}+\bar{f}_{r,i}\right) ^{2}\mathrm{d}r\Bigg ].
\end{eqnarray*}%
Applying Young's inequality ($ab\leq \frac{p-1}{p}a^{\frac{p}{p-1}}+\frac{1}{%
p}b^{p}$ and $ab\leq \frac{p-2}{p}a^{\frac{p}{p-2}}+\frac{2}{p}b^{\frac{p}{2}%
}$ whenever $p>2)$, we obtain%
\begin{eqnarray}
\mathbb{E}\left[ \left\vert S_{t,i}\right\vert ^{p}\right] &\leq &\mathbb{E}%
\Bigg [\displaystyle\int_{0}^{t}\Big (p\left\Vert B\right\Vert _{\infty
}\left\vert S_{r,i}\right\vert ^{p}+p\left\Vert \bar{B}\right\Vert _{\infty
}\sum_{k=1}^{N}\left\vert S_{r,i}\right\vert ^{p-1}\left\vert
S_{r,k}\right\vert +p\left\vert S_{r,i}\right\vert ^{p-1}\left\vert
f_{r,i}\right\vert \Big )\mathrm{d}r  \notag \\
&&+\frac{1}{2}p\left( p-1\right) \int_{0}^{t}\left\vert S_{r,i}\right\vert
^{p-2}\left( \left\Vert D\right\Vert _{\infty }\left\vert S_{r,i}\right\vert
+\left\Vert \bar{D}\right\Vert _{\infty }\sum_{k=1}^{N}\left\vert
S_{r,k}\right\vert +\left\vert \bar{f}_{r,i}\right\vert \right) ^{2}\mathrm{d%
}r\Bigg ]  \notag \\
&\leq &\mathbb{E}\Bigg [\displaystyle\int_{0}^{t}\Big (p\left\Vert
B\right\Vert _{\infty }\left\vert S_{r,i}\right\vert ^{p}+p\left\Vert \bar{B}%
\right\Vert _{\infty }\sum_{k=1}^{N}\left( \frac{p-1}{p}\left\vert
S_{r,i}\right\vert ^{p}+\frac{1}{p}\left\vert S_{r,k}\right\vert ^{p}\right)
\notag \\
&&+p\left( \frac{p-1}{p}\left\vert S_{r,i}\right\vert ^{p}+\frac{1}{p}%
\left\vert f_{r,i}\right\vert ^{p}\right) \Big )\mathrm{d}r  \notag \\
&&+\frac{3}{2}p\left( p-1\right) \int_{0}^{t}\left\vert S_{r,i}\right\vert
^{p-2}\left( \left\Vert D\right\Vert _{\infty }^{2}\left\vert
S_{r,i}\right\vert ^{2}+N\left\Vert \bar{D}\right\Vert _{\infty
}^{2}\sum_{k=1}^{N}\left\vert S_{r,k}\right\vert ^{2}+\left\vert \bar{f}%
_{r,i}\right\vert ^{2}\right) \mathrm{d}r\Bigg ]  \notag \\
&\leq &\mathbb{E}\Bigg [\displaystyle\int_{0}^{t}\Big (p\left\Vert
B\right\Vert _{\infty }\left\vert S_{r,i}\right\vert ^{p}+\left\Vert \bar{B}%
\right\Vert _{\infty }\sum_{k=1}^{N}\left( \left( p-1\right) \left\vert
S_{r,i}\right\vert ^{p}+\left\vert S_{r,k}\right\vert ^{p}\right)  \notag \\
&&+\left( p-1\right) \left\vert S_{r,i}\right\vert ^{p}+\left\vert
f_{r,i}\right\vert ^{p}\Big )\mathrm{d}r  \notag \\
&&+\frac{3}{2}p\left( p-1\right) \int_{0}^{t}\Big (\left\Vert D\right\Vert
_{\infty }^{2}\left\vert S_{r,i}\right\vert ^{p}+N\left\Vert \bar{D}%
\right\Vert _{\infty }^{2}\sum_{k=1}^{N}\left\vert S_{r,i}\right\vert
^{p-2}\left\vert S_{r,k}\right\vert ^{2}  \notag \\
&&+\left\vert S_{r,i}\right\vert ^{p-2}\left\vert \bar{f}_{r,i}\right\vert
^{2}\Big )\mathrm{d}r\Bigg ]  \notag \\
&\leq &\mathbb{E}\Bigg [\displaystyle\int_{0}^{t}\Big (p\left\Vert
B\right\Vert _{\infty }\left\vert S_{r,i}\right\vert ^{p}+\left\Vert \bar{B}%
\right\Vert _{\infty }\sum_{k=1}^{N}\left( \left( p-1\right) \left\vert
S_{r,i}\right\vert ^{p}+\left\vert S_{r,k}\right\vert ^{p}\right)  \notag \\
&&+\left( p-1\right) \left\vert S_{r,i}\right\vert ^{p}+\left\vert
f_{r,i}\right\vert ^{p}\Big )\mathrm{d}r  \notag \\
&&+\frac{3}{2}\left( p-1\right) \int_{0}^{t}\Big (p\left\Vert D\right\Vert
_{\infty }^{2}\left\vert S_{r,i}\right\vert ^{p}+N\left\Vert \bar{D}%
\right\Vert _{\infty }^{2}\sum_{k=1}^{N}\left( \left( p-2\right) \left\vert
S_{r,i}\right\vert ^{p}+2\left\vert S_{r,k}\right\vert ^{p}\right)  \notag \\
&&+\left( p-2\right) \left\vert S_{r,i}\right\vert ^{p}+2\left\vert \bar{f}%
_{r,i}\right\vert ^{p}\Big )\mathrm{d}r\Bigg ]  \notag \\
&\leq &\mathbb{E}\Bigg \{\displaystyle\int_{0}^{t}\Bigg [\Big (p\left\Vert
B\right\Vert _{\infty }+N\left\Vert \bar{B}\right\Vert _{\infty }\left(
p-1\right) +\left( p-1\right) +\frac{3}{2}\left( p-1\right) p\left\Vert
D\right\Vert _{\infty }^{2}  \notag \\
&&+\frac{3}{2}\left( p-1\right) N^{2}\left\Vert \bar{D}\right\Vert _{\infty
}^{2}\left( p-2\right) +\frac{3}{2}\left( p-1\right) \left( p-2\right) \Big )%
\left\vert S_{r,i}\right\vert ^{p}  \notag \\
&&+\left( \left\Vert \bar{B}\right\Vert _{\infty }+3\left( p-1\right)
N\left\Vert \bar{D}\right\Vert _{\infty }^{2}\right)
\sum_{k=1}^{N}\left\vert S_{r,k}\right\vert ^{p}  \notag \\
&&+\left\vert f_{r,i}\right\vert ^{p}+3\left( p-1\right) \left\vert \bar{f}%
_{r,i}\right\vert ^{p}\Big )\Bigg ]\mathrm{d}r\Bigg \}.  \label{ineq3}
\end{eqnarray}%
Summarizing (\ref{ineq3}) over the index $i\in I_{N}$ yields for all $t\in
\lbrack 0;T]$,%
\begin{eqnarray*}
\sum_{k=1}^{N}\mathbb{E}\left[ \left\vert S_{t,k}\right\vert ^{p}\right]
&\leq &\mathbb{E}\Bigg \{\int_{0}^{t}\Bigg [\Big (p\left\Vert B\right\Vert
_{\infty }+N\left\Vert \bar{B}\right\Vert _{\infty }\left( p-1\right)
+\left( p-1\right) +\frac{3}{2}\left( p-1\right) p\left\Vert D\right\Vert
_{\infty }^{2} \\
&&+\frac{3}{2}\left( p-1\right) N^{2}\left\Vert \bar{D}\right\Vert _{\infty
}^{2}\left( p-2\right) +\frac{3}{2}\left( p-1\right) \left( p-2\right) \Big )%
\sum_{k=1}^{N}\left\vert S_{r,k}\right\vert ^{p} \\
&&+N\left( \left\Vert \bar{B}\right\Vert _{\infty }+3\left( p-1\right)
N\left\Vert \bar{D}\right\Vert _{\infty }^{2}\right)
\sum_{k=1}^{N}\left\vert S_{r,k}\right\vert ^{p} \\
&&+\sum_{k=1}^{N}\left( \left\vert f_{r,k}\right\vert ^{p}+3\left(
p-1\right) \left\vert \bar{f}_{r,k}\right\vert ^{p}\right) \Big )\Bigg ]%
\mathrm{d}r\Bigg \} \\
&=&\mathbb{E}\Bigg \{\int_{0}^{t}\Bigg [\Big (p\left\Vert B\right\Vert
_{\infty }+N\left\Vert \bar{B}\right\Vert _{\infty }p+\left( p-1\right) +%
\frac{3}{2}\left( p-1\right) p\left\Vert D\right\Vert _{\infty }^{2} \\
&&+\frac{3}{2}\left( p-1\right) \left( p-2\right) N^{2}\left\Vert \bar{D}%
\right\Vert _{\infty }^{2}+\frac{3}{2}\left( p-1\right) \left( p-2\right) \\
&&+3\left( p-1\right) N^{2}\left\Vert \bar{D}\right\Vert _{\infty }^{2}\Big )%
\sum_{k=1}^{N}\left\vert S_{r,k}\right\vert ^{p} \\
&&+\sum_{k=1}^{N}\left( \left\vert f_{r,k}\right\vert ^{p}+3\left(
p-1\right) \left\vert \bar{f}_{r,k}\right\vert ^{p}\right) \Big )\Bigg ]%
\mathrm{d}r\Bigg \}.
\end{eqnarray*}%
Gronwall's inequality implies that%
\begin{eqnarray*}
\sum_{k=1}^{N}\mathbb{E}\left[ \left\vert S_{t,k}\right\vert ^{p}\right]
&\leq &\mathbb{E}\left[ \left( \int_{0}^{T}\sum_{k=1}^{N}\left( \left\vert
f_{r,k}\right\vert ^{p}+3\left( p-1\right) \left\vert \bar{f}%
_{r,k}\right\vert ^{p}\right) \mathrm{d}r\right) e^{I_{B,D,\bar{B},\bar{D}%
,p}^{3}\cdot T}\right] \\
&=&\mathbb{E}\left[ \left( \sum_{k=1}^{N}\int_{0}^{T}\left\vert
f_{r,k}\right\vert ^{p}+3\left( p-1\right)
\sum_{k=1}^{N}\int_{0}^{T}\left\vert \bar{f}_{r,k}\right\vert ^{p}\mathrm{d}%
r\right) e^{I_{B,D,\bar{B},\bar{D},p}^{3}\cdot T}\right] \\
&=&\left( \sum_{k=1}^{N}\mathbb{E}\left[ \int_{0}^{T}\left\vert
f_{r,k}\right\vert ^{p}\mathrm{d}r\right] +3\left( p-1\right) \sum_{k=1}^{N}%
\mathbb{E}\left[ \int_{0}^{T}\left\vert \bar{f}_{r,k}\right\vert ^{p}\mathrm{%
d}r\right] \right) e^{I_{B,D,\bar{B},\bar{D},p}^{3}\cdot T} \\
&=&\left[ \sum_{k=1}^{N}\left( \left\Vert f_{k}\right\Vert _{\mathcal{H}%
^{p}\left( \mathbb{R}\right) }^{p}+3\left( p-1\right) \left\Vert \bar{f}%
_{k}\right\Vert _{\mathcal{H}^{p}\left( \mathbb{R}\right) }^{p}\right) %
\right] e^{I_{B,D,\bar{B},\bar{D},p}^{3}\cdot T}.
\end{eqnarray*}%
Substituting the above inequality into (\ref{ineq3}) yields for all $t>0$,%
\begin{eqnarray*}
\mathbb{E}\left[ \left\vert S_{t,i}\right\vert ^{p}\right] &\leq &\mathbb{E}%
\Bigg [\int_{0}^{t}\Bigg [\Big (p\left\Vert B\right\Vert _{\infty
}+N\left\Vert \bar{B}\right\Vert _{\infty }\left( p-1\right) +\left(
p-1\right) +\frac{3}{2}\left( p-1\right) p\left\Vert D\right\Vert _{\infty
}^{2} \\
&&+\frac{3}{2}\left( p-1\right) N^{2}\left\Vert \bar{D}\right\Vert _{\infty
}^{2}\left( p-2\right) +\frac{3}{2}\left( p-1\right) \left( p-2\right) \Big )%
\mathbb{E}\left[ \left\vert S_{r,i}\right\vert ^{p}\right] \\
&&+\left( \left\Vert \bar{B}\right\Vert _{\infty }+3\left( p-1\right)
N\left\Vert \bar{D}\right\Vert _{\infty }^{2}\right) \mathbb{E}\left[
\sum_{k=1}^{N}\left\vert S_{r,k}\right\vert ^{p}\right] \\
&&+\mathbb{E}\left[ \left\vert f_{r,i}\right\vert ^{p}+3\left( p-1\right)
\left\vert \bar{f}_{r,i}\right\vert ^{p}\right] \Big )\Bigg ]\mathrm{d}r \\
&\leq &\int_{0}^{t}\Big (p\left\Vert B\right\Vert _{\infty }+N\left\Vert
\bar{B}\right\Vert _{\infty }\left( p-1\right) +\left( p-1\right) +\frac{3}{2%
}\left( p-1\right) p\left\Vert D\right\Vert _{\infty }^{2} \\
&&+\frac{3}{2}\left( p-1\right) N^{2}\left\Vert \bar{D}\right\Vert _{\infty
}^{2}\left( p-2\right) +\frac{3}{2}\left( p-1\right) \left( p-2\right) \Big )%
\mathbb{E}\left[ \left\vert S_{r,i}\right\vert ^{p}\right] \mathrm{d}r \\
&&+\int_{0}^{t}\Bigg \{\left( \left\Vert \bar{B}\right\Vert _{\infty
}+3\left( p-1\right) N\left\Vert \bar{D}\right\Vert _{\infty }^{2}\right) \\
&&\cdot \left[ \sum_{k=1}^{N}\left( \left\Vert f_{k}\right\Vert _{\mathcal{H}%
^{p}\left( \mathbb{R}\right) }^{p}+3\left( p-1\right) \left\Vert \bar{f}%
_{k}\right\Vert _{\mathcal{H}^{p}\left( \mathbb{R}\right) }^{p}\right) %
\right] e^{I_{B,D,\bar{B},\bar{D},p}^{3}\cdot T} \\
&&+\left\vert f_{r,i}\right\vert ^{p}+3\left( p-1\right) \left\vert \bar{f}%
_{r,i}\right\vert ^{p}\Bigg \}\mathrm{d}r
\end{eqnarray*}%
and then we get the desired result.\hfill $\Box $

\paragraph{Proof of Lemma \protect\ref{l4}.}

For convenience of notation, we set $\mathbf{X}=\mathbf{X}^{\phi },$ $%
\mathbf{Y}^{h}=\mathbf{Y}^{\phi ,\phi _{h}^{\prime }}.$ Employing Lemma \ref%
{l3} with $\mathbf{S=Y}^{h},$ $B_{i}\left( t\right) =\left( \partial _{x}%
\bar{b}_{t,i}\right) \left( X_{t,i},\mathbf{X}_{t}\right) ,$ $\bar{B}%
_{i,j}\left( t\right) =\left( \partial _{y_{j}}\bar{b}_{i}\right) \left(
X_{t,i},\mathbf{X}_{t}\right) ,$ $f_{t,i}=\delta _{h,i}\phi _{t,h}^{\prime
}\left( X_{t,i},\mathbf{X}_{t}\right) ,$ $D_{i}\left( t\right) =\left(
\partial _{x}\bar{\sigma}_{i}\right) \left( X_{t,i},\mathbf{X}_{t}\right) ,$
$\bar{D}_{i,j}\left( t\right) =\left( \partial _{y_{j}}\bar{\sigma}%
_{i}\right) \left( X_{t,i},\mathbf{X}_{t}\right) ,$ $\bar{f}_{t,i}=f_{t,i}$
yields that for any $i\in I_{N}$%
\begin{eqnarray*}
\mathbb{E}\left[ \left\vert Y_{t,i}^{h}\right\vert ^{p}\right]  &\leq &\Bigg
[T\left( \left\Vert \bar{B}\right\Vert _{\infty }+3\left( p-1\right)
N\left\Vert \bar{D}\right\Vert _{\infty }^{2}\right) \left( 3p-2\right)
\sum_{k=1}^{N}\left\Vert f_{k}\right\Vert _{\mathcal{H}^{p}\left( \mathbb{R}%
\right) }^{p}e^{I_{B,D,\bar{B},\bar{D},p}^{3}\cdot T} \\
&&+\left( 3p-2\right) \left\Vert f_{i}\right\Vert _{\mathcal{H}^{p}\left(
\mathbb{R}\right) }^{p}\Bigg ]e^{I_{B,D,\bar{B},\bar{D},p}^{4}\cdot T} \\
&\leq &\Bigg [\frac{L_{y}^{\bar{b}}+3\left( p-1\right) \left( L_{y}^{\bar{%
\sigma}}\right) ^{2}}{N}\left( 3p-2\right) TC_{f_{i}}^{p}e^{I_{B,D,\bar{B},%
\bar{D},p}^{3}\cdot T}+\delta _{h,i}\left( 3p-2\right) C_{f_{i}}^{p}\Bigg ]%
e^{I_{B,D,\bar{B},\bar{D},p}^{4}\cdot T}
\end{eqnarray*}
where $C_{f_{i}}^{p}$ is defined in (\ref{cfi})
\begin{eqnarray*}
I_{B,D,\bar{B},\bar{D},p}^{3} &=&p\left\Vert B\right\Vert _{\infty
}+N\left\Vert \bar{B}\right\Vert _{\infty }p+\left( p-1\right) +\frac{3}{2}%
\left( p-1\right) p\left\Vert D\right\Vert _{\infty }^{2} \\
&&+\frac{3}{2}\left( p-1\right) \left( p-2\right) N^{2}\left\Vert \bar{D}%
\right\Vert _{\infty }^{2}+\frac{3}{2}\left( p-1\right) \left( p-2\right) \\
&&+3\left( p-1\right) N^{2}\left\Vert \bar{D}\right\Vert _{\infty }^{2}, \\
I_{B,D,\bar{B},\bar{D},p}^{4} &=&p\left\Vert B\right\Vert _{\infty
}+N\left\Vert \bar{B}\right\Vert _{\infty }\left( p-1\right) +\left(
p-1\right) +\frac{3}{2}\left( p-1\right) p\left\Vert D\right\Vert _{\infty
}^{2} \\
&&+\frac{3}{2}\left( p-1\right) N^{2}\left\Vert \bar{D}\right\Vert _{\infty
}^{2}\left( p-2\right) +\frac{3}{2}\left( p-1\right) \left( p-2\right) .
\end{eqnarray*}%
We now deal with the term $\left\Vert f_{i}\right\Vert _{\mathcal{H}%
^{p}\left( \mathbb{R}\right) }^{p}.$%
\begin{eqnarray*}
\left\Vert f_{i}\right\Vert _{\mathcal{H}^{p}\left( \mathbb{R}\right) }^{p}
&\leq &T\sup_{0\leq t\leq T}\mathbb{E}\left[ \left\vert f_{t,i}\right\vert
^{p}\right] \\
&\leq &\delta _{h,i}\left( 2T\right) ^{p-1}\left[ \left( L^{\phi
_{h}^{\prime }}\right) ^{p}\mathbb{E}\left( 1+\left\vert X_{t,i}\right\vert
\right) ^{p}+\left( \frac{L_{y}^{\phi _{h}^{\prime }}}{N}\right) ^{p}\left(
\sum_{k=1}^{N}\left\vert X_{t,k}\right\vert \right) ^{p}\right] .
\end{eqnarray*}%
From the inequality $\left( \sum_{i=1}^{N}a_{i}\right) ^{p}\leq
N^{p-1}\sum_{i=1}^{N}a_{i}^{p}$ for $a_{i}\geq 0,$ $0\leq i\leq N,$ we
derive that%
\begin{eqnarray*}
\left\Vert f_{i}\right\Vert _{\mathcal{H}^{p}\left( \mathbb{R}\right) }^{p}
&\leq &\delta _{h,i}\left( 2T\right) ^{p-1}\left[ 2^{p-1}\left( L^{\phi
_{h}^{\prime }}\right) ^{p}\mathbb{E}\left( 1+\left\vert X_{t,i}\right\vert
^{p}\right) +\left( L_{y}^{\phi _{h}^{\prime }}\right) ^{p}\frac{%
\sum_{k=1}^{N}\left\vert X_{t,k}\right\vert ^{p}}{N}\right] \\
&\leq &\delta _{h,i}\left( 2T\right) ^{p-1}\left[ 2^{p-1}\left( L^{\phi
_{h}^{\prime }}\right) ^{p}\left( 1+C_{X}^{i,p}\right) +\left( L_{y}^{\phi
_{h}^{\prime }}\right) ^{p}\frac{\sum_{k=1}^{N}C_{X}^{i,p}}{N}\right] \\
&=&\delta _{h,i}\left( 2T\right) ^{p-1}\left[ 2^{p-1}\left( L^{\phi
_{h}^{\prime }}\right) ^{p}\left( 1+C_{X}^{h,p}\right) +\left( L_{y}^{\phi
_{h}^{\prime }}\right) ^{p}\frac{\sum_{k=1}^{N}C_{X}^{i,p}}{N}\right] ,
\end{eqnarray*}%
where the last equality holds due to $\delta _{h,h}=1.$

Put
\begin{equation}
C_{f_{h}}^{p}=\left( 2T\right) ^{p-1}\left[ 2^{p-1}\left( L^{\phi
_{h}^{\prime }}\right) ^{p}\left( 1+C_{X}^{h,p}\right) +\left( L_{y}^{\phi
_{h}^{\prime }}\right) ^{p}\frac{\sum_{k=1}^{N}C_{X}^{i,p}}{N}\right] .
\label{cfi}
\end{equation}%
Analogously, we can get
\begin{eqnarray*}
\sum_{k=1}^{N}\left\Vert f_{k}\right\Vert _{\mathcal{H}^{p}\left( \mathbb{R}%
\right) }^{p} &\leq &\left\Vert f_{h}\right\Vert _{\mathcal{H}^{p}\left(
\mathbb{R}\right) }^{p} \\
&\leq &\left( 2T\right) ^{p-1}\left[ 2^{p-1}\left( L^{\phi _{h}^{\prime
}}\right) ^{p}\left( 1+C_{X}^{h,p}\right) +\left( L_{y}^{\phi _{h}^{\prime
}}\right) ^{p}\frac{\sum_{k=1}^{N}C_{X}^{i,p}}{N}\right] .
\end{eqnarray*}%
where we have used the fact that $\left\Vert B\right\Vert _{\infty }\leq L^{%
\bar{b}},$ $\left\Vert D\right\Vert _{\infty }\leq L^{\bar{\sigma}},$ $%
\left\Vert \bar{B}\right\Vert _{\infty }\leq \frac{L_{y}^{\bar{b}}}{N}%
,\left\Vert \bar{D}\right\Vert _{\infty }\leq \frac{L_{y}^{\bar{\sigma}}}{N}%
. $ Now we consider for all $i,j\in I_{N},$ by virtue of H\"{o}lder's
inequality%
\begin{eqnarray}
\left\Vert Y_{i}^{h}Y_{j}^{\ell }\right\Vert _{\mathcal{H}^{2}\left( \mathbb{%
R}\right) }^{2} &\leq &T\sup_{0\leq t\leq T}\mathbb{E}\left[ \left\vert
Y_{i}^{h}\right\vert ^{4}\right] ^{\frac{1}{2}}\mathbb{E}\left[ \left\vert
Y_{j}^{\ell }\right\vert ^{4}\right] ^{\frac{1}{2}}  \notag \\
&\leq &T\left[ \left( \frac{L_{y}^{\bar{b}}+3\left( L_{y}^{\bar{\sigma}%
}\right) ^{2}}{N}10TC_{f_{h}}^{2}e^{I_{B,D,\bar{B},\bar{D},2}^{3}\cdot
T}+\delta _{h,i}10C_{f_{h}}^{2}\right) e^{I_{B,D,\bar{B},\bar{D},2}^{4}\cdot
T}\right] ^{\frac{1}{2}}\cdot  \notag \\
&&\left[ \left( \frac{L_{y}^{\bar{b}}+3\left( L_{y}^{\bar{\sigma}}\right)
^{2}}{N}10TC_{f_{\ell }}^{2}e^{I_{B,D,\bar{B},\bar{D},2}^{3}\cdot T}+\delta
_{\ell ,j}10C_{f_{\ell }}^{2}\right) e^{I_{B,D,\bar{B},\bar{D},2}^{4}\cdot T}%
\right] ^{\frac{1}{2}}  \notag \\
&\leq &C\Bigg (\delta _{h,i}\delta _{\ell ,j}+\frac{\left( L_{y}^{\bar{b}%
}+3\left( L_{y}^{\bar{\sigma}}\right) ^{2}\right) ^{\frac{1}{2}}}{\sqrt{N}}%
\left( \delta _{h,i}+\delta _{\ell ,j}\right) +\frac{\left( L_{y}^{\bar{b}%
}+3\left( L_{y}^{\bar{\sigma}}\right) ^{2}\right) }{N}\Bigg ).  \label{yy}
\end{eqnarray}%
We thus complete the proof. \hfill $\Box $

\begin{lemma}
\label{f1g1}Suppose that Assumptions \emph{(A1)-(A2) }hold and that there
exists $p\geq 2$ such that $\xi _{i}\in L^{4}\left( \Omega ;\mathbb{R}%
\right) $ for all $i\in I_{N}$. For each $\phi \in \Pi ^{N},$ $h,\ell \in
I_{N}$ satisfying $h\neq \ell $ and $\phi _{h}^{\prime },\phi _{\ell
}^{\prime \prime }\in \Pi ,$ the processes $\mathfrak{f}_{i}^{\phi ,\phi
_{h}^{\prime },\phi _{\ell }^{\prime \prime }}$ and $\mathfrak{g}_{i}^{\phi
,\phi _{h}^{\prime },\phi _{\ell }^{\prime \prime }}$ defined in (\ref{f1})
and (\ref{g1}) admit for all $i\in I_{N}$,%
\begin{eqnarray}
\left\Vert \mathfrak{f}_{i}^{\phi ,\phi _{h}^{\prime },\phi _{\ell }^{\prime
\prime }}\right\Vert _{\mathcal{H}^{2}\left( \mathbb{R}\right) } &\leq &%
\frac{CL_{y}^{\bar{b}}\left( L_{y}^{\bar{b}}+3\left( L_{y}^{\bar{\sigma}%
}\right) ^{2}\right) ^{\frac{1}{4}}}{N^{\frac{1}{4}}}\left( \delta
_{h,i}+\delta _{\ell ,i}\right) +\frac{CL_{y}^{\bar{b}}\left( L_{y}^{\bar{b}%
}+3\left( L_{y}^{\bar{\sigma}}\right) ^{2}\right) ^{\frac{1}{2}}}{N^{\frac{1%
}{2}}}  \notag \\
&&+\delta _{h,i}\frac{C}{N^{\frac{1}{2}}}\Bigg [L^{\phi _{h}^{\prime
}}\left( L_{y}^{\bar{b}}+3\left( L_{y}^{\bar{\sigma}}\right) ^{2}\right) ^{%
\frac{1}{2}}+L_{y}^{\phi _{h}^{\prime }}\Bigg ]  \notag \\
&&+\delta _{\ell ,i}\frac{C}{N^{\frac{1}{2}}}\Bigg [L^{\phi _{\ell }^{\prime
\prime }}\left( L_{y}^{\bar{b}}+3\left( L_{y}^{\bar{\sigma}}\right)
^{2}\right) ^{\frac{1}{2}}+L_{y}^{\phi _{\ell }^{\prime \prime }}\Bigg ],
\label{fineq} \\
\left\Vert \mathfrak{g}_{i}^{\phi ,\phi _{h}^{\prime },\phi _{\ell }^{\prime
\prime }}\right\Vert _{\mathcal{H}^{2}\left( \mathbb{R}\right) } &\leq &%
\frac{CL_{y}^{\bar{\sigma}}\left( L_{y}^{\bar{b}}+3\left( L_{y}^{\bar{\sigma}%
}\right) ^{2}\right) ^{\frac{1}{4}}}{N^{\frac{1}{4}}}\left( \delta
_{h,i}+\delta _{\ell ,i}\right) +\frac{CL_{y}^{\bar{\sigma}}\left( L_{y}^{%
\bar{b}}+3\left( L_{y}^{\bar{\sigma}}\right) ^{2}\right) ^{\frac{1}{2}}}{N^{%
\frac{1}{2}}}  \notag \\
&&+\delta _{h,i}\frac{C}{N^{\frac{1}{2}}}\Bigg [L^{\phi _{h}^{\prime
}}\left( L_{y}^{\bar{b}}+3\left( L_{y}^{\bar{\sigma}}\right) ^{2}\right) ^{%
\frac{1}{2}}+L_{y}^{\phi _{h}^{\prime }}\Bigg ]  \notag \\
&&+\delta _{\ell ,i}\frac{C}{N^{\frac{1}{2}}}\Bigg [L^{\phi _{\ell }^{\prime
\prime }}\left( L_{y}^{\bar{b}}+3\left( L_{y}^{\bar{\sigma}}\right)
^{2}\right) ^{\frac{1}{2}}+L_{y}^{\phi _{\ell }^{\prime \prime }}\Bigg ].
\label{gineq}
\end{eqnarray}
\end{lemma}

\paragraph{Proof.}

For all $t\in \left[ 0,T\right] ,$ we have
\begin{eqnarray}
\mathfrak{f}_{t,i}^{\phi ,\phi _{h}^{\prime },\phi _{\ell }^{\prime \prime
}} &=&\underset{I_{1}}{\underbrace{\partial _{xx}^{2}\left( \bar{b}%
_{t,i}\right) Y_{t,i}^{h}Y_{t,i}^{\ell }}}\underset{I_{2}}{+\underbrace{%
\sum_{j=1}\partial _{xy_{j}}^{2}\left( \bar{b}_{t,i}\right) \left(
X_{t,i}^{\phi },\mathbf{X}_{t}^{\phi }\right) \left(
Y_{t,i}^{h}Y_{t,j}^{\ell }+Y_{t,j}^{h}Y_{t,i}^{\ell }\right) }}  \notag \\
&&+\underset{I_{3}}{\underbrace{\sum_{j,k=1}\partial _{yy}^{2}\left( \bar{b}%
_{t,i}\right) \left( X_{t,i}^{\phi },\mathbf{X}_{t}^{\phi }\right)
Y_{t,j}^{h}Y_{t,k}^{\ell }}}  \notag \\
&&+\underset{I_{4}}{\underbrace{\delta _{h,i}\left( \left( \partial _{x}\phi
_{t,h}^{\prime }\right) \left( X_{t,i}^{\phi },\mathbf{X}_{t}^{\phi }\right)
Y_{t,i}^{\ell }+\left( \mathbf{Y}_{t}^{\ell }\right) ^{\top }\left( \partial
_{y}\phi _{t,h}^{\prime }\right) \left( X_{t,i}^{\phi },\mathbf{X}_{t}^{\phi
}\right) \right) }}  \notag \\
&&+\underset{I_{5}}{\underbrace{\delta _{\ell ,i}\left( \left( \partial
_{x}\phi _{t,\ell }^{\prime \prime }\right) \left( X_{t,i}^{\phi },\mathbf{X}%
_{t}^{\phi }\right) Y_{t,i}^{h}+\left( \mathbf{Y}_{t}^{h}\right) ^{\top
}\left( \partial _{y}\phi _{t,\ell }^{\prime \prime }\right) \left(
X_{t,i}^{\phi },\mathbf{X}_{t}^{\phi }\right) \right) }}.  \label{fti}
\end{eqnarray}%
We will estimate $\left\Vert \mathfrak{f}_{i}^{\phi ,\phi _{h}^{\prime
},\phi _{\ell }^{\prime \prime }}\right\Vert _{\mathcal{H}^{2}\left( \mathbb{%
R}\right) }$ in the following steps. Note that $\delta _{h,i}\delta _{\ell
,i}=0$ as $h\neq \ell .$ Then we get for the term $I_{1}$%
\begin{eqnarray}
\left\Vert \partial _{xx}^{2}\left( \bar{b}_{i}\right) Y_{i}^{h}Y_{i}^{\ell
}\right\Vert _{\mathcal{H}^{2}\left( \mathbb{R}\right) }^{2} &\leq &\left(
L^{\bar{b}}\right) ^{2}\left\Vert Y_{i}^{h}Y_{i}^{\ell }\right\Vert _{%
\mathcal{H}^{2}\left( \mathbb{R}\right) }^{2}  \notag \\
&\leq &C\Bigg [\frac{\left( L_{y}^{\bar{b}}+3\left( L_{y}^{\bar{\sigma}%
}\right) ^{2}\right) ^{\frac{1}{2}}}{\sqrt{N}}\left( \delta _{h,i}+\delta
_{\ell ,i}\right)  \notag \\
&&+\frac{\left( L_{y}^{\bar{b}}+3\left( L_{y}^{\bar{\sigma}}\right)
^{2}\right) }{N}\Bigg ].  \label{f11}
\end{eqnarray}%
Next we deal with $I_{2}$ in (\ref{fti}). The assumption that $\partial
_{xy_{j}}^{2}\left( \bar{b}_{t,i}\right) $ is bounded by $\frac{L_{y}^{\bar{b%
}}}{N}$ and the fact $\left( \sum_{l=1}^{N}c_{l}\right) ^{2}\leq
N\sum_{l=1}^{N}c_{l}^{2}$ for $c_{l}\geq 0,$ $l\in I_{N}$ yield that%
\begin{eqnarray*}
&&\left\Vert \sum_{j=1}\partial _{xy_{j}}^{2}\left( \bar{b}_{i}\right)
\left( X_{i},\mathbf{X}\right) \left( Y_{i}^{h}Y_{j}^{\ell
}+Y_{j}^{h}Y_{i}^{\ell }\right) \right\Vert _{\mathcal{H}^{2}\left( \mathbb{R%
}\right) }^{2} \\
&\leq &\left( \frac{L_{y}^{\bar{b}}}{N}\right) ^{2}\left\Vert
\sum_{j=1}\left( \left\vert Y_{i}^{h}Y_{j}^{\ell }\right\vert +\left\vert
Y_{j}^{h}Y_{i}^{\ell }\right\vert \right) \right\Vert _{\mathcal{H}%
^{2}\left( \mathbb{R}\right) }^{2} \\
&=&\left( \frac{L_{y}^{\bar{b}}}{N}\right) ^{2}\left\Vert \left( \left\vert
Y_{h}^{h}Y_{i}^{\ell }\right\vert +\left\vert Y_{i}^{h}Y_{\ell }^{\ell
}\right\vert +\sum_{j\neq \ell }\left\vert Y_{i}^{h}Y_{j}^{\ell }\right\vert
+\sum_{j\neq h}\left\vert Y_{j}^{h}Y_{i}^{\ell }\right\vert \right)
\right\Vert _{\mathcal{H}^{2}\left( \mathbb{R}\right) }^{2} \\
&\leq &4\left( \frac{L_{y}^{\bar{b}}}{N}\right) ^{2}\Bigg (\left\Vert
Y_{h}^{h}Y_{i}^{\ell }\right\Vert _{\mathcal{H}^{2}\left( \mathbb{R}\right)
}^{2}+\left\Vert Y_{i}^{h}Y_{\ell }^{\ell }\right\Vert _{\mathcal{H}%
^{2}\left( \mathbb{R}\right) }^{2} \\
&&+\left\Vert \sum_{j\neq \ell }\left\vert Y_{i}^{h}Y_{j}^{\ell }\right\vert
\right\Vert _{\mathcal{H}^{2}\left( \mathbb{R}\right) }^{2}+\left\Vert
\sum_{j\neq h}\left\vert Y_{j}^{h}Y_{i}^{\ell }\right\vert \right\Vert _{%
\mathcal{H}^{2}\left( \mathbb{R}\right) }^{2}\Bigg ) \\
&\leq &4\left( \frac{L_{y}^{\bar{b}}}{N}\right) ^{2}\Bigg [\left\Vert
Y_{h}^{h}Y_{i}^{\ell }\right\Vert _{\mathcal{H}^{2}\left( \mathbb{R}\right)
}^{2}+\left\Vert Y_{i}^{h}Y_{\ell }^{\ell }\right\Vert _{\mathcal{H}%
^{2}\left( \mathbb{R}\right) }^{2} \\
&&+\left( N-1\right) \left( \sum_{j\neq \ell }\left\Vert \left\vert
Y_{i}^{h}Y_{j}^{\ell }\right\vert \right\Vert _{\mathcal{H}^{2}\left(
\mathbb{R}\right) }^{2}+\sum_{j\neq h}\left\Vert \left\vert
Y_{j}^{h}Y_{i}^{\ell }\right\vert \right\Vert _{\mathcal{H}^{2}\left(
\mathbb{R}\right) }^{2}\right) \Bigg ]
\end{eqnarray*}
From (\ref{yy}) it follows that%
\begin{eqnarray}
&&\left\Vert \sum_{j=1}\partial _{xy_{j}}^{2}\left( \bar{b}_{i}\right)
\left( X_{i},\mathbf{X}\right) \left( Y_{i}^{h}Y_{j}^{\ell
}+Y_{j}^{h}Y_{i}^{\ell }\right) \right\Vert _{\mathcal{H}^{2}\left( \mathbb{R%
}\right) }^{2}  \notag \\
&\leq &4\left( \frac{L_{y}^{\bar{b}}}{N}\right) ^{2}C\Bigg [\left( \delta
_{\ell ,i}+\frac{\left( L_{y}^{\bar{b}}+3\left( L_{y}^{\bar{\sigma}}\right)
^{2}\right) ^{\frac{1}{2}}}{\sqrt{N}}\left( 1+\delta _{\ell ,i}\right) +%
\frac{\left( L_{y}^{\bar{b}}+3\left( L_{y}^{\bar{\sigma}}\right) ^{2}\right)
}{N}\right)  \notag \\
&&+\left( \delta _{h,i}+\frac{\left( L_{y}^{\bar{b}}+3\left( L_{y}^{\bar{%
\sigma}}\right) ^{2}\right) ^{\frac{1}{2}}}{\sqrt{N}}\left( \delta
_{h,i}+1\right) +\frac{\left( L_{y}^{\bar{b}}+3\left( L_{y}^{\bar{\sigma}%
}\right) ^{2}\right) }{N}\right)  \notag \\
&&+\left( N-1\right) ^{2}\Bigg (\frac{\left( L_{y}^{\bar{b}}+3\left( L_{y}^{%
\bar{\sigma}}\right) ^{2}\right) ^{\frac{1}{2}}}{\sqrt{N}}\delta _{h,i}+%
\frac{\left( L_{y}^{\bar{b}}+3\left( L_{y}^{\bar{\sigma}}\right) ^{2}\right)
}{N}  \notag \\
&&+\frac{\left( L_{y}^{\bar{b}}+3\left( L_{y}^{\bar{\sigma}}\right)
^{2}\right) ^{\frac{1}{2}}}{\sqrt{N}}\delta _{\ell ,i}+\frac{\left( L_{y}^{%
\bar{b}}+3\left( L_{y}^{\bar{\sigma}}\right) ^{2}\right) }{N}\Bigg )\Bigg ]
\notag \\
&\leq &C\left( L_{y}^{\bar{b}}\right) ^{2}\left( \frac{\left( L_{y}^{\bar{b}%
}+3\left( L_{y}^{\bar{\sigma}}\right) ^{2}\right) ^{\frac{1}{2}}}{\sqrt{N}}%
\left( \delta _{h,i}+\delta _{\ell ,i}\right) +\frac{\left( L_{y}^{\bar{b}%
}+3\left( L_{y}^{\bar{\sigma}}\right) ^{2}\right) }{N}\right) .  \label{f22}
\end{eqnarray}%
Note that $h\neq \ell ,$ thus
\begin{eqnarray*}
I_{N}\times I_{N} &=&\left\{ \left( h,\ell \right) \right\} \cup \left\{
\left( h,h\right) \right\} \cup \left\{ \left( \ell ,\ell \right) \right\}
\cup \left\{ \left. \left( h,k\right) \right\vert k\in I_{N}\backslash
\left\{ h,\ell \right\} \right\}  \\
&&\cup \left\{ \left. \left( j,\ell \right) \right\vert j\in I_{N}\backslash
\left\{ h,\ell \right\} \right\} \cup \left\{ \left. \left( j,k\right)
\right\vert j\in I_{N}\backslash \left\{ h\right\} ,k\in I_{N}\backslash
\left\{ \ell \right\} \right\} .
\end{eqnarray*}
Therefore, the term $I_{3}\ $can be decomposed into%
\begin{eqnarray}
&&\sum_{j,k=1}\partial _{yy}^{2}\left( \bar{b}_{t,i}\right) \left(
X_{t,i}^{\phi },\mathbf{X}_{t}^{\phi }\right) Y_{t,j}^{h}Y_{t,k}^{\ell }
\notag \\
&=&\partial _{y_{h}y_{\ell }}^{2}\left( \bar{b}_{t,i}\right) \left(
X_{t,i}^{\phi },\mathbf{X}_{t}^{\phi }\right) Y_{t,h}^{h}Y_{t,\ell }^{\ell
}+\partial _{y_{h}y_{h}}^{2}\left( \bar{b}_{t,i}\right) \left( X_{t,i}^{\phi
},\mathbf{X}_{t}^{\phi }\right) Y_{t,h}^{h}Y_{t,h}^{\ell }  \notag \\
&&+\partial _{y_{\ell }y_{\ell }}^{2}\left( \bar{b}_{t,i}\right) \left(
X_{t,i}^{\phi },\mathbf{X}_{t}^{\phi }\right) Y_{t,\ell }^{h}Y_{t,\ell
}^{\ell }+\sum_{k\in I_{N}\backslash \left\{ h,\ell \right\} }\partial
_{y_{h}y_{k}}^{2}\left( \bar{b}_{t,i}\right) \left( X_{t,i}^{\phi },\mathbf{X%
}_{t}^{\phi }\right) Y_{t,h}^{h}Y_{t,k}^{\ell }  \notag \\
&&+\sum_{j\in I_{N}\backslash \left\{ h,\ell \right\} }\partial
_{y_{j}y_{\ell }}^{2}\left( \bar{b}_{t,i}\right) \left( X_{t,i}^{\phi },%
\mathbf{X}_{t}^{\phi }\right) Y_{t,j}^{h}Y_{t,\ell }^{\ell }  \notag \\
&&+\sum_{j\in I_{N}\backslash \left\{ h\right\} ,k\in I_{N}\backslash
\left\{ \ell \right\} }\partial _{y_{j}y_{k}}^{2}\left( \bar{b}_{t,i}\right)
\left( X_{t,i}^{\phi },\mathbf{X}_{t}^{\phi }\right)
Y_{t,j}^{h}Y_{t,k}^{\ell }.  \label{yyb}
\end{eqnarray}%
It is easy to check that%
\begin{eqnarray}
\left\Vert \partial _{y_{h}y_{\ell }}^{2}\left( \bar{b}_{t,i}\right) \left(
X_{i}^{\phi },\mathbf{X}^{\phi }\right) Y_{h}^{h}Y_{\ell }^{\ell
}\right\Vert _{\mathcal{H}^{2}\left( \mathbb{R}\right) } &\leq &\frac{L_{y}^{%
\bar{b}}}{N^{2}}\left\Vert Y_{h}^{h}Y_{\ell }^{\ell }\right\Vert _{\mathcal{H%
}^{2}\left( \mathbb{R}\right) }\leq \frac{CL_{y}^{\bar{b}}}{N^{2}},
\label{f33} \\
\left\Vert \partial _{y_{h}y_{h}}^{2}\left( \bar{b}_{i}\right) \left(
X_{i}^{\phi },\mathbf{X}^{\phi }\right) Y_{h}^{h}Y_{h}^{\ell }\right\Vert _{%
\mathcal{H}^{2}\left( \mathbb{R}\right) } &\leq &\frac{L_{y}^{\bar{b}}}{N}%
\left\Vert Y_{h}^{h}Y_{h}^{\ell }\right\Vert _{\mathcal{H}^{2}\left( \mathbb{%
R}\right) }\leq \frac{C\left( L_{y}^{\bar{b}}+3\left( L_{y}^{\bar{\sigma}%
}\right) ^{2}\right) ^{\frac{1}{4}}L_{y}^{\bar{b}}}{N^{\frac{5}{4}}},
\label{f44} \\
\left\Vert \partial _{y_{\ell }y_{\ell }}^{2}\left( \bar{b}_{i}\right)
\left( X_{i}^{\phi },\mathbf{X}^{\phi }\right) Y_{\ell }^{h}Y_{\ell }^{\ell
}\right\Vert _{\mathcal{H}^{2}\left( \mathbb{R}\right) } &\leq &\frac{L_{y}^{%
\bar{b}}}{N}\left\Vert Y_{\ell }^{h}Y_{\ell }^{\ell }\right\Vert _{\mathcal{H%
}^{2}\left( \mathbb{R}\right) }\leq \frac{C\left( L_{y}^{\bar{b}}+3\left(
L_{y}^{\bar{\sigma}}\right) ^{2}\right) ^{\frac{1}{4}}L_{y}^{\bar{b}}}{N^{%
\frac{5}{4}}}.  \label{f55}
\end{eqnarray}%
Next
\begin{eqnarray}
&&\left\Vert \sum_{k\in I_{N}\backslash \left\{ h,\ell \right\} }\partial
_{y_{h}y_{k}}^{2}\left( \bar{b}_{i}\right) \left( X_{i}^{\phi },\mathbf{X}%
^{\phi }\right) Y_{h}^{h}Y_{k}^{\ell }\right\Vert _{\mathcal{H}^{2}\left(
\mathbb{R}\right) }+\left\Vert \sum_{j\in I_{N}\backslash \left\{ h,\ell
\right\} }\partial _{y_{j}y_{\ell }}^{2}\left( \bar{b}_{i}\right) \left(
X_{i}^{\phi },\mathbf{X}^{\phi }\right) Y_{j}^{h}Y_{\ell }^{\ell
}\right\Vert _{\mathcal{H}^{2}\left( \mathbb{R}\right) }  \notag \\
&\leq &\frac{L_{y}^{\bar{b}}}{N^{2}}\left( \sum_{k\in I_{N}\backslash
\left\{ h,\ell \right\} }\left\Vert Y_{h}^{h}Y_{k}^{\ell }\right\Vert _{%
\mathcal{H}^{2}\left( \mathbb{R}\right) }+\sum_{j\in I_{N}\backslash \left\{
h,\ell \right\} }\left\Vert Y_{j}^{h}Y_{\ell }^{\ell }\right\Vert _{\mathcal{%
H}^{2}\left( \mathbb{R}\right) }\right)   \notag \\
&\leq &\frac{C\left( L_{y}^{\bar{b}}+3\left( L_{y}^{\bar{\sigma}}\right)
^{2}\right) ^{\frac{1}{4}}L_{y}^{\bar{b}}}{N^{\frac{5}{4}}}.  \label{f66}
\end{eqnarray}
For the last term in (\ref{yyb}),
\begin{eqnarray}
&&\sum_{j\in I_{N}\backslash \left\{ h\right\} ,k\in I_{N}\backslash \left\{
\ell \right\} }\partial _{yy}^{2}\left( \bar{b}_{t,i}\right) \left(
X_{t,i}^{\phi },\mathbf{X}_{t}^{\phi }\right) Y_{t,j}^{h}Y_{t,k}^{\ell }
\notag \\
&=&\sum_{k\in I_{N}\backslash \left\{ \ell \right\} }\partial _{y_{\ell
}y_{k}}^{2}\left( \bar{b}_{t,i}\right) \left( X_{t,i}^{\phi },\mathbf{X}%
_{t}^{\phi }\right) Y_{t,\ell }^{h}Y_{t,k}^{\ell }+\sum_{j\in
I_{N}\backslash \left\{ h,\ell \right\} }\Bigg (\partial
_{y_{j}y_{j}}^{2}\left( \bar{b}_{t,i}\right) \left( X_{t,i}^{\phi },\mathbf{X%
}_{t}^{\phi }\right) Y_{t,j}^{h}Y_{t,j}^{\ell }  \notag \\
&&+\sum_{k\in I_{N}\backslash \left\{ \ell ,j\right\} }\partial
_{y_{j}y_{k}}^{2}\left( \bar{b}_{t,i}\right) \left( X_{t,i}^{\phi },\mathbf{X%
}_{t}^{\phi }\right) Y_{t,j}^{h}Y_{t,k}^{\ell }\Bigg ).  \label{byy}
\end{eqnarray}%
We handle the first and second terms in (\ref{byy}),
\begin{eqnarray}
&&\left\Vert \sum_{k\in I_{N}\backslash \left\{ \ell \right\} }\partial
_{y_{\ell }y_{k}}^{2}\left( \bar{b}_{i}\right) \left( X_{i}^{\phi },\mathbf{X%
}^{\phi }\right) Y_{\ell }^{h}Y_{k}^{\ell }\right\Vert _{\mathcal{H}%
^{2}\left( \mathbb{R}\right) }+\left\Vert \sum_{j\in I_{N}\backslash \left\{
h,\ell \right\} }\partial _{y_{j}y_{j}}^{2}\left( \bar{b}_{i}\right) \left(
X_{i}^{\phi },\mathbf{X}^{\phi }\right) Y_{j}^{h}Y_{j}^{\ell }\right\Vert _{%
\mathcal{H}^{2}\left( \mathbb{R}\right) }  \notag \\
&\leq &\frac{L_{y}^{\bar{b}}}{N^{2}}\sum_{k\in I_{N}\backslash \left\{ \ell
\right\} }\left\Vert Y_{\ell }^{h}Y_{k}^{\ell }\right\Vert _{\mathcal{H}%
^{2}\left( \mathbb{R}\right) }+\frac{L_{y}^{\bar{b}}}{N}\sum_{j\in
I_{N}\backslash \left\{ h,\ell \right\} }\left\Vert Y_{j}^{h}Y_{j}^{\ell
}\right\Vert _{\mathcal{H}^{2}\left( \mathbb{R}\right) }  \notag \\
&\leq &C\frac{L_{y}^{\bar{b}}}{N^{2}}N\frac{\left( L_{y}^{\bar{b}}+3\left(
L_{y}^{\bar{\sigma}}\right) ^{2}\right) ^{\frac{1}{2}}}{N^{\frac{1}{2}}}+C%
\frac{L_{y}^{\bar{b}}}{N}N\frac{\left( L_{y}^{\bar{b}}+3\left( L_{y}^{\bar{%
\sigma}}\right) ^{2}\right) ^{\frac{1}{2}}}{N^{\frac{1}{2}}}  \notag \\
&\leq &\frac{CL_{y}^{\bar{b}}\left( L_{y}^{\bar{b}}+3\left( L_{y}^{\bar{%
\sigma}}\right) ^{2}\right) ^{\frac{1}{2}}}{N^{\frac{1}{2}}}.  \label{f77}
\end{eqnarray}
For the last term in (\ref{byy}),
\begin{eqnarray}
\left\Vert \sum_{j\in I_{N}\backslash \left\{ h,\ell \right\} }\sum_{k\in
I_{N}\backslash \left\{ \ell ,j\right\} }\partial _{y_{j}y_{k}}^{2}\left(
\bar{b}_{i}\right) \left( X_{i}^{\phi },\mathbf{X}^{\phi }\right)
Y_{j}^{h}Y_{k}^{\ell }\right\Vert _{\mathcal{H}^{2}\left( \mathbb{R}\right)
} &\leq &\frac{L_{y}^{\bar{b}}}{N^{2}}\sum_{j\in I_{N}\backslash \left\{
h,\ell \right\} }\sum_{k\in I_{N}\backslash \left\{ \ell ,j\right\}
}\left\Vert Y_{j}^{h}Y_{k}^{\ell }\right\Vert _{\mathcal{H}^{2}\left(
\mathbb{R}\right) }  \notag \\
&\leq &\frac{L_{y}^{\bar{b}}}{N^{2}}N^{2}\frac{\left( L_{y}^{\bar{b}%
}+3\left( L_{y}^{\bar{\sigma}}\right) ^{2}\right) ^{\frac{1}{2}}}{N^{\frac{1%
}{2}}}  \notag \\
&=&\frac{L_{y}^{\bar{b}}\left( L_{y}^{\bar{b}}+3\left( L_{y}^{\bar{\sigma}%
}\right) ^{2}\right) ^{\frac{1}{2}}}{N^{\frac{1}{2}}}.  \label{f88}
\end{eqnarray}
Therefore, from (\ref{f11}), (\ref{f22}), (\ref{f33}), (\ref{f44}), (\ref%
{f55}), (\ref{f66}), (\ref{f77}), (\ref{f88}), we have
\begin{eqnarray}
&&\left\Vert \left(
\begin{array}{c}
Y_{t,i}^{h} \\
\mathbf{Y}_{t}^{h}%
\end{array}%
\right) ^{\top }\left(
\begin{array}{cc}
\partial _{xx}^{2}\left( \bar{b}_{t,i}\right) & \partial _{xy}^{2}\left(
\bar{b}_{t,i}\right) \\
\partial _{yx}^{2}\left( \bar{b}_{t,i}\right) & \partial _{yy}^{2}\left(
\bar{b}_{t,i}\right)%
\end{array}%
\right) \left( X_{t,i}^{\phi },\mathbf{X}_{t}^{\phi }\right) \left(
\begin{array}{c}
Y_{t,i}^{\ell } \\
\mathbf{Y}_{t}^{\ell }%
\end{array}%
\right) \right\Vert _{\mathcal{H}^{2}\left( \mathbb{R}\right) }  \notag \\
&\leq &CL_{y}^{\bar{b}}\left( \frac{\left( L_{y}^{\bar{b}}+3\left( L_{y}^{%
\bar{\sigma}}\right) ^{2}\right) ^{\frac{1}{4}}}{N^{\frac{1}{4}}}\left(
\delta _{h,i}+\delta _{\ell ,i}\right) +\frac{\left( L_{y}^{\bar{b}}+3\left(
L_{y}^{\bar{\sigma}}\right) ^{2}\right) ^{\frac{1}{2}}}{N^{\frac{1}{2}}}%
\right) .  \label{fb1}
\end{eqnarray}%
We now consider the term $\delta _{h,i}\left( \left( \partial _{x}\phi
_{t,h}^{\prime }\right) \left( X_{t,i}^{\phi },\mathbf{X}_{t}^{\phi }\right)
Y_{t,i}^{\ell }+\left( \mathbf{Y}_{t}^{\ell }\right) ^{\top }\left( \partial
_{y}\phi _{t,h}^{\prime }\right) \left( X_{t,i}^{\phi },\mathbf{X}_{t}^{\phi
}\right) \right) $ as follows%
\begin{eqnarray*}
&&\delta _{h,i}\left\Vert \left( \partial _{x}\phi _{h}^{\prime }\right)
\left( X_{i}^{\phi },\mathbf{X}^{\phi }\right) Y_{i}^{\ell }+\left( \mathbf{Y%
}^{\ell }\right) ^{\top }\left( \partial _{y}\phi _{h}^{\prime }\right)
\left( X_{i}^{\phi },\mathbf{X}^{\phi }\right) \right\Vert _{\mathcal{H}%
^{2}\left( \mathbb{R}\right) }^{2} \\
&=&\delta _{h,i}\left\Vert \left( \partial _{x}\phi _{h}^{\prime }\right)
\left( X_{i}^{\phi },\mathbf{X}^{\phi }\right) Y_{i}^{\ell }+\left( \partial
_{y_{\ell }}\phi _{h}^{\prime }\right) \left( X_{i}^{\phi },\mathbf{X}^{\phi
}\right) Y_{\ell }^{\ell }+\sum_{j\neq \ell }\left( \partial _{y_{j}}\phi
_{h}^{\prime }\right) \left( X_{i}^{\phi },\mathbf{X}^{\phi }\right)
Y_{j}^{\ell }\right\Vert _{\mathcal{H}^{2}\left( \mathbb{R}\right) }^{2} \\
&\leq &3\delta _{h,i}\left( \left( L^{\phi _{h}^{\prime }}\right)
^{2}\left\Vert Y_{i}^{\ell }\right\Vert _{\mathcal{H}^{2}\left( \mathbb{R}%
\right) }^{2}+\frac{\left( L_{y}^{\phi _{h}^{\prime }}\right) ^{2}}{N^{2}}%
\left\Vert Y_{\ell }^{\ell }\right\Vert _{\mathcal{H}^{2}\left( \mathbb{R}%
\right) }^{2}+\left\Vert \sum_{j\neq \ell }\left( \partial _{y_{j}}\phi
_{h}^{\prime }\right) \left( X_{i}^{\phi },\mathbf{X}^{\phi }\right)
Y_{j}^{\ell }\right\Vert _{\mathcal{H}^{2}\left( \mathbb{R}\right)
}^{2}\right) \\
&\leq &3\delta _{h,i}\left( \left( L^{\phi _{h}^{\prime }}\right)
^{2}\left\Vert Y_{h}^{\ell }\right\Vert _{\mathcal{H}^{2}\left( \mathbb{R}%
\right) }^{2}+\frac{\left( L_{y}^{\phi _{h}^{\prime }}\right) ^{2}}{N^{2}}%
\left\Vert Y_{\ell }^{\ell }\right\Vert _{\mathcal{H}^{2}\left( \mathbb{R}%
\right) }^{2}+\left( N-1\right) \frac{\left( L_{y}^{\phi _{h}^{\prime
}}\right) ^{2}}{N^{2}}\sum_{j\neq \ell }\left\Vert Y_{j}^{\ell }\right\Vert
_{\mathcal{H}^{2}\left( \mathbb{R}\right) }^{2}\right) \\
&\leq &3T\delta _{h,i}\Bigg [\left( L^{\phi _{h}^{\prime }}\right) ^{2}\frac{%
L_{y}^{\bar{b}}+3\left( L_{y}^{\bar{\sigma}}\right) ^{2}}{N}+\frac{\left(
L_{y}^{\phi _{h}^{\prime }}\right) ^{2}}{N^{2}}\left( \frac{L_{y}^{\bar{b}%
}+3\left( L_{y}^{\bar{\sigma}}\right) ^{2}}{N}4TC_{f_{\ell }}^{p}e^{I_{B,D,%
\bar{B},\bar{D},2}^{3}\cdot T}+4C_{f_{\ell }}^{2}\right) \\
&&+\left( N-1\right) \frac{\left( L_{y}^{\phi _{h}^{\prime }}\right) ^{2}}{%
N^{2}}\sum_{j\neq \ell }\frac{L_{y}^{\bar{b}}+3\left( L_{y}^{\bar{\sigma}%
}\right) ^{2}}{N}4TC_{f_{\ell }}^{p}e^{I_{B,D,\bar{B},\bar{D},p}^{3}\cdot T}%
\Bigg ].
\end{eqnarray*}%
Then%
\begin{eqnarray}
&&\delta _{h,i}\left\Vert \left( \partial _{x}\phi _{h}^{\prime }\right)
\left( X_{i}^{\phi },\mathbf{X}^{\phi }\right) Y_{i}^{\ell }+\left( \mathbf{Y%
}^{\ell }\right) ^{\top }\left( \partial _{y}\phi _{h}^{\prime }\right)
\left( X_{i}^{\phi },\mathbf{X}^{\phi }\right) \right\Vert _{\mathcal{H}%
^{2}\left( \mathbb{R}\right) }^{2}  \notag \\
&\leq &\delta _{h,i}\frac{3T}{N}\Bigg [\left( L^{\phi _{h}^{\prime }}\right)
^{2}\left( L_{y}^{\bar{b}}+3\left( L_{y}^{\bar{\sigma}}\right) ^{2}\right)
+\left( L_{y}^{\phi _{h}^{\prime }}\right) ^{2}\Bigg (\frac{L_{y}^{\bar{b}%
}+3\left( L_{y}^{\bar{\sigma}}\right) ^{2}}{N^{2}}4TC_{f_{\ell
}}^{p}e^{I_{B,D,\bar{B},\bar{D},2}^{3}\cdot T}  \notag \\
&&+\left( L_{y}^{\bar{b}}+3\left( L_{y}^{\bar{\sigma}}\right) ^{2}\right)
4TC_{f_{\ell }}^{p}e^{I_{B,D,\bar{B},\bar{D},p}^{3}\cdot T}+\frac{%
4C_{f_{\ell }}^{2}}{N}\Bigg )\Bigg ]  \notag \\
&\leq &\delta _{h,i}\frac{3T}{N}\Bigg [\left( L^{\phi _{h}^{\prime }}\right)
^{2}\left( L_{y}^{\bar{b}}+3\left( L_{y}^{\bar{\sigma}}\right) ^{2}\right)
\notag \\
&&+\left( L_{y}^{\phi _{h}^{\prime }}\right) ^{2}\Bigg (\left( L_{y}^{\bar{b}%
}+3\left( L_{y}^{\bar{\sigma}}\right) ^{2}\right) 8TC_{f_{\ell
}}^{p}e^{I_{B,D,\bar{B},\bar{D},2}^{3}\cdot T}+4C_{f_{\ell }}^{2}\Bigg )%
\Bigg ]  \notag \\
&\leq &\delta _{h,i}\frac{C}{N}\Bigg [\left( L^{\phi _{h}^{\prime }}\right)
^{2}\left( L_{y}^{\bar{b}}+3\left( L_{y}^{\bar{\sigma}}\right) ^{2}\right)
+\left( L_{y}^{\phi _{h}^{\prime }}\right) ^{2}\Bigg ].  \label{fb2}
\end{eqnarray}%
Analogously,
\begin{eqnarray}
&&\delta _{\ell ,i}\left\Vert \left( \partial _{x}\phi _{t,\ell }^{\prime
\prime }\right) \left( X_{t,i}^{\phi },\mathbf{X}_{t}^{\phi }\right)
Y_{t,i}^{h}+\left( \mathbf{Y}_{t}^{h}\right) ^{\top }\left( \partial
_{y}\phi _{t,\ell }^{\prime \prime }\right) \left( X_{t,i}^{\phi },\mathbf{X}%
_{t}^{\phi }\right) \right\Vert _{\mathcal{H}^{2}\left( \mathbb{R}\right)
}^{2}  \notag \\
&\leq &\delta _{\ell ,i}\frac{C}{N}\Bigg [\left( L^{\phi _{\ell }^{\prime
\prime }}\right) ^{2}\left( L_{y}^{\bar{b}}+3\left( L_{y}^{\bar{\sigma}%
}\right) ^{2}\right) +\left( L_{y}^{\phi _{\ell }^{\prime \prime }}\right)
^{2}\Bigg ].  \label{fb3}
\end{eqnarray}%
From (\ref{fb1}), (\ref{fb2}), (\ref{fb3}), we get the desired result (\ref%
{fineq}). Similarly, we can also derive (\ref{gineq}). The proof is thus
complete.\hfill $\Box $

\paragraph{Proof of Lemma \protect\ref{l5}.}

For convenience of notation, we set $\mathbf{X}=\mathbf{X}^{\phi },$ $%
\mathbf{Y}^{h}=\mathbf{Y}^{\phi ,\phi _{h}^{\prime }}$ and $\mathbf{Z}%
^{h,\ell }=\mathbf{Z}^{\phi ,\phi _{h}^{\prime },\phi _{\ell }^{\prime
\prime }}.$ Employing Lemma \ref{l3} with $\mathbf{S=Z}^{h,\ell },$ for all $%
i\in I_{N}$ and $t\in \left[ 0,T\right] ,$ $B_{i}\left( t\right) =\left(
\partial _{x}\bar{b}_{t,i}\right) \left( X_{t,i},\mathbf{X}_{t}\right) ,$ $%
\bar{B}_{i,j}\left( t\right) =\left( \partial _{y_{j}}\bar{b}_{i}\right)
\left( X_{t,i},\mathbf{X}_{t}\right) ,$ $f_{t,i}=\mathfrak{f}_{t,i}^{\phi
,\phi _{h}^{\prime },\phi _{\ell }^{\prime \prime }},$ $D_{i}\left( t\right)
=\left( \partial _{x}\bar{\sigma}_{i}\right) \left( X_{t,i},\mathbf{X}%
_{t}\right) ,$ $\bar{D}_{i,j}\left( t\right) =\left( \partial _{y_{j}}\bar{%
\sigma}_{i}\right) \left( X_{t,i},\mathbf{X}_{t}\right) ,$ $\bar{f}_{t,i}=%
\mathfrak{g}_{t,i}^{\phi ,\phi _{h}^{\prime },\phi _{\ell }^{\prime \prime
}} $ yields that
\begin{eqnarray*}
\mathbb{E}\left[ \left\vert Z_{t,i}^{h,\ell }\right\vert ^{2}\right]  &\leq &%
\Bigg \{T\left( \left\Vert \bar{B}\right\Vert _{\infty }+3N\left\Vert \bar{D}%
\right\Vert _{\infty }^{2}\right) \cdot \Bigg [\sum_{k=1}^{N}\left(
\left\Vert \mathfrak{f}_{k}^{\phi ,\phi _{h}^{\prime },\phi _{\ell }^{\prime
\prime }}\right\Vert _{\mathcal{H}^{p}\left( \mathbb{R}\right)
}^{2}+3\left\Vert \mathfrak{g}_{k}^{\phi ,\phi _{h}^{\prime },\phi _{\ell
}^{\prime \prime }}\right\Vert _{\mathcal{H}^{p}\left( \mathbb{R}\right)
}^{2}\right) \Bigg ]e^{I_{B,D,\bar{B},\bar{D},2}^{3}\cdot T} \\
&&+\left\Vert \mathfrak{f}_{i}^{\phi ,\phi _{h}^{\prime },\phi _{\ell
}^{\prime \prime }}\right\Vert _{\mathcal{H}^{p}\left( \mathbb{R}\right)
}^{2}+3\left\Vert \mathfrak{g}_{i}^{\phi ,\phi _{h}^{\prime },\phi _{\ell
}^{\prime \prime }}\right\Vert _{\mathcal{H}^{p}\left( \mathbb{R}\right)
}^{2}\Bigg \}e^{I_{B,D,\bar{B},\bar{D},2}^{4}\cdot T},
\end{eqnarray*}
where
\begin{eqnarray*}
I_{B,D,\bar{B},\bar{D},2}^{3} &=&2\left\Vert B\right\Vert _{\infty
}+2N\left\Vert \bar{B}\right\Vert _{\infty }+1+3\left\Vert D\right\Vert
_{\infty }^{2}+3N^{2}\left\Vert \bar{D}\right\Vert _{\infty }^{2}, \\
I_{B,D,\bar{B},\bar{D},2}^{4} &=&2\left\Vert B\right\Vert _{\infty
}+N\left\Vert \bar{B}\right\Vert _{\infty }+1+3\left\Vert D\right\Vert
_{\infty }^{2},
\end{eqnarray*}%
additionally, $\left\Vert B\right\Vert _{\infty }=\max_{i\in
I_{N}}\left\Vert B_{i}\right\Vert _{L^{\infty }},$ $\left\Vert D\right\Vert
_{\infty }=\max_{i\in I_{N}}\left\Vert D_{i}\right\Vert _{L^{\infty }}$, $%
\left\Vert \bar{B}\right\Vert _{\infty }=\max_{i,j\in I_{N}}\left\Vert \bar{B%
}_{i,j}\right\Vert _{L^{\infty }}$ and $\left\Vert \bar{D}\right\Vert
_{\infty }=\max_{i,j\in I_{N}}\left\Vert \bar{D}_{i,j}\right\Vert
_{L^{\infty }}.$

From Lemma \ref{f1g1}, we have%
\begin{eqnarray}
\left\Vert \mathfrak{f}_{i}^{\phi ,\phi _{h}^{\prime },\phi _{\ell }^{\prime
\prime }}\right\Vert _{\mathcal{H}^{2}\left( \mathbb{R}\right) }^{2} &\leq &%
\frac{C\left( L_{y}^{\bar{b}}\right) ^{2}\left( L_{y}^{\bar{b}}+3\left(
L_{y}^{\bar{\sigma}}\right) ^{2}\right) ^{\frac{1}{2}}}{N^{\frac{1}{2}}}%
\left( \delta _{h,i}+\delta _{\ell ,i}\right) +\frac{C\left( L_{y}^{\bar{b}%
}\right) ^{2}\left( L_{y}^{\bar{b}}+3\left( L_{y}^{\bar{\sigma}}\right)
^{2}\right) }{N}  \notag \\
&&+\delta _{h,i}\frac{C}{N}\Bigg [L^{\phi _{h}^{\prime }}\left( L_{y}^{\bar{b%
}}+3\left( L_{y}^{\bar{\sigma}}\right) ^{2}\right) ^{\frac{1}{2}%
}+L_{y}^{\phi _{h}^{\prime }}\Bigg ]^{2}  \notag \\
&&+\delta _{\ell ,i}\frac{C}{N^{{}}}\Bigg [L^{\phi _{\ell }^{\prime \prime
}}\left( L_{y}^{\bar{b}}+3\left( L_{y}^{\bar{\sigma}}\right) ^{2}\right) ^{%
\frac{1}{2}}+L_{y}^{\phi _{\ell }^{\prime \prime }}\Bigg ]^{2}, \\
\left\Vert \mathfrak{g}_{i}^{\phi ,\phi _{h}^{\prime },\phi _{\ell }^{\prime
\prime }}\right\Vert _{\mathcal{H}^{2}\left( \mathbb{R}\right) }^{2} &\leq &%
\frac{C\left( L_{y}^{\bar{\sigma}}\right) ^{2}\left( L_{y}^{\bar{b}}+3\left(
L_{y}^{\bar{\sigma}}\right) ^{2}\right) ^{\frac{1}{2}}}{N^{\frac{1}{2}}}%
\left( \delta _{h,i}+\delta _{\ell ,i}\right) +\frac{C\left( L_{y}^{\bar{%
\sigma}}\right) ^{2}\left( L_{y}^{\bar{b}}+3\left( L_{y}^{\bar{\sigma}%
}\right) ^{2}\right) }{N}  \notag \\
&&+\delta _{h,i}\frac{C}{N}\Bigg [L^{\phi _{h}^{\prime }}\left( L_{y}^{\bar{b%
}}+3\left( L_{y}^{\bar{\sigma}}\right) ^{2}\right) ^{\frac{1}{2}%
}+L_{y}^{\phi _{h}^{\prime }}\Bigg ]^{2}  \notag \\
&&+\delta _{\ell ,i}\frac{C}{N}\Bigg [L^{\phi _{\ell }^{\prime \prime
}}\left( L_{y}^{\bar{b}}+3\left( L_{y}^{\bar{\sigma}}\right) ^{2}\right) ^{%
\frac{1}{2}}+L_{y}^{\phi _{\ell }^{\prime \prime }}\Bigg ]^{2}.
\end{eqnarray}%
Then%
\begin{eqnarray*}
&&\sum_{k=1}^{N}\left( \left\Vert \mathfrak{f}_{k}^{\phi ,\phi _{h}^{\prime
},\phi _{\ell }^{\prime \prime }}\right\Vert _{\mathcal{H}^{p}\left( \mathbb{%
R}\right) }^{2}+3\left\Vert \mathfrak{g}_{k}^{\phi ,\phi _{h}^{\prime },\phi
_{\ell }^{\prime \prime }}\right\Vert _{\mathcal{H}^{p}\left( \mathbb{R}%
\right) }^{2}\right) \\
&=&\sum_{k\in \left\{ h,\ell \right\} }\left( \left\Vert \mathfrak{f}%
_{k}^{\phi ,\phi _{h}^{\prime },\phi _{\ell }^{\prime \prime }}\right\Vert _{%
\mathcal{H}^{p}\left( \mathbb{R}\right) }^{2}+3\left\Vert \mathfrak{g}%
_{k}^{\phi ,\phi _{h}^{\prime },\phi _{\ell }^{\prime \prime }}\right\Vert _{%
\mathcal{H}^{p}\left( \mathbb{R}\right) }^{2}\right) \\
&&+\sum_{k\in I_{N}\backslash \left\{ h,\ell \right\} }\left( \left\Vert
\mathfrak{f}_{k}^{\phi ,\phi _{h}^{\prime },\phi _{\ell }^{\prime \prime
}}\right\Vert _{\mathcal{H}^{p}\left( \mathbb{R}\right) }^{2}+3\left\Vert
\mathfrak{g}_{k}^{\phi ,\phi _{h}^{\prime },\phi _{\ell }^{\prime \prime
}}\right\Vert _{\mathcal{H}^{p}\left( \mathbb{R}\right) }^{2}\right) \\
&\leq &\frac{C\left( L_{y}^{\bar{b}}\right) ^{2}\left( L_{y}^{\bar{b}%
}+3\left( L_{y}^{\bar{\sigma}}\right) ^{2}\right) ^{\frac{1}{2}}}{N^{\frac{1%
}{2}}}+\frac{C\left( L_{y}^{\bar{\sigma}}\right) ^{2}\left( L_{y}^{\bar{b}%
}+3\left( L_{y}^{\bar{\sigma}}\right) ^{2}\right) ^{\frac{1}{2}}}{N^{\frac{1%
}{2}}} \\
&&+\left( N-2\right) \frac{C\left( L_{y}^{\bar{b}}\right) ^{2}\left( L_{y}^{%
\bar{b}}+3\left( L_{y}^{\bar{\sigma}}\right) ^{2}\right) }{N}+\left(
N-2\right) \frac{C\left( L_{y}^{\bar{\sigma}}\right) ^{2}\left( L_{y}^{\bar{b%
}}+3\left( L_{y}^{\bar{\sigma}}\right) ^{2}\right) }{N}.
\end{eqnarray*}%
We then get the desired result. \hfill $\Box $

\paragraph{Proof of Lemma \protect\ref{l6}.}

For convenience of notation, we set $\mathbf{X}=\mathbf{X}^{\phi },$ $%
\mathbf{Y}^{h}=\mathbf{Y}^{\phi ,\phi _{h}^{\prime }}$ and $\mathbf{Z}%
^{h,\ell }=\mathbf{Z}^{\phi ,\phi _{h}^{\prime },\phi _{\ell }^{\prime
\prime }}.$ By lemmata \ref{l3} and \ref{mvt} for all $i\in I_{N}$ and $%
t\in \left[ 0,T\right] ,$ it follows that%
\begin{eqnarray*}
\mathbb{E}\left[ \left\vert u_{t,i}^{\phi }\right\vert ^{2}\right] &=&%
\mathbb{E}\left[ \left\vert \phi _{t,i}\left( X_{t,i}^{\phi },\mathbf{X}%
_{t}^{\phi }\right) \right\vert ^{2}\right] \\
&\leq &\left( L^{\phi }\right) ^{2}\left( 1+\mathbb{E}\left[ \left\vert
X_{t,i}^{\phi }\right\vert ^{2}\right] \right) +\frac{\left( L_{y}^{\phi
}\right) ^{2}}{N}\sum_{j=1}^{N}\mathbb{E}\left[ \left\vert X_{t,j}^{\phi
}\right\vert ^{2}\right] \\
&\leq &C,
\end{eqnarray*}%
from which we obtain (\ref{l61}). Next from Lemma \ref{l21} we have%
\begin{eqnarray*}
\mathbb{E}\left[ \left\vert \upsilon _{t,i}^{\phi ,\phi _{h}^{\prime
}}\right\vert ^{2}\right] &\leq &3\Bigg \{\mathbb{E}\left[ \left\vert \left(
\partial _{x}\phi _{t,i}\right) \left( X_{t,i},\mathbf{X}_{t}\right)
Y_{t,i}^{h}\right\vert ^{2}\right] +\mathbb{E}\left[ \left\vert \left(
\mathbf{Y}_{t}^{h}\right) ^{\top }\left( \partial _{y}\phi _{t,i}\right)
\left( X_{t,i},\mathbf{X}_{t}\right) \right\vert ^{2}\right] \\
&&+\delta _{h,i}\mathbb{E}\left[ \left\vert \phi _{t,h}^{\prime }\left(
X_{t,i},\mathbf{X}_{t}\right) \right\vert ^{2}\right] \Bigg \} \\
&=&3\Bigg \{\mathbb{E}\left[ \left\vert \left( \partial _{x}\phi
_{t,i}\right) \left( X_{t,i},\mathbf{X}_{t}\right) Y_{t,i}^{h}\right\vert
^{2}\right] +\mathbb{E}\left[ \left\vert \sum_{j=1}^{N}\left( \partial
_{y_{j}}\phi _{t,i}\right) \left( X_{t,i},\mathbf{X}_{t}\right)
Y_{t,j}^{h}\right\vert ^{2}\right] \\
&&+\delta _{h,i}\mathbb{E}\left[ \left\vert \phi _{t,h}^{\prime }\left(
X_{t,i},\mathbf{X}_{t}\right) \right\vert ^{2}\right] \Bigg \} \\
&\leq &3\Bigg \{\mathbb{E}\left[ \left\vert \left( \partial _{x}\phi
_{t,i}\right) \left( X_{t,i},\mathbf{X}_{t}\right) Y_{t,i}^{h}\right\vert
^{2}\right] +\frac{\left( L_{y}^{\phi }\right) ^{2}}{N^{2}}\mathbb{E}\Bigg [%
\left\vert Y_{t,h}^{h}\right\vert ^{2}+\left( N-1\right) \sum_{j\neq
h}^{N}\left\vert Y_{t,j}^{h}\right\vert ^{2}\Bigg ] \\
&&+\delta _{h,i}\mathbb{E}\Bigg [\Bigg |L^{\phi _{h}^{\prime }}\left(
1+\left\vert X_{t,i}\right\vert \right) +\frac{L_{y}^{\phi _{h}^{\prime }}}{N%
}\sum_{j=1}^{N}\left\vert X_{t,j}\right\vert \Bigg |^{2}\Bigg ]\Bigg \}.
\end{eqnarray*}%
By Lemma \ref{l4}, we get the desired (\ref{l62}). Now we deal with $\omega
_{i}^{\phi ,\phi _{h}^{\prime },\phi _{\ell }^{\prime \prime }}.$ Similarly,
we consider
\begin{eqnarray}
&&\left\Vert \left(
\begin{array}{c}
Y_{t,i}^{h} \\
\mathbf{Y}_{t}^{h}%
\end{array}%
\right) ^{\top }\left(
\begin{array}{cc}
\partial _{xx}^{2}\phi _{t,i} & \partial _{xy}^{2}\phi _{t,i} \\
\partial _{yx}^{2}\phi _{t,i} & \partial _{yy}^{2}\phi _{t,i}%
\end{array}%
\right) \left( X_{t,i}^{\phi },\mathbf{X}_{t}^{\phi }\right) \left(
\begin{array}{c}
Y_{t,i}^{\ell } \\
\mathbf{Y}_{t}^{\ell }%
\end{array}%
\right) \right\Vert _{\mathcal{H}^{2}\left( \mathbb{R}\right) }  \notag \\
&\leq &C\left( \frac{\left( L_{y}^{\bar{b}}+3\left( L_{y}^{\bar{\sigma}%
}\right) ^{2}\right) ^{\frac{1}{4}}}{N^{\frac{1}{4}}}\left( \delta
_{h,i}+\delta _{\ell ,i}\right) +\frac{\left( L_{y}^{\bar{b}}+3\left( L_{y}^{%
\bar{\sigma}}\right) ^{2}\right) ^{\frac{1}{2}}}{N^{\frac{1}{2}}}\right)
\max \left\{ L_{y}^{\bar{b}},L_{y}^{\phi }\right\} .  \label{w1}
\end{eqnarray}%
Next
\begin{eqnarray}
&&\left\Vert \left( \partial _{x}\phi _{i}\right) \left( X_{i}^{\phi },%
\mathbf{X}^{\phi }\right) Z_{i}^{h,\ell }\right\Vert _{\mathcal{H}^{2}\left(
\mathbb{R}\right) }+\left\Vert \left( \mathbf{Z}^{h,\ell }\right) ^{\top
}\left( \partial _{y}\phi _{i}\right) \left( X_{i}^{\phi },\mathbf{X}^{\phi
}\right) \right\Vert _{\mathcal{H}^{2}\left( \mathbb{R}\right) }  \notag \\
&=&\left\Vert \left( \partial _{x}\phi _{i}\right) \left( X_{i}^{\phi },%
\mathbf{X}^{\phi }\right) Z_{i}^{h,\ell }\right\Vert _{\mathcal{H}^{2}\left(
\mathbb{R}\right) }+\left\Vert \sum_{j=1}^{N}\left( \partial _{y_{j}}\phi
_{i}\right) \left( X_{i}^{\phi },\mathbf{X}^{\phi }\right) Z_{j}^{h,\ell
}\right\Vert _{\mathcal{H}^{2}\left( \mathbb{R}\right) }  \notag \\
&\leq &\left( L^{\phi }\right) ^{2}\left\Vert Z_{i}^{h,\ell }\right\Vert _{%
\mathcal{H}^{2}\left( \mathbb{R}\right) }+\sum_{j\in \left\{ h,\ell \right\}
}^{N}\left\Vert \left( \partial _{y_{j}}\phi _{i}\right) \left( X_{i}^{\phi
},\mathbf{X}^{\phi }\right) Z_{j}^{h,\ell }\right\Vert _{\mathcal{H}%
^{2}\left( \mathbb{R}\right) }  \notag \\
&&+\sum_{j\in I_{N}\backslash \left\{ h,\ell \right\} }^{N}\left\Vert \left(
\partial _{y_{j}}\phi _{i}\right) \left( X_{i}^{\phi },\mathbf{X}^{\phi
}\right) Z_{j}^{h,\ell }\right\Vert _{\mathcal{H}^{2}\left( \mathbb{R}%
\right) }  \notag \\
&\leq &\left( L^{\phi }\right) ^{2}(\frac{CL_{y}^{\bar{b}}\left( L_{y}^{\bar{%
b}}+3\left( L_{y}^{\bar{\sigma}}\right) ^{2}\right) ^{\frac{1}{4}}}{N^{\frac{%
1}{4}}}\left( \delta _{h,i}+\delta _{\ell ,i}\right) +\frac{CL_{y}^{\bar{b}%
}\left( L_{y}^{\bar{b}}+3\left( L_{y}^{\bar{\sigma}}\right) ^{2}\right) ^{%
\frac{1}{2}}}{N^{\frac{1}{2}}}  \notag \\
&&+\frac{CL_{y}^{\bar{\sigma}}\left( L_{y}^{\bar{b}}+3\left( L_{y}^{\bar{%
\sigma}}\right) ^{2}\right) ^{\frac{1}{4}}}{N^{\frac{1}{4}}}\left( \delta
_{h,i}+\delta _{\ell ,i}\right) +\frac{CL_{y}^{\bar{\sigma}}\left( L_{y}^{%
\bar{b}}+3\left( L_{y}^{\bar{\sigma}}\right) ^{2}\right) ^{\frac{1}{2}}}{N^{%
\frac{1}{2}}})  \notag \\
&&+\frac{CL_{y}^{\bar{b}}\left( L_{y}^{\bar{b}}+3\left( L_{y}^{\bar{\sigma}%
}\right) ^{2}\right) ^{\frac{1}{2}}}{N^{\frac{1}{2}}}+\frac{CL_{y}^{\bar{%
\sigma}}\left( L_{y}^{\bar{b}}+3\left( L_{y}^{\bar{\sigma}}\right)
^{2}\right) ^{\frac{1}{2}}}{N^{\frac{1}{2}}}  \notag \\
&&+\left( N-2\right) \frac{L_{y}^{\phi }}{N}\left( \frac{CL_{y}^{\bar{b}%
}\left( L_{y}^{\bar{b}}+3\left( L_{y}^{\bar{\sigma}}\right) ^{2}\right) ^{%
\frac{1}{2}}}{N^{\frac{1}{2}}}+\frac{CL_{y}^{\bar{\sigma}}\left( L_{y}^{\bar{%
b}}+3\left( L_{y}^{\bar{\sigma}}\right) ^{2}\right) ^{\frac{1}{2}}}{N^{\frac{%
1}{2}}}\right)  \notag \\
&\leq &C\Bigg (\frac{L_{y}^{\bar{b}}\left( L_{y}^{\bar{b}}+3\left( L_{y}^{%
\bar{\sigma}}\right) ^{2}\right) ^{\frac{1}{4}}}{N^{\frac{1}{4}}}\left(
\delta _{h,i}+\delta _{\ell ,i}\right) +\frac{L_{y}^{\bar{\sigma}}\left(
L_{y}^{\bar{b}}+3\left( L_{y}^{\bar{\sigma}}\right) ^{2}\right) ^{\frac{1}{4}%
}}{N^{\frac{1}{4}}}\left( \delta _{h,i}+\delta _{\ell ,i}\right)  \notag \\
&&+\frac{L_{y}^{\bar{b}}\left( L_{y}^{\bar{b}}+3\left( L_{y}^{\bar{\sigma}%
}\right) ^{2}\right) ^{\frac{1}{2}}}{N^{\frac{1}{2}}}+\frac{CL_{y}^{\bar{%
\sigma}}\left( L_{y}^{\bar{b}}+3\left( L_{y}^{\bar{\sigma}}\right)
^{2}\right) ^{\frac{1}{2}}}{N^{\frac{1}{2}}}\Bigg ).  \label{w2}
\end{eqnarray}%
Next
\begin{eqnarray}
&&\delta _{h,i}\left\Vert \left( \partial _{x}\phi _{h}^{\prime }\right)
\left( X_{t,i},\mathbf{X}_{t}\right) Y_{t,i}^{\phi ,\phi _{\ell }^{\prime
\prime }}+\left( \mathbf{Y}_{t}^{\phi ,\phi _{\ell }^{\prime \prime
}}\right) ^{\top }\left( \partial _{y}\phi _{h}^{\prime }\right) \left(
X_{t,i},\mathbf{X}_{t}\right) \right\Vert _{\mathcal{H}^{2}\left( \mathbb{R}%
\right) }^{2}  \notag \\
&&+\delta _{\ell ,i}\left\Vert \left( \partial _{x}\phi _{\ell }^{\prime
\prime }\right) \left( X_{t,i},\mathbf{X}_{t}\right) Y_{t,i}^{\phi ,\phi
_{h}^{\prime }}+\left( \mathbf{Y}_{t}^{\phi ,\phi _{h}^{\prime }}\right)
^{\top }\left( \partial _{y}\phi _{\ell }^{\prime \prime }\right) \left(
X_{t,i},\mathbf{X}_{t}\right) \right\Vert _{\mathcal{H}^{2}\left( \mathbb{R}%
\right) }^{2}  \notag \\
&=&\delta _{h,i}\left\Vert \left( \partial _{x}\phi _{h}^{\prime }\right)
\left( X_{t,i}^{\phi },\mathbf{X}_{t}^{\phi }\right) Y_{t,i}^{\phi ,\phi
_{\ell }^{\prime \prime }}+\sum_{k=1}^{N}\left( \partial _{y_{k}}\phi
_{h}^{\prime }\right) \left( X_{t,i},\mathbf{X}_{t}\right) Y_{t,k}^{\phi
,\phi _{\ell }^{\prime \prime }}\right\Vert _{\mathcal{H}^{2}\left( \mathbb{R%
}\right) }^{2}  \notag \\
&&+\delta _{\ell ,i}\left\Vert \left( \partial _{x}\phi _{\ell }^{\prime
\prime }\right) \left( X_{t,i},\mathbf{X}_{t}\right) Y_{t,i}^{\phi ,\phi
_{h}^{\prime }}+\sum_{k=1}^{N}\left( \partial _{y_{k}}\phi _{\ell }^{\prime
\prime }\right) \left( X_{t,i},\mathbf{X}_{t}\right) Y_{t,k}^{\phi ,\phi
_{h}^{\prime }}\right\Vert _{\mathcal{H}^{2}\left( \mathbb{R}\right) }^{2}
\notag \\
&\leq &\delta _{h,i}\left( L^{\phi _{h}^{\prime }}\left\Vert Y_{t,i}^{\phi
,\phi _{\ell }^{\prime \prime }}\right\Vert _{\mathcal{H}^{2}\left( \mathbb{R%
}\right) }^{2}+L_{y}^{\phi _{h}^{\prime }}\left\Vert Y_{t,\ell }^{\phi ,\phi
_{\ell }^{\prime \prime }}\right\Vert _{\mathcal{H}^{2}\left( \mathbb{R}%
\right) }^{2}+L_{y}^{\phi _{h}^{\prime }}\sum_{k\in I_{N}\backslash \left\{
\ell \right\} }\left\Vert Y_{t,k}^{\phi ,\phi _{\ell }^{\prime \prime
}}\right\Vert _{\mathcal{H}^{2}\left( \mathbb{R}\right) }^{2}\right)  \notag
\\
&&+\delta _{\ell ,i}\left( L^{\phi _{\ell }^{\prime \prime }}\left\Vert
Y_{t,i}^{\phi ,\phi _{h}^{\prime }}\right\Vert _{\mathcal{H}^{2}\left(
\mathbb{R}\right) }^{2}+L_{y}^{\phi _{\ell }^{\prime \prime }}\left\Vert
Y_{t,h}^{\phi ,\phi _{h}^{\prime }}\right\Vert _{\mathcal{H}^{2}\left(
\mathbb{R}\right) }^{2}+L_{y}^{\phi _{\ell }^{\prime \prime }}\sum_{k\in
I_{N}\backslash \left\{ h\right\} }\left\Vert Y_{t,k}^{\phi ,\phi
_{h}^{\prime }}\right\Vert _{\mathcal{H}^{2}\left( \mathbb{R}\right)
}^{2}\right)  \notag \\
&\leq &\delta _{h,i}\Bigg (\frac{L^{\phi _{h}^{\prime }}\left( L_{y}^{\bar{b}%
}+3\left( L_{y}^{\bar{\sigma}}\right) ^{2}\right) C}{N}+\delta _{\ell
,i}C+L_{y}^{\phi _{h}^{\prime }}C+L_{y}^{\phi _{h}^{\prime }}C\left(
N-1\right) \frac{L_{y}^{\bar{b}}+3\left( L_{y}^{\bar{\sigma}}\right) ^{2}}{N}%
\Bigg )  \notag \\
&&+\delta _{\ell ,i}\Bigg (\frac{L^{\phi _{\ell }^{\prime \prime }}\left(
L_{y}^{\bar{b}}+3\left( L_{y}^{\bar{\sigma}}\right) ^{2}\right) C}{N}+\delta
_{h,i}C+L_{y}^{\phi _{\ell }^{\prime \prime }}C+L_{y}^{\phi _{\ell }^{\prime
\prime }}C\left( N-1\right) \frac{L_{y}^{\bar{b}}+3\left( L_{y}^{\bar{\sigma}%
}\right) ^{2}}{N}\Bigg ).  \label{w3}
\end{eqnarray}%
From (\ref{w1})-(\ref{w3}), we prove (\ref{l63}) immediately. We thus
complete the proof. \hfill $\Box $

\paragraph{Proof of Theorem \protect\ref{the2}.}

For the case $i=j,$ it is trivial to get $\frac{\delta ^{2}V_{i}}{\delta
\phi _{i}\delta \phi _{j}}\left( \phi ,\phi _{i}^{\prime },\phi _{j}^{\prime
\prime }\right) =\frac{\delta ^{2}V_{j}}{\delta \phi _{j}\delta \phi _{i}}%
\left( \phi ,\phi _{j}^{\prime \prime },\phi _{i}^{\prime }\right) .$ We
explore the case $i\neq j.$ From the expression of $\tilde{C}_{0}^{i,j},$ we
get
\begin{eqnarray*}
&&\left\Vert \partial _{x_{i}x_{j}}^{2}\Delta _{i,j}^{f}\right\Vert
_{L^{\infty }}++\left\Vert \partial _{x_{i}u_{j}}^{2}\Delta _{i,j}^{f}\left(
\mathbf{X}_{t},\mathbf{u}_{t}\right) \right\Vert _{L^{\infty }}+\left\Vert
\partial _{u_{i}x_{j}}^{2}\Delta _{i,j}^{f}\left( \mathbf{X}_{t},\mathbf{u}%
_{t}\right) \right\Vert _{L^{\infty }} \\
&&+\left\Vert \left( \partial _{u_{i}u_{j}}^{2}\Delta _{i,j}^{f}\left(
\mathbf{X}_{t},\mathbf{u}_{t}\right) \right) \right\Vert _{L^{\infty
}}+\left\Vert \partial _{x_{i}x_{j}}^{2}\Delta _{i,j}^{g}\right\Vert
_{L^{\infty }} \\
&&+\Bigg [\sqrt{L_{y}^{\phi _{h}^{\prime }}\left( 1+L_{y}^{\bar{b}}+3\left(
L_{y}^{\bar{\sigma}}\right) ^{2}\right) }+\sqrt{L_{y}^{\phi _{\ell }^{\prime
\prime }}\left( 1+L_{y}^{\bar{b}}+3\left( L_{y}^{\bar{\sigma}}\right)
^{2}\right) }\Bigg ]\cdot \sum_{h\in \left\{ i,j\right\} }^{N}\Big [%
\left\Vert \left( \partial _{u_{h}}\Delta _{i,j}^{f}\right) \left(
0,0\right) \right\Vert _{L^{\infty }} \\
&&+\sum_{k=1}^{N}\left( \left\Vert \partial _{u_{_{h}}x_{k}}\Delta
_{i,j}^{f}\right\Vert _{L^{\infty }}+\left\Vert \partial _{u_{h}u_{k}}\Delta
_{i,j}^{f}\right\Vert _{L^{\infty }}\right) \Big ] \\
&\leq &C\left( \sqrt{L_{y}^{\phi _{h}^{\prime }}\left( 1+L_{y}^{\bar{b}%
}+3\left( L_{y}^{\bar{\sigma}}\right) ^{2}\right) }+\sqrt{L_{y}^{\phi _{\ell
}^{\prime \prime }}\left( 1+L_{y}^{\bar{b}}+3\left( L_{y}^{\bar{\sigma}%
}\right) ^{2}\right) }\right) \left( 2L+\left( N-1\right) \frac{\tilde{L}}{%
N^{2\beta }}\right) +\frac{C\tilde{L}}{N^{2\beta }}.
\end{eqnarray*}%
For $\tilde{C}_{1}^{i,j},$
\begin{eqnarray*}
\tilde{C}_{1}^{i,j} &=&\left( 1+L_{y}^{\bar{b}}+L_{y}^{\bar{\sigma}}\right)
\left( L_{y}^{\bar{b},\bar{\sigma}}\right) ^{\frac{1}{4}}\sum_{h\in \left\{
i,j\right\} }^{N}\Big [\left\Vert \left( \partial _{u_{h}}\Delta
_{i,j}^{f}\right) \left( 0,0\right) \right\Vert _{L^{\infty
}}+\sum_{k=1}^{N}\left( \left\Vert \partial _{u_{_{h}}x_{k}}\Delta
_{i,j}^{f}\right\Vert _{L^{\infty }}+\left\Vert \partial _{u_{h}u_{k}}\Delta
_{i,j}^{f}\right\Vert _{L^{\infty }}\right) \Big ] \\
&\leq &C\left( 1+L_{y}^{\bar{b}}+L_{y}^{\bar{\sigma}}\right) \left( L_{y}^{%
\bar{b},\bar{\sigma}}\right) ^{\frac{1}{4}}\left[ L+\frac{2\tilde{L}}{%
N^{\beta }}+\left( N-1\right) \frac{\tilde{L}}{N^{2\beta }}\right] .
\end{eqnarray*}%
For $\tilde{C}_{2}^{i,j},$%
\begin{eqnarray*}
\tilde{C}_{2}^{i,j} &=&\left( 1+L_{y}^{\bar{b}}+L_{y}^{\bar{\sigma}}\right)
\left( L_{y}^{\bar{b},\bar{\sigma}}\right) ^{\frac{1}{2}}\sum_{h\in
I_{N}\backslash \left\{ i,j\right\} }^{N}\Big [\left\Vert \left( \partial
_{u_{h}}\Delta _{i,j}^{f}\right) \left( 0,0\right) \right\Vert _{L^{\infty
}}+\sum_{k=1}^{N}\left( \left\Vert \partial _{u_{_{h}}x_{k}}\Delta
_{i,j}^{f}\right\Vert _{L^{\infty }}+\left\Vert \partial _{u_{h}u_{k}}\Delta
_{i,j}^{f}\right\Vert _{L^{\infty }}\right) \Big ] \\
&&+\left( L_{y}^{\bar{b},\bar{\sigma}}\right) ^{\frac{1}{2}}\left(
\sum_{\ell \in I_{N}\backslash \left\{ j\right\} }^{N}\left\Vert \partial
_{x_{i}x_{\ell }}^{2}\Delta _{t,i,j}^{f}\right\Vert _{L^{\infty
}}+\sum_{h\in I_{N}\backslash \left\{ i\right\} }^{N}\left\Vert \partial
_{x_{h}x_{j}}^{2}\Delta _{i,j}^{f}\right\Vert _{L^{\infty }}\right)  \\
&&+\left( L_{y}^{\bar{b},\bar{\sigma}}\right) ^{\frac{1}{2}}\sum_{h\in
I_{N}\backslash \left\{ i\right\} }^{N}\left\Vert \partial
_{x_{h}u_{j}}^{2}\Delta _{t,i,j}^{f}\left( \mathbf{X}_{t},\mathbf{u}%
_{t}\right) \right\Vert _{L^{\infty }}+\left[ L_{y}^{\bar{b},\bar{\sigma}%
}\left( \left( L_{y}^{\phi }\right) ^{2}+\left( L^{\phi }\right) ^{2}\right) %
\right] ^{\frac{1}{2}}\sum_{\ell \in I_{N}\backslash \left\{ j\right\}
}^{N}\left\Vert \partial _{x_{i}u_{\ell }}^{2}\Delta _{t,i,j}^{f}\left(
\mathbf{X}_{t},\mathbf{u}_{t}\right) \right\Vert _{L^{\infty }} \\
&&+\left( L_{y}^{\bar{b},\bar{\sigma}}\right) ^{\frac{1}{2}}\sum_{\ell \in
I_{N}\backslash \left\{ j\right\} }^{N}\left\Vert \partial _{u_{i}x_{\ell
}}^{2}\Delta _{t,i,j}^{f}\left( \mathbf{X}_{t},\mathbf{u}_{t}\right)
\right\Vert _{L^{\infty }}+\left[ L_{y}^{\bar{b},\bar{\sigma}}\left( \left(
L_{y}^{\phi }\right) ^{2}+\left( L^{\phi }\right) ^{2}\right) \right] ^{%
\frac{1}{2}}\sum_{h\in I_{N}\backslash \left\{ i\right\} }^{N}\left\Vert
\partial _{u_{h}x_{j}}^{2}\Delta _{t,i,j}^{f}\left( \mathbf{X}_{t},\mathbf{u}%
_{t}\right) \right\Vert _{L^{\infty }} \\
&&+\left( \left( L_{y}^{\phi }\right) ^{2}+\left( L^{\phi }\right)
^{2}\right) ^{\frac{1}{2}}\left( L_{y}^{\bar{b},\bar{\sigma}}\right) ^{\frac{%
1}{2}}\Bigg (\sum_{\ell \in I_{N}\backslash \left\{ j\right\}
}^{N}\left\Vert \partial _{u_{i}u_{\ell }}^{2}\Delta _{t,i,j}^{f}\left(
\mathbf{X}_{t},\mathbf{u}_{t}\right) \right\Vert _{L^{\infty }}+\sum_{h\in
I_{N}\backslash \left\{ i\right\} }^{N}\left\Vert \partial
_{u_{h}u_{j}}^{2}\Delta _{i,j}^{f}\left( \mathbf{X}_{t},\mathbf{u}%
_{t}\right) \right\Vert _{L^{\infty }}\Bigg ) \\
&\leq &\left( 1+L_{y}^{\bar{b}}+L_{y}^{\bar{\sigma}}\right) \left( L_{y}^{%
\bar{b},\bar{\sigma}}\right) ^{\frac{1}{2}}\left[ \left( N-2\right) \frac{%
\tilde{L}}{N^{\beta }}+\left( N-2\right) \left( \frac{\tilde{L}}{N^{\beta }}%
+\left( N-1\right) \frac{\tilde{L}}{N^{2\beta }}\right) \right]  \\
&&+2\left( L_{y}^{\bar{b},\bar{\sigma}}\right) ^{\frac{1}{2}}\left(
N-1\right) \frac{\tilde{L}}{N^{2\beta }} \\
&&+\left( L_{y}^{\bar{b},\bar{\sigma}}\right) ^{\frac{1}{2}}\left[ L+\left(
N-2\right) \frac{\tilde{L}}{N^{2\beta }}\right] +\left[ L_{y}^{\bar{b},\bar{%
\sigma}}\left( \left( L_{y}^{\phi }\right) ^{2}+\left( L^{\phi }\right)
^{2}\right) \right] ^{\frac{1}{2}}\left[ L+\left( N-2\right) \frac{\tilde{L}%
}{N^{2\beta }}\right]  \\
&&+\left( L_{y}^{\bar{b},\bar{\sigma}}\right) ^{\frac{1}{2}}\left[ L+\left(
N-2\right) \frac{\tilde{L}}{N^{2\beta }}\right] +\left[ L_{y}^{\bar{b},\bar{%
\sigma}}\left( \left( L_{y}^{\phi }\right) ^{2}+\left( L^{\phi }\right)
^{2}\right) \right] ^{\frac{1}{2}}\left[ L+\left( N-2\right) \frac{\tilde{L}%
}{N^{2\beta }}\right]  \\
&&+2\left( \left( L_{y}^{\phi }\right) ^{2}+\left( L^{\phi }\right)
^{2}\right) ^{\frac{1}{2}}\left( L_{y}^{\bar{b},\bar{\sigma}}\right) ^{\frac{%
1}{2}}\left[ L+\left( N-2\right) \frac{\tilde{L}}{N^{2\beta }}\right] ,
\end{eqnarray*}%
from which we need $\beta \geq 1.$

For $\tilde{C}_{3}^{i,j},$

\begin{eqnarray*}
\tilde{C}_{3}^{i,j} &=&L_{y}^{\bar{b},\bar{\sigma}}\sum_{\substack{ h\in
I_{N}\backslash \left\{ i\right\}  \\ \ell \in I_{N}\backslash \left\{
j\right\} }}^{N}\left\Vert \partial _{x_{h}x_{\ell }}^{2}\Delta
_{i,j}^{f}\right\Vert _{L^{\infty }}+L_{y}^{\bar{b},\bar{\sigma}}\left(
\left( L_{y}^{\phi }\right) ^{2}+\left( L^{\phi }\right) ^{2}\right) ^{\frac{%
1}{2}}\sum_{\substack{ h\in I_{N}\backslash \left\{ i\right\}  \\ \ell \in
I_{N}\backslash \left\{ j\right\} }}^{N}\left\Vert \left( \partial
_{x_{h}u_{\ell }}^{2}\Delta _{i,j}^{f}\left( \mathbf{X}_{t},\mathbf{u}%
_{t}\right) \right) \right\Vert _{L^{\infty }} \\
&&+L_{y}^{\bar{b},\bar{\sigma}}\left( \left( L_{y}^{\phi }\right)
^{2}+\left( L^{\phi }\right) ^{2}\right) ^{\frac{1}{2}}\sum_{\substack{ h\in
I_{N}\backslash \left\{ i\right\}  \\ \ell \in I_{N}\backslash \left\{
j\right\} }}^{N}\left\Vert \left( \partial _{u_{h}x_{\ell }}^{2}\Delta
_{i,j}^{f}\left( \mathbf{X}_{t},\mathbf{u}_{t}\right) \right) \right\Vert
_{L^{\infty }} \\
&&+\left( \left( L_{y}^{\phi }\right) ^{2}+\left( L^{\phi }\right)
^{2}\right) \left( L_{y}^{\bar{b}}+3\left( L_{y}^{\bar{\sigma}}\right)
^{2}\right) \sum_{\substack{ h\in I_{N}\backslash \left\{ i\right\}  \\ \ell
\in I_{N}\backslash \left\{ j\right\} }}^{N}\left\Vert \partial
_{u_{h}u_{\ell }}^{2}\Delta _{t,i,j}^{f}\left( \mathbf{X}_{t},\mathbf{u}%
_{t}\right) \right\Vert _{L^{\infty }} \\
&&+L_{y}^{\bar{b},\bar{\sigma}}\left( \sum_{\ell \in I_{N}\backslash \left\{
j\right\} }\left\Vert \partial _{x_{i}x_{\ell }}^{2}\Delta
g_{i,j}\right\Vert _{L^{\infty }}+\sum_{h\in I_{N}\backslash \left\{
i\right\} }\left\Vert \partial _{x_{h}x_{j}}^{2}\Delta g_{i,j}\right\Vert
_{L^{\infty }}\right)  \\
&&+\sqrt{\Gamma _{1}}\Bigg [\left( L^{\bar{b}}+L^{\bar{\sigma}}\right)
\left( L_{y}^{\bar{b},\bar{\sigma}}\right) ^{2}+2\left( L_{y}^{\bar{b}%
}+L_{y}^{\bar{\sigma}}\right) \left( 1+L_{y}^{\bar{b},\bar{\sigma}}+\left(
L_{y}^{\bar{b},\bar{\sigma}}\right) ^{2}\right) +2L_{y}^{\bar{b},\bar{\sigma}%
}\left( L^{\phi _{i}^{\prime }}+L_{y}^{\phi _{i}^{\prime }}+L^{\phi
_{j}^{\prime \prime }}+L_{y}^{\phi _{j}^{\prime \prime }}\right) \Bigg ] \\
&&+\sqrt{\Gamma _{1}}\left( L_{y}^{\bar{b}}+L_{y}^{\bar{\sigma}}\right)
\left( 1+L_{y}^{\bar{b},\bar{\sigma}}+\left( L_{y}^{\bar{b},\bar{\sigma}%
}\right) ^{2}\right)  \\
&\leq &\Bigg (L_{y}^{\bar{b},\bar{\sigma}}+L_{y}^{\bar{b},\bar{\sigma}%
}\left( \left( L_{y}^{\phi }\right) ^{2}+\left( L^{\phi }\right) ^{2}\right)
^{\frac{1}{2}}+L_{y}^{\bar{b},\bar{\sigma}}\left( \left( L_{y}^{\phi
}\right) ^{2}+\left( L^{\phi }\right) ^{2}\right) ^{\frac{1}{2}} \\
&&+\left( \left( L_{y}^{\phi }\right) ^{2}+\left( L^{\phi }\right)
^{2}\right) \left( L_{y}^{\bar{b}}+3\left( L_{y}^{\bar{\sigma}}\right)
^{2}\right) \Bigg )\left[ \left( \left( N-1\right) ^{2}-\left( N-2\right)
\right) \frac{\tilde{L}}{N^{2\beta }}+\left( N-2\right) \frac{\tilde{L}}{%
N^{\beta }}\right]  \\
&&+2L_{y}^{\bar{b},\bar{\sigma}}\left[ \left( N-2\right) \frac{\tilde{L}}{%
N^{2\beta }}+\frac{\tilde{L}}{N^{\beta }}\right]  \\
&&+\sqrt{\Gamma _{1}}\Bigg [\left( L^{\bar{b}}+L^{\bar{\sigma}}\right)
\left( L_{y}^{\bar{b},\bar{\sigma}}\right) ^{2}+2\left( L_{y}^{\bar{b}%
}+L_{y}^{\bar{\sigma}}\right) \left( 1+L_{y}^{\bar{b},\bar{\sigma}}+\left(
L_{y}^{\bar{b},\bar{\sigma}}\right) ^{2}\right) +2L_{y}^{\bar{b},\bar{\sigma}%
}\left( L^{\phi _{i}^{\prime }}+L_{y}^{\phi _{i}^{\prime }}+L^{\phi
_{j}^{\prime \prime }}+L_{y}^{\phi _{j}^{\prime \prime }}\right) \Bigg ] \\
&&+\sqrt{\Gamma _{1}}\left( L_{y}^{\bar{b}}+L_{y}^{\bar{\sigma}}\right)
\left( 1+L_{y}^{\bar{b},\bar{\sigma}}+\left( L_{y}^{\bar{b},\bar{\sigma}%
}\right) ^{2}\right) .
\end{eqnarray*}%
For $\tilde{C}_{4}^{i,j},$%
\begin{eqnarray*}
\tilde{C}_{4}^{i,j} &=&\left( L_{y}^{\bar{b},\bar{\sigma}}\right) ^{2}\sum
_{\substack{ h\in I_{N}\backslash \left\{ i\right\}  \\ \ell \in
I_{N}\backslash \left\{ j\right\} }}\left\Vert \partial _{x_{h}x_{\ell
}}^{2}\Delta g_{i,j}\right\Vert _{L^{\infty }}+\sqrt{\Gamma _{1}}\left(
L_{y}^{\bar{b},\bar{\sigma}}\right) ^{2}\left( L_{y}^{\bar{b}}+L_{y}^{\bar{%
\sigma}}\right)  \\
&\leq &\left( L_{y}^{\bar{b},\bar{\sigma}}\right) ^{2}\left[ \left( \left(
N-1\right) ^{2}-\left( N-2\right) \right) \frac{\tilde{L}}{N^{2\beta }}%
+\left( N-2\right) \frac{\tilde{L}}{N^{\beta }}\right] +\sqrt{\Gamma _{1}}%
\left( L_{y}^{\bar{b},\bar{\sigma}}\right) ^{2}\left( L_{y}^{\bar{b}}+L_{y}^{%
\bar{\sigma}}\right) .
\end{eqnarray*}%
At last
\begin{eqnarray*}
\Gamma _{1} &=&\mathbb{E}\Bigg [\sum_{\ell \in I_{N}}\left\vert \left(
\partial _{x_{\ell }}\Delta _{i,j}^{g}\right) \left( 0\right) \right\vert
^{2}+\sum_{\ell ,k\in I_{N}}\left\Vert \partial _{x_{\ell }x_{k}}\Delta
_{i,j}^{g}\right\Vert _{L^{\infty }}^{2}+\sum_{\ell ,k\in I_{N}}\left\Vert
\partial _{x_{\ell }x_{k}}\Delta _{i,j}^{f}\right\Vert _{L^{\infty
}}^{2}+\sum_{\ell ,k\in I_{N}}\left\Vert \partial _{x_{\ell }u_{k}}\Delta
_{i,j}^{f}\right\Vert _{L^{\infty }}^{2} \\
&&+3\int_{0}^{T}\Bigg (\sum_{\ell }\left\vert \left( \partial _{x_{\ell
}}\Delta _{i,j}^{f}\right) \left( t,0,0\right) \right\vert ^{2}\Bigg )%
\mathrm{d}t\Bigg ] \\
&\leq &\left( N-1\right) \frac{\tilde{L}}{N^{\beta }}+N\cdot \frac{\tilde{L}%
}{N^{\beta }}+\left( N^{2}-N\right) \frac{\tilde{L}}{N^{2\beta }}+2L+\left(
N-2\right) \frac{\tilde{L}}{N^{\beta }}+\left( N^{2}-N\right) \frac{\tilde{L}%
}{N^{2\beta }} \\
&&+6TL+3\left( N-2\right) T\frac{\tilde{L}}{N^{\beta }}.
\end{eqnarray*}

\end{document}